\documentclass[11pt]{article}

\usepackage[margin=1in]{geometry}
\usepackage[utf8]{inputenc}
\usepackage[T1]{fontenc}
\usepackage{amsthm}
\usepackage[numbers,sort&compress]{natbib}
\bibpunct[,~]{[}{]}{,}{n}{}{,}
\usepackage{bm}
\usepackage{nicefrac}
\usepackage{microtype}
\usepackage{tikz}
\usetikzlibrary{shapes.geometric,arrows.meta,positioning,fit}
\usepackage{ifthen}
\usepackage{aliascnt}
\usepackage{algorithm}
\usepackage{algorithmic}
\usepackage{authblk}
\usepackage{appendix}
\usepackage{enumitem}
\usepackage{tcolorbox}
\usepackage{multirow}
\usepackage{adjustbox}
\usepackage{times}
\usepackage{toolmp}

\colorlet{ForestGreen}{black}
\colorlet{teal}{black}
\let\argmax\relax
\let\argmin\relax

\newcommand{\av}{{\bf a}}
\newcommand{\bv}{{\bf b}}

\DeclareMathOperator\R{\mathbb{R}}

\def\textiid{i.i.d.\@\xspace}
\newcommand\iid{\ifmmode\text{ i.i.d. } \else \textiid \fi}

\newcommand{\beqs}{\begin{equation*}}
\newcommand{\eeqs}{\end{equation*}}
\newcommand{\beq}{\begin{equation}}
\newcommand{\eeq}{\end{equation}}

\DeclareMathOperator*{\argmax}{arg\,max}
\DeclareMathOperator*{\argmin}{arg\,min}

\newtheorem{thm}{Theorem}[section]

\newaliascnt{lem}{thm}
\newtheorem{lem}[lem]{Lemma}
\aliascntresetthe{lem}
\crefname{lem}{lemma}{lemmas}
\Crefname{lem}{Lemma}{Lemmas}

\newaliascnt{prop}{thm}
\newtheorem{prop}[prop]{Proposition}
\aliascntresetthe{prop}
\crefname{prop}{proposition}{propositions}
\Crefname{prop}{Proposition}{Propositions}

\newaliascnt{corr}{thm}
\newtheorem{corr}[corr]{Corollary}
\aliascntresetthe{corr}
\crefname{corr}{corollary}{corollaries}
\Crefname{corr}{Corollary}{Corollaries}

\newaliascnt{rmk}{thm}
\newtheorem{rmk}[rmk]{Remark}
\aliascntresetthe{rmk}
\crefname{rmk}{remark}{remarks}
\Crefname{rmk}{Remark}{Remarks}

\newaliascnt{mydef}{thm}
\newtheorem{mydef}[mydef]{Definition}
\aliascntresetthe{mydef}
\crefname{mydef}{definition}{definitions}
\Crefname{mydef}{Definition}{Definitions}

\newaliascnt{assumption}{thm}
\newtheorem{assumption}[assumption]{Assumption}
\aliascntresetthe{assumption}
\crefname{assumption}{assumption}{assumptions}
\Crefname{assumption}{Assumption}{Assumptions}

\newaliascnt{step}{thm}

\aliascntresetthe{step}
\crefname{step}{step}{steps}
\Crefname{step}{Step}{Steps}

\newtheorem{example}{Example}
\crefname{example}{example}{examples}
\Crefname{example}{Example}{Examples}

\crefname{appendix}{appendix}{appendices}
\Crefname{appendix}{Appendix}{Appendices}

\newtheorem{result}{Result}

\newcommand{\Leb}{{\texttt{Leb}}}
\newcommand{\SQG}{{\texttt{SQ-G}}}

\newcommand{\II}{\mathrm{I\!I}}
\newcommand{\op}{\mathrm{op}}

\newcommand{\optimalsetlower}{\mathcal{S}}
\newcommand{\SQR}{\texttt{SQ}}

\newcommand{\SQ}{\texttt{SQ-G}}

\newcommand{\hyperobjective}{F_{\max}}

\newcommand{\VaR}{\text{Q}}

\newcommand{\gibbs}{\mu}

\newcommand{\ud}{\mathrm{d}}
\newcommand{\PLcirc}{\texttt{P\L$^\circ$}}
\newcommand{\phivar}{\phi_{\lambda,\delta}}

\newcommand{\dimk}{\texttt{k}}

\title{Superquantile--Gibbs Relaxation for Minima-selection in Bilevel Optimization}

\author[1]{Saeed Masiha}
\author[2]{Zebang Shen}
\author[1]{Negar Kiyavash}
\author[2]{Niao He}

\affil[1]{EPFL School of Management of Technology, Station 5, 1015 Lausanne, Switzerland\\
\texttt{mohammadsaeed.masiha@epfl.ch, negar.kiyavash@epfl.ch}}
\affil[2]{ETH Department of Computer Science, Universit\"atstrasse 6, 8092 Z\"urich, Switzerland\\
\texttt{zebang.shen@inf.ethz.ch, niao.he@inf.ethz.ch}}

\date{}

\hypersetup{
  hidelinks,
  bookmarks=false,
  pdftitle={Superquantile--Gibbs Relaxation for Minima-selection in Bilevel Optimization},
  pdfauthor={Saeed Masiha, Zebang Shen, Negar Kiyavash, and Niao He},
  pdfsubject={Bilevel optimization with nonunique lower-level minimizers},
  pdfkeywords={bilevel optimization, minima selection, local Polyak--Lojasiewicz condition, superquantiles, Gibbs measures, nonsmooth optimization}
}

\begin{document}
\maketitle

\begin{abstract}
Bilevel optimization (BLO) becomes fundamentally more challenging when the
lower-level objective admits multiple minimizers. In contrast to the commonly
studied unique-minimizer setting, this regime introduces two additional
difficulties: (1) evaluating the hyper-objective \(F_{\max}\) requires a
\emph{minima-selection} step over the lower-level solution set; and
(2) \(F_{\max}\) may be discontinuous without additional structure. We address
both issues under a parameter-uniform local Polyak--\L{}ojasiewicz (P\L)
condition on the lower-level objective, denoted by \(\PLcirc\). In contrast to
the global P\L\ condition often assumed in BLO, the P\L\ inequality in
\(\PLcirc\) is imposed only near the local minima. This formulation
accommodates bounded, non-singleton minimizer sets and is motivated by
hyperparameter tuning in over-parameterized learning. We show that
\(F_{\max}\) is Lipschitz continuous and that the lower-level minimizer sets
are connected, compact, embedded submanifolds with a common intrinsic
dimension \(\dimk\).

Importantly, the intrinsic dimension \(\dimk\) determines the explicit
accuracy exponents in our complexity bound. We propose a randomized method
whose output satisfies an \((\mathcal O(\epsilon),\rho)\)-Goldstein stationarity
condition for \(F_{\max}\) in expectation. In the small-accuracy regime, and
suppressing logarithmic and fixed-problem factors, it uses at most
\(\widetilde{\mathcal O}(m^{15+12\dimk}\epsilon^{-2}(\epsilon\rho)^{-12\dimk-14})\) queries to the gradient oracle of the lower-level objective, where \(m\) is
the upper-level dimension. The key ingredient is a
\emph{Superquantile-Gibbs relaxation} that yields an approximate function-evaluation oracle for $F_{\max}$: it turns
minima selection into a sampling problem that can be efficiently addressed using standard Langevin
dynamics. To the best of our knowledge, this is the first work that rigorously tackles the minima selection
step in BLO and explicitly characterizes how the overall complexity of BLO scales with the \emph{intrinsic
dimensionality} of the lower-level problem.
\end{abstract}



\section{Introduction}

We study the pessimistic/optimistic formulation of Bilevel Optimization (BLO), which is a common paradigm in the literature \citep{ye1997exact,ye1997necessary,dempe2007new,liu2021towards,xu2014smoothing}, although alternatives to pessimistic/optimistic formulation for BLO also exist (see e.g., \citep{arbel2022non,jiang2025correspondence,bolte2025bilevel}). The pessimistic formulation of BLO is defined as follows\footnote{All results extend to the optimistic setting where minimizing is considered instead of maximizing $f$ over $x$.\label{footnote}}:
\begin{equation}\label{problem_BO_hyper_grad_ver}
\setlength{\abovedisplayskip}{5pt}
\setlength{\belowdisplayskip}{5pt}
\min_{\theta\in\Theta\subset\R^m}\hyperobjective(\theta) := \left\{\max_{x\in\optimalsetlower(\theta)\subseteq\mathbb{R}^d}f(\theta,x)\right\},
\quad\optimalsetlower(\theta):=\argmin_{x\in\mathbb{R}^{d}}g(\theta,x),
\end{equation}
where $\theta$ and $x$ are the upper- and lower-level variables, $f$ and $g$ are the upper- and lower-level objectives, respectively, and $\Theta$ is a \emph{closed compact convex} subset of $\mathbb{R}^{m}$.
We focus on the setting where the lower-level minima set $\optimalsetlower(\theta)$ is \emph{non-singleton}.
We refer to the maximization of  $f$ over $x \in \optimalsetlower(\theta)$ as \emph{minima-selection} and call $\hyperobjective$ the \emph{hyper-objective}.
BLO captures problems with a nested structure, such as hyper-parameter optimization \citep{feurer2019hyperparameter}, meta-learning \citep{bertinetto2018meta}, and reinforcement learning \citep{hong2023two,liu2021investigating}.
In particular, when the lower-level problem pertains to deep learning, the set $\optimalsetlower(\theta)$ is typically nonconvex and non-singleton \citep{lin2024exploring,kuditipudi2019explaining,nguyen2019connected,garipov2018loss,draxler2018essentially,liang2018understanding,venturi2018spurious,sagun2017empirical}. Therefore, minima-selection becomes important.

\paragraph{Relaxing strong convexity to the global P\L\ condition }
Recent works on BLO have studied nonconvex lower-level problems under the \emph{global} Polyak-Łojasiewicz (P\L) condition \citep{chen2023bilevel,kwon2023penalty,xiao2023generalized,huang2024optimal,shen2025penalty,liu2022bome}. This condition relaxes the classical strong convexity assumption on \(g\), which implies singleton \(\mathcal{S}(\theta)\), to cover  nonconvex and potentially non-singleton optimal solution sets. Another motivation for adopting the global P\L\ condition is to guarantee the continuity of the hyper-objective function \(\hyperobjective\), which might otherwise become discontinuous and render BLO problems generally intractable \citep{kwon2023penalty,chen2023bilevel}. 

All of the aforementioned  work and their results are crucially dependent on the assumption that \(g\) is \(\mathcal{C}^2\) and satisfies the global P\L\ condition.  Under these conditions, it can be proved that the lower-level solution set \(\optimalsetlower(\theta)\) would be limited to either a singleton or an unbounded subset of \(\mathbb{R}^d\) \citep{criscitiello2025}.
More specifically, Huang \cite{huang2024optimal} assume a nondegenerate Hessian of $g$ at global minimizers, implicitly enforcing uniqueness; Liu et al. \cite{liu2022bome} explicitly assume $\mathcal{S}(\theta)$ is a singleton; Xiao et al. \cite{xiao2023generalized} and Shen et al. \cite[Sec.~4]{shen2025penalty} assume additional convexity of the lower-level objective, implying that the minimizer set is a linear subspace; and Kwon et al. \cite{kwon2023penalty} require $g$ to be coercive which again yields a singleton solution\footnote{Kwon et al. \cite{kwon2023penalty} assume the proximal error bound for $g$ (with $\sigma=0$), equivalent to the global P\L\ condition for smooth $g$ \citep[Prop.~C.1]{chen2023bilevel}.}.
The above observations hint at the necessity of a further relaxation to capture richer landscapes.

\paragraph{Our setting: local P\L\ geometry with non-singleton minimizers.}
The limitation above motivates separating the geometry near the lower-level
minima from the behavior of the objective at infinity.  We impose the
\(\PLcirc\) condition of Gong et al.~\citep{gong2024poincare} locally around
the minima, while the tail and smoothness conditions in
\Cref{assum_sample_comp} provide coercivity.  Under these combined
conditions, \(\optimalsetlower(\theta)\) is a connected, compact
\(\mathcal C^2\) embedded submanifold of \(\mathbb R^d\) without boundary.
Thus, unlike the globally P\L\ and coercive setting discussed above, our
framework permits bounded non-singleton minimizer sets and retains a genuine
minima-selection problem.  Such connected sets of minima are also observed
in over-parameterized learning models
\citep{draxler2018essentially,garipov2018loss,nguyen2019connected,kuditipudi2019explaining,lin2024exploring}.

The uniform parametric assumptions additionally make the minimizer mapping
Hausdorff--Lipschitz, and \Cref{lemm_hyp_blo_cont} then shows that
\(\hyperobjective\) is Lipschitz continuous.  They do not, in general, make
\(\hyperobjective\) differentiable; see \Cref{append_diff_hyper_obj}.  We
therefore next specify the nonsmooth solution concept and the approximation
used to address minima-selection.

\paragraph{Solution concept for nonsmooth nonconvex optimization.}
The preceding discussion leaves two algorithmic issues.  First, evaluating
\(F_{\max}(\theta)\) requires maximizing over the generally nonconvex set
\(\optimalsetlower(\theta)\).  Second, even if this minima-selection problem
could be solved exactly, \(F_{\max}\) need not be differentiable.  We address
the second issue first by choosing a tractable stationarity notion for
nonsmooth Lipschitz optimization.  We then address the first by constructing
a superquantile--Gibbs surrogate that provides approximate values of
\(F_{\max}\).
\subsection{Finding a Goldstein Stationary Point of the Hyper-objective}\label{sec:intro_Goldstein_point}
For a nonsmooth Lipschitz objective, Clarke stationarity is the classical
first-order target \citep{clarke1990optimization}.  However, Zhang et
al.~\cite{zhang2020complexity} show that finding an
\(\epsilon\)-Clarke-stationary point in polynomial time is impossible in
general.  Goldstein stationarity introduces a spatial tolerance \(\rho\) by
convexifying Clarke subgradients at points in a \(\rho\)-neighborhood; this
weaker target admits meaningful algorithmic guarantees.  We therefore seek
an \((\epsilon,\rho)\)-Goldstein-stationary point
\citep{goldstein1977optimization}.

For the domain-constrained objective in \eqref{problem_BO_hyper_grad_ver},  akin to \citep{liu2024zeroth}, we measure stationarity by a \textit{generalized }$(\epsilon,\rho)$-Goldstein criterion which is defined for the gradient mapping\footnote{The gradient mapping for $F:\Theta\to \mathbb{R}$ at a point $\theta$ is $\mathcal{G}_{\Theta}(\theta,\nabla F;\eta):=\eta^{-1}(\theta-\text{Proj}_{\Theta}(\theta-\eta\nabla F(\theta)))$.}, denoted by  $(\epsilon,\rho,\eta)$-generalized Goldstein stationary point (see \Cref{def_gen_gold_stat}). 
As $F_{\max}$ may not be differentiable, we propose a zeroth-order (ZO) method to find such a stationary point of $F_{\max}$.

The major obstacle to applying ZO methods to minimize $F_{\max}$ is that exact function values of this hyper-objective are typically unavailable. To address this, we introduce a tractable surrogate, the superquantile-Gibbs $\tilde F_{\SQG}(\theta)$, which we define next and show that it approximates $F_{\max}$ well. Consequently, querying $\tilde F_{\SQG}$ in place of $F_{\max}$  provides the finite-difference gradient estimates for our zeroth-order (ZO) method. We prove that the approximation error incurred by replacing $F_{\max}$ with $\tilde F_{\SQG}$ does not adversely affect the convergence guarantees of our ZO algorithm for obtaining an $(\epsilon,\rho,\eta)$-Goldstein stationary point.

\subsection{Superquantile-Gibbs Relaxation (\SQG)}\label{sec:intro_SQG}
In this part, we introduce the superquantile of the random variable \(f(\theta, X)\) where \(X\) is sampled according to the the Gibbs measure \(\sim \exp(-g(\theta, x)/\lambda)\).
Our idea is that as $\lambda\rightarrow 0$, the Gibbs measure concentrates on the feasible domain $\optimalsetlower(\theta)$, and the superquantile at a high confidence level would approach the maximum value of $f(\theta, \cdot)$ on $\optimalsetlower(\theta)$.
Formally, we define the Superquantile-Gibbs Relaxation (\SQ)-loss using the variational formulation of the superquantile by \citep{rockafellar2000optimization}:
\begin{tcolorbox}[colback=gray!20, colframe=gray!50, boxrule=0pt, arc=0mm, left=0mm, right=0mm, top=0mm, bottom=0mm]
    \begin{equation} \label{eqn:SQ_gibbs_variation}
    \begin{aligned} 
        (\text{\SQ-loss})\quad &\ \tilde F_{\SQ}(\theta) := \min_{\beta\in\mathbb{R}} \phi_{\lambda,\delta}(\theta,\beta) \\
        \text{with }&\ \phi_{\lambda,\delta}(\theta,\beta):=\beta+\frac{1}{\delta}\mathbb{E}_{ X\sim\gibbs^\lambda_\theta}[\max\{f(\theta,X)-\beta,0\}].
    \end{aligned}    
    \end{equation}
\end{tcolorbox}  
\noindent Here $\lambda>0$ and $0<\delta<1$ are parameters controlling the level of \SQG\ relaxation, with $(1-\delta)$ known as the confidence level, and $\gibbs^\lambda_\theta$ denotes the Gibbs measure:
\begin{equation}\label{eq_gibbs}
\setlength{\abovedisplayskip}{5pt}
            \setlength{\belowdisplayskip}{5pt}
\gibbs^\lambda_\theta(x):=\frac{\exp({-{g(\theta,x)}/{\lambda}})}{Z_{\lambda}(\theta)},
\end{equation}
where $Z_{\lambda}(\theta):=\int_{\mathbb{R}^{d}} \exp({-{g(\theta,z)}/{\lambda}})dz$ is a normalizing factor. A key contribution of this project is that we characterize how $\lambda$ and $\delta$ impact the accuracy with which $\tilde F_{\SQ}(\theta)$ approximates $F_{\max}$.

\subsection{Our Contributions} 
Our main contributions are as follows.
\begin{itemize}[leftmargin=*]
\setlength{\itemsep}{0.5pt}
\item 
To allow lower-level optimal sets with rich structures, we replace the
    commonly imposed global P\L\ condition by the parameter-uniform local
\PLcirc\ condition of Gong et
    al.~\citep{gong2024poincare}, together with the regularity
    conditions in \Cref{assum_sample_comp}; in particular, its tail gradient
    error bound and global lower-gradient smoothness yield the quadratic
    growth that supplies coercivity.  Under these combined
    assumptions, \Cref{lemm_hyp_blo_cont} proves Lipschitz continuity of the
    hyper-objective \(\hyperobjective\).
\item We use the Superquantile–Gibbs relaxation (\SQG) to approximate  $\hyperobjective$. This converts the crucial minima-selection step in BLO, originally a nonconvex optimization over a nonconvex set, into sampling from a Gibbs distribution, which we implement efficiently via Langevin dynamics.
    \item Building on the manifold structure of the lower-level optimal set, we derive the following results.

    \begin{result}[Approximation] \label{result_approximation}
The error of the superquantile--Gibbs surrogate has two sources.  Under
    \Cref{mu_PL_assum,assum_secon_fund_form_bound,assum_sample_comp},
    \Cref{rmk_func_approx_cvar_gibbs} proves, uniformly for
    \(\theta\in\Theta\) and for the stated ranges of \(\delta\) and \(\lambda\),
    \[
    \left|F_{\max}(\theta)-\tilde F_{\SQG}(\theta)\right|
    =
    \mathcal O\!\left(
        \delta^{1/\dimk}+\frac{\sqrt\lambda}{\delta}
    \right).
    \]
    Here \(\dimk\) is the common intrinsic dimension of the minimizer
    manifolds, established in \Cref{lemma_common_dimension}.  The term
    \(\delta^{1/\dimk}\) is the error from replacing the maximum on the
    limiting minimizer manifold by a high superquantile
    (\Cref{th_CVaR_unif_approx}).  The term
    \(\sqrt\lambda/\delta\) is the additional finite-temperature error from
    replacing the limiting measure by the Gibbs measure; it follows from
    \Cref{prop_W1_Gibbs_uniform,Theo_func_approx_cvar}.
    Consequently, for sufficiently small \(\epsilon_v>0\), choosing
    \(\delta=\epsilon_v^{\dimk}\) and
    \(\lambda=\epsilon_v^{2(\dimk+1)}\) makes both terms
    \(\mathcal O(\epsilon_v)\).  This is the value-oracle accuracy used in
    the optimization result below.
The explicit constants and admissible ranges of
    \(\delta\) and \(\lambda\) are given in
    \Cref{rmk_func_approx_cvar_gibbs}.
\end{result}
\begin{result}[Optimization] \label{theorem_informal_complexity}
Suppose that
     \Cref{mu_PL_assum,assum_secon_fund_form_bound,assum_sample_comp} hold on a fixed compact thickening of \(\Theta\), as required by
     the perturbed zeroth-order queries, and let \(\dimk\ge1\).  The
     \(\mathcal O(\epsilon_v)\)-accurate surrogate above can then be used in the
     Projected Stochastic Zeroth-Order method with Minima-Selection
     (\Cref{alg: Stochastic-Zeroth-order}).  The method returns a randomized
     iterate \(\widehat\theta\) whose expected generalized-Goldstein residual
     for \(F_{\max}\) is \(\mathcal O(\epsilon)\).  With exact Gibbs samples,
     the required number of samples, and hence evaluations of \(f\), is
     \[
        \widetilde{\mathcal O}\!\left(
        m^{2\dimk+5}\epsilon^{-2}
        (\epsilon\rho)^{-2\dimk-4}
        \right).
     \]

     Replacing exact Gibbs samples by LMC preserves the same residual
     guarantee and the same order of evaluations of \(f\).  The resulting
     number of lower-level gradient evaluations is
     \[
        \widetilde{\mathcal O}\!\left(
        \mathfrak C_{\mathrm{mix}}\,
        d\,m^{15+12\dimk}\epsilon^{-2}
        (\epsilon\rho)^{-12\dimk-14}
        \right).
     \]
     Here \(\widetilde{\mathcal O}\) suppresses logarithmic and fixed-problem
     factors.  The factor \(\mathfrak C_{\mathrm{mix}}\) is uniform in the
     parameter and accuracy variables, but may depend on \(d\).
\end{result}

The precise statements of these results appear in \Cref{rmk_func_approx_cvar_gibbs}, \Cref{thm_cvg_alg_zeroth_goldstein} and \Cref{theorem_gradient_complexity}, respectively.

\item We evaluate our approach on both synthetic and real-world benchmarks and demonstrate the effect of minima-selection. As shown in \Cref{sec:experiment}, our proposed approach outperforms hyper-gradient and penalty-based baselines.
\end{itemize}

\begin{rmk}(Exponential dependence on $\dimk$)\label{rmk_1_2}
The $\epsilon^{-\dimk}$ factor in our results reflects the intrinsic difficulty of the problem.
Since evaluating $F_{\max}(\theta)=\max_{y\in S(\theta)} f(\theta, y)$ is equivalent to solving a nonconvex optimization problem over a nonconvex set, a single $\epsilon$-accurate evaluation of $F_{\max}(\theta)$ yields an $\Omega(\epsilon^{-\dimk})$ complexity lower bound. More specifically, 
consider the instance where $S(\theta) = \mathbb{S}^{\dimk}$ is a hypersphere. Reparameterizing $y$ in polar coordinates yields $F_{\max}(\theta)=\max_{x\in \mathbb{B}_{\dimk}} f(\theta, y(x))$. Here $\mathbb{B}_{\dimk}=[0, \pi]^\dimk$ and the point $y(x)\in \mathbb{S}^{\dimk}$ is given by $y_{i}(x)=\sin(x_{1})\cdots \sin(x_{i-1})\cos(2x_{i})$ for $i\in[\dimk]$ and $y_{\dimk+1}(x)=\sin(x_{1})\cdots \sin(x_{\dimk-1})\sin(2x_{\dimk})$.
It follows from \citep[Section 1.6.2]{nemirovskij1983problem} that the worst-case complexity of finding an $\epsilon$-accurate solution to is $\Omega(\epsilon^{-\dimk})$ for $f(\theta, \cdot)\in \mathcal{C}^1$.
\end{rmk}

\subsection{Related Work}\label{section_related_work}
Bilevel optimization has a long history in mathematical programming; see, e.g., the early bibliographic review by \cite{vicente1994bilevel}, the monograph \cite{dempe2002foundations}, the survey \cite{colson2007overview}, the classical text on mathematical programs with equilibrium constraints (MPECs) by \cite{luo1996mathematical}, and the recent collection \cite{dempe2020bilevel}. Classical approaches often reduced bilevel problems to single-level programs, either via value-function reformulations or by replacing the lower-level problem with its Karush--Kuhn--Tucker (KKT) optimality system. A related value-function-based reduction for bilevel programs with nonconvex lower-level objectives yields single-level formulations with nonsmooth inequality constraints, for which smoothing and augmented-Lagrangian type algorithms were developed in, e.g., \citep{lin2014solving,xu2014smoothing}. The latter route yields MPEC formulations in which complementarity conditions between multipliers and constraint violations appear as explicit equilibrium constraints; see, for example, \cite{ye1995optimality,ye2010new,dempe2012kkt,zemkoho2021comparison}. In practice, these complementarity constraints are frequently handled through mixed-integer or so-called big-\(M\) linearizations; however, their performance can be sensitive to the choice of constants, with overly small values excluding feasible bilevel solutions and overly large values leading to weak relaxations and numerical instability \citep{kleinert2020there}. Complementing these algorithmic considerations, recent work has also highlighted inherent geometric and computational hardness in bilevel programming, even under smoothness, which motivates the search for additional structural assumptions and tailored relaxations \citep{bolte2025geometric}.

\paragraph{Gradient-based bilevel methods with strongly convex lower level.}
Modern large-scale machine-learning applications led to gradient-based methods for bilevel empirical risk minimization and hyperparameter optimization. When the lower-level objective is strongly convex in the inner variable and $f,g$ are smooth, the hyper-objective is differentiable and its gradient can be obtained by implicit differentiation. Representative works in this setting include approximate hyper-gradient methods such as HOAG \citep{pedregosa2016hyperparameter}, stochastic bilevel schemes \citep{ghadimi2018approximation,yang2021provably}, and complexity-refined algorithms for empirical risk minimization \citep{dagreou2024lowerbound,chen2024optimal}, as well as more abstract treatments of optimization layers via automatic differentiation \citep{ablin2020super}, amortized implicit differentiation \citep{arbel2021amortized}, and modular implicit differentiation frameworks \citep{blondel2022efficient,grazzi2020iteration}. In these approaches, computing exact hyper-gradients typically required Hessian or Hessian-inverse information for the lower-level problem, which could be prohibitive in large neural networks. This motivated a substantial line of \emph{Hessian-free} or fully first-order methods that approximated the hyper-gradient using only first-order information, including BOME \citep{liu2022bome}, Hessian/Jacobian-free stochastic methods with improved sample complexity \citep{yang2023hessianfree}, fully first-order stochastic approximation methods such as F$^2$BO \citep{kwon2023fully}, and bilevel stochastic-gradient methods that also accommodate nonlinear lower-level constraints and inexact hyper-gradient computations \citep{giovannelli2025inexact}. 
A different line of work revisits bilevel learning from a functional viewpoint: Petrulionyt{\.e} et al.~\citep{petrulionyte2024functional} treat the lower-level variable as a \emph{prediction function} (rather than a parameter vector) and optimize it over a convex function space. Assuming strong convexity in this functional variable yields a unique lower-level solution and enables implicit differentiation via a functional implicit-function theorem, shifting the strong-convexity requirement from parameters to function space.
These developments made bilevel optimization scalable in settings with a unique lower-level solution, so the minima-selection problem that we study did not arise.

\paragraph{Relaxing strong convexity to global P{\L}.}
A complementary line of work relaxed lower-level strong convexity to a global Polyak--{\L}ojasiewicz (P{\L}) condition and developed bilevel gradient methods under this assumption \citep{chen2023bilevel,kwon2023penalty,xiao2023generalized,huang2024optimal,shen2025penalty,liu2022bome}. To obtain a smooth and well-behaved hyper-objective, these methods typically combined the global P{\L} condition with additional regularity assumptions on the lower-level problem (as discussed earlier), so that minima-selection was effectively avoided. In contrast, our analysis is based on a local P{\L} condition and explicitly exploits the geometry of nontrivial minimizer manifolds. In a complementary direction, Chen et al.~\cite{chen2025setsmoothness} introduce a \emph{set-smoothness} property of the lower-level solution mapping and show that it yields weak convexity/concavity of the induced optimistic/pessimistic hyper-objectives, enabling the computation of approximate Clarke hyper-stationary points. Their analysis relies on assumptions such as a \emph{uniform} error bound together with smoothness, which imply a global P{\L} condition for the lower-level problem. When the lower-level objective is $\mathcal{C}^{2}$ in the lower-level variable, this global P{\L} property forces $\argmin_x g(\theta,x)$ to be either a singleton or an unbounded set \cite{criscitiello2025}. As a consequence, these conditions rule out the bounded non-singleton minimizer manifolds that motivate the minima-selection regime studied in our work.

\paragraph{Alternative formulations of BLO.}
Beyond such regularity-based assumptions (strong convexity or global P{\L}), recent works proposed alternative formulations of BLO that target nonconvex lower-level problems.
For instance, Bolte et al. \cite{bolte2025bilevel} worked under a (piecewise) Morse-type parametric qualification, showed that the graphs of lower-level critical points and local minima are generically decomposed into finitely many smooth manifolds, and used this structure to define branch-wise hyper-objectives and study gradient methods that updated the upper-level variable along individual branches. Arbel and Mairal \cite{arbel2022non} introduced bilevel games with selection maps, where the selection map chooses a specific lower-level critical point for each upper-level parameter (e.g., the limit of a gradient flow from a given initialization). They showed that popular unrolled/iterative-differentiation bilevel algorithms with a finite inner budget approximate equilibria of these selection-based games up to approximation errors, and they extended implicit differentiation to characterize the selection and proposed a correction that removes the finite-budget error. Jiang et al. \cite{jiang2025correspondence} further argued that, for general nonconvex lower levels, the classical hyper-objective may be ill-posed when global minimizers are unattainable; they proposed a correspondence-driven hyper-objective that evaluated the upper-level loss at the output of a prescribed lower-level algorithm, applied Gaussian smoothing to this mapping, and developed a projected stochastic gradient method with  oracle-complexity guarantees for the smoothed problem. 
Another recent line develops \emph{Moreau-envelope} single-level reformulations as an alternative to classical value-function and KKT-based complementarity-constrained reductions.
Bai et al.~\citep{bai2025optimality} derive directional necessary optimality conditions in this framework, Gao et al.~\citep{gao2023moreau} propose a difference-of-weakly-convex reformulation with a sequentially convergent algorithm, and Lu~\citep{Lu25tsp} builds on these ideas to obtain a two-sided smoothed primal--dual method for general \emph{nonconvex} bilevel problems with KKT-type convergence guarantees.
In contrast, our work retains the classical (pessimistic or optimistic) hyper-objective, combines it with a local P{\L} condition to obtain a manifold structure for the lower-level minimizers, and uses a superquantile--Gibbs relaxation to convert the minima-selection into sampling from a Gibbs measure.

\paragraph{Pessimistic bilevel optimization.}
The approaches discussed above treat lower-level non-uniqueness either by enforcing regularity conditions (strong convexity or global P{\L}) or by embedding explicit selection mechanisms (Morse branches, selection maps, or algorithmic correspondences), and thus remain essentially in an optimistic or selection-based setting. By contrast, the \emph{pessimistic} formulation evaluates each upper-level decision against the worst upper-level value over all lower-level minimizers, which is natural when the follower may act adversarially or unreliably \citep{beck2023survey}. This risk-averse viewpoint was studied in mathematical programming, including structural and complexity results for general pessimistic models \citep{wiesemann2013pessimistic}, surveys of modelling and solution schemes \citep{liu2018survey}, and necessary optimality conditions based on generalized differentiation and value-function techniques \citep{dempe2014necessary}. On the algorithmic side, several works constructed single-level relaxations of the pessimistic problem: Zeng and Ye \cite{zeng2020practical} designed a tight relaxation and a computation scheme for (mainly linear) pessimistic bilevel problems; Benchoouk et al. \cite{benchouk2025scholtes} analyzed a Scholtes-type relaxation of the KKT system; and Antoniou et al. \cite{antoniou2024delta} proposed a $\delta$-perturbed formulation that enforced an \emph{approximate} unique follower response and derived error bounds between the perturbed and original pessimistic problems. Pessimistic ideas also appeared in machine learning: Ustun et al. \cite{ustun2024hyperparameter} formulated hyperparameter tuning as a pessimistic bilevel problem over a finite set of tasks and highlighted advantages of this formulation over the optimistic counterpart, and Jiang et al. \cite{jiang2021learning} interpreted adversarial training as a bilevel game in which a learned attacker played the role of the follower. All these approaches, however, either (i) resolved non-uniqueness by modifying the lower-level problem (e.g., $\delta$-perturbations, KKT relaxations) or (ii) specialized to finite or discrete follower responses, rarely addressed the continuous pessimistic setting with a non-singleton solution set.
In contrast, we work directly with the pessimistic
hyper-objective $F_{\max}$ under a local P{\L} condition and model pessimism
through a superquantile risk over Gibbs samples on the minimizer manifolds.
The resulting accuracy exponents depend explicitly on their common intrinsic
dimension $\dimk$; the LMC implementation also carries the ambient-dimension
factors recorded in \Cref{theorem_gradient_complexity}.

\subsection{Organization}
The remainder of the paper is organized as follows.
\Cref{sec:prel} introduces the \PLcirc\ and regularity assumptions and develops
the fixed-parameter minimizer-manifold and limiting-Gibbs preliminaries.
\Cref{sec:4_unif_props} then establishes the parametric results that are uniform
over $\theta\in\Theta$, including Hausdorff--Lipschitz continuity of the
solution mapping, common manifold geometry, and the uniform Wasserstein
and low-temperature Poincar\'e estimates.  Using these uniform ingredients,
\Cref{sec:3} proves the
dimension-dependent superquantile and Gibbs approximation bounds and combines
them into the uniform \SQG\ approximation of $F_{\max}$ in
\Cref{rmk_func_approx_cvar_gibbs}.  Building on that approximation,
\Cref{sec:6} establishes the zeroth-order complexity guarantees in
Theorems~\ref{thm_cvg_alg_zeroth_goldstein}
and~\ref{theorem_gradient_complexity}.
Their accuracy exponents display the common intrinsic
dimension of the lower-level solution manifolds, while the LMC implementation
also retains the dimension-dependent mixing and warm-start factors stated in
\Cref{theorem_gradient_complexity}.  Finally, \Cref{sec:experiment} presents the
numerical study, first on a synthetic dimension-dependent example and then on
a toy problem and a data hyper-cleaning task.

\subsection{Notations}
Unless stated otherwise, $\|\cdot\|$ denotes the Euclidean ($\ell_2$) vector norm and the associated operator norm on matrices and higher-order tensors. For a random variable $Z$, $\mathrm{supp}(Z)$ denotes its support.  For $d \in \mathbb{N}$, $x \in \mathbb{R}^d$, and $r > 0$, we write $ \mathbb{B}_{d}(x;r):= \{y \in \mathbb{R}^d : \|y - x\| \le r\}$
for the closed $d$-dimensional Euclidean ball with center $x$ and radius $r$.
For a Riemannian manifold $M$ and $p\in M$,
$\exp_p:T_pM\to M$ denotes the Riemannian exponential map at $p$.  Thus,
for $v\in T_pM$ and every $t$ for which it is defined,
$\exp_p(tv)$ is the point at time $t$ on the geodesic with initial point
$p$ and initial velocity $v$.
We write $s\uparrow a$ when $s$ approaches $a$ through
values smaller than $a$, and $s\downarrow a$ when $s$ approaches $a$ through
values larger than $a$.
For $p \ge 1$ and probability measures $\mu,\nu$ on $\mathbb{R}^d$ with finite $p$-th moments, the $p$-Wasserstein distance is
\[
\mathbb{W}_p(\mu,\nu):= \left[ \inf_{\pi \in \Pi(\mu,\nu)}\int_{\mathbb{R}^d \times \mathbb{R}^d} \|x-y\|^{p} \,\mathrm{d}\pi(x,y)\right]^{1/p},
\]
where $\Pi(\mu,\nu)$ is the set of all couplings of $\mu$ and $\nu$. 
In particular, $\mathbb{W}_1$ is the 1-Wasserstein distance.
Given a \(\mathcal C^1\) function
\(f:\mathbb R^m\times\mathbb R^d\to\mathbb R\), we write
\(\partial_\theta f(\theta,x)\) and \(\partial_xf(\theta,x)\) for its
gradients with respect to the first and second arguments, respectively.
The diameter of a set $\Theta \subset \mathbb{R}^{d}$, is defined as  $\text{diam}(\Theta):= \sup_{\theta_{1},\theta_{2} \in \Theta} \|\theta_{1} - \theta_{2}\|$. The class of continuous real-valued functions on $\mathbb{R}^{d}$ is denoted by $\mathcal{C}(\mathbb{R}^{d},\mathbb{R})$. We write $f \in \mathcal{C}^{k}$ if $f$ has continuous derivatives up to order $k$. For a set $A \subset \mathbb{R}^d$, we denote its closure by $\overline{A}$.  For a scalar random variable $Z$ and $\delta \in (0,1)$, the $(1-\delta)$-quantile is defined as
\[
    \VaR_{1-\delta}(Z)
    := \inf\{t \in \mathbb{R} : \mathbb{P}[Z > t] \le \delta\}.
\]
For a closed convex set
$\Theta\subset\mathbb R^m$ and $y\in\mathbb R^m$, its Euclidean projection is $\operatorname{Proj}_{\Theta}(y):=\argmin_{z\in\Theta}\|y-z\|.$
We use the notation $\mathcal{O}$ to hide constants, and the notation $\tilde{\mathcal{O}}$ to hide both constants and logarithmic factors. We write $a \lesssim b$ to indicate that $a \le c b$ for some universal positive constant $c > 0$ that is independent of the problem parameters. We write $\mathbb{R}_\alpha(P\Vert Q)$ for the Rényi divergence of order $\alpha\in(0,\infty)\setminus\{1\}$, defined by
\[
\mathbb{R}_\alpha(P\Vert Q)=
\begin{cases}
\displaystyle \frac{1}{\alpha-1}\log \!\int \!\Big(\frac{\mathrm{d}P}{\mathrm{d}Q}\Big)^{\alpha}\,\mathrm{d}Q, & P\ll Q,\\[1ex]
+\infty, & \text{otherwise.}
\end{cases}
\]
We denote $\mathrm{Law}(X)$ as the probability distribution (or law) of a random variable $X$.
 \section{Summary of Assumptions and Preliminaries}\label{sec:prel}
We begin by recalling the local Polyak--Lojasiewicz (\PLcirc) condition,
which describes the geometry of the lower-level objective near its local
minimizers.  The coercivity needed for the global analysis is derived later
from the parametric lower-level assumptions; see the discussion following
\Cref{assum_sample_comp}.

\begin{mydef}\label{def_local_PL}
    Let $\mathcal{M}$ be the collection of all local minima of a function \( g \in \mathcal{C}^{1}(\mathbb{R}^{d}, \mathbb{R}) \).
    \( g\) is locally P\L\ if there exists a constant \(\mu > 0\) such that for any connected component $\mathcal{M}'\subseteq\mathcal{M}$, there exists an open neighborhood \(\mathcal{N}(\mathcal{M}')\) containing \(\mathcal{M}'\) (\(\mathcal{N}(\mathcal{M}') \supset \mathcal{M}'\)) such that 
    \begin{equation}
    \setlength{\abovedisplayskip}{5pt}
            \setlength{\belowdisplayskip}{5pt}
        \forall x \in \mathcal{N}(\mathcal{M}'), \quad g(x) - \min_{x' \in \mathcal{N}(\mathcal{M}')} g(x') \leq (2\mu)^{-1} \|\nabla g(x)\|^{2}.
    \end{equation}
\end{mydef}

\begin{mydef}\label{def_mu_PL_assum}
    A function $g: \R^d \rightarrow \R $ satisfies the \PLcirc\ condition if the following hold.
    \begin{itemize}\setlength{\itemsep}{0.5pt}
        \item $g$ is a $\mathcal{C}^{2}$ and locally P\L\ function, and $g$ is a $\mathcal{C}^{4}$ function in every neighborhood $\mathcal{N}(\mathcal{M}')$\footnote{Gong et al.~\citep{gong2024poincare} require
        \(g\in\mathcal C^3\).
        The pointwise $\mathcal C^4$ clause supports the fixed-parameter estimate
        in \Cref{prop_W1_Gibbs_pointwise}; it does not by itself imply the
        uniform local Hessian-Lipschitz clause of
        \Cref{assum_sample_comp}.}.
        \item Let $\mathcal{N}(\mathcal{M}):= \cup_{\mathcal{M}'}\mathcal{N}(\mathcal{M}')$. For any $x\in \mathbb{R}^{d}\setminus \mathcal{N}(\mathcal{M})$, if $\nabla g(x)=0$, then  $\nabla^{2}g(x)\prec 0$.
        \item The collection of all local minima \(\mathcal{M}\) is contained within a compact set.
    \end{itemize}
\end{mydef}
The global conclusions used below require coercivity in addition to the
\PLcirc\ condition.
\begin{prop}[Connectivity, Proposition~3 in \citep{gong2024poincare}]\label{prop_S_path_connected}
    Let \(d\ge 2\), and let \(g:\mathbb{R}^{d}\to\mathbb{R}\) be coercive and
    satisfy the \PLcirc\ condition.  Then its collection \(\mathcal M\) of
    local minima is connected.
\end{prop}

Throughout the paper, we assume the lower-level ambient dimension satisfies
\(d\ge2\).  Rebjock and Boumal \citep{rebjock2024fast} establish the local
manifold structure associated with the local P\L\ condition;
Gong et al.~\citep{gong2024poincare} obtain the following
global conclusion under
coercivity and the full \PLcirc\ condition.

\begin{prop}[Manifold structure, Corollary~1 in \citep{gong2024poincare}]\label{prop_S_C_2_wihout}
    Let \(d\ge2\).  If \(g:\mathbb R^d\to\mathbb R\) is coercive and satisfies
    the \PLcirc\ condition, then all local minima have the same value and $\mathcal M=\argmin_{x\in\mathbb R^d}g(x)$.
    Moreover, this global minimizer set is a compact \(\mathcal C^2\) embedded
    submanifold of \(\mathbb R^d\) without boundary.
\end{prop}
\begin{rmk}[No-saddle condition]
The second condition in \Cref{def_mu_PL_assum} has two distinct uses.  First,
together with coercivity and the remaining \PLcirc\ clauses, it is an
ingredient in the mountain-pass argument of
Gong et al.~\citep{gong2024poincare}, which establishes the
connectedness conclusion in \Cref{prop_S_path_connected}.  Second, once
\Cref{prop_S_path_connected,prop_S_C_2_wihout} identify the local minima with
the global minimizer set, the local-P\L\ inequality shows that a critical point
inside a local-P\L\ neighborhood is a local, hence global, minimizer, whereas
the no-saddle clause makes every critical point outside those neighborhoods
have a negative-definite Hessian.  This yields the critical-point dichotomy
used in \Cref{prop_parametric_morse_bott_stability}.  Assuming connectedness
directly would replace only the first use; it would not supply the no-saddle
conclusion needed for the second.
\end{rmk}
Combining \Cref{prop_S_path_connected,prop_S_C_2_wihout}, the global minimizer
manifold is connected and hence path-connected.  In particular, any two of
its points can be joined within the manifold by a piecewise \(\mathcal C^1\)
path.

\begin{tcolorbox}[colback=gray!20, colframe=gray!50, boxrule=0pt, arc=0mm, left=0mm, right=0mm, top=0mm, bottom=0mm]
\begin{assumption}\label{mu_PL_assum}
There exists a constant $\mu>0$ such that, for every $\theta\in\Theta$, the
    lower-level objective $g(\theta,\cdot)$ satisfies the \PLcirc\ condition and the
    constant in the local P\L\ inequality of \Cref{def_local_PL} can be chosen to be
    this same $\mu$.  The corresponding local P\L\ neighborhoods may depend on
    $\theta$.
\end{assumption}
\end{tcolorbox}
Following \citep{gong2024poincare}, we impose the next
coordinate-free bound to control the extrinsic curvature of the minimizer
manifolds.  We also use it to establish a nonzero lower bound on the injectivity radius, allowing us to analyze the
volume-radius relationship of geodesic balls, which is crucial for our approximation results.
\begin{tcolorbox}[colback=gray!20, colframe=gray!50, boxrule=0pt, arc=0mm, left=0mm, right=0mm, top=0mm, bottom=0mm]
\begin{assumption}[Second-order fundamental form]\label{assum_secon_fund_form_bound}
    For every $\theta \in \Theta$, the manifold $\mathcal{S}(\theta)$ has a well-defined second fundamental form in the sense of \Cref{def:second-fund-form} in \Cref{append:def}, and its operator norm is uniformly bounded by a constant $\mathsf{C}>0$.
\end{assumption}
\end{tcolorbox}

\noindent We further need the following regularity conditions, which are standard in the BLO literature.
\begin{tcolorbox}[colback=gray!20, colframe=gray!50, boxrule=0pt, arc=0mm, left=0mm, right=0mm, top=0mm, bottom=0mm]
\begin{assumption}\label{assum_sample_comp}
The following statements hold for all $\theta \in \Theta$.
\begin{itemize}[leftmargin=*]\setlength{\itemsep}{0.5pt}
\item The function \(f\) is \(\mathcal C^1\) and
    \(L_{f,1}\)-Lipschitz with respect to \((\theta,x)\). For every
    \(x\in\mathbb R^d\),
    \[
        |f(\theta,x)|\le C_f\|x\|^{\mathsf n_f}+D_f,
    \]
    where \(C_f,D_f\ge0\) and \(\mathsf n_f\ge1\) is an integer.

    \item The lower-level objective $g$ admits a $\mathcal C^2$
    extension to an open neighborhood of $\Theta\times\mathbb R^d$. Its
    lower-level gradient is globally Lipschitz in both arguments, uniformly over
    the other argument: for all \(\theta,\theta'\in\Theta\) and
    \(x,x'\in\mathbb R^d\),
    \begin{align*}
        \|\partial_xg(\theta,x)-\partial_xg(\theta,x')\|
        \le L_{g,xx}\|x-x'\|,\\
        \|\partial_xg(\theta,x)-\partial_xg(\theta',x)\|
        \le L_{g,\theta x}\|\theta-\theta'\|.
    \end{align*}
    We write \(L_{g,2}:=\max\{L_{g,xx},L_{g,\theta x}\}\) when one common bound is sufficient. For all $x\in\mathbb{R}^{d}$, $\|\partial_{\theta}g(\theta,x)\|\le C_{g}\|x\|^{\mathsf{n}_{g}}+D_{g}$ where $C_{g},D_{g}\ge 0$ are constants and $\mathsf{n}_{g}\ge 1$ is an integer.

    \item The Hessian $\partial_x^2 g(\theta, x)$ is continuous w.r.t. $(\theta,x)$.
\item \label{assum_uniform_gibbs_regularity}
    There exist \(r_3>0\) and \(L_{g,3}<\infty\), independent of
    \(\theta\), such that
    \begin{equation}\label{eq_uniform_normal_hessian_lipschitz}
        \|\partial_x^2g(\theta,x)-\partial_x^2g(\theta,y)\|_{\mathrm{op}}
        \le L_{g,3}\|x-y\|
    \end{equation}
    for every \(\theta\in\Theta\) and \(x,y\in\mathbb R^d\) satisfying
    \(\max\{\operatorname{dist}(x,\mathcal S(\theta)),
    \operatorname{dist}(y,\mathcal S(\theta))\}\le r_3\).
\item \emph{Uniform tail gradient error
    bound.} There exist constants
    $\mathsf D_{\mathsf{eb}}>0$ and $\nu_{\mathsf{eb}}>0$, independent of
    $\theta$, such that, for every $\theta\in\Theta$ and every
    $x\in\mathbb R^d$ with $\|x\|\ge \mathsf D_{\mathsf{eb}}$,
    \begin{equation}\label{eq_uniform_tail_gradient_error_bound}
        \|\partial_x g(\theta,x)\|
        \ge \nu_{\mathsf{eb}}
        \operatorname{dist}\bigl(x,\mathcal S(\theta)\bigr).
    \end{equation}
\end{itemize}
\end{assumption}
\end{tcolorbox}
The implication from an error bound to quadratic growth for a function with
Lipschitz-continuous gradient is standard; see
\citep[Theorem~2 and Appendix~A]{karimi2016linear}.  Its argument
applies on the tail because \(\partial_xg\) is globally \(L_{g,xx}\)-Lipschitz and \eqref{eq_uniform_tail_gradient_error_bound}
holds there.  Define the derived constants $\mathsf D_{\mathsf{qg}}:=\mathsf D_{\mathsf{eb}},$ and $\mu_{\mathsf{qg}}:=\nu_{\mathsf{eb}}^2/L_{g,xx}$.
Consequently, for every \(\theta\in\Theta\) and every \(x\in\mathbb R^d\)
with \(\|x\|\ge\mathsf D_{\mathsf{qg}}\),
\begin{equation}\label{eq_tail_EB_implies_QG}
    g(\theta,x)-g^\star(\theta)
    \ge \frac{\nu_{\mathsf{eb}}^2}{2L_{g,xx}}
        \operatorname{dist}^2\bigl(x,\mathcal S(\theta)\bigr)
    =\frac{\mu_{\mathsf{qg}}}{2}
        \operatorname{dist}^2\bigl(x,\mathcal S(\theta)\bigr).
\end{equation}
This is the \emph{derived uniform tail quadratic-growth estimate}.
Since \(\mathcal S(\theta)\) is nonempty and, as a subset of the local-minimum
set, is contained in the compact set supplied by the \(\PLcirc\) condition,
it is bounded.  Hence
\(\operatorname{dist}(x,\mathcal S(\theta))\to\infty\) as
\(\|x\|\to\infty\), and \eqref{eq_tail_EB_implies_QG} shows that
\(g(\theta,\cdot)\) is coercive for every \(\theta\in\Theta\).  Therefore,
\Cref{prop_S_path_connected,prop_S_C_2_wihout} apply under the combined
\Cref{mu_PL_assum,assum_sample_comp}.  This conclusion is pointwise in
\(\theta\); uniform boundedness of
\(\{\mathcal S(\theta)\}_{\theta\in\Theta}\) is established later in
\Cref{prop_parametric_morse_bott_stability}.

The Lipschitz and polynomial-growth conditions in
\Cref{assum_sample_comp} are standard in the literature
\citep{kwon2023penalty,chen2023bilevel}.  Unlike
\cite{kwon2023penalty}, however,
we allow an unbounded lower-level domain; the tail estimate
\eqref{eq_tail_EB_implies_QG} ensures that \(\gibbs_g^\lambda\) is well
defined \citep{hasenpflug2024wasserstein}.
Moreover, whereas \cite{kwon2023penalty} assumes boundedness of \(f\), we use
a polynomial bound.  We also assume a bounded \(\theta\)-domain and a
polynomial bound on \(\partial_\theta g\), which we plan to relax in future
work.

\begin{example}\label{ex:sphere_family}
Fix \(a>0\) and \(d\ge2\), let \(\Theta=[0,1]\), and define
\(s(x):=\sqrt{\|x\|^2+a^2}\),
\(R(\theta):=\sqrt{\theta+1+a^2}\), and
\(g(\theta,x):=(s(x)-R(\theta))^2\).
Its minimizer set is
\(\mathcal S(\theta)
=\{x\in\mathbb R^d:\|x\|=\sqrt{\theta+1}\}\), which is a compact,
connected, non-singleton sphere with radius in
\([1,\sqrt2]\). On the common neighborhood
\(\mathcal N:=\{x:\|x\|>1/2\}\),
\[
    \|\partial_xg(\theta,x)\|^2
    =4g(\theta,x)\frac{\|x\|^2}{\|x\|^2+a^2}
    \ge\frac{4}{1+4a^2}g(\theta,x).
\]
Thus the local P\L\ constant may be chosen uniformly as
\(\mu=(1+4a^2)^{-1}\). The only critical point outside the solution sphere is
\(x=0\), where
\(\partial_{xx}^2g(\theta,0)
=2(1-R(\theta)/a)I\prec0\).
The remaining smoothness, growth, and geometric conditions are verified in
Appendix~\ref{append:verification_sphere_example}. Hence this family satisfies all the
lower-level assumptions above. Taking \(f(\theta,x):=x_1\) also satisfies the
upper-level assumptions and gives
\(\max_{x\in\mathcal S(\theta)}f(\theta,x)=\sqrt{\theta+1}\).
\end{example}

Under Assumption~\ref{assum_sample_comp}, the mappings
$(\theta,x)\mapsto f(\theta,x)$ and $(\theta,x)\mapsto \partial_\theta g(\theta,x)$
are continuous on the compact set $\Theta\times\mathbb{B}_d(0;R)$ with $R>0$.
We define the constants
\begin{align}\label{eq_unif_const_Lip_smooth}
    L_{f,0}(R)
\;:=\;
\sup_{\theta\in\Theta}\ \sup_{x\in\mathbb{B}_d(0;R)}
|f(\theta,x)|,
\quad
L_{g,1}(R)
\;:=\;
\sup_{\theta\in\Theta}\ \sup_{x\in\mathbb{B}_d(0;R)}
\|\partial_\theta g(\theta,x)\|.
\end{align}
By compactness and continuity, these suprema are finite.

\subsection{Implications of the \texorpdfstring{\PLcirc\ }{TEXT}Condition}\label{sec:preliminary}
We review implications of \Cref{mu_PL_assum}, in particular how it helps to quantify the convergence of the Gibbs measures \(\gibbs_{\theta}^\lambda\) to the limiting measure \(\gibbs_{\theta}^0\) under the Wasserstein-1 metric.
In this subsection, $\theta$ is considered fixed, and hence for ease of the notation, we omit it and simply write \(\gibbs^\lambda\), and \(\gibbs^0\).

\subsubsection{Implications of the Embedding Submanifold Structure in \Cref{prop_S_C_2_wihout}} \label{section_implication_PL_circ}

Fix \(\theta\) and write \(\optimalsetlower:=\optimalsetlower(\theta)\). A compact
embedded manifold need not admit one global coordinate chart or one global orthonormal
frame of its normal bundle. We therefore use the coordinate-free normal-bundle
formulation and introduce charts only for local calculations.
\paragraph{Normal bundle and its sections.}
For \(x\in\optimalsetlower\), let \(T_x\optimalsetlower\) be the tangent space and
let \(N_x\optimalsetlower:=(T_x\optimalsetlower)^\perp\) be the normal space. Define
the normal bundle
\[
    N\optimalsetlower
    :=\{(x,v):x\in\optimalsetlower,\ v\in N_x\optimalsetlower\}.
\]
Let
\[
    \pi_N:N\optimalsetlower\to\optimalsetlower,
    \qquad
    \pi_N(x,v):=x,
\]
be the bundle projection. A section of \(N\optimalsetlower\) is a map
\(\widehat s:\optimalsetlower\to N\optimalsetlower\) satisfying
\[
    \pi_N\circ\widehat s=\operatorname{id}_{\optimalsetlower}.
\]
Thus, for every \(x\in\optimalsetlower\), there is a unique vector
\(s(x)\in N_x\optimalsetlower\) such that
\(\widehat s(x)=(x,s(x))\). Conversely, any normal vector field
\(x\mapsto s(x)\in N_x\optimalsetlower\) defines a section by this formula.
We henceforth identify a section with its vector component and write simply
\(s(x)\in N_x\optimalsetlower\).

Let \(U\subseteq\optimalsetlower\) be open. A local \(\mathcal C^1\)
orthonormal normal frame on \(U\) is a collection of \(\mathcal C^1\) maps
\[
    e^{(1)},\ldots,e^{(d-\dimk)}:U\to\mathbb R^d
\]
such that, for every \(x\in U\), the vectors
\(e^{(1)}(x),\ldots,e^{(d-\dimk)}(x)\) form an orthonormal basis of
\(N_x\optimalsetlower\). Every section \(s\) over \(U\) then has the unique
representation
\[
    s(x)=\sum_{j=1}^{d-\dimk}s_j(x)e^{(j)}(x),
    \qquad
    s_j(x)=\langle s(x),e^{(j)}(x)\rangle.
\]
A section \(s\) is continuous, respectively \(\mathcal C^1\), if the map
\(x\mapsto s(x)\in\mathbb R^d\) is continuous, respectively \(\mathcal C^1\),
and \(s(x)\in N_x\optimalsetlower\) for every \(x\). Equivalently, the
coefficient functions \(s_1,\ldots,s_{d-\dimk}\) are continuous, respectively
\(\mathcal C^1\), in every local orthonormal normal frame. Since this definition
uses only local frames, it does not require the normal bundle to admit a global
frame.
The \emph{zero section} is the map
\[
    \widehat 0:\optimalsetlower\to N\optimalsetlower,
    \qquad
    \widehat 0(x):=(x,0),
\]
whose image is
\(\widehat 0(\optimalsetlower)
  =\{(x,0):x\in\optimalsetlower\}\).

\paragraph{Tubular neighborhood.}
The following statement provides the normal coordinates used below.  Its proof
at the stated finite regularity is given in
Appendix~\ref{append:finite_regularity_tubular_neighborhood}; the
\(\mathcal C^\infty\) version is the tubular-neighborhood
theorem~\citep[Theorem~6.24]{lee2003smooth}.

Since \(\optimalsetlower\) is a compact \(\mathcal C^2\) embedded manifold
without boundary, there exists \(\delta_{\mathrm{tub}}>0\) such that, with
\[
    \mathcal V
    :=\bigl\{(u,v)\in N\optimalsetlower:
    \|v\|<\delta_{\mathrm{tub}}\bigr\},
\]
the addition map
\begin{equation}\label{local coordinate decomposition}
    \Phi:\mathcal V\to\Phi(\mathcal V)\subseteq\mathbb R^d,
    \qquad
    \Phi(u,v):=u+v,
\end{equation}
is a \(\mathcal C^1\) diffeomorphism, and \(\Phi(\mathcal V)\) is an open
neighborhood of \(\optimalsetlower\).  Consequently, every
\(y\in\Phi(\mathcal V)\) has a unique representation
\[
    y=u+v,
    \qquad
    u\in\optimalsetlower,
    \quad
    v\in N_u\optimalsetlower,
    \quad
    \|v\|<\delta_{\mathrm{tub}},
\]
and \(u\) and \(v\) depend \(\mathcal C^1\)-smoothly on \(y\).

For local calculations, choose an open subset \(U\subseteq\optimalsetlower\) and
a \(\mathcal C^2\) local parametrization
\[
    \mathsf S:\Gamma\subseteq\mathbb R^\dimk\to U,
\]
where \(\Gamma\) is open. By the local triviality of the normal bundle, after
shrinking \(U\) if necessary, there is a local \(\mathcal C^1\) orthonormal
normal frame \(e^{(1)},\ldots,e^{(d-\dimk)}\) on \(U\). Explicitly, local
triviality means that
\[
    \mathcal T_U:
    U\times\mathbb R^{d-\dimk}
    \to N\optimalsetlower|_U,
    \qquad
    \mathcal T_U(x,t)
    :=
    \left(x,\sum_{j=1}^{d-\dimk}t_je^{(j)}(x)\right),
\]
is a \(\mathcal C^1\) diffeomorphism, where $N\optimalsetlower|_U
    :=\{(x,v)\in N\optimalsetlower:x\in U\}.$ In the resulting manifold and normal coordinates, the tubular map is represented
by

\begin{equation}\label{eq_local_normal_coordinates}
    y(u,t)
    =
    \mathsf S(u)
    +
    \sum_{j=1}^{d-\dimk}t_je^{(j)}(\mathsf S(u)).
\end{equation}
Compactness provides a finite collection of such coordinate neighborhoods, and a
partition of unity combines the corresponding local calculations.

\paragraph{Volume measure.}
The Euclidean inner product induces a Riemannian metric on \(\optimalsetlower\).
If \(\mathsf S:\Gamma\subset\R^\dimk\to\optimalsetlower\) is a local
parametrization, define
\[
    G_{\mathsf S}(u)
    :=
    \bigl[\langle\partial_i\mathsf S(u),\partial_j\mathsf S(u)\rangle\bigr]_{i,j=1}^{\dimk}.
\]
The Riemannian volume measure \(\mathcal M\) is defined locally by
\[
    \mathrm d\mathcal M(\mathsf S(u))
    =\sqrt{\det G_{\mathsf S}(u)}\,\mathrm du.
\]
These local expressions agree on chart overlaps and define \(\mathcal M\) without an
orientability assumption. We suppress the dependence on \(\theta\) when \(\theta\) is
fixed.

\paragraph{Inverse Poincar\'e constant.}
A probability measure \(\mu\) on \(\mathbb R^d\) satisfies a Poincar\'e
inequality with inverse Poincar\'e constant \(C_{\mathrm{PI}}\) if
\[
    \operatorname{Var}_{\mu}(\varphi)
    \le C_{\mathrm{PI}}
    \int_{\mathbb R^d}\|\nabla\varphi\|^2\,\mathrm d\mu
    \qquad\text{for all }\varphi\in\mathbb H^1(\mu),
\]
where \(\mathbb H^1(\mu)\) is the Sobolev space weighted by \(\mu\).
For the Gibbs measure \(\mu_\theta^\lambda\), we denote
its inverse Poincar\'e constant by \(C_{\mathrm{PI}}(\theta,\lambda)\), and
write \(\rho_{\theta,\lambda}:=C_{\mathrm{PI}}(\theta,\lambda)^{-1}\) for
the corresponding spectral-gap constant.
Gong et al.~\citep{gong2024poincare} use the reciprocal
convention: their spectral-gap constant \(\rho_{\theta,\lambda}\) multiplies the
variance on the left-hand side.
Fix \(\theta\) and suppose that the complete lower-level
hypotheses of this paper hold, in the active non-singleton regime
\(\dimk\ge1\).  The short hypothesis comparison at the beginning of
Appendix~\ref{append_proof_uniform_poincare}, followed by the intrinsic tube
comparison proved there, adapts the fixed-potential argument of
Gong et al.~\citep[Theorem~4 in Section~5]{gong2024poincare} and yields constants
\(\lambda_{\mathrm{PI}}(\theta)>0\) and
\(C_{\mathrm{PI}}(\theta)<\infty\) such that
\[
    C_{\mathrm{PI}}(\theta,\lambda)
    =\rho_{\theta,\lambda}^{-1}
    \le C_{\mathrm{PI}}(\theta)
    \qquad
    \bigl(0<\lambda\le\lambda_{\mathrm{PI}}(\theta)\bigr).
\]
Thus the bound is independent of \(\lambda\) on a sufficiently small
fixed-parameter interval, but both the bound and the endpoint may depend on
the fixed \(\theta\).  This result alone gives no
uniform-in-\(\theta\) conclusion.
The required family-level adaptation, including one common temperature
threshold and one common inverse Poincar\'e bound, is proved separately in
\Cref{prop_uniform_poincare}.

\subsubsection{Limiting Gibbs measure}

We identify the weak limit of the Gibbs measure as \(\lambda\downarrow0\).
The mass concentrates near the minimizer manifold, and its distribution along
the manifold is determined by the curvature of \(g\) in the normal directions.
The tail
quadratic-growth estimate \eqref{eq_tail_EB_implies_QG} ensures
that the Gibbs measure is well-defined and that mass away from the manifold is
negligible.
Throughout this subsection, \(\theta\) is fixed.  We write
\(g:=g(\theta,\cdot)\),
\(\optimalsetlower:=\optimalsetlower(\theta)\), and
\(Z_\lambda:=\int_{\mathbb R^d}
e^{-g(x)/\lambda}\,\mathrm dx\).
Using the local normal coordinates \(y(u,t)\) in
\eqref{eq_local_normal_coordinates}, define, for the manifold point
\(a=\mathsf S(u)\),
\begin{equation}\label{eq_limiting_gibbs_weight}
    w(a)
    :=
    \det\!\left(
        \left.\partial_{tt}^2\bigl[g(y(u,t))\bigr]\right|_{t=0}
    \right)^{-1/2}.
\end{equation}
The local P\L--Morse--Bott equivalence
\citep[Corollary~2.17]{rebjock2024fast} ensures that the matrix inside the
determinant is positive definite.  A change of local orthonormal normal frame
conjugates this matrix by an orthogonal matrix, so its determinant, and hence
\(w(a)\), is independent of the selected frame.

\begin{prop}[Laplace limit on a minimizer manifold]
\label{lemma_density_limiting_measure}
Let \Cref{mu_PL_assum} and the
\(\mathcal C^2\)-extension, global lower-gradient Lipschitz,
Hessian-continuity, and uniform tail gradient error-bound clauses of
\Cref{assum_sample_comp} hold.  Then
\(Z_\lambda<\infty\) for every \(\lambda>0\), and there exists a
probability measure \(\gibbs^0\) such that
\(\gibbs^\lambda\Rightarrow\gibbs^0\) as \(\lambda\downarrow0\).
The limiting measure is supported on \(\optimalsetlower\) and,
with \(w\) defined in \eqref{eq_limiting_gibbs_weight}, satisfies
\begin{equation}\label{eqn_density_limiting_measure}
    \gibbs^0(\mathrm da)
    =
    \rho(a)\,\mathrm d\mathcal M(a),
    \qquad
    \rho(a)
    :=
    \frac{w(a)}
    {\displaystyle
     \int_{\optimalsetlower}
     w(b)\,\mathrm d\mathcal M(b)}.
\end{equation}
In particular,
\(\operatorname{supp}(\gibbs^0)=\optimalsetlower\).
\end{prop}
This is the standard manifold Laplace limit.  Its proof in
Appendix~\ref{append:proof_limiting_gibbs_measure} adapts the local-fiber
argument of Hasenpflug et al.~\citep{hasenpflug2024wasserstein}; finitely many
local orthonormal normal frames replace their stronger global-frame
formulation.  The quantitative Wasserstein rate is established separately in
\Cref{prop_W1_Gibbs_pointwise}.

We continue to suppress \(\theta\) in \(\gibbs^0_\theta\), \(\rho_\theta\),
and \(\mathcal M_\theta\) whenever \(\theta\) is fixed.

\paragraph{Wasserstein-1 convergence for fixed \(\theta\).}
The weak convergence in \Cref{lemma_density_limiting_measure} identifies the
limiting Gibbs measure but gives no convergence rate.  For the fixed parameter
\(\theta\), the next proposition strengthens this conclusion to an
\(\mathcal O(\sqrt\lambda)\) rate in Wasserstein-1 distance.
The proof follows the \(p=1\) normal-projection coupling
of Hasenpflug et al.~\citep[Theorem~3.6]{hasenpflug2024wasserstein}.  Our
assumptions provide the same normal Taylor, tubular-Jacobian, Gaussian-tail,
and first-moment estimates used in that proof; finitely many local
orthonormal normal frames replace its global-frame formulation.  The
finite-regularity adaptation is verified in
\Cref{append_valid_assmp_thm_W_1,append_proof_pointwise_W1_Gibbs}.

\begin{prop}[Fixed-\(\theta\) \(\mathbb W_1\) rate]
{\color{teal}\label{prop_W1_Gibbs_pointwise}}
Suppose that
\Cref{mu_PL_assum,assum_secon_fund_form_bound} hold.  Suppose also that the
\(\mathcal C^2\)-extension, global lower-gradient Lipschitz,
Hessian-continuity, and uniform tail gradient error bound in
\Cref{assum_sample_comp} hold.  Then there exist
constants \(\lambda_0>0\) and \(K_0<\infty\) such that
\[
    \mathbb W_1\!\left(
        \gibbs^\lambda,\gibbs^0
    \right)
        \le K_0\lambda^{1/2}
    \qquad
    \text{for all }0<\lambda\le\lambda_0.
\]
\end{prop}

The parameter \(\theta\) is suppressed in \(K_0\) and \(\lambda_0\), as it is
fixed throughout this subsection; these constants need not yet be uniform over
\(\Theta\).  The uniform version is established later in
\Cref{prop_W1_Gibbs_uniform}.
The proof is given in
Appendix~\ref{append_proof_pointwise_W1_Gibbs}.

 \section{Regularity of the Solution Mapping and Geometric Properties}\label{sec:4_unif_props}
The purpose of this section is to turn the
fixed-parameter geometry developed above into estimates that hold uniformly
over \(\theta\in\Theta\).  We first establish parametric
stability of the minimizer manifolds and a common sublevel neighborhood
on which the local P{\L} inequality and the corresponding error bound hold
uniformly.  These results imply Hausdorff--Lipschitz continuity of
\(\theta\mapsto\mathcal S(\theta)\) and, consequently, Lipschitz continuity of
\(F_{\max}\).  We then show that the minimizer manifolds have a common
dimension and derive uniform bounds on their geometry, Riemannian volume, and
limiting Gibbs density.  Finally, under the local Hessian-Lipschitz
clause of \Cref{assum_sample_comp}, we obtain an
\(\mathcal O(\sqrt\lambda)\) Wasserstein--1 bound with constants independent of
\(\theta\).
We conclude the section by proving a uniform low-temperature
Poincar\'e inequality, which supplies the common inverse Poincar\'e constant
and temperature threshold used by the later LMC analysis.
\subsection{Lipschitz Continuity of the Solution Mapping}\label{sec:Lip_cont_sol_mapping}
In this subsection, our goal (Lemma~\ref{lemm_PL_yield_S_theta_lip})
is to show that, under the uniform parametric assumptions in
\Cref{mu_PL_assum,assum_sample_comp}, the set-valued mapping
$\theta \mapsto \mathcal{S}(\theta)$ is globally Lipschitz on $\Theta$ with
respect to the Hausdorff distance.
We then combine this regularity with the Lipschitz continuity of $f(\theta,x)$ to deduce the Lipschitz 
continuity of $F_{\max}$ (Corollary~\ref{lemm_hyp_blo_cont}). 
We begin by recalling the notion of Hausdorff distance.

\begin{mydef}[Hausdorff distance]\label{def_Hays_dist}

Let $\mathcal{S}_1$ and $\mathcal{S}_2$ be two nonempty compact subsets of
$\mathbb{R}^{d}$.  Their Hausdorff distance is
\[
\text{dist}(\mathcal{S}_1,\mathcal{S}_2)
:=
\max\left\{
\sup_{x_{1}\in\mathcal{S}_1}\inf_{x_{2}\in\mathcal{S}_2}\|x_1-x_{2}\|,
\;
\sup_{x_{2}\in\mathcal{S}_2}\inf_{x_{1}\in\mathcal{S}_1}\|x_1-x_{2}\|
\right\}.
\]

\end{mydef}
The local P{\L} inequality is imposed separately for each fixed $\theta$ and
therefore does not, by itself, control how $\mathcal S(\theta)$ varies with
$\theta$.  Under the full assumptions, every critical point of $g(\theta,\cdot)$ is either a global
minimizer or has a negative-definite Hessian and is therefore a strict local
maximum (see \Cref{lem_critical_point_dichotomy}).  
In the following proposition, we fix $\bar\theta\in\Theta$ and describe
$\mathcal S(\theta)$ for $\theta$ sufficiently close to $\bar\theta$.  We show
that $\mathcal S(\theta)$ is obtained by moving each point of
$\mathcal S(\bar\theta)$ in a direction normal to
$\mathcal S(\bar\theta)$ at that point.  The proof first constructs the nearby
manifold by requiring the normal component of the gradient of
$g(\theta,\cdot)$ to vanish.  Positive curvature in the normal directions,
together with the no-saddle condition, is then used to show that every point
on this manifold is a global minimizer.  Finally, connectedness shows that the
constructed manifold is the entire minimizer set $\mathcal S(\theta)$.

\begin{prop}[Parametric Morse--Bott stability]
\label{prop_parametric_morse_bott_stability}
Let \Cref{mu_PL_assum,assum_sample_comp} hold. Fix
$\bar\theta\in\Theta$ and write
$\bar{\mathcal S}:=\mathcal S(\bar\theta)$.  There is a relatively open
neighborhood $V_{\bar\theta}$ of $\bar\theta$ in $\Theta$ and, for each
$\theta\in V_{\bar\theta}$, a $\mathcal C^1$ section $n_\theta$ of the normal
bundle $N\bar{\mathcal S}$ such that $n_{\bar\theta}=0$ and
\[
    \mathcal S(\theta)
    =\bigl\{u+n_\theta(u):u\in\bar{\mathcal S}\bigr\}.
    \]
    The section is taken in a common tubular neighborhood of $\bar{\mathcal S}$,
    and $u\mapsto u+n_\theta(u)$ is an embedding.
    Consequently, there is
$\mathsf D_{\mathcal S}>0$ such that
\[
    \mathcal S(\theta)\subseteq
    \mathbb B_d(0;\mathsf D_{\mathcal S})
    \qquad\forall\theta\in\Theta,
\]
and the solution graph
\begin{equation}\label{eq_solution_graph}
    \mathfrak G
    :=\bigl\{(\theta,x)\in\Theta\times\mathbb R^d:
                   x\in\mathcal S(\theta)\bigr\}
\end{equation}
is compact.
\end{prop}
The proof is given in
Appendix~\ref{append_proof_parametric_morse_bott_stability}.
For the later tail estimates, fix once and for all a working radius
\begin{equation}\label{eq_working_radius_D}
    \mathsf D\ge
    \max\bigl\{\mathsf D_{\mathsf{qg}},\;2\mathsf D_{\mathcal S}\bigr\}.
\end{equation}
Here \(\mathsf D_{\mathsf{qg}}\) is the derived radius defined after
\eqref{eq_tail_EB_implies_QG}.
Then every $\mathcal S(\theta)$ is contained in
$\mathbb B_d(0;\mathsf D)$ and its diameter is at most $\mathsf D$.
Although the P{\L} constant $\mu$ is common to all parameters, the neighborhoods
in the local P{\L} condition may depend on $\theta$.  The main difficulty in the
next lemma is therefore to obtain one neighborhood of the entire solution graph
\(\mathfrak G\)
and one positive sublevel height on which the P{\L} and error-bound estimates
hold for every $\theta\in\Theta$.  The proof combines the finite normal-graph
cover furnished by
\Cref{prop_parametric_morse_bott_stability} with the uniform tail
gradient error bound in \Cref{assum_sample_comp} and its quadratic-growth
consequence \eqref{eq_tail_EB_implies_QG}.
It works on coordinate neighborhoods with local orthonormal normal frames and
does not require a global chart or a global normal frame.

\begin{lem}[Uniform local P{\L} and error bound]
\label{lem_uniform_local_PL_band}
\leavevmode\par\noindent
Under \Cref{mu_PL_assum,assum_sample_comp}, let
$g^\star(\theta):=\min_xg(\theta,x)$ and set $\bar\mu:=\mu^2/(4L_{g,xx})$.
There exist a relatively open neighborhood $\mathcal W$ of
$\mathfrak G$ defined in \eqref{eq_solution_graph} in
$\Theta\times\mathbb R^d$ and a constant $\varepsilon_{\mathrm{PL}}>0$ such that
for every $\theta\in\Theta$ and $x\in\mathbb R^d$,
\begin{equation}\label{eq_uniform_PL_band}
g(\theta,x)-g^\star(\theta)\le\varepsilon_{\mathrm{PL}}
\quad\Longrightarrow\quad
(\theta,x)\in\mathcal W
\end{equation}
and, on this sublevel set,
\begin{align}
g(\theta,x)-g^\star(\theta)
&\le \frac{1}{2\bar\mu}\|\partial_xg(\theta,x)\|^2,
\label{eq_uniform_local_PL}\\
\operatorname{dist}\bigl(x,\mathcal S(\theta)\bigr)
&\le \frac{2}{\mu}\|\partial_xg(\theta,x)\|.
\label{eq_uniform_error_bound}
\end{align}
\end{lem}
The proof of \Cref{lem_uniform_local_PL_band} is given in
Appendix~\ref{append_proof_uniform_local_PL_band}.
The uniform error bound above \eqref{eq_uniform_error_bound}, together with the Lipschitz dependence of
$\partial_xg$ on $\theta$, controls the distance from a minimizer at one
parameter value to the minimizer set at a nearby parameter value.  Applying
this local estimate along line segments in the convex set $\Theta$ yields the
following global Hausdorff bound.
\begin{lem}\label{lemm_PL_yield_S_theta_lip}
Under \Cref{mu_PL_assum,assum_sample_comp}, the minimizer mapping is globally
Lipschitz on the compact convex set $\Theta$ in Hausdorff distance:
\begin{equation}\label{eq_global_Hausdorff_Lipschitz}
\operatorname{dist}\bigl(\mathcal S(\theta_1),\mathcal S(\theta_2)\bigr)
\le
\frac{2L_{g,\theta x}}{\mu}\|\theta_1-\theta_2\|
\qquad\forall\theta_1,\theta_2\in\Theta.
\end{equation}
\end{lem}

The proof of \Cref{lemm_PL_yield_S_theta_lip} is given in
Appendix~\ref{append_proof_Hausdorff_and_hyperobjective_Lipschitz}.
Combining \eqref{eq_global_Hausdorff_Lipschitz} with the Lipschitz continuity of
$f$ gives the following consequence.

\begin{corr}[Lipschitz continuity of the hyper-objective]\label{lemm_hyp_blo_cont}
Under \Cref{mu_PL_assum,assum_sample_comp}, $F_{\max}$ is
\[
    L_{F_{\max},1}
    :=L_{f,1}\left(1+\frac{2L_{g,\theta x}}{\mu}\right)
\]
Lipschitz continuous on $\Theta$.
\end{corr}
The proof is given in
Appendix~\ref{append_proof_Hausdorff_and_hyperobjective_Lipschitz}.

\subsection{Common Dimension and Geometric Properties of \texorpdfstring{$\{\mathcal{S}(\theta)\}_{\theta\in\Theta}$}{TEXT}}\label{subsec:3.2}
By \Cref{prop_parametric_morse_bott_stability}, for every $\bar\theta$ and all
$\theta$ sufficiently close to $\bar\theta$, the manifold $\mathcal S(\theta)$ is
the image of $\mathcal S(\bar\theta)$ under a normal-graph embedding.  Its
dimension is therefore locally constant in $\theta$.  Since $\Theta$ is convex and
hence connected, this local constancy gives a common dimension over all of
$\Theta$.

\begin{lem} \label{lemma_common_dimension}
    Under \Cref{mu_PL_assum,assum_sample_comp}, the manifolds $\{\optimalsetlower(\theta)\}_{\theta\in\Theta}$ admit a common dimension $\dimk$.
\end{lem}
The proof is given in \Cref{append_proof_lemma_common_dimension}. Lemma~\ref{lemma_common_dimension} allows us to work with a single geometric parameter $\dimk$ in all later bounds, rather than letting the dimension depend on $\theta$. In particular, when we study the Gibbs relaxation and the superquantile approximation of $F_{\max}$, the scaling with accuracy is expressed in terms of this common dimension $\dimk$ (see \Cref{rmk_func_approx_cvar_gibbs}).
Throughout the paper, we consider $\dimk\geq 1$, since when $\dimk=0$, $\optimalsetlower(\theta)$ is a singleton and minima-selection is redundant.
Next, we collect basic geometric properties of $\{\mathcal{S}(\theta)\}_{\theta\in\Theta}$ that will be used in our approximation analysis (see \Cref{sec:3}). 
Their proofs combine parametric normal-graph stability and compactness
with the uniform second-fundamental-form bound and the resulting common
ambient tube.
\begin{prop}\label{assum_compact_opt_set}
    Under compactness of \(\Theta\) and \Cref{mu_PL_assum,assum_secon_fund_form_bound,assum_sample_comp}, the following holds for all \(\theta \in \Theta\).
    \begin{itemize}[leftmargin=*]
        \setlength\itemsep{0em}
    \item The working radius $\mathsf D$ fixed in \eqref{eq_working_radius_D}
    satisfies
    $\operatorname{diam}_{\mathbb R^d}(\mathcal S(\theta))\le\mathsf D$.
\item The sectional curvature\footnote{Sectional curvature measures the
    intrinsic curvature associated with a two-dimensional tangent plane; see
    \Cref{def_sec_curv}.} is uniformly bounded: there is a constant
    $\mathsf C_0>0$ such that
    \[
        \bigl|\operatorname{Sec}_{\mathcal S(\theta)}(p,\Pi)\bigr|
        \le \mathsf C_0
    \]
    for every $p\in\mathcal S(\theta)$ and every two-plane
    $\Pi\subset T_p\mathcal S(\theta)$.  This assertion is vacuous when
    $\dimk=1$.
    \item The intrinsic injectivity radius\footnote{The injectivity radius at
    \(p\) is the largest radius on which the Riemannian exponential map at
    \(p\) is a diffeomorphism onto its image; see \Cref{def_injec_rad}.} is
    uniformly bounded below: there is a constant \(r_0>0\) such that
    \(\operatorname{inj}_{\mathcal S(\theta)}(p)\ge r_0\) for every
    \(p\in\mathcal S(\theta)\).
    \end{itemize}
\end{prop}
The proof is given in \Cref{append_proof_assum_compact_opt_set}.
Finally, we show a uniform bound on the limiting density $\rho_\theta$ defined in \cref{eqn_density_limiting_measure}.
\begin{prop}\label{prop_strict_bound_limit_gibbs_density}
Under \Cref{mu_PL_assum,assum_sample_comp}, 
    there exists a constant $\kappa>0$ such that $\kappa \leq \rho_\theta(x) \leq 1 / \kappa$ for all $x\in\optimalsetlower(\theta)$ and $\theta \in \Theta$.
\end{prop}
The proof is given in \Cref{append_proof_prop_strict_bound_limit_gibbs_density}. This result ensures that integration with respect to the limiting Gibbs measure is uniformly comparable to integration with respect to the volume measure on $\mathcal{S}(\theta)$, which is crucial for our approximation analysis in the next section.
\begin{rmk}(On the size of the constant $\kappa$)
We choose
\(
\kappa=\mathsf b(\mu/L_{g,xx})^{(d-\dimk)/2}
\)
in the proof of \Cref{prop_strict_bound_limit_gibbs_density}.  Here \(\mu\)
is the uniform lower bound on the nonzero normal eigenvalues of
\(\partial_x^2g(\theta,x)\), while global \(L_{g,xx}\)-Lipschitz continuity
of the lower gradient gives the upper Hessian bound
\(\|\partial_x^2g(\theta,x)\|_{\mathrm{op}}\le L_{g,xx}\).  The constant
\(\mathsf b>0\) comes from the uniform volume bounds in
\Cref{lem_vol_sol_set_is_bounded}.  Thus a poorly conditioned ratio
\(\mu/L_{g,xx}\) can make \(\kappa\) small; in the idealized isotropic case
\(\mu=L_{g,xx}\), this choice reduces to \(\kappa=\mathsf b\).
\end{rmk}
The pointwise estimate in \Cref{prop_W1_Gibbs_pointwise} allows
its constants to depend on the fixed parameter.  The following proposition
shows that the stability and uniform geometric estimates established above,
together with the uniform local Hessian-Lipschitz condition in
\Cref{assum_sample_comp}, yield one Wasserstein--1 rate that holds
simultaneously for every \(\theta\in\Theta\).

\begin{prop}[Uniform Wasserstein-1 rate]
\label{prop_W1_Gibbs_uniform}\label{lemm_W1_Gibbs_unif}
Suppose that
\Cref{mu_PL_assum,assum_sample_comp,assum_secon_fund_form_bound} hold.  Then
there exist constants
\(K_0<\infty\) and
\(\lambda_0>0\), independent of \(\theta\), such that
\begin{equation}\label{eq_uniform_W1_Gibbs_rate}
    \mathbb W_1\bigl(\gibbs_\theta^\lambda,\gibbs_\theta^0\bigr)
    \le K_0\sqrt\lambda
    \qquad
    (\theta\in\Theta,\ 0<\lambda\le\lambda_0).
\end{equation}
\end{prop}
The proof is given in Appendix~\ref{append_proof_uniform_W1_Gibbs}.
Gong et al.~\citep{gong2024poincare} develop the underlying
Lyapunov and shrinking-tube strategy for one fixed potential.  The following
result is the family-level adaptation:
its appendix proof tracks the analytic and geometric constants uniformly over
the compact parameter family and works intrinsically on the normal bundle,
using only local normal frames.

\begin{prop}[Uniform low-temperature Poincar\'e inequality]
\label{prop_uniform_poincare}
Under Assumptions~\ref{mu_PL_assum}--\ref{assum_sample_comp}, if
\(\dimk\ge1\), then there exist constants
\[
    \lambda_{\mathrm{PI}}^\star>0,
    \qquad
    C_{\mathrm{PI}}^\star<\infty,
\]
independent of \(\theta\), such that, for every \(\theta\in\Theta\), every
\(0<\lambda\le\lambda_{\mathrm{PI}}^\star\), and every
\(\varphi\in\mathbb H^1(\mu_\theta^\lambda)\),
\[
    \operatorname{Var}_{\mu_\theta^\lambda}(\varphi)
    \le
    C_{\mathrm{PI}}^\star
    \int_{\mathbb R^d}\|\nabla\varphi(x)\|^2\,
    \mu_\theta^\lambda(\mathrm dx).
\]
Equivalently, if \(\rho_{\theta,\lambda}\) denotes the spectral-gap
constant, then
\begin{equation}\label{eq_uniform_poincare_gap}
    \rho_{\theta,\lambda}
    \ge
    \rho_{\mathrm{PI}}^\star
    :=(C_{\mathrm{PI}}^\star)^{-1}>0
    \qquad
    \bigl(\theta\in\Theta,\ 0<\lambda\le
    \lambda_{\mathrm{PI}}^\star\bigr).
\end{equation}
\end{prop}
The proof is given in Appendix~\ref{append_proof_uniform_poincare}.

\section{Approximation Guarantees with \SQG\ Relaxation}\label{sec:3}
In this section, we show that the superquantile-Gibbs relaxation
\(\tilde{F}_{\SQG}\), defined in \Cref{eqn:SQ_gibbs_variation}, provides a
uniform approximation to \(F_{\max}\) under
\Cref{mu_PL_assum,assum_sample_comp,assum_secon_fund_form_bound}, for suitable
\(\delta\) and \(\lambda\).  The
uniform Wasserstein estimate in \Cref{prop_W1_Gibbs_uniform} supplies the Gibbs
part of this approximation.  The result justifies using \(\tilde F_{\SQG}\) as
a surrogate for \(F_{\max}\) in the zeroth-order method developed later.  We
establish it through the following steps:
    \begin{enumerate}
\setlength{\itemsep}{0.5pt}
        \item The superquantile loss (\SQR-loss), defined as
\begin{align}\label{eq_def_sqr_loss}
            \phi_{0,\delta}(\theta,\beta)
            &:=\beta+\frac{1}{\delta}
            \mathbb{E}_{X\sim\gibbs^0_\theta}
            [\max\{f(\theta,X)-\beta,0\}],\\
            \tilde F_{\SQR}(\theta)
            &:=\min_{\beta\in\mathbb R}\phi_{0,\delta}(\theta,\beta).
        \end{align}
         closely approximates the hyper-objective \( F_{\max}(\theta) \) point-wisely.  
        \item The difference between \(\tilde F_{\SQR}\) and \(\tilde F_{\SQG}\) is bounded by the Wasserstein-1 distance between \(\gibbs^\lambda_\theta\) and \(\gibbs^{0}_\theta\):  
\(
|\tilde F_{\SQR}(\theta)-\tilde F_{\SQG}(\theta)|\le \mathcal{O}(\delta^{-1} \mathbb{W}_{1}(\gibbs^\lambda_\theta,\gibbs^{0}_\theta)).
\)
    \end{enumerate}  
Below we formally present these results.
\subsection{Approximation of \texorpdfstring{$\hyperobjective$}{TEXT} with \texorpdfstring{$\tilde F_{\SQR}$}{TEXT}}

Recall that \(\widetilde F_{\SQR}(\theta)\) is the upper
\((1-\delta)\)-superquantile of \(f(\theta,X)\) for
\(X\sim\mu_\theta^0\).  Informally, it averages the largest \(\delta\)-fraction
of the values of \(f(\theta,X)\); the variational definition in
\eqref{eq_def_sqr_loss} also treats possible probability mass at the quantile
threshold correctly.  Our goal is to show that this upper-tail average is
close to \(F_{\max}(\theta)\) when \(\delta\) is small.

To see the geometric idea, fix \(\theta\), let \(q\) be the
\((1-\delta)\)-quantile of \(f(\theta,X)\), and choose any maximizer
\(x^\star\in\argmax_{x\in\mathcal S(\theta)}f(\theta,x)\).  The generalized
quantile satisfies
\[
 \mu_\theta^0\{x\in\mathcal S(\theta):f(\theta,x)>q\}\le\delta,
\]
and the variational formula gives
\[
 0\le F_{\max}(\theta)-\widetilde F_{\SQR}(\theta)
 \le F_{\max}(\theta)-q.
\]
If the gap \(F_{\max}(\theta)-q\) were large, the Lipschitz continuity of
\(f(\theta,\cdot)\) would place an intrinsic geodesic ball, with radius
proportional to this gap and centered at \(x^\star\), inside the strict upper
set \(\{f(\theta,\cdot)>q\}\).  Its
probability would therefore be at most \(\delta\).  We convert this probability
bound into a bound on the ball's radius in two steps.  First, the uniform lower
bound on the density of \(\mu_\theta^0\) converts probability into Riemannian
volume.  Second, the small-ball estimate below converts volume into radius:
because the volume of a radius-\(r\) ball is bounded below by a constant
multiple of \(r^{\dimk}\), the radius is at most a constant multiple of
\(\delta^{1/\dimk}\).  The proof of \Cref{th_CVaR_unif_approx} makes this
argument precise without assuming a unique maximizer or an atomless
distribution.

For completeness, let \(\mathsf g_\theta\) denote the Riemannian metric on
\(\mathcal S(\theta)\) induced by the Euclidean inner product, let
\(\operatorname{dist}_{\mathsf g_\theta}\) be its geodesic distance, and let
\(\mathcal M_\theta\) be its Riemannian volume measure.  For
\(p\in\mathcal S(\theta)\), the open intrinsic geodesic ball of radius \(r>0\)
is
\[
    \mathbb B_{\mathsf g_\theta}(p;r)
    :=\bigl\{u\in\mathcal S(\theta):
       \operatorname{dist}_{\mathsf g_\theta}(p,u)<r\bigr\}.
\]
For each fixed \(\theta\) and \(p\), normal coordinates give the familiar
small-radius limit
\[
    \lim_{r\to0^+}
    \frac{\mathcal M_\theta(\mathbb B_{\mathsf g_\theta}(p;r))}{r^{\dimk}}
    =\omega_{\dimk},
\]
where \(\omega_{\dimk}\) is the volume of the Euclidean unit ball in
\(\mathbb R^{\dimk}\).  This pointwise asymptotic is not enough here: the
argument above needs one radius range and two-sided constants that work for
every \(\theta\) and every center \(p\in\mathcal S(\theta)\).  The next lemma
provides exactly this uniform, non-asymptotic comparison.

\begin{lem}\label{lemm_comparison_volume_geod_ball_euc_ball}
Under
\Cref{mu_PL_assum,assum_secon_fund_form_bound,assum_sample_comp}, there are
constants \(\mathsf c_L,\mathsf c_H,c_1>0\), independent of \(\theta\), such
that, for every \(\theta\in\Theta\), \(p\in\mathcal S(\theta)\), and
\[
    0<r<\min\left\{r_0,\frac{c_1}{\sqrt{2\mathsf C_0}}\right\},
\]
the geodesic ball satisfies
\[
    \mathsf c_L\,\omega_{\dimk}r^{\dimk}
    \le
    \mathcal M_\theta\bigl(\mathbb B_{\mathsf g_\theta}(p;r)\bigr)
    \le
    \mathsf c_H\,\omega_{\dimk}r^{\dimk},
\]
where \(\omega_{\dimk}:=\Leb_{\dimk}(\mathbb B_{\dimk}(0;1))\), and
\(r_0,\mathsf C_0\) are given by \Cref{assum_compact_opt_set}.
\end{lem}
The proof is given in
\Cref{append_proof_lemm_VaR_approx}.  
We now state the following key result.
\begin{thm}[$\tilde F_{\SQR}$ approximates $F_{\max}$]\label{th_CVaR_unif_approx}
Recall \(\tilde F_{\SQR}(\theta) := \min_\beta \phi_{0, \delta}(\theta, \beta)\).
Under \Cref{mu_PL_assum,assum_secon_fund_form_bound,assum_sample_comp}, for
every \(\theta\in\Theta\) and every \(\delta\in(0,1)\) satisfying
\[
\delta\le \mathsf{c}_{L}\cdot C_{\dimk}\cdot \kappa\cdot\left(4^{-1}\min\{c_{1}/\sqrt{2\mathsf{C}_{0}},r_{0}\}\right)^{\dimk},
\]
we have
\begin{equation*}
\setlength{\abovedisplayskip}{5pt}
            \setlength{\belowdisplayskip}{5pt}
  |F_{\max}(\theta)-\tilde{F}_{\SQR}(\theta)|\le \frac{4L_{f,1}\delta^{\frac{1}{\dimk}}}{(\kappa\cdot \mathsf{c}_{L}\cdot C_{\dimk})^{\frac{1}{\dimk}}},
\end{equation*}
where $\mathsf{c}_{L}$ and  $c_{1}$ are defined in \Cref{lemm_comparison_volume_geod_ball_euc_ball}, $\kappa$ in \Cref{prop_strict_bound_limit_gibbs_density}, $\dimk$ is the dimension of $\mathcal{S}(\theta)$, $C_{\dimk}=\pi^{\dimk/2}/\Gamma(\dimk/2+1)$, and $r_{0},\mathsf{C}_{0}$ are from \Cref{assum_compact_opt_set}.
\end{thm}

\begin{proof}\label{proof_th_CVaR_unif_approx}
Fix \(\theta\in\Theta\), let \(X\sim\mu_\theta^0\), and set
\[
Z:=f(\theta,X),\qquad
q:=\VaR_{1-\delta}(Z).
\]
The generalized-quantile definition gives
\[
\mathbb P(Z>q)\le\delta\le\mathbb P(Z\ge q).
\]
Indeed, the first inequality follows from right-continuity of
\(t\mapsto\mathbb P(Z>t)\); for every \(t<q\), the definition gives
\(\mathbb P(Z>t)>\delta\), and continuity from above as \(t\uparrow q\)
gives the second inequality.

For every \(\beta\in\mathbb R\), if \(\beta>q\), then
\(\phi_{0,\delta}(\theta,\beta)\ge\beta>q\), whereas if \(\beta\le q\),
then
\[
\phi_{0,\delta}(\theta,\beta)
\ge
\beta+\frac{q-\beta}{\delta}\,\mathbb P(Z\ge q)
\ge q.
\]
Moreover, \(Z\le F_{\max}(\theta)\) almost surely, so
\(\phi_{0,\delta}(\theta,F_{\max}(\theta))=F_{\max}(\theta)\).  Since \(Z\) is bounded and
\(0<\delta<1\), the continuous variational objective is coercive in
\(\beta\), and hence its minimum is attained.  Consequently,
\begin{equation}\label{eq_22}
0
\le F_{\max}(\theta)-\widetilde F_{\SQR}(\theta)
\le F_{\max}(\theta)-q.
\end{equation}

Set
\[
r_\delta
:=
\left(
\frac{\delta}{\kappa\mathsf c_L C_{\dimk}}
\right)^{1/\dimk}.
\]
The assumed bound on \(\delta\) gives
\[
r_\delta
\le\frac14\min\left\{r_0,\frac{c_1}{\sqrt{2\mathsf C_0}}\right\},
\qquad
2r_\delta
<\min\left\{r_0,\frac{c_1}{\sqrt{2\mathsf C_0}}\right\}.
\]
If \(L_{f,1}=0\), then \(f(\theta,\cdot)\) is constant on
\(\mathcal S(\theta)\), so \(q=F_{\max}(\theta)\) and the result follows from
\eqref{eq_22}.  Suppose therefore that \(L_{f,1}>0\), and assume for a
contradiction that
\[
F_{\max}(\theta)-q>4L_{f,1}r_\delta.
\]
Choose any
\(x^\star\in\argmax_{x\in\mathcal S(\theta)}f(\theta,x)\).  For every
\(x\in\mathbb B_{\mathsf g_\theta}(x^\star;2r_\delta)\),
\begin{align*}
\|x-x^\star\|
&\le\operatorname{dist}_{\mathsf g_\theta}(x,x^\star)
\le2r_\delta,\\
f(\theta,x)
&\ge F_{\max}(\theta)-L_{f,1}\|x-x^\star\|\\
&\ge F_{\max}(\theta)-2L_{f,1}r_\delta
>q.
\end{align*}
Thus this geodesic ball is contained in the strict upper set
\(\{x\in\mathcal S(\theta):f(\theta,x)>q\}\).  Consequently,
\begin{align*}
\delta
&\ge \mu_\theta^0\{f(\theta,\cdot)>q\}\\
&\ge \mu_\theta^0\!\left(
\mathbb B_{\mathsf g_\theta}(x^\star;2r_\delta)
\right)\\
&\ge \kappa\,
\mathcal M_\theta\!\left(
\mathbb B_{\mathsf g_\theta}(x^\star;2r_\delta)
\right)\\
&\ge \kappa\mathsf c_L C_{\dimk}(2r_\delta)^{\dimk}
=2^{\dimk}\delta,
\end{align*}
The first inequality is
$\mathbb P(Z>q)\le\delta$, because $X\sim\mu_\theta^0$ and
$Z=f(\theta,X)$.  The second follows from the preceding Lipschitz estimate,
which proves
$\mathbb B_{\mathsf g_\theta}(x^\star;2r_\delta)
\subseteq\{x\in\mathcal S(\theta):f(\theta,x)>q\}$.  For the third, the density of
$\mu_\theta^0$ with respect to $\mathcal M_\theta$ is at least $\kappa$ by
\Cref{prop_strict_bound_limit_gibbs_density}.  The fourth is the lower
geodesic-ball volume bound in
\Cref{lemm_comparison_volume_geod_ball_euc_ball}.  Finally, the equality uses
$r_\delta^{\dimk}=\delta/(\kappa\mathsf c_L C_{\dimk})$.
Since \(\dimk\ge1\) and \(\delta>0\), the resulting
inequality \(\delta\ge2^{\dimk}\delta\) is impossible.
Hence
\[
F_{\max}(\theta)-q
\le
4L_{f,1}
\left(
\frac{\delta}{\kappa\mathsf c_L C_{\dimk}}
\right)^{1/\dimk}.
\]
Combining this estimate with \eqref{eq_22} proves the result.
\end{proof}

\begin{rmk}(Sharpness of the polynomial exponent \(1/\dimk\))
The polynomial exponent \(1/\dimk\) in \Cref{th_CVaR_unif_approx} is dictated by the $\dimk$-dimensional volume growth of small balls on $\mathcal{S}(\theta)$ and cannot, in general, be increased. Indeed, fix
$\theta$ and a point $x^{\ast}(\theta)\in \mathcal{S}(\theta)$, and consider
on \(\mathcal S(\theta)\) the toy function
$f(\theta,x):=-\operatorname{dist}_{\mathsf g_\theta}
(x,x^{\ast}(\theta))$. Clearly, $F_{\max}(\theta)=0$. For $r>0$, the set
\(A_r(\theta):=\{x\in\mathcal S(\theta):-f(\theta,x)<r\}\)
is exactly \(\mathbb B_{\mathsf g_\theta}(x^\ast(\theta);r)\).
Let \(r=-\VaR_{1-\delta}(f(\theta,Y))\), where
\(Y\sim\mu_\theta^0\) and \(\delta\) is sufficiently small.
By \Cref{prop_strict_bound_limit_gibbs_density} and \Cref{lemm_comparison_volume_geod_ball_euc_ball}, we have
\begin{align}
    \delta
    &=\mu_\theta^0(A_r(\theta))
      =\mu_\theta^0\bigl(
        \mathbb B_{\mathsf g_\theta}(x^\ast(\theta);r)\bigr)\nonumber\\
    &\le\kappa^{-1}\,
      \mathcal M_\theta\bigl(
        \mathbb B_{\mathsf g_\theta}(x^\ast(\theta);r)\bigr)
      \le\kappa^{-1}\mathsf c_H C_{\dimk}r^{\dimk}.
\end{align}
Then $\delta^{1/\dimk}\lesssim r=F_{\max}(\theta)-\VaR_{1-\delta}(f(\theta,Y))$. Thus one cannot, in general, replace the exponent
$1/\dimk$ in \Cref{th_CVaR_unif_approx} by a larger polynomial exponent.
Although this conical toy function is not
\(\mathcal C^1\) at \(x^\ast(\theta)\), Appendix~
\ref{appendix_sharpness_sqr_rate} gives a \(\mathcal C^1\), globally
Lipschitz modification and proves that no larger polynomial exponent of
\(\delta\) is valid under the current assumptions.
\end{rmk}

\subsection{Approximation of \text{\SQR-Loss} with \text{\SQG-Loss}}
For each fixed \(\beta\), the variational objectives
\(\phi_{0,\delta}\) and \(\phi_{\lambda,\delta}\) differ only through the
underlying measures.  The map
\(x\mapsto(f(\theta,x)-\beta)_+\) is Lipschitz with a constant independent of
\(\beta\), so the coupling characterization of Wasserstein--1 distance gives
a uniform perturbation bound that passes directly to the two minimum values.
\begin{lem}\label{Theo_func_approx_cvar}
Under \Cref{mu_PL_assum,assum_sample_comp}, for every
\(\theta\in\Theta\), \(\lambda>0\), and \(\delta\in(0,1)\), the measures
\(\mu_\theta^\lambda\) and \(\mu_\theta^0\) have finite first moments, the
minima defining \(\tilde F_{\SQG}(\theta)\) and
\(\tilde F_{\SQR}(\theta)\) are attained, and
\begin{equation}\label{eq_cvar_upper_bound_w_1}
   \bigl|\tilde{F}_{\SQG}(\theta)-\tilde{F}_{\SQR}(\theta)\bigr|
   \le \frac{L_{f,1}}{\delta}\,
   \mathbb{W}_{1}\!\left(\mu_{\theta}^{\lambda},\mu_{\theta}^{0}\right).
\end{equation}
\end{lem}
The proof of \Cref{Theo_func_approx_cvar} is given in \Cref{append_proof_Theo_func_approx_cvar}.
As a consequence of \Cref{th_CVaR_unif_approx},
\Cref{Theo_func_approx_cvar}, and \Cref{prop_W1_Gibbs_uniform},
$\tilde{F}_{\SQG}(\theta)$ approximates $F_{\max}(\theta)$ uniformly when
$\delta$ and $\lambda$ are chosen appropriately.
\begin{tcolorbox}[colback=gray!20, colframe=gray!50, boxrule=0pt, arc=0mm, left=0mm, right=0mm, top=0mm, bottom=0mm]
\begin{thm}\label{rmk_func_approx_cvar_gibbs}
Under
\Cref{mu_PL_assum,assum_secon_fund_form_bound,assum_sample_comp}, let
\(K_0\) and \(\lambda_0\) be the constants from
\Cref{prop_W1_Gibbs_uniform}; in particular, they are independent of
\(\theta\).  Then, for every \(\theta\in\Theta\), every
\(\delta\in(0,1)\) satisfying
\[
    \delta<
    \mathsf c_L C_{\dimk}\kappa
    \left(4^{-1}\min\{c_1/\sqrt{2\mathsf C_0},r_0\}\right)^{\dimk},
\]
and every \(0<\lambda\le\lambda_0\),
\begin{equation}\label{eq_func_approx_cvar_gibbs}
    \setlength{\abovedisplayskip}{5pt}
            \setlength{\belowdisplayskip}{5pt}
     \left|F_{\max}(\theta)-\tilde{F}_{\SQG}(\theta)\right|\le \frac{4L_{f,1}\delta^{\frac{1}{\dimk}}}{(\kappa\mathsf c_L C_{\dimk})^{\frac{1}{\dimk}}}+\frac{L_{f,1}}{\delta}\cdot K_{0}\lambda^{\frac{1}{2}}.
\end{equation}
\end{thm}
\end{tcolorbox}
The bound in \eqref{eq_func_approx_cvar_gibbs} makes the role of the manifold dimension $\dimk$ explicit: the bias term scales as $\delta^{1/\dimk}$, so larger $\dimk$ requires smaller $\delta$ to reach the same accuracy. This bound also yields a simple rule for choosing the parameters to reach any target accuracy: pick $\delta=\epsilon_{v}^{\dimk}$ and $\lambda=\epsilon_{v}^{2(\dimk+1)}$, where
\[
0<\epsilon_v<
\min\left\{1,
(\mathsf c_L C_{\dimk}\kappa)^{1/\dimk}
\,4^{-1}\min\{c_1/\sqrt{2\mathsf C_0},r_0\}\right\},
\qquad
\epsilon_v\le\lambda_0^{1/(2(\dimk+1))},
\]
    then
    \[
    \left|F_{\max}(\theta)-\tilde{F}_{\SQG}(\theta)\right|=\mathcal{O}(\epsilon_{v}).
    \]
This provides a clear, dimension-aware recipe for setting $(\delta,\lambda)$ that we use in \Cref{sec:6}.

 \section{Oracle Complexity of Optimizing $F_{\max}$}\label{sec:6}
In this section, we develop a zeroth-order algorithm to find a Goldstein stationary point of the 
hyper-objective $F_{\max}$. Zeroth-order methods require pointwise evaluations of the objective, 
but direct evaluation of $F_{\max}$ is intractable because it involves an inner maximization over 
the lower-level solution set. To obtain a tractable objective oracle, we replace $F_{\max}$ by its 
Gibbs–superquantile surrogate $\tilde{F}_{\SQG}$ defined in \eqref{eqn:SQ_gibbs_variation}, which 
can be estimated from Gibbs samples. As shown in \Cref{rmk_func_approx_cvar_gibbs}, suitable choices of 
$\lambda$ and $\delta$ ensure that $\tilde{F}_{\SQG}$ uniformly approximates $F_{\max}$. We then 
query $\tilde{F}_{\SQG}$ at multiple points to construct a finite-difference estimator of a proxy 
gradient of $F_{\max}$, and use this estimator within a zeroth-order scheme to obtain a Goldstein 
stationary point of the original hyper-objective.
For the zeroth-order results in this section, fix \(\bar\rho>0\) and set
\[
    \Theta^{+\bar\rho}
    :=\{\vartheta\in\mathbb R^m:
          \operatorname{dist}(\vartheta,\Theta)\le\bar\rho\}
     =\Theta+\bar\rho\,\overline{\mathbb B_m(0;1)}.
\]
This set is compact and convex.  We impose the explicit algorithmic
neighborhood condition that
\Cref{mu_PL_assum,assum_secon_fund_form_bound,assum_sample_comp} hold, with
the same type of uniform constants, when the parameter set is replaced by
\(\Theta^{+\bar\rho}\).  This condition is not inferred from their validity
on \(\Theta\).  All uniform constants used by the algorithms below are chosen
by applying the preceding results to \(\Theta^{+\bar\rho}\); the optimization
domain and the projected iterates remain in \(\Theta\).

To evaluate $\tilde{F}_{\SQG}(\theta)$, we use its variational representation in \eqref{eqn:SQ_gibbs_variation}, which involves solving the inner minimization problem $\min_{\beta \in \mathbb{R}} \phi_{\lambda,\delta}(\theta, \beta)$. 
Since $\phi_{\lambda,\delta}(\theta,\cdot)$ is convex but possibly nonsmooth, a natural choice is to apply the subgradient method\footnote{Because the problem is one-dimensional, the bisection method is also applicable, provided the domain is bounded.}.
To ensure convergence, the nonsmooth optimization must be restricted to a \emph{bounded} interval. 
In \Cref{lem_3}, we show that the optimal set $\argmin_{\beta}\phi_{\lambda,\delta}(\theta,\beta)$ lies within $[-\mathrm{B},\mathrm{B}]$, so that
\[
\tilde{F}_{\SQG}(\theta) = \min_{\beta \in [-\mathrm{B}, \mathrm{B}]} \phi_{\lambda,\delta}(\theta, \beta).
\]

Given that (i) $\phi_{\lambda,\delta}(\theta,\cdot)$ is convex and
\(\sigma\)-Lipschitz, where \(\sigma:=1+\delta^{-1}\), (ii) the feasible set
is \([ -\mathrm B,\mathrm B]\), and (iii) the sample subgradient is
conditionally valid and bounded by \(\sigma\), we use constant-step PSGD
with
\(\eta=2\mathrm B/(\sigma\sqrt L)\) when \(\mathrm B>0\).
The self-contained projection--martingale argument in
\Cref{append_proofs_lem_5_2}, which is consistent with
Lan et al.~\citep[Proposition~2]{lan2012validation}, gives,
for every \(\hat\epsilon>0\) and
\(0<\Delta\le1/2\),
\begin{align}\label{eq_whp_gaurante_PSGD}
\mathbb P\!\left(
    \phi_{\lambda,\delta}(\theta,\hat\beta)
    -\min_{\beta}\phi_{\lambda,\delta}(\theta,\beta)>\hat\epsilon
    \,\middle|\,\theta
\right)\le\Delta
\end{align}
provided
\[
L\ge
64\,\mathrm B^2\sigma^2\hat\epsilon^{-2}
\log^2(2/\Delta).
\]
If \(\mathrm B=0\), the feasible interval is the singleton \(\{0\}\) and
the PSGD optimization error is zero.

Building on the above discussion, we now use the PSGD output to construct a tractable zeroth-order oracle for $F_{\max}$. For any query $\theta$, we evaluate the scalar quantity $\phi_{\lambda,\delta}(\theta,\hat{\beta})$, where $\hat{\beta}$ is the approximate minimizer returned by PSGD, and regard this as a noisy evaluation of $\widetilde{F}_{\SQG}(\theta)$ and hence of $F_{\max}(\theta)$. Querying this oracle at carefully chosen perturbations of $\theta$ yields a multi-point finite-difference estimator of a proxy gradient, which drives the projected zeroth-order scheme in \Cref{alg: Stochastic-Zeroth-order}.

Below, we analyze in detail the convergence of \Cref{alg:PSGD} (inner minimization of the Gibbs–superquantile) and \Cref{alg: Stochastic-Zeroth-order} (outer minimization of the hyper-objective).

\subsection{Inner Minimization of Gibbs-superquantile}
In this subsection, we first show that all global minimizers of \(\phi_{\lambda,\delta}(\theta,\cdot)\) lie in a common bounded interval. This boundedness allows us to run PSGD on a compact domain and invoke standard convergence guarantees. We then use this result to justify treating the noisy  values \(\phi_{\lambda,\delta}(\theta,\hat{\beta})\), where \(\hat{\beta}\) is the  PSGD output, as approximate evaluations of \(F_{\max}(\theta)\).

We first briefly introduce the PSGD (\Cref{alg:PSGD}) to solve $\min_{\beta\in[-\mathrm{B},\mathrm{B}]}\phi_{\lambda,\delta}(\theta,\beta)$. 
Specifically, for a given \(\theta\in\Theta^{+\bar\rho}\), we have
\begin{align}
    \beta^{l+1}\leftarrow\text{Proj}_{[-\mathrm{B},\mathrm{B}]}(\beta^{l}-\eta\,\tilde\partial_{\beta}\phi_{\lambda,\delta}(\theta,\beta^{l};\tilde{X}^{(l)})),
\end{align}
where \( \tilde\partial_{\beta} \phi_{\lambda,\delta}(\theta,\beta;\tilde{X}^{(l)}) \) represents a sub-gradient of the noisy sample of $\phi_{\lambda,\delta}(\theta,\beta)$ as follows:
\begin{align}\label{eq_def_nois_phi}
    \phi_{\lambda,\delta}(\theta,\beta;\tilde{X}^{(l)}):=\beta+\frac{1}{\delta}\max\{f(\theta,\tilde{X}^{(l)})-\beta,0\},
\end{align}
At equality we select the hinge subgradient zero, so the sample
subgradient used below is
\(1-\delta^{-1}\mathbf 1_{\{f(\theta,\widetilde X^{(l)})>\beta\}}\).
The samples \( \{ \tilde{X}^{(l)} \}_{l=0}^{L-1} \) are independent and identically distributed (i.i.d.) draws from \( \mu_{\theta}^{\lambda} \), and $\mathrm{B}$ is nonnegative. The PSGD output is the pre-update average
\(L^{-1}\sum_{l=0}^{L-1}\beta^{l}\).

We next show in \Cref{lem_3} that the minimizers 
\(\beta^{\star}(\theta)\in\arg\min_{\beta}\phi_{\lambda,\delta}(\theta,\beta)\) are uniformly bounded.
For the upper bound on $\beta^\ast(\theta)$, recall the definition of $\tilde{F}_{\SQR}(\theta)$ in \eqref{eq_def_sqr_loss}; then
\begin{align}
\beta^{\ast}(\theta)
&\le \beta^{\ast}(\theta)+\tfrac{1}{\delta}\,\mathbb{E}_{X\sim \mu_{\theta}^{\lambda}}\!\big[\max\{f(\theta,X)-\beta^{\ast}(\theta),0\}\big]
= \phi_{\lambda,\delta}(\theta,\beta^{\ast}(\theta)) \nonumber\\
&= \tilde{F}_{\SQG}(\theta)
\le \tilde{F}_{\SQR}(\theta) + \big|\tilde{F}_{\SQG}(\theta)-\tilde{F}_{\SQR}(\theta)\big| \nonumber\\
&\le \tilde{F}_{\SQR}(\theta) + \tfrac{L_{f,1}}{\delta}\,K_{0}\,\lambda^{1/2} \label{eq_202020}\\
&\le L_{f,0}(\mathsf{D}) + \tfrac{L_{f,1}}{\delta}\,K_{0}\,\lambda^{1/2}, \label{eq_212121}
\end{align}
where \eqref{eq_202020} applies
\Cref{prop_W1_Gibbs_uniform} to \Cref{Theo_func_approx_cvar} for
\(0<\lambda\le\lambda_0\), and \eqref{eq_212121} uses
\[
\tilde{F}_{\SQR}(\theta)\le \max_{x\in\mathcal{S}(\theta)}f(\theta,x)\le \max_{x\in \mathbb{B}_{d}(0;\mathsf{D})} f(\theta,x)\le L_{f,0}(\mathsf{D}),
\]
with $L_{f,0}(\cdot)$ defined in \eqref{eq_unif_const_Lip_smooth}
and $\mathsf D$ fixed so that
$\mathcal S(\theta)\subseteq\mathbb B_d(0;\mathsf D)$ by
\Cref{prop_parametric_morse_bott_stability} and
\eqref{eq_working_radius_D}.
Throughout this section, \(K_0\) and \(\lambda_0\) denote
the common constants furnished by \Cref{prop_W1_Gibbs_uniform} on the compact
parameter set under consideration.  In the zeroth-order results this set is
the thickening \(\Theta^{+\bar\rho}\); in particular, neither constant depends
on the query parameter.

A key part of the proof is a uniform lower bound on $\beta^\ast(\theta)$: since $\beta^{*}(\theta)$ is a $(1-\delta)$-quantile, we have 
$\mathbb{P}\{f(\theta,X)\le \beta^{\ast}(\theta)\}\ge 1-\delta$ for $X\sim \mu^\lambda_\theta$.
To turn this probability statement into a lower bound on
\(\beta^\ast(\theta)\), it suffices to show that the corresponding quantile
event contains a point in one fixed ball.  For this purpose,
\Cref{lem_gibbs_mass_outside} bounds the probability that a Gibbs sample lies
outside that ball: there are constants
\(\mathsf W_0<\infty\), \(c_{\mathrm{tail}}>0\), and
\(\lambda_{\mathrm{tail},0}>0\), independent of \(\theta\), such that
\[
\mu_\theta^\lambda\!\left(
    \mathbb R^d\setminus\mathbb B_d(0;2\mathsf D)
\right)
\le
\mathsf W_0\lambda^{1-d/2}
\exp\!\left(-\frac{c_{\mathrm{tail}}}{\lambda}\right)
\]
for every \(\theta\in\Theta\) and
\(0<\lambda\le\lambda_{\mathrm{tail},0}\).  The exponential decay dominates
the polynomial factor \(\lambda^{1-d/2}\), so the right-hand side converges to
zero as \(\lambda\to0^+\).  We may therefore choose
\(0<\lambda_{\mathrm{mass}}\le\lambda_{\mathrm{tail},0}\) so that
\[
    \mu_\theta^\lambda\!\left(
        \mathbb B_d(0;2\mathsf D)
    \right)
    \ge\frac34
\]
for every \(\theta\in\Theta\) and
\(0<\lambda\le\lambda_{\mathrm{mass}}\).

Combining this estimate with \(\delta\le1/2\) yields
\[
\mu^\lambda_\theta\!\left(
    \{x:f(\theta,x)\le\beta^\ast(\theta)\}
    \cap\mathbb B_d(0;2\mathsf D)
\right)
\ge\frac14.
\]
Hence there exists \(x_0\in\mathbb B_d(0;2\mathsf D)\) such that
\(f(\theta,x_0)\le\beta^\ast(\theta)\), and consequently
\[
    \beta^\ast(\theta)\ge-L_{f,0}(2\mathsf D).
\]
The upper estimate \eqref{eq_212121} additionally requires
\(\lambda\le\lambda_0\).  We therefore set
\[
    \lambda_\beta
    :=\min\{\lambda_0,\lambda_{\mathrm{mass}}\}.
\]
For every \(0<\lambda\le\lambda_\beta\), both the upper and lower bounds are
available, giving
\(\beta^\ast(\theta)\in[-\mathrm B,\mathrm B]\).
We now state the result; its complete proof is deferred to
\Cref{append_proof_lemm_lb_beta_min_rsq}.
\begin{lem}\label{lem_3}
Suppose that
    \Cref{mu_PL_assum,assum_secon_fund_form_bound,assum_sample_comp} hold.
    Then, for every \(\theta\in\Theta\), every \(0<\delta\le1/2\), every
    \(0<\lambda\le\lambda_\beta\), and every $\beta^\star\in\argmin_{\beta\in\mathbb R}
        \phi_{\lambda,\delta}(\theta,\beta),$ one has
    \begin{align}\label{eq_bound_beta_minimizer}
       |\beta^\star|
       &\le \mathrm B(\delta,\lambda):=\max\left\{
       L_{f,0}(2\mathsf D),
       L_{f,0}(\mathsf D)
       +\frac{L_{f,1}K_0}{\delta}\sqrt\lambda
       \right\}.
    \end{align}
    The quantity \(\mathrm B(\delta,\lambda)\) is independent of \(\theta\),
    although it depends on the displayed values of \(\delta\) and \(\lambda\).
\end{lem}
Combining the existing convergence analysis of PSGD with \Cref{lem_3}, 
we can apply the guarantee in \eqref{eq_whp_gaurante_PSGD} to the PSGD output $\hat{\beta}$.
 We then define the following function as an accessible approximation of $\phi_{\lambda,\delta}(\theta,\hat{\beta})$:
\begin{align}\label{eq_def_tilde_psi}
    \tilde{\psi}(\theta)
    :=\phi_{\lambda,\delta}(\theta,\hat{\beta};\hat{X}_{I})
    =\hat{\beta}+\frac{1}{\delta|I|}\sum_{i\in I}
      \max\{f(\theta,\hat{X}^{(i)})-\hat{\beta},0\}.
\end{align}
Here
$\hat{X}_{I}=\{\hat{X}^{(i)}\}_{i\in I}$ consists of $|I|$ i.i.d. samples
from $\gibbs^{\lambda}_{\theta}$ that are independent of every
sample used to construct $\hat\beta$.

\begin{algorithm}[t]
\caption{$\mathrm{PSGD}(\theta,\beta^{0};(\lambda,\delta))$}
\textbf{Input: $\theta$, $\beta^{0}\in[-\mathrm B,\mathrm B]$,
$\lambda>0$, $\delta$, and $L\ge1$}
    \begin{algorithmic}[1]
        \IF{$\mathrm B=0$}
        \STATE $\hat\beta=0$ and terminate
        \ELSE
        \STATE $\sigma=1+\delta^{-1}$ and
        $\eta=2\mathrm B/(\sigma\sqrt L)$
        \FOR{$l=0$ to $L-1$}
        \STATE draw $\tilde X^{(l)}\sim\gibbs^\lambda_\theta$
        independently of the preceding samples and iterates
        \STATE $\beta^{l+1}=\operatorname{Proj}_{[-\mathrm{B},\mathrm{B}]}
        (\beta^{l}-\eta\,\tilde \partial_{\beta}
        \phivar(\theta,\beta^{l};\tilde{X}^{(l)}))$
        \ENDFOR
        \STATE $\hat\beta=L^{-1}\sum_{l=0}^{L-1}\beta^{l}$
        \ENDIF
\end{algorithmic}\label{alg:PSGD}
\end{algorithm}
The next lemma shows that the quantity \( \tilde{\psi}(\theta) \) defined in \eqref{eq_def_tilde_psi} closely approximates \(\tilde{F}_{\SQG}(\theta)\). The proof is deferred to \Cref{append_proofs_lem_5_2}. Before stating it, we briefly explain the role of the high-probability guarantee in~\eqref{eq_whp_gaurante_PSGD}. The convergence analysis of the outer zeroth-order algorithm is often based on an $L^2$-type control of the oracle error, i.e., $\mathbb{E}[|\tilde{\psi}(\theta)-F_{\max}(\theta)|^{2}\mid\theta]$, rather than just on its expectation.
An in-expectation convergence result for PSGD would only control the mean optimality gap and does not, in general, translate into a sufficiently tight second-moment bound. In contrast, the high-probability estimate in~\eqref{eq_whp_gaurante_PSGD} provides a tail bound on the PSGD optimization error; by choosing the target accuracy $\hat{\epsilon}$ and the failure probability $\Delta$ appropriately, this tail bound can be converted into the required $L^2$ control on the oracle error used in Lemma~\ref{lem_5_2}.
\begin{lem}\label{lem_5_2}

Suppose that
\Cref{mu_PL_assum,assum_secon_fund_form_bound,assum_sample_comp} hold and
\(\dimk\ge1\).  Let
\(0<\epsilon_v\le2^{-1/\dimk}\),
\(\delta:=\epsilon_v^{\dimk}\),
\(\lambda_{\mathrm{mom}}:=
\min\{\lambda_\beta,\lambda_{\mathrm{tail},2\mathsf n_f}\}\), and
\(0<\lambda\le\lambda_{\mathrm{mom}}\).  Then
\(0<\delta\le1/2\).  Fix \(\theta\in\Theta\), let
\(\mathrm B=\mathrm B(\delta,\lambda)\) be given by \Cref{lem_3}, and take
\(\beta^0\in[-\mathrm B,\mathrm B]\).  Set
\(\sigma:=1+\delta^{-1}\) and
\(\Delta:=\min\{1/2,
\epsilon_v^2/(32\sigma^2\mathrm B^2+\epsilon_v^2)\}\).
If \(\mathrm B>0\), run \Cref{alg:PSGD} with
\(L\) independent training samples and with
\[
L\ge512\,\mathrm B^2\sigma^2\epsilon_v^{-2}\log^2(2/\Delta).
\]
If \(\mathrm B=0\), use the singleton convention in \Cref{alg:PSGD}.
Let \(I\) be a nonempty finite index set and let the
validation batch \(\widehat X_I\) consist of i.i.d. samples from
\(\mu_\theta^\lambda\), independent of the complete PSGD trajectory.  If
\[
|I|\ge4\widetilde C_f\epsilon_v^{-2-2\dimk},
\]
then
\[
\mathbb E\!\left[
|\widetilde\psi(\theta)-\widetilde F_{\SQG}(\theta)|^2
\,\middle|\,\theta\right]\le\epsilon_v^2.
\]
Here \(\widetilde C_f\) is independent of
\(\theta,\delta,\lambda,\epsilon_v\) on this range.
Moreover, a successful integer \(L\) can be chosen with
\[
L=\mathcal O\!\left(1+
\dimk^2\mathrm B^2\epsilon_v^{-2-2\dimk}
\bigl(1+\log^2(\epsilon_v^{-1})\bigr)\right).
\]
\end{lem}

We now show that \( \tilde{\psi}(\theta) \) provides a good approximation to
\( F_{\max}(\theta) \). Applying \Cref{lem_5_2} together with
\Cref{rmk_func_approx_cvar_gibbs} and
\(|a+b|^2\le2|a|^2+2|b|^2\) yields
\begin{align}
  \MoveEqLeft[2]\mathbb{E}\!\left[
  |\tilde{\psi}(\theta)-F_{\max}(\theta)|^2\,\middle|\,\theta\right]
  \nonumber\\
  &\le2\mathbb{E}\!\left[
  |\tilde{\psi}(\theta)-\tilde{F}_{\SQG}(\theta)|^2
  \,\middle|\,\theta\right]
  +2|\tilde{F}_{\SQG}(\theta)-F_{\max}(\theta)|^2 \nonumber\\
  &\le2\mathbb{E}\!\left[
  |\phi_{\lambda,\delta}(\theta,\hat{\beta};\hat{X}_{I})
   -\tilde{F}_{\SQG}(\theta)|^2\,\middle|\,\theta\right]\nonumber\\
  &\quad+\mathcal{O}\bigl(\delta^{2/\dimk}+\lambda\delta^{-2}\bigr)
  \nonumber\\
  &= \mathcal{O}(\epsilon_{v}^{2}),\label{eq_append_21}
\end{align}
where, in the final bound, we choose
\(\delta=\epsilon_v^{\dimk}\) and
\(\lambda=\epsilon_v^{2(\dimk+1)}\), with
\begingroup

\[
\begin{gathered}
\epsilon_{\mathrm{geo}}
:=
(\kappa\mathsf c_L C_{\dimk})^{1/\dimk}
\,4^{-1}\min\left\{
\frac{c_1}{\sqrt{2\mathsf C_0}},r_0
\right\},
\\
0<\epsilon_v<\min\{2^{-1/\dimk},1,\epsilon_{\mathrm{geo}}\},
\qquad
\epsilon_v\le
\lambda_{\mathrm{mom}}^{1/(2(\dimk+1))}.
\end{gathered}
\]
The Equation~\eqref{eq_append_21} will serve as a key ingredient in our later convergence analysis.

\subsection{Projected Zeroth-order Algorithm with Minima-Selection}

In this subsection, we analyze the projected stochastic zeroth-order
method with minima-selection (PSZO-MinSel) introduced in
Algorithm~\ref{alg: Stochastic-Zeroth-order}. Our goal is to quantify how many
oracle calls suffice to produce an iterate whose expected generalized-Goldstein
residual for $F_{\max}$ is of order $\epsilon$.
\begin{algorithm}[t]
\caption{Projected Stochastic Zeroth-order with Minima-selection (PSZO-MinSel)}\label{alg: Stochastic-Zeroth-order}
\textbf{Input:} $\beta_{0}\in[-\mathrm B,\mathrm B]$,
$\theta_{0}\in\Theta$, $\lambda>0$, $0<\delta\le1/2$,
$0<\rho<\bar\rho$, an outer stepsize $\eta>0$, and positive integers
$N,b_{u},L,|I|$.
\begin{algorithmic}[1]
        \FOR{$n$ = $0$ to $N-1$}
\STATE Draw fresh directions
         \(u_{n,1},\ldots,u_{n,b_u}\) independently and uniformly from
         \(\mathbb S^{m-1}\), using randomness independent of everything
         generated before iteration \(n\).
         \STATE $\beta_{n,t}^{+}=\mathrm{PSGD}(\theta_{n}+\rho u_{n,t},\beta_{0};(\lambda,\delta)),$ $\forall t\in[1:b_{u}]$.
         \STATE $\beta_{n,t}^{-}=\mathrm{PSGD}(\theta_{n}-\rho u_{n,t},\beta_{0};(\lambda,\delta)),$ $\forall t\in[1:b_{u}]$.
         \STATE Use fresh, mutually independent PSGD
         sample families for every \(t\) and sign in the preceding calls.
\STATE For each \(t\), draw fresh batches
         $\widehat X_{I,n,t}^{+}$ and $\widehat X_{I,n,t}^{-}$, each containing
         $|I|$ i.i.d. samples, from
         $\mu_{\theta_n+\rho u_{n,t}}^\lambda$ and
         $\mu_{\theta_n-\rho u_{n,t}}^\lambda$, respectively. Use distinct sampling randomness for every \(t\) and sign, shared with neither earlier draws nor any PSGD sample family.
        \STATE $\tilde{\psi}(\theta_{n}+\rho u_{n,t})=\phi_{\lambda,\delta}(\theta_{n}+\rho u_{n,t},\beta_{n,t}^{+};\widehat X_{I,n,t}^{+})$, $\forall t\in[1:b_u]$.
        \STATE$\tilde{\psi}(\theta_{n}-\rho u_{n,t})=\phi_{\lambda,\delta}(\theta_{n}-\rho u_{n,t},\beta_{n,t}^{-};\widehat X_{I,n,t}^{-})$, $\forall t\in[1:b_u]$.
        \STATE $\tilde{\mathsf{g}}_{\rho}(\theta_{n})=\frac{m}{2\rho}\cdot\frac{1}{b_{u}}\sum_{t=1}^{b_{u}}\left[\tilde\psi(\theta_{n}+\rho u_{n,t})-\tilde\psi(\theta_{n}-\rho u_{n,t})\right]u_{n,t}$
	        \STATE $\theta_{n+1}=\text{Proj}_{\Theta}(\theta_{n}-\eta\tilde{\mathsf{g}}_{\rho}(\theta_{n}))$
	        \ENDFOR
\STATE Draw
	        $J\sim\operatorname{Unif}\{0,\ldots,N-1\}$ independently of the
	        complete algorithmic history.
	        \STATE \textbf{Return}
	        $\widehat\theta:=\theta_J$.
	    \end{algorithmic}
\end{algorithm}
In \Cref{alg: Stochastic-Zeroth-order}, we use the stochastic
zeroth-order update
\[
\theta_{n+1} \leftarrow \operatorname{Proj}_{\Theta}\!\left(
\theta_{n}-\eta\,\tilde{\mathsf g}_{\rho}(\theta_{n})\right).
\]
where $\tilde{\mathsf g}_{\rho}(\theta_{n})
:= b_{u}^{-1}\sum_{t=1}^{b_{u}}\tilde{\mathsf g}_{\rho}(\theta_{n},u_{n,t}),$
and $\tilde{\mathsf g}_{\rho}(\theta_{n},u_{n,t})
:= \frac{m}{2\rho}\big[\tilde{\psi}(\theta_{n}+\rho u_{n,t})-\tilde{\psi}(\theta_{n}-\rho u_{n,t})\big]\,u_{n,t}$.
Here $\tilde{\psi}$ is defined in \eqref{eq_def_tilde_psi},
$\{u_{n,t}\}_{t=1}^{b_{u}}$ are i.i.d. draws from the uniform distribution on
the unit sphere $\mathbb{S}^{m-1}$, and the outer step is the common constant
$\eta$ specified in \Cref{thm_cvg_alg_zeroth_goldstein}.
The complete procedure is summarized in Algorithm~\ref{alg: Stochastic-Zeroth-order}.
Because $\theta_n\in\Theta$, $\|u_{n,t}\|=1$, and
$0<\rho<\bar\rho$, the whole ball $\mathbb B_m(\theta_n;\rho)$ lies in
$\operatorname{int}(\Theta^{+\bar\rho})$.  In particular, every query
$\theta_n\pm\rho u_{n,t}$ belongs to $\Theta^{+\bar\rho}$, so
\Cref{lem_5_2} applies with the generic compact parameter set replaced by
$\Theta^{+\bar\rho}$.  Each validation batch has the Gibbs law at the same
perturbed parameter and is independent of the training samples that form its
PSGD output.  The resulting two-point quantity is an approximate oracle for
the gradient of the ball-smoothed $F_{\max}$; its error is controlled in mean
square and it is not assumed to be unbiased.  It drives the projected update
in $\Theta$.
To state the result, we recall the ambient Goldstein
$\rho$-subdifferential and the projected generalized-Goldstein criterion.

\begin{mydef}[\cite{goldstein1977optimization}]\label{def_gold_subdiff}
For $\theta\in\mathbb R^m$ and $\rho\ge0$, write
\(
\mathbb B_m(\theta;\rho)
:=\{v\in\mathbb R^m:\|v-\theta\|\le\rho\}.
\)
Suppose that \(F\) is defined and locally Lipschitz on an open neighborhood
of \(\mathbb B_m(\theta;\rho)\), and let \(\partial F(v)\) denote its
ambient Clarke subdifferential there.  The ambient
Goldstein $\rho$-subdifferential is
\[
    \partial_\rho F(\theta)
    :=\operatorname{conv}\!\left(
       \bigcup_{v\in\mathbb B_m(\theta;\rho)}\partial F(v)
    \right).
\]
\end{mydef}
\begin{mydef}[\cite{liu2024zeroth}]\label{def_gen_gold_stat}
At a point \(\theta\in\Theta\) for which
\(\partial_\rho F(\theta)\) in \Cref{def_gold_subdiff} is defined,
$\theta$ is an $(\epsilon,\rho,\eta)$-generalized Goldstein
stationary point (GGSP) of the constrained problem if
\[
    \min_{g\in\partial_\rho F(\theta)}
    \bigl\|\mathcal G_\Theta(\theta,g;\eta)\bigr\|
    \le\epsilon,
\]
where, for $g\in\mathbb R^m$ and $\eta>0$,
\[
    \mathcal G_\Theta(\theta,g;\eta)
    :=\eta^{-1}\!\left[
       \theta-\operatorname{Proj}_\Theta(\theta-\eta g)
    \right].
\]
\end{mydef}
Our main result shows that \Cref{alg: Stochastic-Zeroth-order} produces an
iterate whose expected generalized-Goldstein residual for $F_{\max}$ is of
order $\epsilon$.  The error from replacing $F_{\max}$ by $\tilde\psi$ is
controlled by \eqref{eq_append_21}; the complete proof is deferred to
\Cref{append_proof_thm_cvg_alg_zeroth_goldstein}.

\begin{tcolorbox}[colback=gray!20, colframe=gray!50, boxrule=0pt, arc=0mm, left=0mm, right=0mm, top=0mm, bottom=0mm]
\begin{thm}\label{thm_cvg_alg_zeroth_goldstein}
Let $L_F:=L_{f,1}(1+2L_{g,\theta x}/\mu)$.  Suppose that
    \Cref{mu_PL_assum,assum_secon_fund_form_bound,assum_sample_comp}
    hold on $\Theta^{+\bar\rho}$, with $\dimk\ge1$, $0<\rho<\bar\rho$, and all
    uniform constants chosen on that set. Set
    \(\epsilon_v:=\rho\epsilon/m\),
    \(\delta:=\epsilon_v^{\dimk}\), and
    \(\lambda:=\epsilon_v^{2(\dimk+1)}\), and suppose
$0<\epsilon_v<\min\{2^{-1/\dimk},1,\epsilon_{\mathrm{geo}}\}$ and
    $\epsilon_v\le\lambda_{\mathrm{mom}}^{1/(2(\dimk+1))}$.

    If $L_F=0$, then $F_{\max}$ is constant on the connected thickening; one
    may return \(\theta_0\), and the residual below is zero for every
    $\eta>0$.  Suppose $L_F>0$ and
    set
    \(\eta:=\frac{\rho}{\sqrt m\,L_F}\) and
    \(\mathrm B:=\mathrm B(\delta,\lambda)\).  There exist finite
    constants $C_{\mathrm{res}},C_N,C_u,$ and $C_{\mathrm{PSGD}}>0$,
    depending only on the fixed problem data and $\bar\rho$, such that the
    following lower bounds suffice:
    \[
    \begin{aligned}
      N&\ge C_Nm\epsilon^{-2}\rho^{-2},\\
      b_u&\ge C_um^2\epsilon^{-2},\\
      |I|&\ge4\widetilde C_fm^{2\dimk+2}
          (\rho\epsilon)^{-2\dimk-2},\\
      L&\ge\max\left\{1,
          C_{\mathrm{PSGD}}\dimk^2m^{2\dimk+2}\mathrm B^2
          (\rho\epsilon)^{-2\dimk-2}
          \left[1+\log^2\!\left(\frac{m}{\rho\epsilon}\right)\right]
          \right\}.
    \end{aligned}
    \]
    Let $\{\theta_i\}_{i=0}^{N}$ be generated by
    \Cref{alg: Stochastic-Zeroth-order} with this common outer stepsize.
    Then
    \begin{align}\label{eq_120}
        \min_{0\le i\le N-1}\mathbb{E}\!\left[
        \min_{g\in\partial_\rho F_{\max}(\theta_i)}
        \bigl\|\mathcal G_\Theta(\theta_i,g;\eta)\bigr\|
        \right]
        \le C_{\mathrm{res}}\epsilon,
    \end{align}
    and the randomized iterate returned by the algorithm satisfies
    \[
      \mathbb E\!\left[
      \min_{g\in\partial_\rho F_{\max}(\widehat\theta)}
      \bigl\|\mathcal G_\Theta(\widehat\theta,g;\eta)\bigr\|
      \right]
      \le C_{\mathrm{res}}\epsilon.
    \]
\end{thm}
\end{tcolorbox}
Based on \Cref{thm_cvg_alg_zeroth_goldstein}, choose, for the complexity
bound,
\[
\begin{aligned}
N&=\left\lceil C_Nm\epsilon^{-2}\rho^{-2}\right\rceil, \quad b_u=\left\lceil C_um^2\epsilon^{-2}\right\rceil,\\
|I|&=\max\left\{1,\left\lceil4\widetilde C_fm^{2\dimk+2}
      (\rho\epsilon)^{-2\dimk-2}\right\rceil\right\},\\
L&=\max\Biggl\{1,
\Biggl\lceil C_{\mathrm{PSGD}}\dimk^2m^{2\dimk+2}\mathrm B^2
(\rho\epsilon)^{-2\dimk-2}\times\left[1+\log^2\!\left(\frac{m}{\rho\epsilon}\right)\right]
\Biggr\rceil\Biggr\}.
\end{aligned}
\]
Writing \(Q_f\) for the total number of Gibbs samples,
the two signs give the ceiling-aware bound
\begin{align}\label{eq_0021}
Q_f
&\le2N b_u(L+|I|)\nonumber\\
&\le2\left(1+C_Nm\epsilon^{-2}\rho^{-2}\right)
\left(1+C_um^2\epsilon^{-2}\right)\nonumber\\
&\quad\times\left\{
2+C_{\mathrm{PSGD}}\dimk^2m^{2\dimk+2}\mathrm B^2
(\rho\epsilon)^{-2\dimk-2}
\left[1+\log^2\!\left(\frac{m}{\rho\epsilon}\right)\right]
\right.\nonumber\\
&\hspace{39mm}\left.
+4\widetilde C_fm^{2\dimk+2}
(\rho\epsilon)^{-2\dimk-2}
\right\}.
\end{align}
In the small-accuracy regime, it reduces to the
previously advertised leading order
\[
\mathcal O\!\left(
\frac{m^{2\dimk+5}\left\{
\dimk^2\mathrm B^2\left[1+\log^2\!\left(\frac{m}{\rho\epsilon}\right)\right]
+\widetilde C_f\right\}}
{\epsilon^{2\dimk+6}\rho^{2\dimk+4}}
\right).
\]
\endgroup

 \subsection{Gradient Oracle Complexity under LMC-based Implementation}
\label{section_gradient_oracle_complexity}
In the previous subsections, we analyze the sampling-oracle
complexity required to produce an iterate whose expected generalized-Goldstein
residual is of order \(\epsilon\).
For a more concrete complexity result, we analyze the lower-level gradient
oracle complexity when Langevin Monte Carlo (LMC) implements the Gibbs
sampling oracle.  The Langevin dynamics associated with \Cref{eq_gibbs} is
\begin{equation}\label{eqn_Langevin_dynamics}
    \ud X(t)=-\partial_xg(\theta,X(t))\ud t+\sqrt{2\lambda}\ud W(t),
\end{equation}
where \(W\) is a \(d\)-dimensional Brownian motion.  Under the coercivity and
regularity conditions used below, the Gibbs measure is its unique invariant
law; see, for example, \citep{chiang1987diffusion}.  Its Euler--Maruyama
discretization is the unadjusted LMC iteration
\begin{equation}\label{eqn:langevin}
    \widehat X^{k+1}
    =\widehat X^k-h_k\partial_xg(\theta,\widehat X^k)
      +\sqrt{2\lambda h_k}\,\xi_k,
\end{equation}
where the \(\xi_k\)'s are independent standard Gaussian vectors.  We write
\(\widehat\mu_n:=\operatorname{Law}(\widehat X^n)\).
Our idea is to replace the ideal Gibbs oracle in \Cref{alg:PSGD,alg: Stochastic-Zeroth-order} by the output of LMC algorithm.
While the analyses of this more concrete implementation closely follow the ones in the previous two subsections, we further need to ensure the following two points hold:
\begin{itemize}
    \item The discretization step and the number of LMC iterations must be chosen
jointly so that \(\widehat\mu_n\) is sufficiently close to the target Gibbs
measure \(\mu_\theta^\lambda\) in Wasserstein-1 distance;
    \item Conditioned on the first point, replacing Gibbs samples with LMC samples introduces only a controlled bias in the Gibbs–superquantile estimator.
\end{itemize}
Here ``sufficiently close'' accounts for the
factor \(\delta^{-1}\) in the superquantile objective: the induced error in
the expected hinge loss must be of order \(\delta\epsilon_v\), so that its
contribution after multiplication by \(\delta^{-1}\) is of order
\(\epsilon_v\).  We make this joint choice uniformly on the compact parameter
thickening used by the zeroth-order algorithm.

To control the error caused by replacing exact Gibbs samples with LMC
samples, we convert a R\'enyi-divergence bound for
\(\widehat\mu_n\) relative to \(\mu_\theta^\lambda\) into a Wasserstein
bound between these two measures.  We first verify that the corresponding
\(\mathbb W_2\) distance is finite.  The derived tail quadratic-growth
estimate gives Gaussian tails, and hence a finite second moment, for
\(\mu_\theta^\lambda\).  The LMC initialization is also Gaussian, while the
globally Lipschitz drift has at most linear growth; induction in
\eqref{eqn:langevin} therefore shows that every finite LMC iterate, and hence
\(\widehat\mu_n\), has a finite second moment.  We may now apply
\cite[Theorem~1.1]{liu2020poincare} to obtain
\[
    \mathbb W_1(\widehat\mu_n,\mu_\theta^\lambda)
    \le \mathbb W_2(\widehat\mu_n,\mu_\theta^\lambda)
    \le 2C_{\mathrm{PI}}^{1/2}
       \bigl(e^{\mathbb R_2(\widehat\mu_n\Vert\mu_\theta^\lambda)}-1\bigr)^{1/2},
\]
where \(C_{\mathrm{PI}}\) is the inverse Poincar\'e constant.
We next record the \(q=2\) Poincar\'e-branch guarantee needed below.  It
adapts the proof of \cite[Theorem~6.2.9]{chewi2024logconcave} by applying
the source's change-of-measure lemma at order three.  This adaptation is
proved in \Cref{append_proof_prop_chewi2024logconcave_1}.
\begin{prop}[Poincar\'e-branch LMC bound]
\label{prop_chewi2024logconcave_1}
Let \(\pi(\ud x)\propto e^{-G(x)}\ud x\) on \(\mathbb R^d\).  Assume that
\(\nabla G\) is \(\beta\)-Lipschitz for some \(\beta\ge1\), and that \(\pi\)
has spectral gap at least \(0<\alpha\le1\), namely
\(\operatorname{Var}_\pi(\varphi)\le
\alpha^{-1}\int\|\nabla\varphi\|^2\ud\pi\) for every
\(\varphi\in H^1(\pi)\).  For a constant step \(\bar h>0\), let
\(\widehat\mu_j:=\operatorname{Law}(Y_j)\) for the recursion
\(Y_{j+1}=Y_j-\bar h\nabla G(Y_j)+\sqrt{2\bar h}\,\xi_j\), where the
\(\xi_j\)'s are independent standard Gaussian vectors, and suppose that
\(\mathbb R_3(\widehat\mu_0\Vert\pi)\le R_3<\infty\).

There exist universal constants
\(\varepsilon_{\mathrm{src}}\in(0,1/2]\) and
\(C_T,c_h,C_{\mathrm{Ch}}>0\) such that, for every
\(0<\varepsilon_{\mathrm R}\le\varepsilon_{\mathrm{src}}\), the following
joint choice is valid.  Set
\(T:=C_T\alpha^{-1}[R_3+\log(1/\varepsilon_{\mathrm R})]\), define
\[
 h_0:=c_h\min\left\{
 \frac{1}{4\beta^2T},
 \frac{1}{2\beta^{3/2}\sqrt T},
 \frac{\varepsilon_{\mathrm R}^2}{2\beta^2dT}
 \right\}.
\]
Choose \(n:=\lceil T/h_0\rceil\) and \(\bar h:=T/n\).  Then
\(n\bar h=T\), \(0<\bar h\le h_0\),
\(\mathbb R_2(\widehat\mu_n\Vert\pi)\le\varepsilon_{\mathrm R}^2\), and
\[
 n\le C_{\mathrm{Ch}}\,
 \beta^2\alpha^{-2}d\,\varepsilon_{\mathrm R}^{-2}
 \left[R_3^2+\log^2(1/\varepsilon_{\mathrm R})\right].
\]
\end{prop}
The three entries in \(h_0\) are, respectively, the two auxiliary
discretization restrictions and the accuracy-dependent restriction.  The
rounding above is part of the construction: it preserves the exact horizon
\(n\bar h=T\) while ensuring \(\bar h\le h_0\).  Thus the proposition does
not assert that an arbitrary larger number of iterations works with the same
step; a longer run requires the step to be selected again consistently.

We now construct one set of LMC parameters that works for every query
\(\vartheta\in\Theta^{+\bar\rho}\).  \Cref{prop_chewi2024logconcave_1}
requires common bounds on three quantities:
the smoothness of the standard potential, the spectral gap of the target
measure, and the R\'enyi divergence of the initialization.  Replacing
\(L_{g,xx}\) and \(C_{\mathrm{PI}}^\star\) by their maxima with \(1\), if
necessary, preserves their respective upper-bound properties; hence we take
\(L_{g,xx}\ge1\) and \(C_{\mathrm{PI}}^\star\ge1\).

Let
\(\lambda_{\mathrm{adm}}:=
\min\{1,\lambda_{\mathrm{mom}},\lambda_{\mathrm{PI}}^\star\}\),
\(\alpha_\star:=(C_{\mathrm{PI}}^\star)^{-1}\), and
\(\beta_\lambda:=L_{g,xx}/\lambda\).  The restriction
\(0<\lambda\le\lambda_{\mathrm{adm}}\) places us simultaneously in the
moment and uniform-Poincar\'e temperature ranges.  For the Gibbs target,
\(G_{\vartheta,\lambda}:=g(\vartheta,\cdot)/\lambda\) is the standard
potential.  Its gradient is \(\beta_\lambda\)-Lipschitz, while
\Cref{prop_uniform_poincare} gives the common spectral-gap lower bound
\(\alpha_\star\).

It remains to control the initial R\'enyi divergence.  Fix
\(0<c<3/(2L_{g,2})\) and, for
\(0<\lambda\le\lambda_{\mathrm{adm}}\), define
\[
 \mathcal R_3(\lambda)
 :=\sup_{\vartheta\in\Theta^{+\bar\rho}}
 \mathbb R_3\!\left(
   \mathcal N(0,c\lambda I_d)\,\middle\|\,\mu_\vartheta^\lambda
 \right).
\]
The Gaussian warm-start estimate in
\Cref{lem_gaussian_lmc_initialization} yields
\[
 \lambda\mathcal R_3(\lambda)\le\mathfrak R_3^\star
 \qquad(0<\lambda\le\lambda_{\mathrm{adm}}),
\]
where \(\mathfrak R_3^\star<\infty\) is defined in
\eqref{eq_uniform_lmc_warm_factor}.  Thus the scaled warm-start divergence is
bounded uniformly in \(\vartheta\) and \(\lambda\).

The construction in \Cref{prop_chewi2024logconcave_1} can now be
made uniform over \(\Theta^{+\bar\rho}\).  Replace the parameter-dependent
initial divergence by its common upper bound \(\mathcal R_3(\lambda)\), and
choose the common integration time
\[
T_\star
:=
C_T\alpha_\star^{-1}
\left[
\mathcal R_3(\lambda)
+\log(1/\varepsilon_{\mathrm R})
\right].
\]
Next, define the largest step allowed by the three restrictions in
\Cref{prop_chewi2024logconcave_1} by
\[
h_{\max}
:=
c_h\min\left\{
\frac{1}{4\beta_\lambda^2T_\star},
\frac{1}{2\beta_\lambda^{3/2}\sqrt{T_\star}},
\frac{\varepsilon_{\mathrm R}^2}
     {2\beta_\lambda^2dT_\star}
\right\}.
\]

Choose the step and iteration count jointly by
\begin{equation}\label{eq_uniform_lmc_joint_policy}
 n_\star(\lambda,\varepsilon_{\mathrm R})
 :=\left\lceil\frac{T_\star}{h_{\max}}\right\rceil,
 \qquad
 \bar h_\star(\lambda,\varepsilon_{\mathrm R})
 :=\frac{T_\star}{n_\star(\lambda,\varepsilon_{\mathrm R})}.
\end{equation}
Then \(\bar h_\star\le h_{\max}\) and
\(n_\star\bar h_\star=T_\star\).  Because every actual smoothness,
inverse-gap, and initial-divergence value is bounded by
\(\beta_\lambda\), \(\alpha_\star^{-1}\), and
\(\mathcal R_3(\lambda)\), respectively, the same pair works for every
\(\vartheta\in\Theta^{+\bar\rho}\), with every chain initialized from
\(\mathcal N(0,c\lambda I_d)\).  In the physical-temperature
iteration \eqref{eqn:langevin}, this standard-potential step corresponds to
the constant choice \(h_k=\bar h_\star/\lambda\).

Set
\(\varepsilon_{\mathrm{Ch}}
:=\min\{\varepsilon_{\mathrm{src}},1/2\}\).
For every
\(0<\varepsilon_{\mathrm R}\le\varepsilon_{\mathrm{Ch}}\),
\Cref{prop_chewi2024logconcave_1} applies with the common bounds
\(\beta_\lambda\), \(\alpha_\star\), and
\(\mathcal R_3(\lambda)\).  The joint policy
\eqref{eq_uniform_lmc_joint_policy} therefore satisfies, simultaneously for
every \(\vartheta\in\Theta^{+\bar\rho}\),
\[
 \mathbb R_2(\widehat\mu_{n_\star,\vartheta}
                 \Vert\mu_\vartheta^\lambda)
 \le\varepsilon_{\mathrm R}^2.
\]
For completeness, the ceiling in this specific common policy gives
\[
\begin{aligned}
 n_\star
 &\le 1+\frac{1}{c_h}\max\left\{
 4\beta_\lambda^2T_\star^2,
 2\beta_\lambda^{3/2}T_\star^{3/2},
 2\beta_\lambda^2dT_\star^2\varepsilon_{\mathrm R}^{-2}
 \right\} \\
 &\le C\beta_\lambda^2\alpha_\star^{-2}d
 \varepsilon_{\mathrm R}^{-2}
 \left[\mathcal R_3(\lambda)^2
       +\log^2(1/\varepsilon_{\mathrm R})\right].
\end{aligned}
\]
Here the second inequality uses the same normalizations as
\Cref{prop_chewi2024logconcave_1}, together with the definition of
\(T_\star\).  Substituting
\(\beta_\lambda=L_{g,xx}/\lambda\) and
\(\alpha_\star^{-1}=C_{\mathrm{PI}}^\star\) yields
\begin{equation}\label{eq_uniform_lmc_mixing_bound}
 n_\star(\lambda,\varepsilon_{\mathrm R})
 \le C_{\mathrm{mix}}d(C_{\mathrm{PI}}^\star)^2L_{g,xx}^2
 \lambda^{-2}\varepsilon_{\mathrm R}^{-2}
 \left[\mathcal R_3(\lambda)^2
       +\log^2(1/\varepsilon_{\mathrm R})\right]
\end{equation}
for a universal constant \(C_{\mathrm{mix}}\) independent of
\(\vartheta,\lambda,\varepsilon_{\mathrm R}\).

For the oracle construction below, take
\begin{equation}\label{eq_uniform_lmc_accuracy_threshold}
 \varepsilon_{\mathrm R}:=\delta\epsilon_v
 =\epsilon_v^{\dimk+1}\le\varepsilon_{\mathrm{Ch}},
 \qquad\text{equivalently}\qquad
 \epsilon_v\le
 \varepsilon_{\mathrm{Ch}}^{1/(\dimk+1)}.
\end{equation}
The coupling \(\lambda=\epsilon_v^{2(\dimk+1)}\) gives
\(\lambda=\varepsilon_{\mathrm R}^2\).  Since
\(\lambda\mathcal R_3(\lambda)\le\mathfrak R_3^\star\) and
\[
 \lambda\log(1/\varepsilon_{\mathrm R})
 =\frac{\lambda}{2}\log(1/\lambda)\le\frac{1}{2e},
\]
\eqref{eq_uniform_lmc_mixing_bound} yields
\[
 n_\star\le\mathfrak C_{\mathrm{mix}}d\,
 \epsilon_v^{-10(\dimk+1)},
 \qquad
 \mathfrak C_{\mathrm{mix}}
 :=C_{\mathrm{mix}}(C_{\mathrm{PI}}^\star)^2L_{g,xx}^2
   [1+(\mathfrak R_3^\star)^2].
\]
The constant \(\mathfrak C_{\mathrm{mix}}\) is uniform in the parameter,
temperature, and accuracy variables.  It may, however, depend on \(d\), the
minimizer-manifold geometry, and the fixed problem constants.  Thus the
display preserves the \(\epsilon_v^{-10(\dimk+1)}\) exponent.

Finally, for \(\varepsilon_{\mathrm R}\le1/2\), the
transportation inequality of \cite[Theorem~1.1]{liu2020poincare}, together
with \Cref{prop_uniform_poincare}, implies
\[
 \sup_{\vartheta\in\Theta^{+\bar\rho}}
 \mathbb W_2(\widehat\mu_{n_\star,\vartheta},
             \mu_\vartheta^\lambda)
 \le C_W(C_{\mathrm{PI}}^\star)^{1/2}\varepsilon_{\mathrm R}
 =O(\delta\epsilon_v),
\]
with a universal \(C_W\).  This is the distributional estimate used in the
LMC oracle proof.

Analogously to \Cref{lem_5_2} for exact Gibbs samples,
\Cref{lemm_E_1} in \Cref{append_comp_complex_gibbs} shows that this
Wasserstein estimate gives the LMC-based value oracle the mean-square
accuracy required by the outer zeroth-order analysis.  The resulting fully
quantified complexity statement is as follows.

\begin{thm}\label{theorem_gradient_complexity}
Fix \(\bar\rho>0\), and suppose that
\Cref{mu_PL_assum,assum_secon_fund_form_bound,assum_sample_comp} hold
uniformly on \(\Theta^{+\bar\rho}\).  Let
\(L_F:=L_{f,1}(1+2L_{g,\theta x}/\mu)\), let \(0<\rho<\bar\rho\), and set
\(\epsilon_v:=\rho\epsilon/m\), \(\delta:=\epsilon_v^{\dimk}\), and
\(\lambda:=\epsilon_v^{2(\dimk+1)}\).  Also set
\(\mathrm B:=\mathrm B(\delta,\lambda)\) as in \Cref{lem_3}, and choose
\(\theta_0\in\Theta\) and \(\beta_0\in[-\mathrm B,\mathrm B]\).
Assume \(\dimk\ge1\), the admissible range of
\Cref{thm_cvg_alg_zeroth_goldstein}, and, in addition,
\[
 \epsilon_v\le\min\left\{
 2^{-1/\dimk},
 \varepsilon_{\mathrm{Ch}}^{1/(\dimk+1)}
 \right\},
 \qquad
 0<\lambda\le\lambda_{\mathrm{adm}}.
\]
If \(L_F=0\), return \(\widehat\theta:=\theta_0\)
without making an oracle call.  Then the generalized-Goldstein residual is
zero for every \(\eta>0\), and the theorem is immediate.  Henceforth suppose
\(L_F>0\).
Fix \(c\in(0,3/(2L_{g,2}))\).  Initialize every LMC chain from
\(\mathcal N(0,c\lambda I_d)\), and run it using the common joint policy
\eqref{eq_uniform_lmc_joint_policy} with R\'enyi target
\(\varepsilon_{\mathrm R}=\delta\epsilon_v\).  Let
\(\eta:=\rho/(\sqrt mL_F)\).  There exist finite constants
\(C_N,C_u,C_{\mathrm{PSGD}},C_{\mathrm{res}}^{\mathrm{LMC}}>0\), depending
only on the fixed problem data on the thickening, \(\bar\rho\), and the fixed
initialization scale \(c\), and independent of the accuracy and algorithmic
parameters.  The residual constant \(C_{\mathrm{res}}^{\mathrm{LMC}}\) may
inherit dependence on \(d\) through the uniform Poincar\'e and LMC-oracle
constants.  Let \(\widetilde C_f\) be the uniform validation-moment constant
from \Cref{lemm_E_1}.  Choose
\[
\begin{aligned}
 N&:=\left\lceil C_Nm\epsilon^{-2}\rho^{-2}\right\rceil,
 &
 b_u&:=\left\lceil C_um^2\epsilon^{-2}\right\rceil,\\
|I|&:=\max\left\{1,\left\lceil
 4\widetilde C_fm^{2\dimk+2}(\rho\epsilon)^{-2\dimk-2}
 \right\rceil\right\},
\end{aligned}
\]
and
\[
 L:=\max\left\{1,\left\lceil
 C_{\mathrm{PSGD}}\dimk^2m^{2\dimk+2}\mathrm B^2
 (\rho\epsilon)^{-2\dimk-2}
 \left[1+\log^2\!\left(\frac{m}{\rho\epsilon}\right)\right]
 \right\rceil\right\}.
\]
Run \Cref{alg:zroth-order-LMC} with these choices and the initialization
specified above.  It returns \(\widehat\theta=\theta_J\).  The index
\(J\) is independent of the algorithmic history and is uniformly
distributed on \(\{0,\ldots,N-1\}\).  Moreover,
\[
 \mathbb E\!\left[
 \min_{g\in\partial_\rho F_{\max}(\widehat\theta)}
 \|\mathcal G_\Theta(\widehat\theta,g;\eta)\|
 \right]
 \le C_{\mathrm{res}}^{\mathrm{LMC}}\epsilon .
\]

Let \(Q_f\) be the total number of evaluations of
\(f\), and let \(Q_g\) be the total number of evaluations of the lower-level
gradient \(\partial_x g\) made by the LMC chains.
For \(L_F=0\), one has \(Q_f=Q_g=0\).  In the case
\(L_F>0\), one also has \(\mathrm B>0\), and
\[
Q_f=2Nb_u(L+|I|),
\qquad
Q_g=n_\star Q_f.
\]
More precisely, if
\[
 C_Nm\epsilon^{-2}\rho^{-2}\ge1,
 \qquad C_um^2\epsilon^{-2}\ge1,
\]
and
\[
 m^{2\dimk+2}(\rho\epsilon)^{-2\dimk-2}
 \left\{
 4\widetilde C_f+C_{\mathrm{PSGD}}\dimk^2\mathrm B^2
 \left[1+\log^2\!\left(\frac{m}{\rho\epsilon}\right)\right]
 \right\}\ge1,
\]
then the ceiling-aware bounds in the proof simplify to
\[
 Q_f=\mathcal O\!\left(
 \frac{m^{2\dimk+5}
 \left\{\dimk^2\mathrm B^2
 [1+\log^2(m/(\rho\epsilon))]+\widetilde C_f\right\}}
 {\rho^{2\dimk+4}\epsilon^{2\dimk+6}}
 \right)
\]
and
\[
 Q_g=\mathcal O\!\left(
 \mathfrak C_{\mathrm{mix}}
 \frac{d\,m^{15+12\dimk}
 \left\{\dimk^2\mathrm B^2
 [1+\log^2(m/(\rho\epsilon))]+\widetilde C_f\right\}}
 {\epsilon^{12\dimk+16}\rho^{12\dimk+14}}
 \right).
\]
\end{thm}
The proof of \Cref{theorem_gradient_complexity} is given in
\Cref{append_e_3}.
 \section{Experiments}\label{sec:experiment}
Before detailing the experiments, we outline our adaptation of \Cref{alg: Stochastic-Zeroth-order} to the \emph{optimistic} formulation (see \cref{footnote}), the standard in hyperparameter optimization.
We will compare this optimistic variant against existing baselines in both synthetic and data hyper-cleaning setups.
 To approximate
\(
F_{\min}(\theta) \coloneqq \min_{x \in \mathcal{S}(\theta)} f(\theta, x),
\)
we use the lower-tail superquantile
\[
\tilde{F}_{\mathrm{LSQ}}(\theta) \coloneqq \max_{\beta}\Bigl\{\, \beta - \delta^{-1}\,\mathbb{E}\bigl[\max\{\beta - f(\theta, X),\, 0\}\bigr] \Bigr\}.
\]
In our implementation of \Cref{alg: Stochastic-Zeroth-order}, we replace \(\phi_{\lambda,\delta}(\theta,\beta)\) with the surrogate \( \psi_{\lambda,\delta}(\theta,\beta) \coloneqq \beta - \delta^{-1}\,\mathbb{E}_{\mu_{g}^{\lambda}(\theta)}\!\left[\max\{\beta - f(\theta, X),\, 0\}\right].\)
Additionally, in \Cref{alg:PSGD} we change the projected gradient \emph{descent} step in line~2 to a projected gradient \emph{ascent} step.

\paragraph{Compute budget and reporting protocol.}
Our method is zeroth-order and each outer update invokes a Gibbs--superquantile oracle implemented via PSGD and (parallel-chain) Langevin Monte Carlo, so the wall-clock cost per outer iteration can be higher than that of gradient-based baselines.
Accordingly, in the main text we report \emph{iteration-based} results to illustrate algorithmic behavior and minima-selection effects, and in Appendix~\ref{app:fixed_runtime} we additionally provide \emph{compute-matched} comparisons under a fixed wall-clock time budget (including per-iteration timings for all methods).

Experiments are run on a single machine with \texttt{CPU: Apple M1 Pro (10 cores)}, \texttt{RAM: 16 GB} with no discrete GPU.
Code to reproduce all experiments is released at \href{https://github.com/msmasiha1997/ZSGD-MinSel.git}{GitHub}. Toy-example methods are implemented in JAX, while data hyper-cleaning methods are implemented in PyTorch.

\subsection{Toy Example}
\paragraph{Optimistic formulation:} We evaluate our algorithm along with three more algorithms in the literature on the optimistic BLO problem as follows
\[
\min_{\theta\in[-\pi,\pi]}\;\min_{x\in \mathcal{S}(\theta)}
\;f(\theta,x),
\]
where $f(\theta,x)
=2\|x\|+\langle x,u(\theta)\rangle$ and

\begin{equation}\label{eq_def_g_k}
    g_{\dimk}(\theta,x)
    = \frac{1}{4}\Bigl(\sum_{i=1}^{\dimk+1} x_{i}^{2} - \theta^{2}\Bigr)^{2}
      - \frac{1}{2}\sum_{i=1}^{\dimk+1}x_i^2
      + \frac{1}{2}\sum_{i=\dimk+2}^{d}x_i^2.
\end{equation}
The last sum is zero when \(\dimk=d-1\). The positive quadratic term in the inactive coordinates is necessary to make
the objective bounded below and to obtain the stated minimizer sphere when \(\dimk<d-1\). Here
$u(\theta)=[\cos\theta,\sin\theta,0,\dots,0]^\top\in\R^d$.
The solution set is the \(\dimk\)-dimensional sphere of radius \(\sqrt{1+\theta^{2}}\) embedded in \(\mathbb{R}^{d}\) along the first \(\dimk+1\) coordinates:
\[
\mathcal{S}(\theta)= \Bigl\{\, x \in \mathbb{R}^{d} : \sum_{i=1}^{\dimk+1} x_{i}^{2} = 1+\theta^{2},\; x_{\dimk+2}=\cdots=x_{d}=0 \,\Bigr\}.
\]

So the induced hyper-objective
\[
F_{\min}(\theta)=\sqrt{1+\theta^2}
\]
has the unique global minimizer $\theta^*=0$.

Throughout the experiments, we compare the following algorithms. The fixed-intrinsic-dimension rows reported below
were produced by the original implementation, which omitted the stabilizing last term in \eqref{eq_def_g_k}; those rows
must therefore be rerun before they can be interpreted as evidence for the corrected model.
\begin{itemize}
    \item \textbf{PSZO-MinSel (\Cref{alg: Stochastic-Zeroth-order})}\, adapted to $\min_{x\in\mathcal{S}(\theta)}f(\theta,x)$ instead of pessimistic formulation.
    \item \textbf{V-PBGD} \citep{shen2025penalty}: joint gradient descent on the penalized objective $\Phi(\theta,x)=f(\theta,x)+\gamma(g(\theta,x)-g^{*}(\theta))$,
    \item \textbf{G-PBGD} \citep{shen2025penalty}: joint gradient descent on the penalized objective $\Phi(\theta,x)=f(\theta,x)+\tfrac\gamma2\|\nabla_x g(\theta,x)\|^2$,
    \item \textbf{RHG} \citep{franceschi2017forward}: truncated reverse‐mode hypergradient with re‐initialization of $\theta$ every outer step,
    \item \textbf{IAPTT-GM} \citep{liu2021towards}: unroll $T$ inner steps $x_{k+1}=x_k-\beta\nabla_x g(\theta,x_k)$ from $x_0=z$, pick the pessimistic truncation $k_\star=\arg\max_k F(\theta,x_k)$, and update $\theta,z$ by gradient descent on the truncated objective $F(\theta,x_{k_\star})$.
\end{itemize}
In this part, all methods are implemented in JAX with identical initialization $\theta_0=1.0$, $x_0\sim\mathcal N(0,I_d)$, and random seed.
\begin{figure}[t]
    \centering        \includegraphics[width=0.7\textwidth]{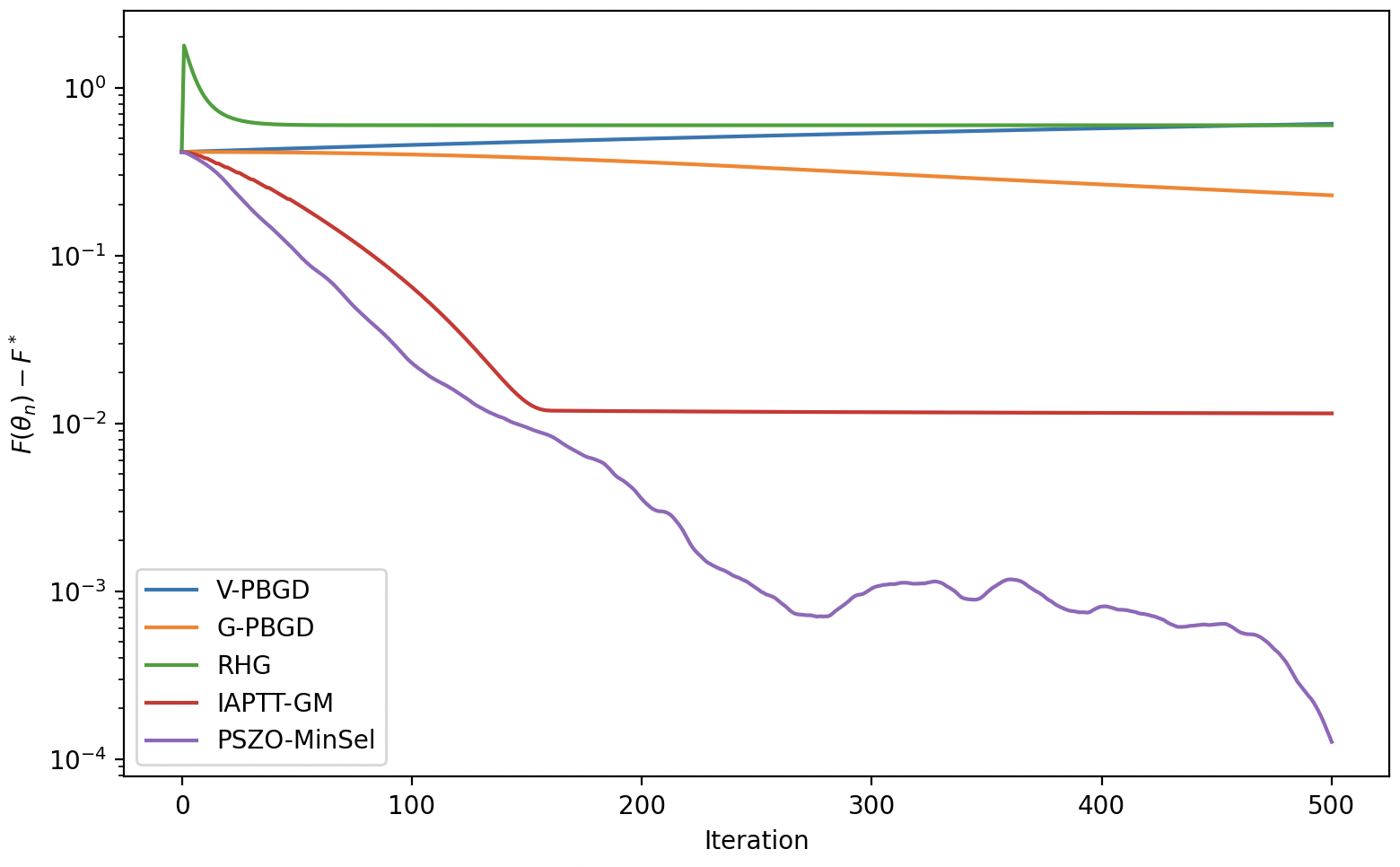}
\label{fig:image1}
    \caption{Comparison of $F(\theta_n)-F^*$ versus iteration for $d=2$.  Our algorithm converges faster with the chosen hyperparameters $\delta=0.1$ and $\lambda=0.01$ in optimistic formulation, while V-PBGD, G-PBGD, RHG, and IAPTT-GM obtain higher values.}
    \label{fig:objective}
\end{figure}

Figure~\ref{fig:objective} plots $F(\theta_n)-F^*$ over iterations for $d=2$. The zeroth-order method rapidly converges toward zero, while the other methods tend to plateau or oscillate around suboptimal objective values, primarily due to the absence of a minima-selection phase. 

\subsubsection{Dimension Sweep and Minima‐selection}
\begin{table}[t]
\centering
\begin{tabular}{c|ccccc}
\hline
$d$ & V-PBGD & G-PBGD & RHG & IAPTT-GM & \textbf{PSZO-MinSel} \\ \hline
2   & 0.297  & 0.160 & 0.304  & 0.152 & $\mathbf{0.009\pm 0.001}$  \\
5   & 0.273  &  0.160 & 0.302 & 0.151 & $\mathbf{0.015\pm 0.012}$ \\
10  & 0.169 & 0.153 & 0.303 & 0.152 & $\mathbf{0.023\pm 0.017}$   \\
20  & 0.113  & 0.174 & 0.302 & 0.154 & $\mathbf{0.043\pm 0.025}$ \\
\hline
\end{tabular}
\caption{Best-so-far absolute error $|\hat{\theta}_N-\theta^*|$ after $N=500$ outer iterations, where
$\hat{\theta}_N$ is the selected iterate defined by the best objective value along the trajectory.
PSZO-MinSel reports mean $\pm$ 95\% CI over 10 runs.}
\label{tab:dim_sweep}
\end{table}
In this experiment, we examine the setting where the intrinsic dimension scales with the ambient dimension, specifically setting $\dimk = d - 1$. This simplifies the lower-level objective to $g(\theta,x)=\frac{1}{4}(\|x\|^{2}-\theta^{2})^{2}-\frac{1}{2}\|x\|^{2}$. \Cref{tab:dim_sweep} reports the best-so-far absolute error $|\hat{\theta}_N-\theta^*|$ after $N=500$ outer iterations across varying dimensions $d$,
where $\hat{\theta}_N$ denotes the selected iterate with the smallest achieved hyper-objective value along the run. Despite this scaling, our algorithm consistently outperforms the three baseline methods, achieving a significantly closer approximation to the global minimizer \(\theta^{\ast} = 0\). For robustness, we report the mean and confidence intervals over $10$ independent runs for our randomized method, compared to single runs for the deterministic baselines. Baselines are deterministic given initialization; repeated runs match up to numerical tolerance.

\subsubsection{Constant Manifold Dimension}
We next report the original fixed-dimension experiment, which set $\dimk=1$ while varying the ambient dimension $d$.
The implementation used to produce \Cref{tab:dim_fixed} omitted the stabilizing positive quadratic term in the inactive
coordinates of \eqref{eq_def_g_k}. For \(\dimk<d-1\), that implemented objective is unbounded below and does not have
the claimed sphere of minimizers. The values are retained for provenance, but the experiment must be rerun with the
corrected objective before it can test dependence on intrinsic rather than ambient dimension.

\begin{table}[t]
\centering
\begin{tabular}{c| ccccc}
\hline
$d$ & V-PBGD & G-PBGD & RHG & IAPTT-GM & \textbf{PSZO-MinSel}\\ \hline
5   & 0.495 &  0.313 & 0.302 & 0.119 & $\mathbf{0.015\pm 0.014}$ \\
10  & 0.499  & 0.313  & 0.316 & 0.114 & $\mathbf{0.018\pm 0.012}$ \\
20  & 0.511 & 0.315  & 0.307 & 0.121 & $\mathbf{0.018\pm 0.011}$\\
30  &  0.532 & 0.313 & 0.311 & 0.121 & $\mathbf{0.02\pm 0.02} $ \\ \hline
\end{tabular}
\caption{Best-so-far absolute error $|\hat{\theta}_N-\theta^*|$ after $N=500$ outer iterations, with a fixed intrinsic dimension $\dimk=1$.}

\label{tab:dim_fixed}
\end{table}

\paragraph{Pessimistic formulation.}
To show that our method is able to handle both the optimistic and pessimistic settings, we solve the previous synthetic problem also in the latter setup.
To our knowledge, existing hypergradient and penalty-based methods do not model the minimizer-selection step. They either (i) assume the lower-level solution set is a singleton or (ii) return points satisfying KKT conditions, without guaranteeing that these correspond to a selected minimizer of the original bilevel objective. Thus, the connection between their outputs and the true BLO solution is unclear. We also observe that their algorithms are \emph{identical} under the optimistic and pessimistic formulations; when the lower-level solution is non-unique, comparing them in both settings is hence not meaningful. Because these methods are designed for the optimistic case, we benchmark against them only in the optimistic setting. For the pessimistic formulation, we report results only for our method.

We now study the following pessimistic toy problem
\begin{align}\label{eq_010121}
    \min_{\theta\in[-\pi,\pi]}\;\max_{x\in \mathcal{S}(\theta)} \; f(\theta,x),
\end{align}
with \(f\) and \(g_{\dimk}\) identical to the previous setup. The induced hyper-objective
\(
F_{\max}(\theta)=3\sqrt{1+\theta^2}
\)
admits the unique global minimizer \(\theta^*=0\).

We also report the original fixed-\(\dimk\) pessimistic run in \Cref{tab:dim_fixed_pess}. The \(g_{d-1}\) column uses
the well-defined original objective, whereas the \(g_1\) column was produced by the unstabilized implementation
described above and must be rerun before it can support a claim about intrinsic dimension.

\begin{table}[t]
\centering
\begin{tabular}{c|cc}
\hline
\(d\) & LL: \(g_{d-1}(\theta,x)\) \eqref{eq_def_g_k} & LL: \(g_{1}(\theta,x)\) \eqref{eq_def_g_k} \\
\hline
5   & \(\mathbf{0.011\pm 0.010}\) & \(\mathbf{0.021\pm 0.011}\) \\
10  & \(\mathbf{0.025\pm 0.020}\) & \(\mathbf{0.023\pm 0.024}\) \\
20  & \(\mathbf{0.030\pm 0.025}\) & \(\mathbf{0.026\pm 0.023}\) \\
30  & \(\mathbf{0.045\pm 0.029}\) & \(\mathbf{0.027\pm 0.021}\) \\
\hline
\end{tabular}
\caption{Best-so-far absolute error \(|\hat{\theta}_N-\theta^*|\) after \(N=500\) outer iterations, evaluated using \Cref{alg: Stochastic-Zeroth-order} in the pessimistic setting.
The first column corresponds to \eqref{eq_010121} with lower-level \(g_{d-1}\) from \eqref{eq_def_g_k}; the second uses the same problem with lower-level \(g_{1}\) from \eqref{eq_def_g_k}.}
\label{tab:dim_fixed_pess}
\end{table}
\subsection{Data Hyper-cleaning}
In this section, we evaluate our method on the data hyper-cleaning benchmark~\citep{franceschi2017forward,shaban2019truncated}, a canonical BLO task widely used in the machine learning community.
Given a set of corrupted\footnote{A sample is \emph{clean} if its label equals the true class; otherwise it is \emph{corrupted}. In MNIST, we corrupt by flipping the label of a fraction \(p\) of training points uniformly to a different (incorrect) class. $p$ is the pollute-rate.
} training samples $\mathcal{D}_{\text{tr}}=\{(d_i^{\text{tr}},\,\ell_i^{\text{tr}})\}_{i=1}^N$, a \emph{small} clean validation set $\mathcal{D}_{\text{val}}=\{(d_i^{\text{val}},\,\ell_i^{\text{val}})\}_{i=1}^M$, and a clean test set, the goal is to learn a data cleaner that down-weights mislabeled training points so as to improve generalization on unseen data. Let $\theta$ denote the data-cleaner parameters (producing per-example weights) and $x$ the model parameters. We use a sigmoid parameterization for the weights,
\(\omega_i(\theta)\;=\;(1+\exp(-\theta_i))^{-1},\ (i=1,\dots,N),\) and constrain $\|\theta\|\le B$ to avoid exploding weights. The task is cast as the optimistic bilevel problem in the following:
\begin{equation}
\label{eq:hypercleaning}
\begin{aligned}
\min_{\theta,\,x}\quad
& -\frac{1}{M}\sum_{i=1}^{M} \log P\!\left(\ell_i^{\text{val}}\mid d_i^{\text{val}};\,x\right) \\
\text{s.t.}\quad
& \|\theta\|\le B,\;\;
x \in \arg\min_{x'\in\mathbb{R}^{d_y}}
\Bigg[
-\frac{1}{N}\sum_{i=1}^{N} \omega_i(\theta)\,\log P\!\left(\ell_i^{\text{tr}}\mid d_i^{\text{tr}};\,x'\right)
+\frac{\eta}{2}\,\|x'\|^2
\Bigg].
\end{aligned}
\end{equation}

\paragraph{Setup.}
We use MNIST with \(N=5{,}000\) training, \(M=100\) validation, and \(10{,}000\) test examples, and uniformly corrupt a fraction \(p\in\{0.4,0.6,0.8\}\) of the training labels. In practice, the corrupted/imprecise data are widely available, but clean/precise data are scarce. Thus we intentionally adopt a very small validation set. Much of the literature uses roughly equal training/validation splits, in which case the effect of the corruption rate in the uncleaned training data during training is often not apparent. We compare our zeroth-order method (\Cref{alg: Stochastic-Zeroth-order}) against V-PBGD and G-PBGD~\citep{shen2025penalty}, RHG~\citep{franceschi2017forward}, and IAPTT-GM~\citep{liu2021towards}, using both a linear model and a two-layer MLP which are standard Neural network models in data hyper-cleaning tasks.

\paragraph{Results.}
\Cref{tab:hypercleaning-results} reports solution quality (test accuracy and F1 score\footnote{F1-score is defined as follows:
We take \emph{clean} samples as the positive class and \emph{corrupted} samples as the negative class. Each sample \(i\) is assigned a score \(s_i=\sigma(\theta_i)\). We predict ``clean'' for the top \((1-p)\)-fraction by \(s_i\). We compute
\(
\text{Precision}=\frac{\mathrm{TP}}{\mathrm{TP}+\mathrm{FP}},\
\text{Recall}=\frac{\mathrm{TP}}{\mathrm{TP}+\mathrm{FN}},
\)
where \(\mathrm{TP}\), \(\mathrm{FP}\), and \(\mathrm{FN}\) denote true positives, false positives, and false negatives, respectively. The F1 score is the harmonic mean of Precision and Recall,
\(
\text{F1}=\frac{2\,\text{Precision}\cdot\text{Recall}}{\text{Precision}+\text{Recall}}.
\)
Note that this F1 measures the \emph{cleaner’s} quality, not the classifier’s accuracy.
}, averaged over $10$ runs with $95\%$ confidence intervals). Our zeroth-order method (see \Cref{alg: Stochastic-Zeroth-order}) adapted to the optimistic formulation attains higher accuracy and F1-score in three different pollute rates ($p=40\%,60\%,80\%$).
Our method’s high performance stems from the \emph{minima-selection} step: by choosing the iterate with the lowest validation loss along the trajectory, it produces higher-quality cleaner parameters $\theta$, yields better model parameters $x$, and ultimately improves generalization on the held-out test set. Appendix~\ref{app:fixed_runtime} additionally reports compute-matched comparisons under a fixed wall-clock time budget.

\begin{table}[ht]
\centering
\caption{Data hyper-cleaning on MNIST across three noise phases (pollute rate $p\in\{40,60,80\}\%$): test accuracy (\%) and F1 score (\%) over $10$ runs ($\pm$ 95\% CI). F1 reflects data-cleaner quality.}
\label{tab:hypercleaning-results}
\tiny
\resizebox{\textwidth}{!}{
\begin{tabular}{cccccccc}
\toprule
\multirow{2}{*}{$p$} & \multirow{2}{*}{Method} & \multicolumn{3}{c}{Linear model} & \multicolumn{3}{c}{2-layer MLP} \\
\cmidrule(lr){3-4}\cmidrule(lr){5-6}
& & Test acc. & F1 & Test acc.  & F1 \\
\midrule
\multirow{5}{*}{40} 
& RHG     & $55.64 \pm 0.45$ & $36.5\pm 0.4$ & $54.64 \pm 0.29$  & $36.4\pm 0.2$ \\
& IAPTT-GM               & $60.4\pm 0.22$  & $54.3\pm0.4$ & $59.37\pm 0.25$  & $52.5\pm0.4$ \\
& G-PBGD               & $73.89\pm 0.13$  & $65.9\pm 0.3$ & $72.54\pm 0.32$  & $65.8\pm 0.3$  \\
& V-PBGD                 & $65.34\pm 0.33$ &  $59.2\pm 0.2$ & $64.77\pm 0.2$ &  $57.2\pm 0.2$ \\
& \textbf{PSZO-MinSel} &$\mathbf{82.85\pm 0.32}$ & $\mathbf{72.1\pm 0.4}$ & $\mathbf{85.31\pm 0.24}$ & $\mathbf{73.2\pm 0.5}$ \\
\cmidrule(lr){1-8}
\multirow{5}{*}{60} 
& RHG                    & $51.24 \pm 0.19$ & $29.8\pm 0.1$ & $52.94 \pm 0.42$  & $28.5\pm 0.5$ \\
& IAPTT-GM               & $59.52\pm 0.38$ & $43.2\pm0.5$ & $58.72\pm 0.35$ & $44.2\pm0.3$ \\
& G-PBGD                 & $64.56\pm 0.52$ & $50.2\pm 0.5$ & $62.65\pm 0.42$ & $51.7\pm0.4$ \\
& V-PBGD                 & $63.14\pm 0.43$ & $48.7\pm 0.5$ & $62.63\pm 0.53$ & $50.4\pm0.4$ \\
& \textbf{PSZO-MinSel} & $\mathbf{79.50\pm 0.18}$ & $\mathbf{55.2\pm 0.1}$ & $\mathbf{80.70\pm 0.28}$ & $\mathbf{55.9\pm 0.2}$ \\
\cmidrule(lr){1-8}
\multirow{5}{*}{80} 
& RHG                    & $48.14 \pm 0.22$ & $21.5\pm 0.2$ & $50.32 \pm 0.34$  & $21.9\pm 0.3$ \\
& IAPTT-GM               & $56.94\pm 0.44$   & $27.5\pm0.3$ & $58.24\pm 0.36$ & $26.1\pm 0.3$ \\
& G-PBGD                 & $46.84\pm 0.32$  & $28.9\pm 0.3$ & $47.70\pm 0.44$  & $27.3\pm0.3$ \\
& V-PBGD                 & $61.23\pm0.32$  & $33.5\pm 0.3$ & $60.42\pm0.42$  & $34.9\pm 0.4$ \\
& \textbf{PSZO-MinSel} & $\mathbf{69.45\pm 0.45}$ & $\mathbf{35.2\pm 0.4}$ & $\mathbf{70.27\pm 0.48}$ & $\mathbf{44.1\pm 0.3}$ \\
\bottomrule
\end{tabular}
}
\end{table}

\section{Conclusion and Future Work}
\label{sec:conclusion}
This paper studies minima-selection in bilevel optimization and proposes a method that is both theoretically grounded and implementable. Under the \PLcirc\ condition, exploiting the manifold structure of the lower-level solution set, allows to introduce a dimension-aware viewpoint: the difficulty of BLO with minima-selection is governed by the intrinsic dimension of the minimizer manifold, rather than the ambient dimension of the lower-level variable. The proposed superquantile--Gibbs (\SQG) relaxation turns minima-selection into a sampling-based procedure with an explicit selection rule. Overall, the results suggest that minima-selection can be handled in a structured way by combining geometric information (the minimizer manifold), stochastic approximation (Gibbs measures), and risk-type criteria (superquantiles).

Looking ahead, several directions appear particularly promising. On the algorithmic side, a natural next step is to move from zeroth-order schemes to \emph{hyper-gradient methods}, in which one first selects a point on the minimizer manifold (using ideas similar to the Gibbs-based selection studied here) to maximize\footnote{For the optimistic formulation, one can also minimize the upper-level function over the lower-level solution manifold.} the upper-level objective, and then differentiates through this choice via implicit differentiation. In parallel, it is important to improve the efficiency of the Gibbs sampling oracle. Our current framework relies on a standard Langevin Monte Carlo implementation; geometry-aware Langevin samplers tailored to the identified manifold structure could sharpen the dependence on both the intrinsic dimension $\dimk$ and the ambient dimension $d$. Such developments would make the proposed PSZO-MinSel algorithm more practical for large-scale applications.
 
\appendix
\newpage
\section{Definitions and Technical Lemmas}
This appendix collects the background components used throughout the paper.
We first fix notation and recall the Riemannian notions used in
\Cref{sec:4_unif_props}.  We then prove the finite-regularity
tubular-neighborhood statement used in the preliminaries and the Laplace limit
in \Cref{lemma_density_limiting_measure}.  We then briefly summarize standard
differentiability results for bilevel hyper-objectives.
Finally, we compare our assumptions with those of
\citep[Theorem~3.6]{hasenpflug2024wasserstein} and identify the quantitative
adaptation still required for the later Wasserstein rate.

\subsection{Additional Notations}
\label[appendix]{append_additional_notations}
For a scalar \(a\), let \(a_+:=\max\{a,0\}\); for a vector,
\((\cdot)_+\) is applied componentwise.  For a set or event \(A\),
\(\mathbf 1_A\) denotes its indicator.  If \(T\) is measurable and \(\nu\)
is a measure, \(T_\#\nu\) denotes the pushforward measure, defined by
\((T_\#\nu)(B):=\nu(T^{-1}(B))\), and \(\nu|_A\) denotes the restriction of
\(\nu\) to \(A\).  For self-adjoint matrices or operators,
\(A\preceq B\) means that \(B-A\) is positive semidefinite; \(A\succeq B\)
has the reverse meaning.  The notation \(O_\theta(\cdot)\) allows its hidden
constant to depend on the fixed parameter \(\theta\), but not on the
asymptotic variable.  Finally, \(\lambda_1(M)\) denotes the first nonzero
Laplace--Beltrami eigenvalue of a compact connected manifold \(M\), while
\(\lambda_1^{\mathrm N}(\Omega)\) denotes the first nonzero Neumann
Laplacian eigenvalue of a domain \(\Omega\).
{\color{ForestGreen}For a \(\mathcal C^1\) map
\(\Phi:M\to N\) between finite-dimensional manifolds, we write
\[
    D\Phi(x)=D\Phi_x:T_xM\longrightarrow T_{\Phi(x)}N
\]
for its differential at \(x\).  For \(v\in T_xM\), it is characterized by
\[
    D\Phi_x[v]
    :=\left.\frac{\ud}{\ud s}\right|_{s=0}(\Phi\circ c)(s),
\]
where \(c\) is any \(\mathcal C^1\) curve with \(c(0)=x\) and
\(\dot c(0)=v\); the value is independent of the choice of such a curve.
When \(M\) is a vector space, we use the canonical
identification \(T_xM\simeq M\).  Subscripts indicate the differentiated
variable; for example, \(D_bF\) is the differential of \(F\) with respect to
its \(b\)-argument.}
For a linear subspace \(V\subseteq\mathbb R^d\), \(P_V\) denotes the
orthogonal projection onto \(V\). If \(A:\mathbb R^d\to\mathbb R^d\) is
linear, then
\[
    P_VA\big|_V:V\to V,
    \qquad
    v\longmapsto P_V(Av),
\]
means that the input of \(A\) is restricted to \(V\) and its output is
projected onto \(V\). Thus,
\(P_{N_p\mathcal S(\theta)}\partial_x^2g(\theta,p)
\big|_{N_p\mathcal S(\theta)}\) is the normal-to-normal part of the Hessian
at \(p\). If the Hessian maps \(N_p\mathcal S(\theta)\) into itself, the
projection is redundant, and this operator is simply
\(\partial_x^2g(\theta,p)\big|_{N_p\mathcal S(\theta)}\).

\subsection{Riemannian Notions}\label[appendix]{append:def}
In this subsection we fix notation for basic Riemannian objects and record definitions that will be used repeatedly in later proofs.
Any compact manifold $\mathcal{S}$ of class \(\mathcal{C}^3\) (indeed \(\mathcal{C}^1\)) without boundary admits a Riemannian metric\footnote{A fundamental theorem in differential geometry states that \textit{every paracompact manifold \(M\) of class \(\mathcal{C}^k\) with \(k \ge 1\) can be endowed with a Riemannian metric}.}.

Given a Riemannian metric \(\mathsf g\) on \(\mathcal S\), the length of a
piecewise \(\mathcal C^1\) curve \(\Gamma:[a,b]\to\mathcal S\) is
\[
    \ell_{\mathsf g}(\Gamma)
    :=\int_a^b
      \sqrt{\mathsf g_{\Gamma(t)}
      \bigl(\dot\Gamma(t),\dot\Gamma(t)\bigr)}\,\ud t.
\]
For \(p,q\in\mathcal S\), the associated Riemannian distance is
\begin{align}\label{def_geod_dist}
    \mathrm{dist}_{\mathsf g}(p,q)
    :=\inf\bigl\{\ell_{\mathsf g}(\Gamma):
       {}&\Gamma:[a,b]\to\mathcal S\text{ is piecewise }\mathcal C^1,\notag\\[-2pt]
       &\Gamma(a)=p,\ \Gamma(b)=q\bigr\}.
\end{align}
We use the convention that this infimum is \(+\infty\) when
    \(p\) and \(q\) lie in different connected components.  If they lie in the
    same component of the compact manifold, that component is complete, and
    the Hopf--Rinow theorem gives a minimizing geodesic between them.  In all
    applications below, \(\mathcal S=\mathcal S(\theta)\) is connected, so
    this condition is automatically satisfied.

\begin{mydef}[Second fundamental form]\label{def:second-fund-form}
Let \(\mathcal S\subset\mathbb R^d\) be a \(\mathcal C^2\) embedded
submanifold with the Euclidean-induced metric.  For \(x\in\mathcal S\) and
\(u,v\in T_x\mathcal S\), extend \(v\) to a local \(\mathcal C^1\) tangent
vector field \(V\) on \(\mathcal S\).
Choose a \(\mathcal C^1\) curve \(\gamma:(-\varepsilon,\varepsilon)
\to\mathcal S\) satisfying \(\gamma(0)=x\) and \(\dot\gamma(0)=u\).  The
directional derivative of \(V\) at \(x\) in the direction \(u\) is
\[
 D_uV(x):=\left.\frac{\ud}{\ud t}V(\gamma(t))\right|_{t=0}
          =DV(x)[u]\in\mathbb R^d.
\]
It depends only on \(x\), \(u\), and \(V\), and not on the chosen curve.
At the fixed point \(x\), every \(w\in\mathbb R^d\) has the unique
orthogonal decomposition
\[
 w=w^\top+w^\perp,
 \qquad w^\top\in T_x\mathcal S,
 \qquad w^\perp\in N_x\mathcal S:=(T_x\mathcal S)^\perp.
\]
Thus \(w^\perp\) denotes the orthogonal projection of \(w\) onto the normal
space \(N_x\mathcal S\).  The second fundamental form is
\[
 \II_{\mathcal S}(x)[u,v]
 :=\bigl(D_uV(x)\bigr)^\perp\in N_x\mathcal S,
\]
namely, the normal component of the change of the tangent vector field \(V\)
in the direction \(u\).
The normal component is
independent of the chosen tangent extension \(V\).  Thus
\(\II_{\mathcal S}(x)\) is a well-defined symmetric bilinear map
\(T_x\mathcal S\times T_x\mathcal S\to N_x\mathcal S\).

We retain the diagonal operator-norm convention
\[
 \|\II_{\mathcal S}(x)\|_{\op}
 :=\sup_{\substack{u\in T_x\mathcal S\\\|u\|=1}}
       \|\II_{\mathcal S}(x)[u,u]\|.
\]
In particular, if \(\gamma\) is a unit-speed geodesic in \(\mathcal S\), then
\[
 \|\ddot\gamma(t)\|
 =\|\II_{\mathcal S}(\gamma(t))
       [\dot\gamma(t),\dot\gamma(t)]\|
 \le\|\II_{\mathcal S}(\gamma(t))\|_{\op}.
\]
\Cref{assum_secon_fund_form_bound} requires this norm to be bounded by one
constant \(\mathsf C\) independent of \(\theta\).  See also
\cite[Appendix~D]{gong2024poincare}.
\end{mydef}

\begin{mydef}[Reach]\label{def_reach}
Let \(M\subseteq\mathbb R^d\) be a nonempty closed set.  For
\(x\in\mathbb R^d\), let
\[
    \Pi_M(x)
    :=\mathop{\arg\min}_{u\in M}\|x-u\|
\]
denote the set of nearest points in \(M\) to \(x\).  The \emph{reach} of
\(M\) is an ambient Euclidean quantity.  Call \(r>0\) admissible if
\(\Pi_M(x)\) is a singleton for every \(x\) satisfying
\(\operatorname{dist}(x,M)<r\).  Then
\[
    \operatorname{reach}(M)
    :=\sup\bigl(\{0\}\cup
       \{r>0:r\text{ is admissible}\}\bigr)
    \in[0,+\infty].
\]
Thus \(\operatorname{reach}(M)\ge r\) means that every point whose Euclidean
distance from \(M\) is strictly less than \(r\) has a unique nearest point
on \(M\).  Reach should not be confused with the intrinsic injectivity
radius defined next.
\end{mydef}

\begin{mydef}[Injectivity radius]\label{def_injec_rad}
For \(p\in\mathcal S\), let
\[
    \mathbb B_{T_p\mathcal S}(0;r)
    :=\{v\in T_p\mathcal S:\|v\|<r\}.
\]
The pointwise injectivity radius is
\[
    \operatorname{inj}_{\mathcal S}(p)
    :=
    \sup\left\{
        r>0:
        \left.\exp_p\right|_{\mathbb B_{T_p\mathcal S}(0;r)}
        \text{ is a diffeomorphism onto its image}
    \right\}.
\]
The restriction can fail to be a diffeomorphism either because the
exponential map loses injectivity at the cut locus or because its differential
becomes singular at a conjugate point.  The injectivity radius of the manifold
is
\[
    \operatorname{inj}(\mathcal S)
    :=\inf_{p\in\mathcal S}\operatorname{inj}_{\mathcal S}(p).
\]
\end{mydef}

\begin{mydef}[Sectional Curvature]\label{def_sec_curv}
    In a Riemannian manifold \((\mathcal{S},\mathsf{g})\), the sectional curvature measures how curved the manifold is in each two-dimensional ``slice" (or ``plane") of the tangent space at a point.  Concretely 1) Pick a point \(p \in \mathcal{S}\). 2) Pick a 2-dimensional subspace \(A_{p} \subset T_{p}\mathcal{S}\).  (In other words, \(A_{p}\) is a plane in the tangent space at \(p\).)
The sectional curvature \(K(p,A_{p})\) is defined using the \emph{Riemann curvature tensor}\footnote{In a Riemannian manifold \((\mathcal{S}, \mathsf{g})\) with Levi-Civita connection \(\nabla\), the Riemann curvature tensor is a \((1,3)\)-tensor \(R\) that measures how much the covariant derivative \(\nabla\) fails to commute.} \(R\).  If \(X\) and \(Y\) are any two linearly independent vectors spanning \(A_{p}\), then
\[
  K(p,A_{p}) 
  = \frac{\mathsf{g}(R(X,Y)Y,X)}{\mathsf{g}(X,X)\mathsf{g}(Y,Y)-\mathsf{g}(X,Y)^{2}}.
\]
Equivalently, if you choose an orthonormal basis \(\{\,e_{1}, e_{2}\}\) of \(A_{p}\) with respect to \(\mathsf{g}\) (so \(\mathsf{g}(e_{i}, e_{j}) = \delta_{ij}\)), then the denominator becomes \(1\), and one obtains the simpler formula
\[
  K(p,A_{p})
  = \mathsf{g}(R(e_{1}, e_{2})\,e_{2},\,e_{1}).
\]
\end{mydef}

\subsection{Finite-regularity tubular neighborhood}
\label[appendix]{append:finite_regularity_tubular_neighborhood}
We prove the tubular-neighborhood statement preceding
\eqref{local coordinate decomposition}.  Write
\(M:=\optimalsetlower\), let \(k:=\dim M\), and let
\(q:=d-k\).  Because \(M\) is a \(\mathcal C^2\) embedded manifold without
boundary, its normal bundle
\[
    NM:=\{(u,v):u\in M,\ v\in N_uM\}
\]
is a \(\mathcal C^1\) manifold of dimension \(k+q=d\).  The addition map
\[
    \Phi:NM\to\mathbb R^d,
    \qquad
    \Phi(u,v):=u+v,
\]
is therefore \(\mathcal C^1\).

We first establish local invertibility along the zero section.  Fix \(x\in M\).
At \((x,0)\), a tangent vector to the total space \(NM\) is uniquely
represented by a pair \((\xi,\eta)\), with \(\xi\in T_xM\) describing motion
of the base point and \(\eta\in N_xM\) describing motion in the normal fiber.
Thus
\[
    T_{(x,0)}NM\cong T_xM\oplus N_xM,
\]
and differentiation of the addition map gives
\[
    D\Phi(x,0)[\xi,\eta]=\xi+\eta.
\]
The orthogonal decomposition \(\mathbb R^d=T_xM\oplus N_xM\) shows that this
derivative is a linear isomorphism.  The \(\mathcal C^1\) inverse-function
theorem in local coordinates therefore provides neighborhoods
\(\mathcal O_x\subset NM\) of \((x,0)\) and \(W_x\subset\mathbb R^d\) of
\(x\) such that \(\Phi:\mathcal O_x\to W_x\) is a \(\mathcal C^1\)
diffeomorphism.

We next choose one radius that works along all of \(M\).  Invertibility is an
open condition on linear maps, and \(D\Phi\) is continuous.  Hence, in a local
normal-bundle trivialization about each \(x\in M\), there are an open
neighborhood \(U_x\subset M\) and a number \(r_x>0\) such that
\(D\Phi(u,v)\) is invertible whenever \(u\in U_x\) and
\(v\in N_uM\) satisfies \(\|v\|<r_x\).  Choose a finite subcover
\(U_{x_1},\ldots,U_{x_N}\) of the compact manifold \(M\), and set
\[
    r_0:=\min_{1\le i\le N}r_{x_i}>0.
\]
Then \(D\Phi(u,v)\) is invertible for every \(u\in M\) and every
\(v\in N_uM\) with \(\|v\|<r_0\).

It remains to ensure that normal fibers based at different points do not meet
inside a sufficiently thin tube.  Suppose, to the contrary, that the
restriction of \(\Phi\) is noninjective for every radius in \((0,r_0)\).
For each sufficiently large integer \(j\), choose distinct pairs
\((x_j,v_j),(x'_j,v'_j)\in NM\) such that
\[
    \|v_j\|,\|v'_j\|<1/j,
    \qquad
    x_j+v_j=x'_j+v'_j.
\]
After passing to subsequences, compactness of \(M\) gives
\(x_j\to x\in M\) and \(x'_j\to x'\in M\).  Since
\(v_j,v'_j\to0\), the displayed equality implies \(x=x'\).  Consequently,
both sequences of pairs converge to \((x,0)\).  They eventually lie in the
neighborhood \(\mathcal O_x\), contradicting the injectivity of \(\Phi\) on
\(\mathcal O_x\).  Thus there exists \(\delta_{\mathrm{tub}}\in(0,r_0)\)
such that \(\Phi\) is injective on
\[
    \mathcal V
    :=\bigl\{(u,v)\in NM:\|v\|<\delta_{\mathrm{tub}}\bigr\}.
\]

The derivative of \(\Phi\) is invertible throughout \(\mathcal V\), so
\(\Phi\) is a local \(\mathcal C^1\) diffeomorphism there and its image is
open.  Since \(\Phi\) is also injective, its local inverses agree with its
set-theoretic inverse; hence
\[
    \Phi:\mathcal V\to\Phi(\mathcal V)
\]
is a \(\mathcal C^1\) diffeomorphism.  Moreover,
\(\Phi(u,0)=u\) for every \(u\in M\), so \(\Phi(\mathcal V)\) is an open
neighborhood of \(M\).  Finally, the inverse map assigns to each
\(y\in\Phi(\mathcal V)\) its unique pair \((u,v)\); composing this inverse
with the bundle projection and the ambient vector component shows that both
\(u\) and \(v\) depend \(\mathcal C^1\)-smoothly on \(y\).  This proves the
statement.

\subsection{Proof of \Cref{lemma_density_limiting_measure}}
\label[appendix]{append:proof_limiting_gibbs_measure}
\begin{proof}
Fix \(\theta\in\Theta\) and suppress this fixed parameter throughout the
proof.  Write
\[
 \begin{aligned}
 g^\star&:=g^\star(\theta), &
 g&:=g(\theta,\cdot)-g^\star,\\
 \mathcal S&:=\mathcal S(\theta), &
 k&:=\dim\mathcal S, &
 q&:=d-k.
 \end{aligned}
\]
Then \(g\ge0\), \(\{g=0\}=\mathcal S\), and this additive shift does not
change the normalized Gibbs measure.  Let \(\mathcal M\)
denote the Riemannian volume measure on \(\mathcal S\).  In the
notation of Hasenpflug et al.~\citep{hasenpflug2024wasserstein}, we take
\(\ell:=g\), \(\pi_0\equiv1\), and inverse temperature
\(t:=\lambda^{-1}\).  All constants below are fixed-parameter constants; no
parameter-uniform conclusion is asserted.

Their Assumptions~2.2 and~3.1 impose a smooth orientable minimizer manifold
with a global smooth orthonormal normal frame, whereas \(\mathcal S\) here
may be nonorientable and is only \(\mathcal C^2\).  Thus their theorem is not
applied verbatim.  We instead adapt the manifold-Laplace part of its proof on
the intrinsic normal bundle, using finitely many local orthonormal normal
frames.

First, \(\mathcal S\) is nonempty and compact by
\Cref{prop_S_path_connected,prop_S_C_2_wihout}.  Define
\[
 B:=\sup_{u\in\mathcal S}\|u\|<\infty.
\]
Since
\(\operatorname{dist}(x,\mathcal S)\ge\|x\|-B\),
\eqref{eq_tail_EB_implies_QG} gives
\(g(x)\ge c\|x\|^2\) outside a sufficiently large ball for some \(c>0\).
Hence
\[
 \int_{\mathbb R^d}e^{-g(x)/\lambda}\,\mathrm dx<\infty
 \quad(\lambda>0),
 \qquad
 Z_\lambda
 =e^{-g^\star/\lambda}
  \int_{\mathbb R^d}e^{-g(x)/\lambda}\,\mathrm dx<\infty.
\]
Choose a sufficiently small fixed tube \(\mathcal U\) about
\(\mathcal S\), with radius \(r\) furnished by
Appendix~\ref{append:finite_regularity_tubular_neighborhood} and reduced
below as necessary.  Choose a ball outside which the quadratic lower bound is
positive.  The part of \(\mathcal U^c\) inside that ball is compact
and disjoint from \(\{g=0\}=\mathcal S\), so continuity, together with the
quadratic tail, gives
\(\varepsilon_{\mathcal U}:=\inf_{\mathcal U^c}g>0\).  Consequently, for
\(0<\lambda\le1\),
\[
 \int_{\mathcal U^c}e^{-g/\lambda}\,\mathrm dx
 \le e^{-(\lambda^{-1}-1)\varepsilon_{\mathcal U}}
     \int_{\mathbb R^d}e^{-g}\,\mathrm dx
 =o(\lambda^{q/2}).
\]
This is the off-tube factorization used in the proof of
Hasenpflug et al.~\citep[Theorem~3.6, Step~2, equations~(29)--(30)]
{hasenpflug2024wasserstein}.

For this choice, normal addition is a \(\mathcal C^1\) diffeomorphism onto
\(\mathcal U\), and
\[
 \mathrm dx=J(u,v)\,\mathrm dv\,\mathrm d\mathcal M(u),
 \qquad J(u,0)=1.
\]
Moreover, \(J(u,v)\to1\) uniformly in \(u\) as \(\|v\|\to0\).
This is the intrinsic, finite-regularity counterpart of the tubular identity
in their equation~(11).  Cover \(\mathcal S\) by finitely many local
\(\mathcal C^1\) orthonormal normal frames \(E_\ell\), and set locally
\[
 A_\ell(u):=
 E_\ell(u)^\top\partial_x^2g(u)E_\ell(u).
\]
The local P\L--Morse--Bott conclusion and global
\(L_{g,xx}\)-smoothness yield
\(\mu I_q\preceq A_\ell(u)\preceq L_{g,xx}I_q\).  Indeed, for a normal
eigenpair \((\alpha,v)\), inserting the expansions along \(u+sv\) into
\(g\le(2\mu)^{-1}\|\partial_xg\|^2\) gives \(\alpha\ge\mu\).
On overlaps, \(\widetilde E_\ell=E_\ell Q\), with \(Q\in O(q)\), and hence
\(\widetilde A_\ell=Q^\top A_\ell Q\).  Orthogonal changes preserve fiber
Lebesgue measure, the Gaussian integral, and \(\det A_\ell\).  Therefore the
weight \(w=\det(A_\ell)^{-1/2}\) in
\eqref{eq_limiting_gibbs_weight} is globally well defined; no orientability
or trivialization of the normal bundle is used.

Define intrinsically
\[
\begin{split}
 A_\lambda(u)&:=
 \int_{\substack{v\in N_u\mathcal S\\\|v\|<r}}
 e^{-g(u+v)/\lambda}J(u,v)\,\mathrm dv,\\
 A_{\lambda,\varphi}(u)&:=
 \int_{\substack{v\in N_u\mathcal S\\\|v\|<r}}
 \varphi(u+v)e^{-g(u+v)/\lambda}J(u,v)\,\mathrm dv
\end{split}
\]
for bounded continuous \(\varphi:\mathbb R^d\to\mathbb R\).  Continuity of
the Hessian on the compact closed tube gives, uniformly in \(u\),
\[
 g(u+v)
 =\tfrac12\langle v,\partial_x^2g(u)v\rangle+o(\|v\|^2),
 \qquad g(u+v)\ge\tfrac\mu4\|v\|^2
\]
after reducing \(r\).  In a local frame, set
\(v=\sqrt\lambda E_\ell(u)z\).  On bounded \(z\)-sets the integrands
converge uniformly; on their complements they are dominated by
\(C e^{-\mu\|z\|^2/4}\).  The Jacobian converges uniformly to one, and
\(\varphi\) is uniformly continuous on the closed tube.  Splitting into a
bounded \(z\)-set and its Gaussian tail therefore gives, uniformly for
\(u\in\mathcal S\),
\[
\begin{split}
 \lambda^{-q/2}A_\lambda(u)
 &\longrightarrow(2\pi)^{q/2}w(u),\\
 \lambda^{-q/2}A_{\lambda,\varphi}(u)
 &\longrightarrow(2\pi)^{q/2}\varphi(u)w(u).
\end{split}
\]
This is the qualitative local-frame version of
Hasenpflug et al.~\citep[Lemma~3.17, equation~(24), and Appendix~A.2]
{hasenpflug2024wasserstein}; their stronger quantitative remainder is not
used.  Integrating over compact \(\mathcal S\) and adding the off-tube term
therefore gives, for every bounded continuous \(\varphi\),
\[
 \int_{\mathbb R^d}\varphi e^{-g/\lambda}\,\mathrm dx
 =\lambda^{q/2}\left[(2\pi)^{q/2}
   \int_{\mathcal S}\varphi w\,\mathrm d\mathcal M+o(1)\right].
\]
Taking \(\varphi=1\) gives the corresponding normalizer expansion.  Although
the source writes \(n\in\mathbb N\), equation~(11), the rescaling
\(v=n^{-1/2}z\), and the off-tube factor use it only as a positive scalar;
the argument therefore holds for every real \(t=\lambda^{-1}\to\infty\).
Dividing by the \(\varphi=1\) expansion yields
\[
 \int\varphi\,\mathrm d\gibbs^\lambda
 \longrightarrow
 \frac{\displaystyle\int_{\mathcal S}\varphi(u)w(u)\,
       \mathrm d\mathcal M(u)}
      {\displaystyle\int_{\mathcal S}w(u)\,
       \mathrm d\mathcal M(u)},
\]
which is \(\gibbs^\lambda\Rightarrow\gibbs^0\) with the density in
\eqref{eqn_density_limiting_measure}.  Finally, \(w\) is continuous and
strictly positive because the normal Hessian is continuous and positive
definite.  Riemannian volume has support \(\mathcal S\), so the positive
density implies \(\operatorname{supp}(\gibbs^0)=\mathcal S\).
\end{proof}

\subsection{Differentiability of Hyper-objective}\label[appendix]{append_diff_hyper_obj}
The authors of \cite{kwon2023penalty,chen2023bilevel} reformulate BLO via the classical value-function approach and assume that a P\L-type condition (or an error bound) holds for the penalized objective $\sigma f+g$, uniformly over $\sigma\in[0,\sigma_{0}]$ where $\sigma_{0}>0$ is a constant. More precisely, Kwon et al.~\cite{kwon2023penalty} impose a proximal error bound for $\sigma f+g$, whereas Chen et al.~\cite{chen2023bilevel} assume a global P\L\ condition for $\sigma f+g$; Chen et al.~\cite[Proposition C.1]{chen2023bilevel} show these two assumptions are equivalent.

Under this global regularity, both works establish a Danskin-type result implying that the hyper-objective $\hyperobjective$ is differentiable. However, their conclusion also entails that every choice $x\in\optimalsetlower(\theta)$ yields the same hyper-gradient $\nabla \hyperobjective(\theta)$, so minima-selection becomes redundant. Finally, requiring a uniform P\L\ condition (or error bound) for $\sigma f+g$ is quite stringent and may fail even when $f$ and $g$ individually satisfy P\L-type properties.

\subsection{Relation to \citep[Theorem~3.6]{hasenpflug2024wasserstein}}\label[appendix]{append_valid_assmp_thm_W_1}
Hasenpflug et al.~\citep[Assumption~2.2 and
Theorem~3.6]{hasenpflug2024wasserstein} formulate their result for smooth
orientable manifolds equipped with global smooth orthonormal normal frames.
Those hypotheses are stronger than ours, so their theorem is not directly
applicable.  For \(p=1\), however, every estimate used in their proof remains
valid here after the global-frame calculation is performed in finitely many
local orthonormal normal frames.  The fixed-\(\theta\) verification is given
in \Cref{append_proof_pointwise_W1_Gibbs}; the separate proof of
\Cref{prop_W1_Gibbs_uniform} tracks the corresponding constants uniformly in
the parameter.

\subsection{Proof of \Cref{prop_W1_Gibbs_pointwise}}
\label[appendix]{append_proof_pointwise_W1_Gibbs}
\begin{proof}
Fix \(\theta\) throughout and write
\[
 g:=g(\theta,\cdot),\qquad
 g^\star:=g^\star(\theta),\qquad
 \mathcal S:=\mathcal S(\theta),\qquad
 k:=\dim\mathcal S,\qquad q:=d-k.
\]
In the notation of Hasenpflug et al.~\citep{hasenpflug2024wasserstein}, take
the potential \(\ell:=g-g^\star\), Lebesgue reference measure
\(\pi_0\equiv1\), and \(t:=\lambda^{-1}\).  The shift leaves the Gibbs measure
and normal Hessian unchanged, and their one-component limit in
equations~(13)--(14) is exactly \(\gibbs^0\) in
\eqref{eqn_density_limiting_measure}.  Constants below may depend on
\(\theta\).

We adapt, rather than directly apply, their Theorem~3.6.  First consider the
analytic hypotheses used in its \(p=1\) proof.  By
\Cref{prop_S_path_connected,prop_S_C_2_wihout}, \(\mathcal S\) is a nonempty
compact connected embedded \(\mathcal C^2\) manifold without boundary; thus
the manifold content of clause (M) used in the proof is present, while its
stronger frame formulation is treated below.  If
\(B_\theta:=\sup_{u\in\mathcal S}\|u\|\), then
\eqref{eq_tail_EB_implies_QG} and
\(\operatorname{dist}(x,\mathcal S)\ge\|x\|-B_\theta\) give, outside a
sufficiently large ball,
\(g(x)-g^\star\ge c_\theta\|x\|^2\).  Consequently
\(\int_{\mathbb R^d}(1+\|x\|)
e^{-[g(x)-g^\star]/\lambda}\,\mathrm dx<\infty\) for
every \(\lambda>0\), giving clause (I).  For any sufficiently small fixed
tube \(\mathcal U_\theta\), continuity,
\(\mathcal S=\{x:g(x)=g^\star\}\), and the same tail bound give
\(\inf_{x\in\mathcal U_\theta^c}[g(x)-g^\star]>0\), which is clause (T).
Factoring
\(e^{-(\lambda^{-1}-1)
\inf_{x\in\mathcal U_\theta^c}[g(x)-g^\star]}\) and using
integrability at \(\lambda=1\) further gives
\[
 \int_{\mathcal U_\theta^c}(1+\|x\|)
 e^{-[g(x)-g^\star]/\lambda}\,\mathrm dx
 \le C_\theta e^{-c_\theta/\lambda}
 \qquad(0<\lambda\le1).
\]

For clauses (S) and (P), choose the tube with compact closure in the
neighborhood where \(g\) is \(\mathcal C^4\).  Compactness bounds \(D_x^3g\)
there, so Taylor's formula gives
\[
 \left|g(u+v)-g^\star
 -\frac12\langle v,\partial_x^2g(u)v\rangle\right|
 \le C_\theta\|v\|^3
 \quad
 (u\in\mathcal S,\ v\in N_u\mathcal S,\ \|v\|\ \text{small}).
\]
The local P\L--Morse--Bott conclusion and the normalization in
\Cref{def_local_PL}, together with global \(L_{g,xx}\)-smoothness, give the
exact normal bounds
\[
 \mu I_q
 \preceq
 \left.\partial_x^2g(u)\right|_{N_u\mathcal S}
 \preceq L_{g,xx}I_q.
\]
The lower bound is \(\mu I_q\), with no factor-of-two loss; this calculation
also appears in the proof of \Cref{lemma_density_limiting_measure}.

It remains to replace the stronger geometry in their Assumption~2.2.  The
finite-regularity tubular result in
Appendix~\ref{append:finite_regularity_tubular_neighborhood} supplies
\(r_\theta>0\) for which normal addition is a \(\mathcal C^1\)
diffeomorphism from the radius-\(r_\theta\) normal disk bundle onto
\(\mathcal U_\theta\).  Its base point is the unique nearest-point projection
\(\pi_\theta\), and intrinsically
\[
 \mathrm dx
 =J_\theta(u,v)\,\mathrm dv\,\mathrm d\mathcal M_\theta(u),
 \qquad
 J_\theta(u,0)=1,\qquad
 |J_\theta(u,v)-1|\le C_\theta\|v\|.
\]
Here \(\mathrm dv\) is Euclidean measure on \(N_u\mathcal S\).  To perform
the fiber calculation, use a finite cover by local \(\mathcal C^1\)
orthonormal normal frames.  There
\(J_\theta(u,v)=\det(I-A_v)\), so boundedness of the second fundamental form
gives the displayed estimate.  On overlaps the frames differ orthogonally;
thus fiber measure is preserved, normal Hessians are conjugate, and their
determinant, Gaussian integral, and the weight \(w\) are invariant.  A finite
partition of unity assembles the estimates without orientability or a
trivial normal bundle.

For clarity, define the intrinsic fiber mass
\[
 A_\lambda(u)
 :=
 \int_{\substack{v\in N_u\mathcal S\\\|v\|<r_\theta}}
 e^{-[g(u+v)-g^\star]/\lambda}J_\theta(u,v)\,\mathrm dv.
\]
The rescaling \(v=\sqrt\lambda z\) used in the proof of
Hasenpflug et al.~\citep[Lemma~3.17, equations~(23)--(24)]
{hasenpflug2024wasserstein} now applies in each local frame.  The cubic
remainder and \(J_\theta(u,v)=1+O_\theta(\|v\|)\) contribute, after
rescaling, at most
\(C_\theta\sqrt\lambda(\|z\|+\|z\|^3)e^{-c_\theta\|z\|^2}\), and the
truncated Gaussian tail is \(O_\theta(e^{-c_\theta/\lambda})\).
Thus their calculation, with this first-order Jacobian remainder in place of
the smoother global-frame amplitude, yields
\[
\begin{aligned}
 \sup_{u\in\mathcal S}
 \left|\lambda^{-q/2}A_\lambda(u)-(2\pi)^{q/2}w(u)\right|
 &\le C_\theta\sqrt\lambda,\\
 \int_{\mathcal U_\theta}\|x-\pi_\theta(x)\|
 e^{-[g(x)-g^\star]/\lambda}\,\mathrm dx
 &\le C_\theta\lambda^{(q+1)/2}.
\end{aligned}
\]
Integrating the first estimate and adding the exponentially small complement
gives
\[
 \widehat Z_\lambda
 :=\int_{\mathbb R^d}e^{-[g(x)-g^\star]/\lambda}\,\mathrm dx
 =
 \lambda^{q/2}
 \left[
 (2\pi)^{q/2}\int_{\mathcal S}w\,\mathrm d\mathcal M_\theta
 +O_\theta(\sqrt\lambda)
 \right],
\]
and hence \(\widehat Z_\lambda\ge c_\theta\lambda^{q/2}\) for all sufficiently
small \(\lambda\).  This is the one-manifold version of the normalizer step
in their Remark~3.18, equations~(25)--(26).  The normalized fiber-density
comparison and the compactness of \(\mathcal S\), exactly as in their
Lemma~3.19 and Lemma~A.23, then give, for
\(\sigma_\lambda:=
(\pi_\theta)_\#(\gibbs^\lambda|_{\mathcal U_\theta})/
\gibbs^\lambda(\mathcal U_\theta)\),
\[
 \mathbb W_1(\sigma_\lambda,\gibbs^0)
 \le C_\theta\sqrt\lambda.
\]

The source states \(n\in\mathbb N\), but throughout Lemmas~3.17 and~3.19 and
equations~(27)--(30), \(n\) occurs only in scalar powers, \(e^{-n\ell}\),
Gaussian rescalings, and \(e^{-(n-1)c}\).  Its only finite-integer minimum,
in equation~(26), extends a large-\(n\) bound to small \(n\) and is irrelevant
here.  Thus the proof holds for every real \(t\ge t_0\), or
\(0<\lambda\le\lambda_0:=t_0^{-1}\).

Finally, the coupling in the proof of their Theorem~3.6,
equations~(27)--(30), has an intrinsic one-manifold form.  If
\(a_\lambda:=\gibbs^\lambda(\mathcal U_\theta)\), then
\[
 \gamma_\lambda
 :=
 (\operatorname{Id},\pi_\theta)_\#
       (\gibbs^\lambda|_{\mathcal U_\theta})
 +(\gibbs^\lambda|_{\mathcal U_\theta^c})\otimes\sigma_\lambda
\]
has first marginal exactly \(\gibbs^\lambda\) and second marginal exactly
\(\sigma_\lambda\): the two contributions to the latter are
\(a_\lambda\sigma_\lambda\) and
\((1-a_\lambda)\sigma_\lambda\).  The normal first-moment estimate bounds its
tube cost by \(O_\theta(\sqrt\lambda)\) after normalization, and the weighted
off-tube estimate bounds its remaining cost exponentially.  Couple
\(\sigma_\lambda\) to \(\gibbs^0\) as in Lemma~A.23 and glue the two
couplings.  The resulting coupling has exactly the required marginals
\(\gibbs^\lambda\) and \(\gibbs^0\), and its cost is at most
\(K_0\sqrt\lambda\).  Hence
\[
 \mathbb W_1(\gibbs_\theta^\lambda,\gibbs_\theta^0)
 \le K_0\sqrt\lambda
 \qquad(0<\lambda\le\lambda_0).
\]
Both \(K_0\) and \(\lambda_0\) may depend on the fixed parameter.  No
uniform-in-\(\theta\) claim, no uniform Hessian-Lipschitz bound, and no later
uniform result are used.
\end{proof}

\subsection{Verification of \Cref{ex:sphere_family}}
\label[appendix]{append:verification_sphere_example}

Let
\[
    \begin{gathered}
        \rho(\theta):=\sqrt{\theta+1},
        \qquad
        s(x):=\sqrt{\|x\|^2+a^2},\\
        R(\theta):=\sqrt{\rho(\theta)^2+a^2},
        \qquad
        R_{\max}:=\sqrt{a^2+2}.
    \end{gathered}
\]
Because \(a>0\), the function in \Cref{ex:sphere_family} admits a
\(\mathcal C^\infty\) extension to an open neighborhood of
\(\Theta\times\mathbb R^d\). Its lower gradient and Hessian are
\[
    \partial_xg(\theta,x)
    =2\left(1-\frac{R(\theta)}{s(x)}\right)x,
    \qquad
    \partial_{xx}^2g(\theta,x)
    =2\left(1-\frac{R(\theta)}{s(x)}\right)I
      +\frac{2R(\theta)}{s(x)^3}xx^\top.
\]
Since \(s(x)\ge a\) and \(R(\theta)\le R_{\max}\),
\[
    \|\partial_{xx}^2g(\theta,x)\|_{\mathrm{op}}
    \le 2+\frac{4R_{\max}}{a}.
\]
Thus \(\partial_xg\) is globally Lipschitz in \(x\), uniformly in \(\theta\).
Differentiating the displayed Hessian and using
\(s(x)\ge a\), \(\|x\|/s(x)\le1\), and \(R(\theta)\le R_{\max}\) shows that
\[
    \sup_{\theta\in\Theta,\,x\in\mathbb R^d}
    \|D_x^3g(\theta,x)\|_{\mathrm{op}}<\infty.
\]
The mean-value formula
therefore gives the Hessian-Lipschitz bound in
\eqref{eq_uniform_normal_hessian_lipschitz}, in fact globally in \(x\).
Moreover,
\[
    \partial_\theta\partial_xg(\theta,x)
    =-\frac{x}{R(\theta)s(x)},
    \qquad
    \|\partial_\theta\partial_xg(\theta,x)\|
    \le\frac{1}{\sqrt{a^2+1}},
\]
which gives the required global Lipschitz bound in \(\theta\), uniformly in
\(x\). Finally,
\[
    \partial_\theta g(\theta,x)
    =1-\frac{s(x)}{R(\theta)},
    \qquad
    |\partial_\theta g(\theta,x)|
    \le 1+\frac{\|x\|+a}{\sqrt{a^2+1}},
\]
so the polynomial-growth condition holds with \(\mathsf n_g=1\).

The gradient formula shows that the only critical points are the solution
sphere and \(x=0\). The local P\L\ estimate and the negative-definite Hessian at
\(0\) are established directly in \Cref{ex:sphere_family}; hence
\Cref{mu_PL_assum} holds with one common local P\L\ constant.

For the uniform tail gradient error bound in
\Cref{assum_sample_comp}, the displayed gradient formula and
\(s(x)^2-R(\theta)^2=\|x\|^2-\rho(\theta)^2\) give
\[
    \|\partial_xg(\theta,x)\|
    =2\frac{\|x\|}{s(x)}
      \frac{\|x\|+\rho(\theta)}{s(x)+R(\theta)}
      \bigl|\|x\|-\rho(\theta)\bigr|.
\]
Moreover,
\(\bigl|\|x\|-\rho(\theta)\bigr|
=\operatorname{dist}(x,\mathcal S(\theta))\).  For \(\|x\|\ge1\),
\[
    \frac{\|x\|}{s(x)}
    =\frac{\|x\|}{\sqrt{\|x\|^2+a^2}}
    \ge\frac{1}{\sqrt{1+a^2}},
\]
while \(s(x)\le \|x\|+a\), \(R(\theta)\le\rho(\theta)+a\), and
\(\|x\|+\rho(\theta)\ge2\) imply
\[
    s(x)+R(\theta)
    \le \|x\|+\rho(\theta)+2a
    \le(1+a)\bigl(\|x\|+\rho(\theta)\bigr).
\]
Consequently,
\[
    \|\partial_xg(\theta,x)\|
    \ge
    \frac{2}{(1+a)\sqrt{1+a^2}}
    \operatorname{dist}\bigl(x,\mathcal S(\theta)\bigr)
    \qquad(\|x\|\ge1).
\]
Thus the example satisfies \eqref{eq_uniform_tail_gradient_error_bound} with
\[
    \mathsf D_{\mathsf{eb}}=1,
    \qquad
    \nu_{\mathsf{eb}}
    =\frac{2}{(1+a)\sqrt{1+a^2}}.
\]
The derived tail quadratic-growth estimate needed below then follows from
\eqref{eq_tail_EB_implies_QG}.
The second fundamental form of a sphere
of radius \(\rho(\theta)\) has norm \(1/\rho(\theta)\le1\), giving the required
uniform geometric bound.
Finally, \(f(\theta,x)=x_1\) is globally Lipschitz and
satisfies the polynomial-growth condition in \Cref{assum_sample_comp}.

\section{Proofs of Results in \Cref{sec:4_unif_props}}\label[appendix]{append_proofs_sec_4}
This appendix contains deferred proofs for results in \Cref{sec:4_unif_props} concerning the structure and regularity of the minimizer sets $\{\mathcal{S}(\theta)\}_{\theta\in\Theta}$.
\subsection{Proof of \Cref{prop_parametric_morse_bott_stability}}
\label[appendix]{append_proof_parametric_morse_bott_stability}

\begin{lem}[Critical-point dichotomy]\label{lem_critical_point_dichotomy}
Under \Cref{mu_PL_assum,assum_sample_comp}, for every $\theta\in\Theta$ and
$x\in\mathbb R^d$,
\[
    \partial_xg(\theta,x)=0
    \quad\Longrightarrow\quad
    x\in\mathcal S(\theta)
    \quad\text{or}\quad
    \partial_{xx}^2g(\theta,x)\prec0.
\]
\end{lem}

\begin{proof}

Fix $\theta\in\Theta$ and abbreviate $g_\theta:=g(\theta,\cdot)$.  Let
$\mathcal M_\theta$ be the collection of all local minima of $g_\theta$.
By \Cref{mu_PL_assum}, $g_\theta$ satisfies the \PLcirc\ condition, while the
bounded-local-minimum argument and the derived estimate
\eqref{eq_tail_EB_implies_QG} show that $g_\theta$ is coercive.  Hence
\Cref{prop_S_path_connected,prop_S_C_2_wihout} apply:
$\mathcal M_\theta$ is connected,
$\mathcal M_\theta=\mathcal S(\theta)$, and this set is a compact
$\mathcal C^2$ embedded manifold without boundary.

Now let $x$ be a critical point of $g_\theta$.  If $x$ belongs to one of the
local-P{\L} neighborhoods $\mathcal N(\mathcal M')$ from
\Cref{def_local_PL}, then the local P{\L} inequality at $x$ gives
\[
    g_\theta(x)-\min_{z\in\mathcal N(\mathcal M')}g_\theta(z)
    \le \frac{1}{2\mu}\|\nabla g_\theta(x)\|^2=0.
\]
Thus $x$ minimizes $g_\theta$ on an open neighborhood and is a local, hence
global, minimizer.  If $x$ lies outside the union of those neighborhoods, the
second clause of \Cref{def_mu_PL_assum} gives
$\nabla^2g_\theta(x)\prec0$.  This proves the dichotomy.
\end{proof}

\begin{lem}[Continuous dependence of contractive fixed points]
\label{lem_continuous_contractive_fixed_point}
Let $P$ be a metric space, let
$B:=\overline{\mathbb B}_q(0;R)\subset\mathbb R^q$, and let
$\mathcal T:P\times B\to B$ be continuous.  Suppose that there exists
$\kappa\in[0,1)$ such that
\[
    \|\mathcal T(p,z)-\mathcal T(p,z')\|
    \le\kappa\|z-z'\|
    \qquad
    (p\in P,\ z,z'\in B).
\]
Then, for every $p\in P$, the map $\mathcal T(p,\cdot)$ has a unique fixed
point $z(p)\in B$, and the map $p\mapsto z(p)$ is continuous.
\end{lem}

\begin{proof}
Existence and uniqueness follow from the Banach fixed-point theorem because
$B$ is complete.  Fix $p_0\in P$.  The fixed-point identities and the
contraction estimate give
\begin{align*}
    \|z(p)-z(p_0)\|
    &\le
    \|\mathcal T(p,z(p))-\mathcal T(p,z(p_0))\|\\
    &\quad+
    \|\mathcal T(p,z(p_0))-\mathcal T(p_0,z(p_0))\|\\
    &\le
    \kappa\|z(p)-z(p_0)\|
    +\|\mathcal T(p,z(p_0))-\mathcal T(p_0,z(p_0))\|.
\end{align*}
Therefore,
\[
    \|z(p)-z(p_0)\|
    \le
    \frac{
        \|\mathcal T(p,z(p_0))-\mathcal T(p_0,z(p_0))\|
    }{1-\kappa}.
\]
The right-hand side converges to zero as $p\to p_0$ by continuity of
$\mathcal T$.  Hence $z(p)\to z(p_0)$.
\end{proof}

\begin{proof}[Proof of \Cref{prop_parametric_morse_bott_stability}]
Fix $\bar\theta\in\Theta$, put
$M:=\mathcal S(\bar\theta)$, and let $k:=\dim M$.  Throughout this proof,
parameter neighborhoods are understood relative to $\Theta$: they have the
form $\Theta\cap O$, where $O\subseteq\mathbb R^m$ is open.  Likewise,
$\theta\to\bar\theta$ always means that $\theta\in\Theta$ and
$\|\theta-\bar\theta\|\to0$.

The proof has two parts.  Steps 1--4 construct a $\mathcal C^1$-small normal
graph over $M$ whose points satisfy the normal criticality equation.  Steps
5--7 identify this graph with $\mathcal S(\theta)$ and establish the global
boundedness and compactness conclusions.  Whenever necessary, we replace
$U_{\bar\theta}$ by a smaller relative neighborhood without changing its
notation.

\emph{Step 1: Normal curvature at the reference manifold.}
By the local
P{\L}--Morse--Bott equivalence for $\mathcal C^2$ functions
\citep[Corollary~2.17]{rebjock2024fast}, for every $u\in M$,
\begin{equation}\label{eq_append_MB_kernel}
    \ker H_{\bar\theta}(u)=T_uM,
    \qquad
    H_{\bar\theta}(u)|_{N_uM}\succ0,
    \qquad
    H_\theta(x):=\partial_{xx}^2g(\theta,x).
\end{equation}
More precisely, the implication from $\mu$-P{\L} to $\mu$-Morse--Bott in
\citep[Corollary~2.17]{rebjock2024fast} gives the quantitative bound
\begin{equation}\label{eq_append_normal_gap_at_base}
    H_{\bar\theta}(u)|_{N_uM}\succeq\mu I
    \qquad\forall u\in M.
\end{equation}
For completeness, the constant can also be verified directly.  If
$v\in N_uM$ is a unit eigenvector with eigenvalue $\lambda>0$,
then, as $t\to0$,
\[
g(\bar\theta,u+tv)-g^\star(\bar\theta)
=\tfrac12\lambda t^2+o(t^2),
\qquad
\|\partial_xg(\bar\theta,u+tv)\|^2
=\lambda^2t^2+o(t^2).
\]
Substitution in the local P{\L} inequality and division by $t^2/2$ give
$\lambda\le\lambda^2/\mu$, and therefore $\lambda\ge\mu$.

\emph{Step 2: Uniform solution of the normal equation.}
Use the normal bundle $NM$ and the tubular map
$\tau(u,n):=u+n$ from \eqref{local coordinate decomposition}.  On its domain,
define
\[
    \widetilde g_\theta(u,n):=g(\theta,\tau(u,n)),
    \qquad
    \mathcal F_\theta(u,n)
    :=
    P_{N_uM}\partial_xg(\theta,\tau(u,n))
    \in N_uM,
\]
where $P_{N_uM}$ is the orthogonal projection onto $N_uM$.  Thus
$\mathcal F_\theta(u,n)$ is the component of the gradient in the normal space
at $u$.  Equivalently, for every $\eta\in N_uM$,
\[
    D_n\widetilde g_\theta(u,n)[\eta]
    =
    \left\langle\mathcal F_\theta(u,n),\eta\right\rangle.
\]
Hence $\mathcal F_\theta(u,n)=0$ exactly when the directional derivative of
$\widetilde g_\theta$ vanishes in every normal direction.  Moreover,
$\mathcal F_{\bar\theta}(u,0)=0$, and
\[
    D_n\mathcal F_{\bar\theta}(u,0)
    =
    H_{\bar\theta}(u)\big|_{N_uM}.
\]

We now solve $\mathcal F_\theta(u,n)=0$ simultaneously for every $u\in M$ and
every $\theta$ sufficiently close to $\bar\theta$.  Choose finitely many
coordinate neighborhoods $\mathcal O_j$ of $M$, each equipped with a local
orthonormal normal frame, and compact sets $K_j\subset\mathcal O_j$ whose
relative interiors cover $M$.  By \eqref{eq_append_normal_gap_at_base}, the
normal Hessian is invertible and
\begin{equation}\label{eq_append_normal_inverse_bound}
    \left\|
        \left(H_{\bar\theta}(u)\big|_{N_uM}\right)^{-1}
    \right\|_{\mathrm{op}}
    \le\frac{1}{\mu}.
\end{equation}
For a fixed $j$, the identity
$D_n\mathcal F_{\bar\theta}(u,0)
=H_{\bar\theta}(u)|_{N_uM}$ holds for every $u\in K_j$.  Joint continuity of
$D_n\mathcal F_\theta(u,n)$ and compactness of $K_j$ therefore give a normal-fiber
radius $r_{\mathrm n,j}>0$ and a relative neighborhood
$U_j\subseteq\Theta$ of $\bar\theta$ such
that
\begin{equation}\label{eq_append_chart_derivative_bound}
    \left\|
        D_n\mathcal F_\theta(u,n)
        -H_{\bar\theta}(u)\big|_{N_uM}
    \right\|_{\mathrm{op}}
    \le\frac{\mu}{2}
\end{equation}
whenever $u\in K_j$, $\theta\in U_j$, and $n\in N_uM$ satisfies
$\|n\|\le r_{\mathrm n,j}$.  The normal-fiber radius is also chosen small enough
that all these points
belong to the domain of the tubular map.

Because the cover is finite, the choices
\begin{equation}\label{eq_append_common_normal_radius_and_parameter_neighborhood}
    r_{\mathrm n}:=\min_j r_{\mathrm n,j}>0,
    \qquad
    U_{\bar\theta}:=\bigcap_j U_j
\end{equation}
make \eqref{eq_append_chart_derivative_bound} valid on every $K_j$ whenever
$\theta\in U_{\bar\theta}$ and $\|n\|\le r_{\mathrm n}$.  We next control the residual at
$n=0$.  Since $u\in M=\mathcal S(\bar\theta)$,
$\partial_xg(\bar\theta,u)=0$.  Therefore, the definition of
$\mathcal F_\theta$, the norm-one property of the orthogonal projection, and
the Lipschitz dependence of $\partial_xg$ on $\theta$ in
\Cref{assum_sample_comp} give
\[
\begin{aligned}
    \|\mathcal F_\theta(u,0)\|
    &=\|\mathcal F_\theta(u,0)-\mathcal F_{\bar\theta}(u,0)\|\\
    &=\left\|
        P_{N_uM}\bigl(
            \partial_xg(\theta,u)-\partial_xg(\bar\theta,u)
        \bigr)
      \right\|\\
    &\le L_{g,\theta x}\|\theta-\bar\theta\|.
\end{aligned}
\]
If $L_{g,\theta x}>0$, define the parameter-space radius
\[
    \delta_{\bar\theta}
    :=\frac{\mu r_{\mathrm n}}{4L_{g,\theta x}}>0
\]
and replace $U_{\bar\theta}$ by its intersection with
$\{\theta\in\Theta:\|\theta-\bar\theta\|<\delta_{\bar\theta}\}$.  If
$L_{g,\theta x}=0$, no further shrinking is needed.  In either case,
\begin{equation}\label{eq_append_base_residual_bound}
    \|\mathcal F_\theta(u,0)\|\le\frac{\mu r_{\mathrm n}}{4}
    \qquad
    (u\in M,\ \theta\in U_{\bar\theta}).
\end{equation}
For $u\in M$ and $\theta\in U_{\bar\theta}$, define
\[
    T_{\theta,u}(n)
    :=n-
    \left(H_{\bar\theta}(u)\big|_{N_uM}\right)^{-1}
    \mathcal F_\theta(u,n)
\]
on $\{n\in N_uM:\|n\|\le r_{\mathrm n}\}$.  A vector $n$ is a fixed point of
$T_{\theta,u}$ exactly when $\mathcal F_\theta(u,n)=0$.  We introduce this map
to solve the latter equation by the Banach fixed-point theorem; multiplication
by the inverse normal Hessian makes its derivative small near
$(\bar\theta,u,0)$.  Indeed, its derivative with respect to $n$
satisfies
\[
    D_nT_{\theta,u}(n)
    =
    \left(H_{\bar\theta}(u)\big|_{N_uM}\right)^{-1}
    \left[
        H_{\bar\theta}(u)\big|_{N_uM}
        -D_n\mathcal F_\theta(u,n)
    \right],
\]
and hence \eqref{eq_append_normal_inverse_bound} and
\eqref{eq_append_chart_derivative_bound} give
$\|D_nT_{\theta,u}(n)\|_{\mathrm{op}}\le1/2$.  Likewise,
\eqref{eq_append_normal_inverse_bound} and
\eqref{eq_append_base_residual_bound} give
\begin{align*}
    \|T_{\theta,u}(0)\|
    &\le\frac{r_{\mathrm n}}{4},\\
    \|T_{\theta,u}(n)\|
    &\le\|T_{\theta,u}(n)-T_{\theta,u}(0)\|
       +\|T_{\theta,u}(0)\|\\
    &\le\frac12\|n\|+\frac{r_{\mathrm n}}{4}
     \le\frac{3r_{\mathrm n}}{4}
\end{align*}
whenever $\|n\|\le r_{\mathrm n}$.  Thus $T_{\theta,u}$ is a $1/2$-contraction
mapping $\{n\in N_uM:\|n\|\le r_{\mathrm n}\}$ into itself.
In each local orthonormal normal frame, $T_{\theta,u}$ is represented by a
jointly continuous family of self-maps of the fixed ball
$\overline{\mathbb B}_{d-k}(0;r_{\mathrm n})$, with common contraction factor
$1/2$.  Therefore, \Cref{lem_continuous_contractive_fixed_point} gives a unique
fixed point $n_\theta(u)$ whose coefficient vector in that frame is continuous
in $(\theta,u)$.  Thus $(\theta,u)\mapsto n_\theta(u)$ is a continuous map into
the normal bundle $NM$ and satisfies
\begin{equation}\label{eq_append_normal_equation}
    \mathcal F_\theta(u,n_\theta(u))=0
    \qquad
    (u\in M,\ \theta\in U_{\bar\theta}).
\end{equation}
At $\theta=\bar\theta$, the point $n=0$ is a fixed point; uniqueness therefore
gives $n_{\bar\theta}(u)=0$ for every $u\in M$.

\emph{Step 3: $\mathcal C^1$-regularity and convergence of the section.}
To establish $\mathcal C^1$-regularity, fix one of the coordinate
neighborhoods above, choose a $\mathcal C^2$ parametrization
$\mathsf S_j:\Gamma_j\to\mathcal O_j$, and set
$A_j:=\mathsf S_j^{-1}(K_j)$.  Thus $A_j$ is compact.  For
$a\in\Gamma_j$, write $u=\mathsf S_j(a)$ and let
\[
    E_j(a)
    :=
    \begin{bmatrix}
        e_j^{(1)}(\mathsf S_j(a))&\cdots&
        e_j^{(d-k)}(\mathsf S_j(a))
    \end{bmatrix},
\]
and represent the fixed point as
\[
    n_\theta(\mathsf S_j(a))
    =E_j(a)t_{\theta,j}(a),
    \qquad
    t_{\theta,j}(a)\in\mathbb R^{d-k}.
\]
In these coordinates, \eqref{eq_append_normal_equation} becomes
\[
    \mathbf F_j(\theta,a,t_{\theta,j}(a))=0,
    \qquad
    \mathbf F_j(\theta,a,t)
    :=E_j(a)^\top
      \partial_xg\bigl(\theta,\mathsf S_j(a)+E_j(a)t\bigr).
\]
The $\mathcal C^2$ extension of $g$, the $\mathcal C^2$ regularity of
$\mathsf S_j$, and the $\mathcal C^1$ regularity of $E_j$ make
$\mathbf F_j$ a $\mathcal C^1$ map on an open set in $(\theta,a,t)$.  Moreover,
\[
    D_t\mathbf F_j(\theta,a,t)
    =E_j(a)^\top
      H_\theta\bigl(\mathsf S_j(a)+E_j(a)t\bigr)E_j(a).
\]
At $t=t_{\theta,j}(a)$, \eqref{eq_append_normal_gap_at_base} and
\eqref{eq_append_chart_derivative_bound} therefore give
\[
    D_t\mathbf F_j(\theta,a,t_{\theta,j}(a))
    \succeq\frac{\mu}{2}I_{d-k},
    \qquad
    \left\|
        \bigl[D_t\mathbf F_j(\theta,a,t_{\theta,j}(a))\bigr]^{-1}
    \right\|_{\mathrm{op}}
    \le\frac{2}{\mu}.
\]
Thus the derivative with respect to $t$ is invertible.  The implicit-function
theorem, applied at each solution, gives a local $\mathcal C^1$ solution map
$(\theta,a)\mapsto t_{\theta,j}(a)$.  Its derivative with respect to the
manifold coordinate is
\begin{equation}\label{eq_append_section_derivative}
    D_at_{\theta,j}(a)
    =-\bigl[
        D_t\mathbf F_j(\theta,a,t_{\theta,j}(a))
      \bigr]^{-1}
      D_a\mathbf F_j(\theta,a,t_{\theta,j}(a)).
\end{equation}

The vector $n_\theta(u)\in N_uM$ was already defined, independently of any
frame, as the unique fixed point of $T_{\theta,u}$.  The maps
$t_{\theta,j}$ are only its coordinate representations.  Hence, on the overlap
of two coordinate neighborhoods, both representations describe the same
vector $n_\theta(u)$.  The local $\mathcal C^1$ representations therefore show
that $u\mapsto n_\theta(u)$ is a globally defined $\mathcal C^1$ section of
$NM$.

We next prove the quantitative convergence
\begin{equation}\label{eq_parametric_normal_graph_C1_convergence}
    \|n_\theta\|_{\mathcal C^1(M;\mathbb R^d)}
    :=
    \sup_{u\in M}\|n_\theta(u)\|
    +
    \sup_{\substack{u\in M,\ v\in T_uM\\ \|v\|=1}}
        \|Dn_\theta(u)[v]\|
    \longrightarrow0
    \qquad(\theta\to\bar\theta).
\end{equation}
We first control the first supremum.  Since $n_{\bar\theta}=0$, we have
$t_{\bar\theta,j}(a)=0$ for every $a\in A_j$.  Joint continuity of
$(\theta,a)\mapsto t_{\theta,j}(a)$ and compactness of $A_j$ imply
\[
    \sup_{a\in A_j}\|t_{\theta,j}(a)\|\longrightarrow0
    \qquad\text{as }\theta\to\bar\theta.
\]
Indeed, otherwise there would exist $\varepsilon>0$, a sequence
$\theta_\ell\to\bar\theta$, and points $a_\ell\in A_j$ such that
$\|t_{\theta_\ell,j}(a_\ell)\|\ge\varepsilon$.  Compactness gives a subsequence
with $a_\ell\to a_\star\in A_j$, whereas joint continuity gives
$t_{\theta_\ell,j}(a_\ell)\to t_{\bar\theta,j}(a_\star)=0$, a contradiction.
The cover is finite, so this convergence holds simultaneously for every $j$.
Given $n_\theta(u)
    =E_j(a)t_{\theta,j}(a)$, since the columns of $E_j(a)$ are orthonormal and the sets $K_j$ cover $M$, it
follows that
\begin{equation}\label{eq_append_section_C0_convergence}
    \sup_{u\in M}\|n_\theta(u)\|\longrightarrow0.
\end{equation}

We next control the derivative term.  The identity
$\mathbf F_j(\bar\theta,a,0)=0$ for every $a$ gives
$D_a\mathbf F_j(\bar\theta,a,0)=0$.  The uniform convergence just proved,
joint continuity of $D_a\mathbf F_j$, and compactness of $A_j$ imply
\[
    \sup_{a\in A_j}
    \bigl\|D_a\mathbf F_j(\theta,a,t_{\theta,j}(a))\bigr\|_{\mathrm{op}}
    \longrightarrow0.
\]
The uniform inverse bound and \eqref{eq_append_section_derivative} therefore
give
\[
    \sup_{a\in A_j}\|D_at_{\theta,j}(a)\|_{\mathrm{op}}
    \longrightarrow0.
\]
Moreover,
\[
    D_a\bigl[n_\theta\circ\mathsf S_j\bigr](a)
    =D_aE_j(a)t_{\theta,j}(a)+E_j(a)D_at_{\theta,j}(a).
\]
The frames and their first derivatives, as well as the inverses of the chart
derivatives onto the tangent spaces, are bounded on the finitely many sets
$A_j$.  Consequently,
\begin{equation}\label{eq_append_section_derivative_convergence}
    \sup_{\substack{u\in M,\ v\in T_uM\\ \|v\|=1}}
        \|Dn_\theta(u)[v]\|
    \longrightarrow0.
\end{equation}
Combining \eqref{eq_append_section_C0_convergence} and
\eqref{eq_append_section_derivative_convergence}
proves \eqref{eq_parametric_normal_graph_C1_convergence}.

\emph{Step 4: Geometry of the normal graph.}
Recall from
\eqref{eq_append_common_normal_radius_and_parameter_neighborhood} that
$U_{\bar\theta}$ is a relative parameter neighborhood on which the common
normal-fiber radius $r_{\mathrm n}$ is valid; the intervening steps have only
shrunk this neighborhood.  Choose $\rho_{\mathrm n}\in(0,r_{\mathrm n})$.
By \eqref{eq_append_section_C0_convergence} and
\eqref{eq_append_section_derivative_convergence}, after shrinking
$U_{\bar\theta}$ once more, we may assume that
\begin{equation}\label{eq_append_section_smallness}
    \sup_{u\in M}\|n_\theta(u)\|<\rho_{\mathrm n},
    \qquad
    \sup_{\substack{u\in M,\ v\in T_uM\\ \|v\|=1}}
        \|Dn_\theta(u)[v]\|<\frac12
    \qquad(\theta\in U_{\bar\theta}).
\end{equation}
The first inequality keeps the graph
$\{(u,n_\theta(u)):u\in M\}$ inside the domain of the tubular
diffeomorphism.  The graph map $u\mapsto(u,n_\theta(u))$ is a $\mathcal C^1$
embedding into $NM$, and composition with the tubular diffeomorphism shows
that
\[
    \Phi_\theta:M\to\mathbb R^d,
    \qquad \Phi_\theta(u):=u+n_\theta(u),
\]
is an embedding.  For $v\in T_uM$, the second inequality in
\eqref{eq_append_section_smallness} gives
\[
\begin{aligned}
    \bigl\|P_{T_uM}D\Phi_\theta(u)[v]\bigr\|
    &=\bigl\|v+P_{T_uM}Dn_\theta(u)[v]\bigr\|\\
    &\ge\frac12\|v\|.
\end{aligned}
\]
Therefore
$D\Phi_\theta(u)[T_uM]\cap N_uM=\{0\}$.
To see why the intersection is trivial, let
$w\in D\Phi_\theta(u)[T_uM]\cap N_uM$.  Then
$w=D\Phi_\theta(u)[v]$ for some $v\in T_uM$, while $w\in N_uM$ gives
$P_{T_uM}w=0$.  The preceding inequality yields
\[
    0
    =\bigl\|P_{T_uM}w\bigr\|
    =\bigl\|P_{T_uM}D\Phi_\theta(u)[v]\bigr\|
    \ge\frac12\|v\|,
\]
so $v=0$ and hence $w=0$. These spaces have dimensions $k$ and
$d-k$, respectively, so their direct sum is $\mathbb R^d$, which gives the
first relation below.  The second follows from
\eqref{eq_append_normal_gap_at_base} and
\eqref{eq_append_chart_derivative_bound}, evaluated at
$n=n_\theta(u)$:
\begin{equation}\label{eq_append_transverse_splitting}
    D\Phi_\theta(u)[T_uM]\oplus N_uM=\mathbb R^d,
    \qquad
    z^\top H_\theta(\Phi_\theta(u))z
       \ge\frac{\mu}{2}\|z\|^2
       \quad(z\in N_uM).
\end{equation}
Finally, joint continuity of $H_\theta(x)$, compactness of $M$, and
\eqref{eq_append_section_C0_convergence} allow us to shrink
$U_{\bar\theta}$ again so that
\begin{equation}\label{eq_append_graph_hessian_comparison}
    \bigl\|H_\theta(\Phi_\theta(u))-H_{\bar\theta}(u)\bigr\|_{\mathrm{op}}
    <\frac{\mu}{4}
    \qquad(u\in M,\ \theta\in U_{\bar\theta}).
\end{equation}

\emph{Step 5: The constructed graph consists of global minimizers.}
We first prove that $\Phi_\theta(M)\subseteq\mathcal S(\theta)$.  Define
$h_\theta:=g(\theta,\cdot)\circ\Phi_\theta$ on compact $M$, and let
$u_\star$ maximize $h_\theta$.
Since $M$ has no boundary, the first-order necessary condition gives
\[
    Dh_\theta(u_\star)[v]=0
    \qquad\forall v\in T_{u_\star}M.
\]
By the chain rule and the definition of $h_\theta$, this means
\[
    \left\langle
        \partial_xg(\theta,\Phi_\theta(u_\star)),
        D\Phi_\theta(u_\star)[v]
    \right\rangle
    =0
    \qquad\forall v\in T_{u_\star}M.
\]
Thus $\partial_xg(\theta,\Phi_\theta(u_\star))$ is orthogonal to
$D\Phi_\theta(u_\star)[T_{u_\star}M]$.  On the other hand,
\eqref{eq_append_normal_equation} and the definition of $\mathcal F_\theta$
give
\[
    P_{N_{u_\star}M}
    \partial_xg(\theta,\Phi_\theta(u_\star))=0,
\]
so the same gradient is orthogonal to $N_{u_\star}M$.  The direct-sum relation
in \eqref{eq_append_transverse_splitting} shows that these two subspaces span
$\mathbb R^d$.  Therefore,
\[
    \partial_xg(\theta,\Phi_\theta(u_\star))=0.
\]
Here $q:=d-k=d-\dim M\ge1$.  Indeed, a nonempty full-dimensional embedded submanifold of
$\mathbb R^d$ is open in $\mathbb R^d$ by the inverse-function theorem, whereas
$M$ is nonempty and compact; hence $k=d$ is impossible.  Choose a nonzero
$z\in N_{u_\star}M$.  The curvature estimate in
\eqref{eq_append_transverse_splitting} gives
\[
    z^\top H_\theta(\Phi_\theta(u_\star))z
    \ge\frac{\mu}{2}\|z\|^2>0.
\]
Thus the Hessian at this critical point is not negative definite.  The
critical-point dichotomy in \Cref{lem_critical_point_dichotomy} therefore
implies that $\Phi_\theta(u_\star)$ is a global minimizer.
Hence
\[
    g^\star(\theta)\le h_\theta(u)
    \le h_\theta(u_\star)=g^\star(\theta)
    \qquad\forall u\in M,
\]
so $h_\theta\equiv g^\star(\theta)$.  In particular,
$\Phi_\theta(M)\subseteq\mathcal S(\theta)$.  Since $h_\theta$ is constant, the
preceding tangent-plus-normal argument applies at every $u\in M$ and gives
$\partial_xg(\theta,\Phi_\theta(u))=0$.

\emph{Step 6: The graph is the entire minimizer set.}
We now prove the reverse inclusion.
Differentiating this last identity along $M$ yields
\begin{equation}\label{eq_append_tangent_in_kernel}
    D\Phi_\theta(u)[T_uM]
    \subseteq\ker H_\theta(\Phi_\theta(u)).
\end{equation}
The left-hand side has dimension $k$, whereas positivity of the Hessian on the
$(d-k)$-dimensional space $N_uM$ in
\eqref{eq_append_transverse_splitting} implies that the kernel has dimension at
most $k$.  Thus equality holds in \eqref{eq_append_tangent_in_kernel}.  Since
$\Phi_\theta(u)$ is a minimizer, the Morse--Bott characterization also gives
\[
    T_{\Phi_\theta(u)}\mathcal S(\theta)
    =\ker H_\theta(\Phi_\theta(u))
    =D\Phi_\theta(u)[T_uM].
\]
Put $y:=\Phi_\theta(u)$.  Since $\Phi_\theta$ is an embedding,
\[
    T_y\Phi_\theta(M)
    =D\Phi_\theta(u)[T_uM]
    =T_{\Phi_\theta(u)}\mathcal S(\theta).
\]
Thus the differential of the inclusion
$\iota:\Phi_\theta(M)\hookrightarrow\mathcal S(\theta)$ is the identity
isomorphism from $T_y\Phi_\theta(M)$ to $T_y\mathcal S(\theta)$.  The
inverse-function theorem on manifolds shows that $\iota$ is a local
diffeomorphism at every point.  Hence $\Phi_\theta(M)$ is relatively open in
$\mathcal S(\theta)$.
It is
compact, hence closed in $\mathcal S(\theta)$, and it is nonempty.  The set
$\mathcal S(\theta)$ is connected by \Cref{prop_S_path_connected}; consequently,
\begin{equation}\label{eq_append_solution_normal_graph}
    \mathcal S(\theta)=\Phi_\theta(M)
    \qquad\forall\theta\in U_{\bar\theta}.
\end{equation}
Thus the local assertion of the proposition holds with
$V_{\bar\theta}:=U_{\bar\theta}$.

\emph{Step 7: Uniform boundedness and compactness.}
The neighborhoods $U_{\bar\theta}$ cover compact $\Theta$.  Choose finitely many
of them, with base manifolds $M_i=\mathcal S(\theta_i)$, and use the finite
shrinking lemma to choose relative open sets $V_i$ which still cover $\Theta$ and
satisfy $\overline V_i\subset U_{\theta_i}$.  Formula
\eqref{eq_append_solution_normal_graph}, compactness of
each $M_i$, and continuity of the finitely many sections give a constant
$\mathsf D_{\mathcal S}$ with
$\mathcal S(\theta)\subseteq\mathbb B_d(0;\mathsf D_{\mathcal S})$ for all
$\theta$.  Finally, the
solution graph \(\mathfrak G\) defined in
\eqref{eq_solution_graph} is closed: if
$\theta_j\to\theta$, $x_j\in\mathcal S(\theta_j)$, and $x_j\to x$, then for every
$y\in\mathbb R^d$,
\[
    g(\theta_j,x_j)\le g(\theta_j,y)
    \quad\Longrightarrow\quad
    g(\theta,x)\le g(\theta,y).
\]
Thus $x\in\mathcal S(\theta)$.  The graph is closed and bounded in
$\Theta\times\mathbb R^d$, and is therefore compact.
\end{proof}

\subsection{Proof of \Cref{lem_uniform_local_PL_band}}
\label[appendix]{append_proof_uniform_local_PL_band}
\begin{proof}
We first construct a uniform neighborhood of the solution graph \(\mathfrak G\) defined in \eqref{eq_solution_graph} on which the
two estimates hold, and then show that one common sublevel set lies in that
neighborhood.

\emph{Uniform transverse tube.}
For each reference parameter $\vartheta\in\Theta$,
\Cref{prop_parametric_morse_bott_stability} provides a relatively open
neighborhood $U_\vartheta$ of $\vartheta$ such that, for every
$\theta\in U_\vartheta$, the manifold $\mathcal S(\theta)$ is a normal graph
over $\mathcal S(\vartheta)$.  Since the sets $U_\vartheta$ cover the compact
set $\Theta$, finitely many reference parameters
$\theta_1,\ldots,\theta_N\in\Theta$ can be selected so that
\[
    \Theta=\bigcup_{i=1}^N U_i,
    \qquad
    U_i:=U_{\theta_i}.
\]
We may then choose relatively open sets $V_i\subseteq\Theta$ such that
\[
    \Theta=\bigcup_{i=1}^N V_i,
    \qquad
    \overline V_i\subset U_i,
\]
where the closure is taken in $\Theta$.  The smaller sets $V_i$ will allow us
to obtain uniform estimates on their compact closures.

For each $i$, set $M_i:=\mathcal S(\theta_i)$.  The same proposition gives, for
every $\theta\in U_i$, a normal section $n_{i,\theta}$ satisfying
$n_{i,\theta}(u)\in N_uM_i$.  Define
\[
    \Phi_{i,\theta}(u):=u+n_{i,\theta}(u),
    \qquad u\in M_i,
\]
so that $\Phi_{i,\theta}(M_i)=\mathcal S(\theta)$.  The construction also gives
\begin{equation}\label{eq_uniform_band_transverse_splitting}
    D\Phi_{i,\theta}(u)[T_uM_i]\oplus N_uM_i=\mathbb R^d.
\end{equation}
Set $H_\theta(x):=\partial_{xx}^2g(\theta,x)$.  The normal-Hessian estimate in
\Cref{sec:preliminary} gives
\begin{equation}\label{eq_uniform_band_base_gap}
    z^\top H_{\theta_i}(u)z\ge\mu\|z\|^2
    \qquad (u\in M_i,\ z\in N_uM_i).
\end{equation}
Applied with reference parameter $\bar\theta=\theta_i$, estimate
\eqref{eq_append_graph_hessian_comparison} from the proof of
\Cref{prop_parametric_morse_bott_stability} states that, after choosing $U_i$
small enough, the Hessian at the nearby minimizer
$\Phi_{i,\theta}(u)\in\mathcal S(\theta)$ is uniformly close to the Hessian at
the corresponding reference minimizer $u\in M_i=\mathcal S(\theta_i)$.  Since
$\overline V_i\subset U_i$, this gives
\begin{equation}\label{eq_uniform_band_graph_base_comparison}
    \bigl\|H_\theta(\Phi_{i,\theta}(u))-H_{\theta_i}(u)\bigr\|_{\mathrm{op}}
    \le\frac{\mu}{4}
    \qquad(\theta\in\overline V_i,\ u\in M_i).
\end{equation}
To control normal displacements, consider
\[
    \mathcal D_i
    :=\left\{
        (\theta,u,w):
        \theta\in\overline V_i,\ u\in M_i,\ w\in N_uM_i,\
        \|w\|\le1
    \right\}.
    \]
    The set \(\overline V_i\times M_i\) is compact, and \(\mathcal D_i\) is the
    closed unit-disk bundle of the pullback normal bundle over this base;
    hence \(\mathcal D_i\) is compact.  No global normal frame is used.  On the
    explicitly parameterized compact set
    \(\mathcal D_i\times[0,1]_t\times[0,1]_s\), consider
    \[
        (\theta,u,w,t,s)
        \longmapsto
        \bigl\|H_\theta(\Phi_{i,\theta}(u)+st w)
                 -H_\theta(\Phi_{i,\theta}(u))\bigr\|_{\mathrm{op}}
    \]
    which is continuous and vanishes when \(s=0\).  Uniform continuity therefore permits
    the choice of $r_i\in(0,1]$ such that
\[
    \bigl\|H_\theta(\Phi_{i,\theta}(u)+tz)
             -H_\theta(\Phi_{i,\theta}(u))\bigr\|_{\mathrm{op}}
    \le\frac{\mu}{4}
\]
whenever $\theta\in\overline V_i$, $u\in M_i$, $z\in N_uM_i$,
    $\|z\|\le r_i$, and $t\in[0,1]$; for $z\ne0$, take
    $s=\|z\|$ and $w=z/\|z\|$, while $z=0$ is immediate.  Adding this estimate to
    \eqref{eq_uniform_band_graph_base_comparison} yields, for the same variables,
\begin{equation}\label{eq_uniform_band_hessian_comparison}
    \bigl\|H_\theta\bigl(\Phi_{i,\theta}(u)+tz\bigr)
          -H_{\theta_i}(u)\bigr\|_{\mathrm{op}}
    \le\frac{\mu}{2}.
\end{equation}
Any subsequent decrease of $r_i$ preserves this estimate.
Consequently, \eqref{eq_uniform_band_base_gap} implies
\begin{equation}\label{eq_uniform_band_transverse_gap}
    z^\top H_\theta\bigl(\Phi_{i,\theta}(u)+tz\bigr)z
    \ge\frac{\mu}{2}\|z\|^2.
\end{equation}

For each $i$, define
\begin{align*}
\mathcal E_i
&:=\bigl\{(\theta,u,z):\theta\in V_i,\ u\in M_i,\
z\in N_uM_i,\ \|z\|<r_i\bigr\},\\
\Psi_i(\theta,u,z)
&:=\bigl(\theta,\Phi_{i,\theta}(u)+z\bigr).
\end{align*}
We will show that
\begin{equation}\label{eq_uniform_tube_definition}
    \mathcal W:=\bigcup_{i=1}^N\Psi_i(\mathcal E_i)
\end{equation}
is a relatively open neighborhood of the solution graph $\mathfrak G$.  We
first prove relative openness and then verify that
$\mathfrak G\subseteq\mathcal W$.

Fix $i$.  For each $\theta\in V_i$, consider the spatial map
\[
    \psi_{i,\theta}(u,z):=\Phi_{i,\theta}(u)+z,
    \qquad (u,z)\in NM_i,\quad \|z\|<r_i.
\]
At $(u,0)$, the tangent space of the normal bundle has the canonical
identification $T_{(u,0)}NM_i=T_uM_i\oplus N_uM_i$, and
\[
    D\psi_{i,\theta}(u,0)[v,w]
    =D\Phi_{i,\theta}(u)[v]+w.
\]
The first term ranges over $D\Phi_{i,\theta}(u)[T_uM_i]$, and
\eqref{eq_uniform_band_transverse_splitting} shows that this space and $N_uM_i$
form a direct sum equal to $\mathbb R^d$.  Hence this derivative is a linear
isomorphism.

Choose finitely many local $\mathcal C^2$ charts
$\mathsf S_{i,\ell}:\Gamma_{i,\ell}\to M_i$ with local $\mathcal C^1$
orthonormal normal frames $E_{i,\ell}(a)$, and compact sets
$A_{i,\ell}\Subset\Gamma_{i,\ell}$ such that the sets
$\mathsf S_{i,\ell}(\operatorname{int}A_{i,\ell})$ cover $M_i$.  In these
coordinates, define
\[
\sigma_{i,\ell}(\theta,a,t)
:=\sigma_{\min}\!\left(
    D_{(a,t)}
    [\Phi_{i,\theta}(\mathsf S_{i,\ell}(a))+E_{i,\ell}(a)t]
\right).
\]
This function is continuous and is positive at $t=0$ by the preceding
direct-sum calculation.  Compactness gives
\[
    c_{i,\ell}
    :=\min_{(\theta,a)\in\overline V_i\times A_{i,\ell}}
      \sigma_{i,\ell}(\theta,a,0)>0.
\]
Continuity and compactness of this zero section give
$\rho_{i,\ell}>0$ such that
$\sigma_{i,\ell}(\theta,a,t)\ge c_{i,\ell}/2$ whenever
$(\theta,a)\in\overline V_i\times A_{i,\ell}$ and
$\|t\|<\rho_{i,\ell}$.  Taking the minimum over the finitely many charts and
decreasing $r_i$ accordingly makes the spatial derivative invertible for every
$\theta\in V_i$, $u\in M_i$, and $z\in N_uM_i$ with $\|z\|<r_i$; here
$\|E_{i,\ell}(a)t\|=\|t\|$ because the frame is orthonormal.  This decrease preserves
\eqref{eq_uniform_band_hessian_comparison}.

We now include the parameter variable in the inverse-function argument.  Near
any $(\theta,u,z)\in\mathcal E_i$, use one of the preceding charts and write
$u=\mathsf S_{i,\ell}(a)$ and $z=E_{i,\ell}(a)t$.  The local map is
\[
    \widehat\Psi_{i,\ell}(\vartheta,a,t)
    :=\Bigl(
        \vartheta,
        \Phi_{i,\vartheta}(\mathsf S_{i,\ell}(a))
        +E_{i,\ell}(a)t
      \Bigr).
\]
The implicit-function construction in the proof of
\Cref{prop_parametric_morse_bott_stability} gives joint local
$\mathcal C^1$ dependence of the normal graph on $(\vartheta,a)$.  Because
\Cref{assum_sample_comp} supplies a $\mathcal C^2$ extension of $g$ to an open
neighborhood of $\Theta\times\mathbb R^d$, this construction and hence
$\widehat\Psi_{i,\ell}$ are defined and $\mathcal C^1$ on an ambient open
parameter neighborhood.  On $\Theta$, the local extension agrees with
$\Phi_{i,\theta}$ by uniqueness of the normal solution.  Its derivative has
the block-triangular form
\[
    D\widehat\Psi_{i,\ell}
    =
    \begin{bmatrix}
        I & 0\\
        * & D_{(a,t)}
        [\Phi_{i,\vartheta}(\mathsf S_{i,\ell}(a))+E_{i,\ell}(a)t]
    \end{bmatrix}.
\]
The spatial block is invertible by the choice of $r_i$, so the full derivative
is invertible.  The ordinary $\mathcal C^1$ inverse-function theorem gives an
open image in the ambient parameter space.  Since $V_i$ is relatively open,
the ambient parameter neighborhood can be chosen so that its intersection
with $\Theta$ lies in $V_i$; moreover, the map and its local inverse preserve
the first coordinate.  Restricting the open image to
$\Theta\times\mathbb R^d$ therefore gives a relative neighborhood contained in
$\Psi_i(\mathcal E_i)$.  Thus $\Psi_i(\mathcal E_i)$ is relatively open, and so
is their finite union $\mathcal W$.  This argument is local and uses neither a
global normal frame nor global injectivity of $\Psi_i$.

It remains to show that $\mathcal W$ contains the solution graph.  Let
$(\theta,p)\in\mathfrak G$.  Since the sets $V_i$ cover $\Theta$, choose $i$
such that $\theta\in V_i$.  Moreover,
\[
    p\in\mathcal S(\theta)=\Phi_{i,\theta}(M_i),
\]
so there is some $u\in M_i$ such that $p=\Phi_{i,\theta}(u)$.  Therefore
$(\theta,u,0)\in\mathcal E_i$ and
\[
    (\theta,p)
    =\Psi_i(\theta,u,0)
    \in\Psi_i(\mathcal E_i)
    \subseteq\mathcal W.
\]
Thus $\mathfrak G\subseteq\mathcal W$, proving that $\mathcal W$ is a
relatively open neighborhood of $\mathfrak G$.

Fix $(\theta,x)\in\mathcal W$ and choose a representation
$p:=\Phi_{i,\theta}(u)\in\mathcal S(\theta)$ and $x=p+z$ as above.  Since
$\partial_xg(\theta,p)=0$, Taylor's formula, the global
$L_{g,xx}$-Lipschitz continuity of the lower gradient, and
\eqref{eq_uniform_band_transverse_gap} give
\begin{align}
g(\theta,x)-g^\star(\theta)
&=\int_0^1(1-t)
  z^\top H_\theta(p+tz)z\,\mathrm dt
 \le\frac{L_{g,xx}}{2}\|z\|^2,
\label{eq_uniform_gap_upper_z}\\
\langle\partial_xg(\theta,x),z\rangle
&=\int_0^1 z^\top H_\theta(p+tz)z\,\mathrm dt
 \ge\frac{\mu}{2}\|z\|^2.
\label{eq_uniform_grad_lower_z}
\end{align}
By the Cauchy--Schwarz inequality (with the conclusion immediate when $z=0$),
\begin{equation}\label{eq_uniform_EB_tube}
\operatorname{dist}\bigl(x,\mathcal S(\theta)\bigr)
\le\|z\|
\le\frac{2}{\mu}\|\partial_xg(\theta,x)\|.
\end{equation}
Combining \eqref{eq_uniform_gap_upper_z} and
\eqref{eq_uniform_grad_lower_z} also yields
\[
g(\theta,x)-g^\star(\theta)
\le\frac{2L_{g,xx}}{\mu^2}\|\partial_xg(\theta,x)\|^2
=\frac{1}{2\bar\mu}\|\partial_xg(\theta,x)\|^2,
\qquad
\bar\mu=\frac{\mu^2}{4L_{g,xx}}.
\]
Thus both claimed estimates hold throughout $\mathcal W$.

\emph{A common sublevel set.}
It remains to place one common sublevel set inside $\mathcal W$.  The uniform
solution bound in \Cref{prop_parametric_morse_bott_stability} gives
\[
g^\star(\theta)
=\min_{x\in\overline{\mathbb B}_d(0;\mathsf D_{\mathcal S}+1)}g(\theta,x).
\]
The minimum is taken over a fixed compact ball, so $g^\star$ is continuous on
$\Theta$.  Set
\[
R:=\max\left\{
\mathsf D_{\mathsf{qg}},\,
\mathsf D_{\mathcal S}+\sqrt{2/\mu_{\mathsf{qg}}}
\right\}.
\]
For every $\theta\in\Theta$ and $\|x\|\ge R$,
the derived estimate \eqref{eq_tail_EB_implies_QG} and
$\operatorname{dist}(x,\mathcal S(\theta))
\ge\|x\|-\mathsf D_{\mathcal S}$ imply
\begin{equation}\label{eq_uniform_gap_outside_R}
g(\theta,x)-g^\star(\theta)
\ge\frac{\mu_{\mathsf{qg}}}{2}
   (\|x\|-\mathsf D_{\mathcal S})^2
\ge1.
\end{equation}
The set
\[
K:=\bigl(\Theta\times\overline{\mathbb B}_d(0;R)\bigr)\setminus\mathcal W
\]
is compact and disjoint from $\mathfrak G$.  The continuous function
$(\theta,x)\mapsto g(\theta,x)-g^\star(\theta)$ is therefore strictly positive
on $K$.  If $K\ne\varnothing$, define
\[
    c:=\min_{(\theta,x)\in K}
       \bigl(g(\theta,x)-g^\star(\theta)\bigr)>0;
\]
if $K=\varnothing$, set $c:=1$.  With
\begin{equation}\label{eq_uniform_epsilon_PL}
    \varepsilon_{\mathrm{PL}}:=\frac12\min\{1,c\},
\end{equation}
\eqref{eq_uniform_gap_outside_R} and the definition of $c$ show that
\[
g(\theta,x)-g^\star(\theta)\le\varepsilon_{\mathrm{PL}}
\quad\Longrightarrow\quad
(\theta,x)\in\mathcal W.
\]
The estimates already proved on $\mathcal W$ complete the proof.
\end{proof}

\subsection{Proofs of \Cref{lemm_PL_yield_S_theta_lip,lemm_hyp_blo_cont}}
\label[appendix]{append_proof_Hausdorff_and_hyperobjective_Lipschitz}

\begin{proof}[Proof of \Cref{lemm_PL_yield_S_theta_lip}]
Let
\[
    M_{\mathcal S}
    :=L_{g,1}(\mathsf D_{\mathcal S})
    =\sup_{\theta\in\Theta}
       \sup_{x\in\mathbb B_d(0;\mathsf D_{\mathcal S})}
       \|\partial_\theta g(\theta,x)\|<\infty.
\]
Fix $\theta_1,\theta_2\in\Theta$, take
$x_i\in\mathcal S(\theta_i)$, and write
$\delta:=\|\theta_1-\theta_2\|$.  Since $\Theta$ is convex and both minimizers
belong to the displayed ball, the fundamental theorem of calculus along the
segment $[\theta_1,\theta_2]$, together with optimality of $x_2$ for
$g(\theta_2,\cdot)$, gives
\begin{align}
0&\le g(\theta_1,x_2)-g^\star(\theta_1)\nonumber\\
&=g(\theta_1,x_2)-g(\theta_2,x_2)
  +g(\theta_2,x_2)-g(\theta_2,x_1)
  +g(\theta_2,x_1)-g(\theta_1,x_1)\nonumber\\
&\le2M_{\mathcal S}\delta.
\label{eq_append_uniform_objective_comparison}
\end{align}
If $M_{\mathcal S}>0$, put
$\gamma:=\varepsilon_{\mathrm{PL}}/(2M_{\mathcal S})$.  Whenever
$\delta\le\gamma$, \eqref{eq_append_uniform_objective_comparison} places $x_2$
in the common sublevel set in \eqref{eq_uniform_PL_band} for
\(\theta_1\), where \eqref{eq_uniform_error_bound} applies.  Since
$\partial_xg(\theta_2,x_2)=0$, \eqref{eq_uniform_error_bound} and the
$\theta$-Lipschitz property of the lower gradient yield
\[
\operatorname{dist}(x_2,\mathcal S(\theta_1))
\le\frac{2}{\mu}\|\partial_xg(\theta_1,x_2)\|
\le\frac{2L_{g,\theta x}}{\mu}\delta.
\]
Taking the supremum over $x_2$ and interchanging $\theta_1,\theta_2$ proves the
same Hausdorff bound for every pair with distance at most $\gamma$.
If $M_{\mathcal S}=0$, then for any $x_2\in\mathcal S(\theta_2)$,
\eqref{eq_append_uniform_objective_comparison} gives
\[
    0\le g(\theta_1,x_2)-g^\star(\theta_1)\le0.
\]
Thus $x_2\in\mathcal S(\theta_1)$, so
$\mathcal S(\theta_2)\subseteq\mathcal S(\theta_1)$.  Interchanging
$\theta_1$ and $\theta_2$ gives the reverse inclusion.  Therefore
\[
    \mathcal S(\theta_1)=\mathcal S(\theta_2),
\]
and their Hausdorff distance is zero.  The claimed estimate consequently holds
without any restriction on $\delta$.

For an arbitrary pair with $M_{\mathcal S}>0$, choose
$N:=\max\{1,\lceil\delta/\gamma\rceil\}$ and
\[
    \theta_j:=\theta_1+\frac{j}{N}(\theta_2-\theta_1),
    \qquad j=0,\ldots,N.
\]
Convexity keeps every $\theta_j$ in $\Theta$, and successive points are at most
$\gamma$ apart.  The triangle inequality for Hausdorff distance therefore gives
\[
\operatorname{dist}(\mathcal S(\theta_1),\mathcal S(\theta_2))
\le\sum_{j=1}^N
\operatorname{dist}(\mathcal S(\theta_{j-1}),\mathcal S(\theta_j))
\le\frac{2L_{g,\theta x}}{\mu}\delta.
\]
This is \eqref{eq_global_Hausdorff_Lipschitz}.
\end{proof}

\begin{proof}[Proof of \Cref{lemm_hyp_blo_cont}]
Choose $x_1\in\argmax_{x\in\mathcal S(\theta_1)}f(\theta_1,x)$ and choose
$x_2\in\mathcal S(\theta_2)$ nearest to $x_1$.  Then
\begin{align*}
F_{\max}(\theta_1)-F_{\max}(\theta_2)
&\le f(\theta_1,x_1)-f(\theta_2,x_2)\\
&\le L_{f,1}\bigl(\|\theta_1-\theta_2\|+\|x_1-x_2\|\bigr)\\
&\le L_{f,1}\left(1+\frac{2L_{g,\theta x}}{\mu}\right)
      \|\theta_1-\theta_2\|.
\end{align*}
Interchanging $\theta_1$ and $\theta_2$ proves the asserted absolute-value
bound.
\end{proof}

\subsection{Proof of \Cref{lemma_common_dimension}}\label[appendix]{append_proof_lemma_common_dimension}
\begin{proof}
Fix $\bar\theta\in\Theta$.  By
\Cref{prop_parametric_morse_bott_stability}, every $\mathcal S(\theta)$ with
$\theta$ in a relative neighborhood of $\bar\theta$ is the image of
$\mathcal S(\bar\theta)$ under the embedding
$u\mapsto u+n_\theta(u)$; the proof of that proposition shows that this image is
the entire minimizer manifold.  Hence
\[
    \dim\mathcal S(\theta)=\dim\mathcal S(\bar\theta)
\]
throughout that neighborhood.  Thus
$\theta\mapsto\dim\mathcal S(\theta)$ is locally constant.  Since the standing
parameter set $\Theta$ is convex and therefore connected, a locally constant
integer-valued function on $\Theta$ is constant.  This gives the common dimension
$\dimk$.
\end{proof}

\subsection{Uniform ambient tubular geometry for
\Cref{assum_compact_opt_set,prop_W1_Gibbs_uniform}}
\label[appendix]{append_proof_uniform_ambient_tube}
The lemma below is first used in \Cref{assum_compact_opt_set} to obtain a
uniform reach and hence a uniform intrinsic injectivity radius.  The same
ambient tube and its change-of-variables formula are subsequently used in
\Cref{prop_W1_Gibbs_uniform}.  Its proof is independent of the conclusions of
\Cref{assum_compact_opt_set}.  The radius constructed here concerns
nearest-point projection in the ambient Euclidean space and is distinct from
the intrinsic injectivity radius derived from it.
\begin{lem}[Uniform ambient tubular geometry]
\label{lem_uniform_ambient_tube}
Under
\Cref{mu_PL_assum,assum_sample_comp,assum_secon_fund_form_bound}, there exist
constants \(r_{\mathrm{tub}}>0\) and \(C_J<\infty\), independent of \(\theta\),
such that, for every \(\theta\in\Theta\), the addition map
\[
    \Phi_\theta(u,n):=u+n
\]
is a \(\mathcal C^1\) diffeomorphism from
\[
\bigl\{(u,n):u\in\mathcal S(\theta),\ n\in N_u\mathcal S(\theta),\
                 \|n\|<r_{\mathrm{tub}}\bigr\}
\]
onto
\[
    \mathcal U_\theta
    :=\bigl\{x:\operatorname{dist}(x,\mathcal S(\theta))
                       <r_{\mathrm{tub}}\bigr\}.
\]
Thus every \(x\in\mathcal U_\theta\) has a unique nearest point
\(\pi_\theta(x)\in\mathcal S(\theta)\).  For every integrable \(\varphi\),
\begin{equation}\label{eq_uniform_tubular_volume_formula}
\int_{\mathcal U_\theta}\varphi(x)\,\mathrm dx
=
\int_{\mathcal S(\theta)}
\int_{\substack{n\in N_u\mathcal S(\theta)\\
                 \|n\|<r_{\mathrm{tub}}}}
\varphi(u+n)J_\theta(u,n)\,\mathrm dn\,\mathrm d\mathcal M_\theta(u),
\end{equation}
where \(\mathrm dn\) is Euclidean Lebesgue measure on the normal fiber and
\begin{equation}\label{eq_uniform_tubular_jacobian}
    J_\theta(u,0)=1,
    \qquad
    |J_\theta(u,n)-1|\le C_J\|n\|.
\end{equation}
\end{lem}

\begin{proof}
If \(\dimk=0\), connectedness makes every \(\mathcal S(\theta)\) a singleton.
The addition map is then a translation from its normal fiber
\(\mathbb R^d\); the conclusion holds, for example, with
\(r_{\mathrm{tub}}=1\) and \(J_\theta\equiv1\).  We therefore assume
\(\dimk\ge1\).

We first obtain a radius on a neighborhood of one reference parameter.  Fix
\(\theta_i\in\Theta\), put \(M_i:=\mathcal S(\theta_i)\), and write the
normal-graph parametrizations supplied by
\Cref{prop_parametric_morse_bott_stability} as
\[
    F_{i,\theta}(u):=u+n_{i,\theta}(u),
    \qquad u\in M_i.
\]
By \eqref{eq_parametric_normal_graph_C1_convergence}, after shrinking the
relative parameter neighborhood \(U_i\) of \(\theta_i\),
\begin{equation}\label{eq_uniform_tube_graph_smallness}
\sup_{u\in M_i}\|Dn_{i,\theta}(u)\|_{\mathrm{op}}\le\frac12
\end{equation}
for every \(\theta\in U_i\).  The same convergence permits
\(\sup_{u\in M_i}\|n_{i,\theta}(u)\|\) to be made smaller than any positive
number fixed below.  Here and below the operator norm of the derivative is
restricted to \(T_uM_i\).

We prove that short normal segments based at distinct points of
\(\mathcal S(\theta)\) cannot meet, uniformly for \(\theta\in U_i\).  Let
\(\mathsf C\) be the common second-fundamental-form bound in
\Cref{assum_secon_fund_form_bound}.  Since \(M_i\) is connected and
positive-dimensional, choose
\(0<\delta_i<\operatorname{diam}_{M_i}(M_i)\) so small that
\(3\mathsf C\delta_i/2\le1/2\).  The set
\[
    P_i:=\{(u,v)\in M_i\times M_i:
                  \operatorname{dist}_{M_i}(u,v)\ge\delta_i\}
\]
is compact and contains no diagonal point.  Hence
\[
    s_i:=\min_{(u,v)\in P_i}\|u-v\|>0.
\]
Shrink \(U_i\) further so that
\(\sup_{u\in M_i}\|n_{i,\theta}(u)\|\le s_i/8\).  If
\((u,v)\in P_i\), then
\begin{equation}\label{eq_uniform_tube_far_separation}
\|F_{i,\theta}(u)-F_{i,\theta}(v)\|
\ge s_i-2(s_i/8)=\frac{3s_i}{4}.
\end{equation}

Consider next \(u\ne v\) with
\(\operatorname{dist}_{M_i}(u,v)<\delta_i\), and put
\(x=F_{i,\theta}(u)\), \(y=F_{i,\theta}(v)\).  Apply
\(F_{i,\theta}\) to curves in \(M_i\) whose lengths decrease to
\(\operatorname{dist}_{M_i}(u,v)\).  By
\eqref{eq_uniform_tube_graph_smallness}, their image lengths are at most
three halves of their original lengths.  Thus the intrinsic distance
\(\ell:=\operatorname{dist}_{\mathcal S(\theta)}(x,y)\) satisfies
\(\mathsf C\ell\le1/2\).  Let
\(\gamma:[0,\ell]\to\mathcal S(\theta)\) be a unit-speed minimizing
geodesic from \(x\) to \(y\).  Such a curve exists by the direct method on the
compact manifold.  For completeness, in a \(\mathcal C^2\) chart the induced
metric coefficients are \(\mathcal C^1\).  A constant-speed energy minimizer
satisfies the Euler--Lagrange equation weakly, which makes its metric momentum
absolutely continuous.  Uniform positive definiteness and the \(\mathcal C^1\)
metric coefficients then make its velocity absolutely continuous and bootstrap
the curve to a \(\mathcal C^2\) solution of the geodesic equation.
Its ambient acceleration is
\(\mathrm{II}_\theta(\dot\gamma,\dot\gamma)\), so
\(\|\ddot\gamma\|\le\mathsf C\).  Therefore
\begin{equation}\label{eq_uniform_tube_chord_lower}
\|y-x\|
\ge \ell-\int_0^\ell\mathsf C s\,\mathrm ds
\ge\frac{3\ell}{4}.
\end{equation}
For \(n\in N_x\mathcal S(\theta)\) and
\(m\in N_y\mathcal S(\theta)\), integrate
\(\|\dot\gamma(s)-\dot\gamma(0)\|\le\mathsf C s\) from the first endpoint
and the analogous estimate backward from the second endpoint.  Orthogonality
at the endpoints gives
\begin{align}
|\langle n,y-x\rangle|
&\le\frac{\mathsf C}{2}\|n\|\ell^2,
&
|\langle m,y-x\rangle|
&\le\frac{\mathsf C}{2}\|m\|\ell^2.
\label{eq_uniform_tube_normal_chord}
\end{align}
If \(x+n=y+m\), then \(y-x=n-m\).  Consequently,
\eqref{eq_uniform_tube_chord_lower}--\eqref{eq_uniform_tube_normal_chord}
imply
\[
    \frac{9}{16}\ell^2
    \le\|y-x\|^2
    =\langle n-m,y-x\rangle
    \le\frac{\mathsf C}{2}(\|n\|+\|m\|)\ell^2.
\]
This is impossible when
\(\|n\|,\|m\|<1/[4(1+\mathsf C)]\), because the last coefficient is then
strictly smaller than \(1/4\).  For a far pair,
\eqref{eq_uniform_tube_far_separation} and \(x+n=y+m\) would instead give
\(3s_i/4\le\|x-y\|=\|n-m\|<s_i/2\).  Thus the addition map is injective on
normal vectors of length less than
\begin{equation}\label{eq_uniform_tube_local_radius}
    r_i<\min\left\{\frac{s_i}{4},\frac{1}{4(1+\mathsf C)}\right\},
\end{equation}
uniformly for \(\theta\in U_i\).

We next verify local invertibility using the actual normal spaces of
\(\mathcal S(\theta)\).  For \(n\in N_u\mathcal S(\theta)\), define the shape
operator \(A_n^\theta:T_u\mathcal S(\theta)\to T_u\mathcal S(\theta)\) by
\[
    \langle A_n^\theta v,w\rangle
    :=\langle\mathrm{II}_\theta(u)[v,w],n\rangle
    \qquad
    \bigl(v,w\in T_u\mathcal S(\theta)\bigr).
\]
Choose a local orthonormal normal frame near \(u\), and let \(\widetilde n\)
be the local normal field obtained by holding the fiber coordinates of \(n\)
fixed.  The Weingarten formula is
\[
    D_v\widetilde n(u)
    =-A_n^\theta v+\nabla_v^\perp\widetilde n,
    \qquad
    \nabla_v^\perp\widetilde n
    :=P_{N_u\mathcal S(\theta)}D_v\widetilde n(u).
\]
Thus, for a base variation \(v\in T_u\mathcal S(\theta)\) and a vertical
variation \(\xi\in N_u\mathcal S(\theta)\), differentiation of
\(\Phi_\theta(u,n)=u+n\) in these local bundle coordinates gives
\[
    D\Phi_\theta(u,n)[v,\xi]
    =(I-A_n^\theta)v+\nabla_v^\perp\widetilde n+\xi.
\]
The first term is tangent at \(u\), whereas the last two terms are normal.
Consequently, relative to the tangent/normal splitting at \(u\), the
differential has block form
\begin{equation}\label{eq_uniform_tube_differential}
    D\Phi_\theta(u,n)
    =
    \begin{bmatrix}
       I-A_n^\theta&0\\
       *&I
    \end{bmatrix},
\end{equation}
where the lower-left block represents the normal-connection term
\(v\mapsto\nabla_v^\perp\widetilde n\).
The defining identity and the assumed second-fundamental-form bound give the
following estimate.  Indeed,
polarization of the symmetric bilinear form \(\mathrm{II}_\theta\) shows that
\(\|\mathrm{II}_\theta(v,w)\|\le\mathsf C\) for unit tangent vectors
\(v,w\).  Hence
\begin{equation}\label{eq_uniform_shape_operator_bound}
    \|A_n^\theta\|_{\mathrm{op}}
    =\sup_{\|v\|=\|w\|=1}
       |\langle\mathrm{II}_\theta(v,w),n\rangle|
    \le\mathsf C\|n\|.
\end{equation}
Thus \(I-A_n^\theta\) is invertible for the radii in
\eqref{eq_uniform_tube_local_radius}.  The \(\mathcal C^1\)
inverse-function theorem makes \(\Phi_\theta\) a local diffeomorphism, and
the injectivity proved above makes it a diffeomorphism onto its image.

Repeating this construction at every reference parameter gives a relative
open cover of \(\Theta\).  Compactness selects reference parameters
\(\theta_1,\ldots,\theta_N\), their neighborhoods \(U_i\), and their radii
\(r_i\).  By the finite shrinking lemma for the compact metric space
\(\Theta\), choose relatively open \(V_i\) with
\(\overline V_i\subset U_i\) that still cover \(\Theta\), and set
\(r_{\mathrm{tub}}:=\min_i r_i>0\), decreasing the radii once more if
necessary.  The resulting conclusion holds for every \(\theta\in\Theta\).
Its image is exactly \(\mathcal U_\theta\).  Indeed,
\(\operatorname{dist}(u+n,\mathcal S(\theta))\le\|n\|<r_{\mathrm{tub}}\).
Conversely, compactness of \(\mathcal S(\theta)\) gives a nearest point \(u\)
to any \(x\in\mathcal U_\theta\); first-order optimality of the squared
distance on the manifold gives \(x-u\in N_u\mathcal S(\theta)\), and
\(\|x-u\|<r_{\mathrm{tub}}\).  If \(v\) were another nearest point, the same
first-order condition would give \(x-v\in N_v\mathcal S(\theta)\) with
\(\|x-v\|<r_{\mathrm{tub}}\).  These would be two normal representations of
\(x\), contradicting injectivity.  Thus every point at distance less than
\(r_{\mathrm{tub}}\) has a unique nearest point.  By \Cref{def_reach},
\[
    \operatorname{reach}\bigl(\mathcal S(\theta)\bigr)
    \ge r_{\mathrm{tub}}
    \qquad(\theta\in\Theta).
\]

Finally, \eqref{eq_uniform_tube_differential} and the change-of-variables
formula give, relative to fiber Lebesgue measure and Riemannian volume,
\begin{equation}\label{eq_uniform_tube_jacobian_shape}
    J_\theta(u,n)=\det(I-A_n^\theta)>0.
\end{equation}
If \(\alpha_1,\ldots,\alpha_{\dimk}\) are the eigenvalues of
\(A_n^\theta\), then \(|\alpha_j|\le\mathsf C\|n\|\).  Hence
\[
|J_\theta(u,n)-1|
\le \dimk\,\mathsf C(1+\mathsf C r_{\mathrm{tub}})^{\dimk-1}\|n\|,
\]
which proves \eqref{eq_uniform_tubular_jacobian}.  Formula
\eqref{eq_uniform_tubular_volume_formula} follows by applying the local
change-of-variables formula in finitely many manifold charts with local
orthonormal normal frames and summing with a partition of unity.  Orthogonal
changes of normal frame preserve fiber Lebesgue measure and the determinant,
so neither orientability nor a global normal frame is used.
\end{proof}

\subsection{Proof of \Cref{assum_compact_opt_set}}\label[appendix]{append_proof_assum_compact_opt_set}
\begin{proof}\label{proof_assum_compact_opt_set}
We first record that the Riemannian objects used below have sufficient
finite regularity.  Fix \(\theta\in\Theta\) and
\(p\in\mathcal S(\theta)\).  Since the
minimizer manifold is already known to be \(\mathcal C^2\), it has near \(p\)
a \(\mathcal C^2\) graph representation
\[
    p+a+\psi(a),
    \qquad
    a\in T_p\mathcal S(\theta),\quad
    \psi(a)\in N_p\mathcal S(\theta).
\]
Define
\[
    F(a,b)
    :=P_{N_p\mathcal S(\theta)}
      \partial_xg(\theta,p+a+b),
    \qquad
    (a,b)\in T_p\mathcal S(\theta)\times N_p\mathcal S(\theta).
\]
The pointwise \(\mathcal C^4\) clause in the \(\PLcirc\) condition makes
\(F\) a \(\mathcal C^3\) map near \((0,0)\).  Moreover,
\[
    D_bF(0,0)
    =P_{N_p\mathcal S(\theta)}\partial_x^2g(\theta,p)
      \big|_{N_p\mathcal S(\theta)}.
\]
Indeed, the local P\L--Morse--Bott conclusion gives
\[
    \ker\partial_x^2g(\theta,p)=T_p\mathcal S(\theta).
\]
Symmetry of the Hessian therefore makes \(N_p\mathcal S(\theta)\) invariant,
so \(D_bF(0,0)\) is the restriction of \(\partial_x^2g(\theta,p)\) to
\(N_p\mathcal S(\theta)\), and this restriction is positive definite.  The
\(\mathcal C^3\) implicit-function theorem consequently gives a unique local
\(\mathcal C^3\) graph \(b=\phi(a)\) satisfying
\(F(a,\phi(a))=0\).  Every point of the existing minimizer graph satisfies
the same equation, so local uniqueness gives \(\psi=\phi\) after shrinking
the graph domain.  Thus each fixed \(\mathcal S(\theta)\) is locally
\(\mathcal C^3\), and its Euclidean-induced metric is \(\mathcal C^2\).  No
uniform \(\mathcal C^3\) estimate in \(\theta\) is asserted or used.

\textbf{Property 1.}
By \Cref{prop_parametric_morse_bott_stability},
\(\mathcal S(\theta)\subseteq
\mathbb B_d(0;\mathsf D_{\mathcal S})\) for every \(\theta\in\Theta\).
Therefore
\[
    \operatorname{diam}_{\mathbb R^d}\mathcal S(\theta)
    \le 2\mathsf D_{\mathcal S}
    \le \mathsf D,
\]
where the last inequality is the choice in \eqref{eq_working_radius_D}.

\textbf{Property 2.}
We use the following pointwise form of the Gauss equation.
\begin{lem}[Gauss equation for Euclidean submanifolds]
\label{lem_sectional_curv}
If \(M\subset\mathbb R^d\) is an embedded submanifold with a
\(\mathcal C^2\) induced metric and \(X,Y\in T_pM\) are orthonormal, then
\[
    \operatorname{Sec}_M(p,\operatorname{span}\{X,Y\})
    =\langle\mathrm{II}(X,X),\mathrm{II}(Y,Y)\rangle
     -\|\mathrm{II}(X,Y)\|^2.
\]
\end{lem}
This identity follows by decomposing the ambient derivative as
\(D_XY=\nabla_XY+\mathrm{II}(X,Y)\) in the definition of the curvature
tensor and using that the ambient Euclidean curvature vanishes.

The diagonal norm convention in \Cref{def:second-fund-form} gives
\(\|\mathrm{II}(v,v)\|\le\mathsf C\|v\|^2\).  For orthonormal \(X,Y\),
polarization yields
\[
\mathrm{II}(X,Y)
=\frac14\bigl[\mathrm{II}(X+Y,X+Y)
               -\mathrm{II}(X-Y,X-Y)\bigr],
\]
and hence
\[
    \|\mathrm{II}(X,Y)\|
    \le\frac{\mathsf C}{4}
       \bigl(\|X+Y\|^2+\|X-Y\|^2\bigr)
    =\mathsf C.
\]
Together with \Cref{lem_sectional_curv}, this gives
\[
    -2\mathsf C^2
    \le \operatorname{Sec}_{\mathcal S(\theta)}
       (p,\operatorname{span}\{X,Y\})
    \le \mathsf C^2.
\]
Thus Property~2 holds with
\(\mathsf C_0:=2\mathsf C^2\), independently of \(\theta,p\), and the
two-plane.  The assertion is vacuous when \(\dimk<2\), and the argument uses
neither orientability nor a normal frame.

\textbf{Property 3.}
We use the following reach-to-injectivity estimate.
\begin{prop}[Reach and intrinsic injectivity radius]
\label{prop_inject_lower_bound}
Let \(M\subset\mathbb R^d\) be a closed \(\mathcal C^2\) embedded
submanifold with reach \(\tau_M>0\).  Then
\[
    \operatorname{inj}_M(p)
    \ge \operatorname{inj}(M)
    \ge \pi\tau_M
    \qquad(p\in M),
\]
where \(\operatorname{inj}(M):=\inf_{z\in M}\operatorname{inj}_M(z)\).
The second inequality is \citep[Proposition~A.1(ii)]{aamari2019estimating};
the first follows from the displayed definition of \(\operatorname{inj}(M)\).
\end{prop}

The common-tube result \Cref{lem_uniform_ambient_tube}, proved independently
in \Cref{append_proof_uniform_ambient_tube}, gives, for every \(\theta\), a
unique nearest point in \(\mathcal S(\theta)\) throughout the open distance
tube of one common radius \(r_{\mathrm{tub}}>0\).  Its proof uses only
\Cref{prop_parametric_morse_bott_stability}, compactness of \(\Theta\), the
\(\mathcal C^1\)-small normal graphs, and
\Cref{assum_secon_fund_form_bound}; it does not use the curvature or
injectivity-radius conclusions of the present proposition.  By
\Cref{def_reach},
\[
    \operatorname{reach}\bigl(\mathcal S(\theta)\bigr)
    \ge r_{\mathrm{tub}}
    \qquad(\theta\in\Theta).
\]
Each \(\mathcal S(\theta)\) is compact and hence closed, so
\Cref{prop_inject_lower_bound} applies and gives
\[
    \operatorname{inj}_{\mathcal S(\theta)}(p)
    \ge r_0,
    \qquad
    r_0:=\pi r_{\mathrm{tub}}>0,
\]
for every \(\theta\in\Theta\) and \(p\in\mathcal S(\theta)\).  The logical
chain is unique nearest-point tube, then positive reach, then intrinsic
injectivity radius; these three notions are not being identified.
\end{proof}

\subsection{Uniform volume bounds}
\label[appendix]{append_proof_uniform_manifold_volume}
\begin{lem}[Uniform volume bounds]\label{lem_vol_sol_set_is_bounded}
Under \Cref{mu_PL_assum,assum_sample_comp} and compactness of \(\Theta\), there
exists \(0<\mathsf b<1\) such that
\[
    \mathsf b
    \le \mathcal M_\theta\bigl(\mathcal S(\theta)\bigr)
    \le \mathsf b^{-1}
    \qquad\forall\theta\in\Theta.
\]
\end{lem}
\begin{proof}
Fix \(\bar\theta\in\Theta\), and put
\(M_{\bar\theta}:=\mathcal S(\bar\theta)\).  By
\Cref{prop_parametric_morse_bott_stability}, there is a relative neighborhood
\(U_{\bar\theta}\) of \(\bar\theta\) such that, for every
\(\theta\in U_{\bar\theta}\),
\[
    F_{\bar\theta,\theta}(u)
    :=u+n_{\bar\theta,\theta}(u)
\]
is an injective \(\mathcal C^1\) embedding of \(M_{\bar\theta}\) onto the
entire manifold \(\mathcal S(\theta)\).  By
\eqref{eq_parametric_normal_graph_C1_convergence}, as recorded more directly
in \eqref{eq_append_section_smallness}, we may first shrink this particular
relative neighborhood, while retaining \(\bar\theta\), so that
\[
    \sup_{\substack{u\in M_{\bar\theta},\ v\in T_uM_{\bar\theta}\\
                    \|v\|=1}}
    \|Dn_{\bar\theta,\theta}(u)[v]\|
    \le\frac12
    \qquad(\theta\in U_{\bar\theta}).
\]
The resulting relative neighborhoods still cover \(\Theta\), because every
\(\bar\theta\) belongs to its own neighborhood.  Compactness of \(\Theta\)
now gives a finite subcover with reference parameters
\(\theta_1,\ldots,\theta_N\).  Write
\[
    M_i:=\mathcal S(\theta_i),
    \qquad
    F_{i,\theta}(u):=u+n_{i,\theta}(u)
\]
on the corresponding member of this finite cover.

The same normal-graph representation shows that
\(\theta\mapsto\dim\mathcal S(\theta)\) is locally constant.  Since the
standing set \(\Theta\) is convex and hence connected, all these manifolds
have one common dimension, denoted by \(\dimk\).

For \(u\in M_i\) and \(v\in T_uM_i\), differentiation in the ambient space
gives
\[
    DF_{i,\theta}(u)[v]
    =v+Dn_{i,\theta}(u)[v].
\]
The triangle and reverse-triangle inequalities therefore give
\[
    \frac12\|v\|
    \le\|DF_{i,\theta}(u)[v]\|
    \le\frac32\|v\|.
\]
Set
\[
    A_{i,\theta}(u)
    :=DF_{i,\theta}(u)\big|_{T_uM_i},
\]
and define its tangential Jacobian by
\[
    \operatorname{Jac}_{\dimk}F_{i,\theta}(u)
    :=
    \sqrt{\det\!\left(A_{i,\theta}(u)^*A_{i,\theta}(u)\right)}.
\]
This is the product of the singular values of \(A_{i,\theta}(u)\).  The last
two-sided norm estimate places each singular value in \([1/2,3/2]\), and hence
\[
    2^{-\dimk}
    \le\operatorname{Jac}_{\dimk}F_{i,\theta}(u)
    \le(3/2)^{\dimk}.
\]

Let \(\mathcal M_i\) be the Euclidean-induced Riemannian volume measure on
\(M_i\).  We briefly justify the required area formula.  If
\(\mathsf S:\Gamma\to M_i\) is a local parametrization, then
\[
    \mathrm d\mathcal M_i(\mathsf S(a))
    =\sqrt{\det\!\left(D\mathsf S(a)^\top D\mathsf S(a)\right)}\,\mathrm da,
\]
whereas the volume density of the image parametrization
\(F_{i,\theta}\circ\mathsf S\) is
\[
    \sqrt{\det\!\left(
        D(F_{i,\theta}\circ\mathsf S)(a)^\top
        D(F_{i,\theta}\circ\mathsf S)(a)
    \right)}\,\mathrm da.
\]
The ratio of these two densities is
\(\operatorname{Jac}_{\dimk}F_{i,\theta}(\mathsf S(a))\).  Applying this
local change of variables on a finite atlas and summing with a partition of
unity gives
\[
    \mathcal M_\theta\bigl(\mathcal S(\theta)\bigr)
    =\int_{M_i}\operatorname{Jac}_{\dimk}F_{i,\theta}(u)
       \,\mathrm d\mathcal M_i(u).
\]
The multiplicity is one because \(F_{i,\theta}\) is injective and its image is
the whole of \(\mathcal S(\theta)\).

For every \(i\), compactness of \(M_i\) and a finite atlas give
\(\mathcal M_i(M_i)<\infty\), while nonemptiness and a local chart with
strictly positive metric density give \(\mathcal M_i(M_i)>0\).  If
\(\dimk=0\), Riemannian volume is counting measure, so a nonempty compact
zero-dimensional manifold has finite positive volume; under connectedness it
is a singleton.  Consequently,
\[
    L
    :=\min_{1\le i\le N}2^{-\dimk}\mathcal M_i(M_i)>0,
    \qquad
    U
    :=\max_{1\le i\le N}(3/2)^{\dimk}\mathcal M_i(M_i)<\infty,
\]
and \(L\le U\).  Define
\[
    \mathsf b:=\frac12\min\{1,L,U^{-1}\}.
\]
Then \(0<\mathsf b<1\).  Given \(\theta\in\Theta\), choose a covering index
\(i\).  The Jacobian bounds and the area formula yield
\[
    \mathsf b
    \le L
    \le\mathcal M_\theta\bigl(\mathcal S(\theta)\bigr)
    \le U
    \le\mathsf b^{-1}.
\]
The constant \(\mathsf b\) depends only on the finite reference family and
its volumes, not on \(\theta\).  Finally, the proof uses Riemannian densities
and local coordinate changes of variables, so it requires neither
orientability nor a global chart or global normal frame.
\end{proof}

\subsection{Proof of \Cref{prop_strict_bound_limit_gibbs_density}}\label[appendix]{append_proof_prop_strict_bound_limit_gibbs_density}
To establish the uniform bounds on the limiting density $\rho_\theta$ in~\eqref{eqn_density_limiting_measure}, it is necessary to control the $\theta$-dependent normalizing factor, which involves the Riemannian volume $\mathcal{M}\!\left(\mathcal{S}(\theta)\right)$. The required uniform volume bounds are supplied by \Cref{lem_vol_sol_set_is_bounded}.
\begin{lem}\label{lemma_unif_bound_A_theta}
Let $\Theta$ be compact and let $\mathcal A(\theta)$ be a nonempty compact set
for every $\theta\in\Theta$.  If $\theta\mapsto\mathcal A(\theta)$ is
continuous in Hausdorff distance, then the family
$\{\mathcal A(\theta)\}_{\theta\in\Theta}$ is uniformly bounded.
\end{lem}
\begin{proof}
Define
\[
    M(\theta):=\max_{x\in\mathcal A(\theta)}\|x\|.
\]
For $x\in\mathcal A(\theta_2)$, compactness provides
$y\in\mathcal A(\theta_1)$ with
$\|x-y\|\le d_{\mathrm H}(\mathcal A(\theta_1),\mathcal A(\theta_2))$.
Thus
\[
    \|x\|\le M(\theta_1)
    +d_{\mathrm H}\bigl(\mathcal A(\theta_1),\mathcal A(\theta_2)\bigr).
\]
Taking the maximum over $x$ and then interchanging $\theta_1,\theta_2$ gives
\[
    |M(\theta_1)-M(\theta_2)|
    \le d_{\mathrm H}\bigl(\mathcal A(\theta_1),\mathcal A(\theta_2)\bigr).
\]
Hence $M$ is continuous.  Its maximum on compact $\Theta$ is finite and
provides a common bounding radius.
\end{proof}
Now we provide the proof of \Cref{prop_strict_bound_limit_gibbs_density}:
\begin{proof}
Put \(q:=d-\dimk\).  Fix \(\theta\in\Theta\),
\(u\in\mathcal S(\theta)\), and choose any local orthonormal frame of the
normal bundle near \(u\), with column matrix \(E(u)\in\mathbb R^{d\times q}\).
The normal-coordinate Hessian at \(u\) is
\[
    A_\theta(u):=E(u)^\top\partial_x^2g(\theta,u)E(u).
\]
The local P\L--Morse--Bott argument, with the normalization in
\Cref{def_local_PL}, gives
\[
    A_\theta(u)\succeq\mu I_q.
\]
There is no parameter-dependent upper bound to construct.  Indeed, global
\(L_{g,xx}\)-Lipschitz continuity of
\(x\mapsto\partial_xg(\theta,x)\), together with differentiability, implies
\[
    \|\partial_x^2g(\theta,x)\|_{\mathrm{op}}\le L_{g,xx}
    \qquad(\theta\in\Theta,\ x\in\mathbb R^d).
\]
Since \(E(u)\) is an isometry onto \(N_u\mathcal S(\theta)\), it follows that
\[
    \mu I_q\preceq A_\theta(u)\preceq L_{g,xx}I_q.
\]
If another local orthonormal normal frame is used, then
\(\widetilde E(u)=E(u)Q(u)\) for an orthogonal matrix \(Q(u)\), and
\(\widetilde A_\theta(u)=Q(u)^\top A_\theta(u)Q(u)\).  Thus
\(\det A_\theta(u)\), and all the following bounds, are independent of the
local frame; no global normal frame or orientability is required.

Write
\[
    w_\theta(u):=\det A_\theta(u)^{-1/2},
    \qquad
    V_\theta:=\mathcal M_\theta(\mathcal S(\theta)).
\]
From the preceding eigenvalue bounds,
\[
    L_{g,xx}^{-q/2}\le w_\theta(u)\le\mu^{-q/2}.
\]
Consequently, the denominator in \eqref{eqn_density_limiting_measure}
satisfies
\[
    L_{g,xx}^{-q/2}V_\theta
    \le
    \int_{\mathcal S(\theta)}w_\theta(v)\,\mathrm d\mathcal M_\theta(v)
    \le
    \mu^{-q/2}V_\theta.
\]
Since \Cref{lem_vol_sol_set_is_bounded} gives
\(\mathsf b\le V_\theta\le\mathsf b^{-1}\), the density
\(\rho_\theta=w_\theta/\int w_\theta\,\mathrm d\mathcal M_\theta\) obeys
\[
    \mathsf b\left(\frac{\mu}{L_{g,xx}}\right)^{q/2}
    \le \rho_\theta(u)\le
    \mathsf b^{-1}\left(\frac{L_{g,xx}}{\mu}\right)^{q/2}.
\]
Thus the proposition holds with the explicit, parameter-independent choice
\[
    \kappa
    :=\mathsf b\left(\frac{\mu}{L_{g,xx}}\right)^{(d-\dimk)/2}>0.
\]
Because \(0<\mathsf b<1\) and \(\mu\le L_{g,xx}\), this choice also satisfies
\(\kappa\le1\).
\end{proof}

\subsection{Proof of \Cref{prop_W1_Gibbs_uniform}}
\label[appendix]{append_proof_uniform_W1_Gibbs}
We first establish the uniform off-tube mass and first-moment estimates used
in the normalizer expansion and the final coupling.
Let \(\mathcal W\subseteq\Theta\times\mathbb R^d\) be the
relatively open neighborhood of the solution graph constructed in the proof
of \Cref{lem_uniform_local_PL_band}; see
\eqref{eq_uniform_tube_definition}.  For each \(\theta\in\Theta\), write
\[
    h_\theta(x):=g(\theta,x)-g^\star(\theta),
    \qquad
    \mathcal W_\theta:=\{x:(\theta,x)\in\mathcal W\}.
\]
The transverse estimate established in the proof of
\Cref{lem_uniform_local_PL_band} gives the following auxiliary localization
result.

\begin{corr}[Uniform Gibbs localization]
\label{corr_uniform_gibbs_localization}
Under \Cref{mu_PL_assum,assum_sample_comp},
\begin{equation}\label{eq_uniform_energy_gap}
    h_\theta(x)
    \ge \frac{\mu}{4}
    \operatorname{dist}^2\bigl(x,\mathcal S(\theta)\bigr)
    \quad(x\in\mathcal W_\theta),
    \qquad
    h_\theta(x)>\varepsilon_{\mathrm{PL}}
    \quad(x\notin\mathcal W_\theta).
\end{equation}
Consequently, for every \(r>0\), with
\[
    a_r:=\min\left\{\varepsilon_{\mathrm{PL}},\frac{\mu r^2}{4}\right\},
\]
one has
\begin{equation}\label{eq_uniform_energy_gap_at_distance}
    \operatorname{dist}\bigl(x,\mathcal S(\theta)\bigr)\ge r
    \quad\Longrightarrow\quad
    h_\theta(x)\ge a_r
    \qquad(\theta\in\Theta).
\end{equation}
Moreover, there is a constant \(C_{\mathrm{loc}}<\infty\), independent of
\(\theta,r,\lambda\), such that, for \(0<\lambda\le1\),
\begin{equation}\label{eq_uniform_weighted_localization}
\sup_{\theta\in\Theta}
\int_{\{\operatorname{dist}(x,\mathcal S(\theta))\ge r\}}
(1+\|x\|)\,\gibbs_\theta^\lambda(\mathrm dx)
\le
C_{\mathrm{loc}}\lambda^{-d/2}
\exp\!\left(-\frac{a_r}{2\lambda}\right).
\end{equation}
In particular,
\[
\lim_{\lambda\downarrow0}\sup_{\theta\in\Theta}
\gibbs_\theta^\lambda
\bigl\{x:\operatorname{dist}(x,\mathcal S(\theta))\ge r\bigr\}=0.
\]
The same conclusion holds with \(1+\|x\|\) replaced by any fixed polynomial
weight \(1+\|x\|^s\), with a constant that may depend on \(s\).
\end{corr}

The ambient-ball lower bound used to prove this corollary
has scale \(\lambda^{d/2}\).  This deliberately crude scale suffices for
localization.  In the remainder of the proof of
\Cref{prop_W1_Gibbs_uniform}, the normal-fiber Laplace expansion yields the
sharp scale \(\lambda^{(d-\dimk)/2}\).

\begin{proof}
Fix \(\theta\in\Theta\).  If
\(x\in\mathcal W_\theta\), the construction in the proof of
\Cref{lem_uniform_local_PL_band} provides an index \(i\), a point
\(p=\Phi_{i,\theta}(u)\in\mathcal S(\theta)\), and
\(z\in N_uM_i\) such that \(x=p+z\).  By
\eqref{eq_uniform_band_hessian_comparison},
\[
z^\top\partial_x^2g(\theta,p+tz)z\ge\frac{\mu}{2}\|z\|^2
\qquad(0\le t\le1).
\]
Since \(g(\theta,p)=g^\star(\theta)\) and
\(\partial_xg(\theta,p)=0\), the integral Taylor identity gives
\[
\begin{aligned}
h_\theta(x)
&=\int_0^1(1-t)
   z^\top\partial_x^2g(\theta,p+tz)z\,\mathrm dt\\
&\ge \frac{\mu}{2}\|z\|^2\int_0^1(1-t)\,\mathrm dt
=\frac{\mu}{4}\|z\|^2
\ge\frac{\mu}{4}\operatorname{dist}^2
       \bigl(x,\mathcal S(\theta)\bigr).
\end{aligned}
\]
This proves the first inequality in \eqref{eq_uniform_energy_gap}.  The
implication \eqref{eq_uniform_PL_band}, by contraposition, gives
\(h_\theta(x)>\varepsilon_{\mathrm{PL}}\) when
\(x\notin\mathcal W_\theta\).  Splitting according to membership in
\(\mathcal W_\theta\) proves \eqref{eq_uniform_energy_gap_at_distance}.

Set
\[
    \widehat Z_\lambda(\theta)
    :=\int_{\mathbb R^d}e^{-h_\theta(x)/\lambda}\,\mathrm dx.
\]
Choose any \(u_\theta\in\mathcal S(\theta)\).  The global
\(L_{g,xx}\)-Lipschitz continuity of \(\partial_xg(\theta,\cdot)\), together
with \(h_\theta(u_\theta)=0\) and
\(\partial_xg(\theta,u_\theta)=0\), yields
\[
    0\le h_\theta(u_\theta+y)\le\frac{L_{g,xx}}2\|y\|^2.
\]
Hence, for \(0<\lambda\le1\),
\begin{equation}\label{eq_uniform_localization_crude_normalizer}
\begin{aligned}
\widehat Z_\lambda(\theta)
&\ge
\int_{\mathbb B_d(u_\theta;\sqrt\lambda)}
e^{-h_\theta(x)/\lambda}\,\mathrm dx\\
&\ge
e^{-L_{g,xx}/2}\operatorname{vol}\bigl(\mathbb B_d(0;1)\bigr)
\lambda^{d/2}
=:c_Z\lambda^{d/2}.
\end{aligned}
\end{equation}

We next bound a shifted exponential moment.  By
\eqref{eq_working_radius_D}, every solution set lies in
\(\mathbb B_d(0;\mathsf D)\), and
\eqref{eq_tail_EB_implies_QG} applies outside that ball.  Thus, when
\(\|x\|\ge2\mathsf D\),
\[
\operatorname{dist}\bigl(x,\mathcal S(\theta)\bigr)
\ge\|x\|-\mathsf D\ge\frac12\|x\|,
\qquad
h_\theta(x)\ge\frac{\mu_{\mathsf{qg}}}{8}\|x\|^2.
\]
It follows that
\begin{equation}\label{eq_uniform_localization_exponential_moment}
\sup_{\substack{\theta\in\Theta\\0<\lambda\le1}}
\int_{\mathbb R^d}(1+\|x\|)
e^{-h_\theta(x)/(2\lambda)}\,\mathrm dx
\le C_{\mathrm{mom}}<\infty.
\end{equation}
Indeed, on \(\mathbb B_d(0;2\mathsf D)\) the integrand is bounded by
\(1+2\mathsf D\), while outside this ball it is bounded by
\((1+\|x\|)e^{-\mu_{\mathsf{qg}}\|x\|^2/16}\), an integrable function
independent of \(\theta\) and \(\lambda\).

On the set in \eqref{eq_uniform_weighted_localization},
\[
e^{-h_\theta(x)/\lambda}
\le e^{-a_r/(2\lambda)}e^{-h_\theta(x)/(2\lambda)}.
\]
Dividing the integral of this inequality by
\eqref{eq_uniform_localization_crude_normalizer} and using
\eqref{eq_uniform_localization_exponential_moment} proves
\eqref{eq_uniform_weighted_localization} with
\(C_{\mathrm{loc}}=C_{\mathrm{mom}}/c_Z\).  Its right-hand side tends to zero
because exponential decay dominates \(\lambda^{-d/2}\).  Replacing
\(1+\|x\|\) throughout by \(1+\|x\|^s\) proves the final assertion for each
fixed \(s\ge0\).
\end{proof}
\begin{proof}[Proof of \Cref{prop_W1_Gibbs_uniform}]
Put \(q:=d-\dimk\).
For \(u\in\mathcal S(\theta)\), let
\(H_\theta^\perp(u)\) be the self-adjoint operator on
\(N_u\mathcal S(\theta)\) determined by
\(\langle v,H_\theta^\perp(u)w\rangle
=\langle v,\partial_x^2g(\theta,u)w\rangle\).
In a local orthonormal normal frame with column matrix \(E(u)\), this operator
is represented by
\(E(u)^\top\partial_x^2g(\theta,u)E(u)\).  An orthogonal change of frame
conjugates this matrix and therefore preserves its determinant, which agrees
with the determinant in \eqref{eq_limiting_gibbs_weight}.
The local
P\L--Morse--Bott result, with the normalization in \Cref{def_local_PL}, gives
the lower normal-curvature bound with the same constant \(\mu\).  Global
\(L_{g,xx}\)-Lipschitz continuity of the lower gradient gives the upper bound.
Thus
\begin{equation}\label{eq_uniform_normal_hessian_bounds}
    \mu I_q\preceq H_\theta^\perp(u)
    \preceq L_{g,xx}I_q
    \qquad
    (\theta\in\Theta,\ u\in\mathcal S(\theta)).
\end{equation}
There is no factor of two in the lower bound: applying the P\L\ inequality to
a second-order expansion along a normal eigenvector gives
\(\alpha^2/(2\mu)\ge\alpha/2\), hence every positive normal eigenvalue
\(\alpha\) is at least \(\mu\).

For \(n\in N_u\mathcal S(\theta)\), define
\[
    Q_\theta(u,n):=\frac12
       \langle n,H_\theta^\perp(u)n\rangle.
\]
Since \(h_\theta(u)=0\) and \(\partial_xg(\theta,u)=0\), Taylor's integral
identity applies.  Moreover, if \(\|n\|\le r_3\) and \(t\in[0,1]\), then both
\(u\) and \(u+tn\) lie within distance \(r_3\) of \(\mathcal S(\theta)\), so
\eqref{eq_uniform_normal_hessian_lipschitz} gives
\(
\|\partial_x^2g(\theta,u+tn)-\partial_x^2g(\theta,u)\|_{\mathrm{op}}
\le L_{g,3}t\|n\|
\).  Hence
\begin{align}
h_\theta(u+n)-Q_\theta(u,n)
&=\int_0^1(1-t)
  \bigl\langle n,
  [\partial_x^2g(\theta,u+tn)-\partial_x^2g(\theta,u)]n
  \bigr\rangle\,\mathrm dt,
\nonumber\\
|h_\theta(u+n)-Q_\theta(u,n)|
&\le\frac{L_{g,3}}6\|n\|^3.
\label{eq_uniform_cubic_remainder}
\end{align}
By \Cref{lem_uniform_ambient_tube}, choose
\begin{equation}\label{eq_uniform_fiber_radius}
0<r_\ast\le\min\{r_{\mathrm{tub}},r_3\}
\quad\text{so that}\quad
\frac{L_{g,3}r_\ast}{6}\le\frac{\mu}{4}.
\end{equation}
Then \eqref{eq_uniform_normal_hessian_bounds} and
\eqref{eq_uniform_cubic_remainder} imply the common lower bound
\begin{equation}\label{eq_uniform_fiber_quadratic_lower}
    h_\theta(u+n)\ge\frac{\mu}{4}\|n\|^2
    \qquad(\|n\|<r_\ast).
\end{equation}

For \(u\in\mathcal S(\theta)\), let
\begin{equation}\label{eq_uniform_fiber_integral}
A_{\lambda,\theta}(u)
:=\int_{\substack{n\in N_u\mathcal S(\theta)\\\|n\|<r_\ast}}
e^{-h_\theta(u+n)/\lambda}J_\theta(u,n)\,\mathrm dn.
\end{equation}
Set \(c_q:=(2\pi)^{q/2}\) and
\(w_\theta(u):=\det H_\theta^\perp(u)^{-1/2}\).  We claim that there is
\(C_A<\infty\), independent of \((\theta,u)\), such that
\begin{equation}\label{eq_uniform_fiber_expansion}
\left|\lambda^{-q/2}A_{\lambda,\theta}(u)-c_qw_\theta(u)\right|
\le C_A\sqrt\lambda
\qquad(0<\lambda\le\lambda_A)
\end{equation}
for one common \(\lambda_A>0\).

To prove the claim, set \(n=\sqrt\lambda z\) and write
\(b_\theta(u,z):=\frac12\langle z,H_\theta^\perp(u)z\rangle\).  Then
\[
\lambda^{-q/2}A_{\lambda,\theta}(u)
=\int_{\|z\|<r_\ast/\sqrt\lambda}
e^{-h_\theta(u+\sqrt\lambda z)/\lambda}
J_\theta(u,\sqrt\lambda z)\,\mathrm dz,
\]
and \eqref{eq_uniform_cubic_remainder} gives
\begin{equation}\label{eq_uniform_rescaled_cubic_error}
\left|\lambda^{-1}h_\theta(u+\sqrt\lambda z)
      -b_\theta(u,z)\right|
\le\frac{L_{g,3}}6\sqrt\lambda\|z\|^3.
\end{equation}
For nonnegative \(a,b\),
\(|e^{-a}-e^{-b}|\le|a-b|(e^{-a}+e^{-b})\).  Combining this inequality with
\eqref{eq_uniform_tubular_jacobian},
\eqref{eq_uniform_normal_hessian_bounds}, and
\eqref{eq_uniform_fiber_quadratic_lower}, the absolute difference between the
rescaled integrand and \(e^{-b_\theta(u,z)}\) is bounded by
\begin{equation}\label{eq_uniform_fiber_integrable_majorant}
C\sqrt\lambda(\|z\|+\|z\|^3)
\left(e^{-\mu\|z\|^2/4}+e^{-\mu\|z\|^2/2}\right),
\end{equation}
where \(C\) is independent of \(\theta,u,\lambda\).  The right-hand side is
integrable on \(\mathbb R^q\).  The Gaussian integral is
\[
\int_{\mathbb R^q}e^{-b_\theta(u,z)}\,\mathrm dz
=c_qw_\theta(u),
\]
and the omitted Gaussian tail
\(\{\|z\|\ge r_\ast/\sqrt\lambda\}\) is bounded by
\(C e^{-c/\lambda}\), hence by \(C\sqrt\lambda\) for all sufficiently small
\(\lambda\).  Integrating \eqref{eq_uniform_fiber_integrable_majorant} proves
\eqref{eq_uniform_fiber_expansion}, uniformly in \((\theta,u)\).

Define
\begin{equation}\label{eq_uniform_leading_normalizer}
    W_\theta:=\int_{\mathcal S(\theta)}w_\theta(u)
                    \,\mathrm d\mathcal M_\theta(u).
\end{equation}
The eigenvalue bounds \eqref{eq_uniform_normal_hessian_bounds} and
\Cref{lem_vol_sol_set_is_bounded} give
\begin{equation}\label{eq_uniform_leading_normalizer_bounds}
    \mathsf b L_{g,xx}^{-q/2}
    \le W_\theta
    \le\mathsf b^{-1}\mu^{-q/2}.
\end{equation}
Integrating \eqref{eq_uniform_fiber_expansion} over
\(\mathcal S(\theta)\) therefore gives the tube contribution
\[
\lambda^{q/2}\bigl[c_qW_\theta+O(\sqrt\lambda)\bigr],
\]
uniformly in \(\theta\).  On the complement of the radius-\(r_\ast\) tube,
\Cref{corr_uniform_gibbs_localization}, together with
\eqref{eq_uniform_energy_gap_at_distance} and
\eqref{eq_uniform_localization_exponential_moment} give the unnormalized
bound
\[
\int_{\{\operatorname{dist}(x,\mathcal S(\theta))\ge r_\ast\}}
e^{-h_\theta(x)/\lambda}\,\mathrm dx
\le C_{\mathrm{mom}}e^{-a_{r_\ast}/(2\lambda)}.
\]
This is \(\lambda^{q/2}O(\sqrt\lambda)\) uniformly for small \(\lambda\).
Consequently, the shifted partition function satisfies
\begin{equation}\label{eq_uniform_sharp_partition_function}
\widehat Z_\lambda(\theta)
=\lambda^{q/2}\left[c_qW_\theta+O(\sqrt\lambda)\right]
\end{equation}
uniformly in \(\theta\).  After decreasing a common \(\lambda_0\),
\begin{equation}\label{eq_uniform_sharp_normalizer_lower}
    \widehat Z_\lambda(\theta)\ge c_0\lambda^{q/2}
    \qquad(\theta\in\Theta,\ 0<\lambda\le\lambda_0)
\end{equation}
for some \(c_0>0\).

Let \(\mathcal U_\theta^\ast\) be the radius-\(r_\ast\) tube, put
\(p_{\lambda,\theta}:=\gibbs_\theta^\lambda(\mathcal U_\theta^\ast)\), and
let \(\nu_{\lambda,\theta}\) be the Gibbs measure conditioned on this tube.
Its nearest-point projection
\(\sigma_{\lambda,\theta}:=(\pi_\theta)_\#\nu_{\lambda,\theta}\) has density
\begin{equation}\label{eq_uniform_projected_density}
\frac{\mathrm d\sigma_{\lambda,\theta}}{\mathrm d\mathcal M_\theta}(u)
=\frac{A_{\lambda,\theta}(u)}
       {\displaystyle\int_{\mathcal S(\theta)}A_{\lambda,\theta}(v)
                    \,\mathrm d\mathcal M_\theta(v)}.
\end{equation}
Write
\(\rho_\theta=w_\theta/W_\theta\).  By
\eqref{eq_uniform_fiber_expansion}, the \(L^1\)-norm of the error in
\(\lambda^{-q/2}A_{\lambda,\theta}=c_qw_\theta+e_{\lambda,\theta}\) is at
most \(C_A\mathsf b^{-1}\sqrt\lambda\).  The denominator
\(c_qW_\theta\) is bounded away from zero by
\eqref{eq_uniform_leading_normalizer_bounds}.  The elementary normalization
estimate
\[
\left\|\frac{b+e}{\int(b+e)}-\frac{b}{\int b}\right\|_{L^1}
\le\frac{2\|e\|_{L^1}}{\int b-\|e\|_{L^1}}
\]
therefore yields
\begin{equation}\label{eq_uniform_projected_density_L1}
\left\|
\frac{\mathrm d\sigma_{\lambda,\theta}}{\mathrm d\mathcal M_\theta}
-\rho_\theta
\right\|_{L^1(\mathcal M_\theta)}
\le C_\rho\sqrt\lambda
\end{equation}
uniformly in \(\theta\).  Since
\(\operatorname{diam}\mathcal S(\theta)\le\mathsf D\), total variation and a
maximal coupling give
\begin{equation}\label{eq_uniform_projected_W1}
\mathbb W_1(\sigma_{\lambda,\theta},\gibbs_\theta^0)
\le\frac{\mathsf D}{2}C_\rho\sqrt\lambda.
\end{equation}

The cost of projection in the normal direction is of the same order.  Indeed,
\eqref{eq_uniform_fiber_quadratic_lower}, the uniform bound on
\(J_\theta\), and the volume upper bound imply
\begin{equation}\label{eq_uniform_normal_first_moment}
\int_{\mathcal S(\theta)}
\int_{\substack{n\in N_u\mathcal S(\theta)\\\|n\|<r_\ast}}
\|n\|e^{-h_\theta(u+n)/\lambda}J_\theta(u,n)
\,\mathrm dn\,\mathrm d\mathcal M_\theta(u)
\le C_N\lambda^{(q+1)/2}.
\end{equation}
Division by \eqref{eq_uniform_sharp_normalizer_lower} gives a uniform
\(O(\sqrt\lambda)\) contribution.

It remains to account for all mass outside the tube.  Choose any
\(u_\theta^0\in\mathcal S(\theta)\) and define
\[
T_\theta(x):=
\begin{cases}
\pi_\theta(x),&x\in\mathcal U_\theta^\ast,\\
u_\theta^0,&x\notin\mathcal U_\theta^\ast.
\end{cases}
\]
Then \((\operatorname{Id},T_\theta)_\#\gibbs_\theta^\lambda\) is a coupling
whose second marginal is exactly
\begin{equation}\label{eq_uniform_complete_image_measure}
\eta_{\lambda,\theta}
:=(T_\theta)_\#\gibbs_\theta^\lambda
=p_{\lambda,\theta}\sigma_{\lambda,\theta}
 +(1-p_{\lambda,\theta})\delta_{u_\theta^0}.
\end{equation}
The tube part of this coupling is bounded by
\eqref{eq_uniform_normal_first_moment} divided by the normalizer.  For the
off-tube part, \(\|x-u_\theta^0\|\le\|x\|+\mathsf D\), and
\Cref{corr_uniform_gibbs_localization}, through
\eqref{eq_uniform_weighted_localization}, makes its cost exponentially small,
uniformly in \(\theta\).  Hence
\[
    \mathbb W_1(\gibbs_\theta^\lambda,\eta_{\lambda,\theta})
    \le C\sqrt\lambda.
\]
Finally, mix a coupling between \(\sigma_{\lambda,\theta}\) and
\(\gibbs_\theta^0\) with the coupling
\(\delta_{u_\theta^0}\otimes\gibbs_\theta^0\).  This gives
\[
\mathbb W_1(\eta_{\lambda,\theta},\gibbs_\theta^0)
\le p_{\lambda,\theta}
      \mathbb W_1(\sigma_{\lambda,\theta},\gibbs_\theta^0)
 +(1-p_{\lambda,\theta})\mathsf D
\le C\sqrt\lambda,
\]
where the last term is again controlled by
\eqref{eq_uniform_weighted_localization}.  Both couplings have the stated
probability marginals; no subprobability is left unmatched.  The triangle
inequality proves \eqref{eq_uniform_W1_Gibbs_rate} with constants \(K_0\) and
\(\lambda_0\) independent of \(\theta\).
\end{proof}

\subsection{Proof of \Cref{prop_uniform_poincare}}
\label[appendix]{append_proof_uniform_poincare}
The argument below is a parameter-uniform adaptation of
the Lyapunov and shrinking-tube construction of Gong et
al.~\citep{gong2024poincare}: it derives common constants for the whole
parameter family rather than treating their fixed-potential theorem as a
uniform result.
\begin{proof}
For each fixed \(\theta\), our assumptions provide the four inputs used by
Gong et al.~\citep[Assumptions~1--4]{gong2024poincare}.  First,
\Cref{mu_PL_assum} gives the required local regularity and the source's local
P\L\ inequality with constant \(2\mu\).  Second, its no-saddle clause gives a
negative-definite Hessian at every stationary point outside the local-minimum
neighborhoods.  Third, \eqref{eq_uniform_tail_gradient_error_bound} is the
source's tail error bound, while global lower-gradient smoothness gives
\(\|\partial_x^2g\|_{\mathrm{op}}\le L_{g,xx}\) and hence
\(|\Delta_xg|\le dL_{g,xx}\); the derived tail quadratic growth gives
coercivity.  Finally, in the regime \(\dimk\ge1\),
\Cref{prop_S_path_connected,prop_S_C_2_wihout,assum_secon_fund_form_bound}
give a compact connected boundaryless
\(\mathcal C^2\) minimizer manifold with bounded second fundamental form.
Finite local orthonormal normal frames identify this coordinate-free bound
with the local shape-matrix bound in the source.  The proof below establishes
the global tube comparison intrinsically on the normal bundle and tracks all
of these analytic and geometric constants uniformly in \(\theta\); it does
not invoke the source's fixed-potential theorem as a uniform black box.
The source uses the spectral-gap convention, so our inverse Poincar\'e
constant is \(C_{\mathrm{PI}}=\rho^{-1}\).

Write
\[
    h_\theta(x):=g(\theta,x)-g^\star(\theta),
    \qquad k:=\dimk,
    \qquad q:=d-k.
\]
Here \(k\ge1\) by hypothesis and \(q\ge1\), because a nonempty compact
boundaryless embedded submanifold of \(\mathbb R^d\) cannot have ambient
dimension \(d\).
All constants constructed below are independent of \(\theta\).
The corresponding fixed-potential steps are
\citep[Lemma~3, Lemma~4, Theorem~2, Proposition~6, and Theorem~4]
{gong2024poincare}; each required estimate is established below with constants
common to the parameter family.

\paragraph{A common minimum neighborhood.}
The function
\((\theta,x)\mapsto\operatorname{dist}(x,\mathcal S(\theta))\) is continuous:
this follows from the Hausdorff estimate
\eqref{eq_global_Hausdorff_Lipschitz}.  The solution graph
\(\mathfrak G\) defined in \eqref{eq_solution_graph} is
compact and is contained in the relatively open set \(\mathcal W\) from
\Cref{lem_uniform_local_PL_band}.  A contradiction argument using these two
facts gives \(r'>0\) such that
\[
 \bigl\{(\theta,x):
       \operatorname{dist}(x,\mathcal S(\theta))\le r'\bigr\}
 \subseteq\mathcal W.
\]
Indeed, otherwise there would be
\((\theta_j,x_j)\notin\mathcal W\) with
\(\operatorname{dist}(x_j,\mathcal S(\theta_j))\to0\).  Choosing nearest
points \(u_j\in\mathcal S(\theta_j)\), compactness of \(\mathfrak G\) would
give, along a subsequence,
\((\theta_j,u_j)\to(\bar\theta,\bar u)\in\mathfrak G\), and hence
\((\theta_j,x_j)\to(\bar\theta,\bar u)\in\mathcal W\), a contradiction to
relative openness of \(\mathcal W\).

Choose
\begin{equation}\label{eq_uniform_poincare_minimum_radius}
 0<r_{\mathcal S}<\min\{r',r_{\mathrm{tub}}\},
 \qquad
 \frac{L_{g,xx}}2r_{\mathcal S}^2
 \le\varepsilon_{\mathrm{PL}},
\end{equation}
and put
\[
 \mathcal U_{\mathcal S}
 :=\bigl\{(\theta,x):
   \operatorname{dist}(x,\mathcal S(\theta))<r_{\mathcal S}\bigr\}.
\]
For \((\theta,x)\) in the corresponding closed tube, choose the nearest
point \(u\in\mathcal S(\theta)\).  Global lower-gradient smoothness and
\(\partial_xg(\theta,u)=0\) imply
\[
 0\le h_\theta(x)
 \le\frac{L_{g,xx}}2\|x-u\|^2
 \le\frac{L_{g,xx}}2r_{\mathcal S}^2
 \le\varepsilon_{\mathrm{PL}}.
\]
Thus \eqref{eq_uniform_error_bound} yields the common local error bound
\begin{equation}\label{eq_uniform_poincare_local_EB}
 \|\partial_xg(\theta,x)\|
 \ge\nu_{\mathrm{loc}}
       \operatorname{dist}(x,\mathcal S(\theta)),
 \qquad
 \nu_{\mathrm{loc}}:=\frac\mu2,
 \qquad ((\theta,x)\in\mathcal U_{\mathcal S}).
\end{equation}
In the notation of the cited source, the distance from the boundary of this
minimum neighborhood to the minimizer manifold is therefore uniformly
\(\delta_0=r_{\mathcal S}\).

\paragraph{Uniform critical-region constants.}
Let \(\mathsf D_{\mathcal S}\) be the common solution bound furnished by
\Cref{prop_parametric_morse_bott_stability}, and set
\[
 R_0:=\max\{\mathsf D_{\mathsf{eb}},\mathsf D_{\mathcal S},1\}.
\]
The tail error bound shows that a critical point with \(\|x\|\ge R_0\)
belongs to \(\mathcal S(\theta)\).  We next isolate all remaining
nonminimum critical points in the closed set
\[
 \mathfrak X:=
 \bigl\{(\theta,x):\theta\in\Theta,\ \|x\|\le R_0,\
       \partial_xg(\theta,x)=0,\
       h_\theta(x)\ge\varepsilon_{\mathrm{PL}}\bigr\}.
\]
Here \(h_\theta(x)\) is continuous jointly in \((\theta,x)\).  To see the
only non-immediate point, if \(\theta_j\to\theta\) and
\(u_j\in\mathcal S(\theta_j)\), compactness of the solution graph and
continuity of \(g\) show
\(g^\star(\theta_j)=g(\theta_j,u_j)\to g^\star(\theta)\).  Hence
\(\mathfrak X\) is a closed subset of the compact set
\(\Theta\times\overline{\mathbb B_d(0;R_0)}\), and is compact.

Every nonminimum critical point lies in \(\mathfrak X\).  Indeed, the tail
argument places it in the fixed ball, while a critical point with
\(h_\theta(x)\le\varepsilon_{\mathrm{PL}}\) satisfies
\eqref{eq_uniform_local_PL} and hence \(h_\theta(x)=0\).  Moreover, a
stationary point in one of the local-minimum P\L\ neighborhoods is a local
minimum and therefore, by
\Cref{prop_S_path_connected,prop_S_C_2_wihout}, belongs to
\(\mathcal S(\theta)\).  Thus every point of \(\mathfrak X\) is in the
region where the no-saddle clause applies, and
\[
 \partial_x^2g(\theta,x)\prec0
 \qquad ((\theta,x)\in\mathfrak X).
\]

If \(\mathfrak X\ne\varnothing\), continuity of the largest Hessian
eigenvalue on this compact set gives
\[
 \alpha_-:=-\max_{(\theta,x)\in\mathfrak X}
       \lambda_{\max}(\partial_x^2g(\theta,x))>0.
\]
Joint Hessian continuity then supplies a relatively open neighborhood
\(\mathcal U_-\) of \(\mathfrak X\) and, after a harmless reduction, a
constant
\[
 0<\mu_-\le\min\left\{\frac{\alpha_-}{2},
                         \frac{M_\Delta}{2d}\right\}
 \quad\text{such that}\quad
 \partial_x^2g(\theta,x)\preceq-\mu_-I
 \quad((\theta,x)\in\mathcal U_-).
\]
Here \(M_\Delta:=dL_{g,xx}>0\).  If \(\mathfrak X=\varnothing\), take
\(\mathcal U_-:=\varnothing\) and any bookkeeping value
\(0<\mu_-\le M_\Delta/(2d)\).  The stronger reduction
\(d\mu_-\le M_\Delta/2\) is used below to match the compact-region drift
constant in the source.

The relatively open sets \(\mathcal U_{\mathcal S}\) and \(\mathcal U_-\)
contain every critical point in
\(\Theta\times\overline{\mathbb B_d(0;2R_0)}\).  Therefore
\[
 \mathcal K:=
 \bigl(\Theta\times\overline{\mathbb B_d(0;2R_0)}\bigr)
 \setminus(\mathcal U_{\mathcal S}\cup\mathcal U_-)
\]
is compact and contains no critical point.  If \(\mathcal K\ne\varnothing\),
set
\[
 g_0:=\min_{(\theta,x)\in\mathcal K}
       \|\partial_xg(\theta,x)\|>0;
\]
if it is empty, take any \(g_0>0\), since the corresponding region is
absent.  This proves the required compact-region gradient gap without
inferring uniformity merely from pointwise positivity.

Global lower-gradient smoothness and \(g\in\mathcal C^2\) give, everywhere,
\[
 \|\partial_x^2g(\theta,x)\|_{\mathrm{op}}\le L_{g,xx},
 \qquad
 |\Delta_xh_\theta(x)|\le M_\Delta:=dL_{g,xx}.
\]
We use the common notation
\[
 L:=L_{g,xx},\qquad
 C_\star:=\frac{4M_\Delta}{\nu_{\mathrm{loc}}^2}
          =\frac{16dL_{g,xx}}{\mu^2}.
\]
Since \(R_0\ge1\), the source's tail Laplacian bound also holds in the form
\(|\Delta h_\theta(x)|\le M_\Delta\|x\|^2\) for \(\|x\|\ge R_0\).  No third
derivative or uniform \(\mathcal C^3\) bound is used in these estimates.

\paragraph{A uniform spectral gap on the minimizer manifolds.}
Use the finite normal-graph cover constructed in the proof of
\Cref{lem_vol_sol_set_is_bounded}.  For a covering index \(i\), write
\[
 F_{i,\theta}:M_i\longrightarrow\mathcal S(\theta),
 \qquad M_i:=\mathcal S(\theta_i),
\]
where all tangent singular values of \(DF_{i,\theta}\) lie in
\([a,b]=[1/2,3/2]\), and its tangential Jacobian lies in
\([a^k,b^k]\).  For \(\varphi\in H^1(\mathcal S(\theta))\) and
\(\psi:=\varphi\circ F_{i,\theta}\), change of variables gives
\begin{align*}
 \int_{\mathcal S(\theta)}\|\nabla\varphi\|^2
       \,\mathrm d\mathcal M_\theta
 &\ge a^kb^{-2}\int_{M_i}\|\nabla\psi\|^2
       \,\mathrm d\mathcal M_i,\\
 \inf_{c\in\mathbb R}\int_{\mathcal S(\theta)}|\varphi-c|^2
       \,\mathrm d\mathcal M_\theta
 &\le b^k\inf_{c\in\mathbb R}\int_{M_i}|\psi-c|^2
       \,\mathrm d\mathcal M_i.
\end{align*}
The variational characterization of the first nonzero Laplace--Beltrami
eigenvalue consequently yields
\[
 \lambda_1(\mathcal S(\theta))
 \ge\frac{a^k}{b^{k+2}}\lambda_1(M_i)
 =\frac{4}{3^{k+2}}\lambda_1(M_i).
\]
Each reference manifold is compact, connected, and positive-dimensional.
Consequently, \(\lambda_1(M_i)>0\).  Taking the minimum over the finite
family gives
\begin{equation}\label{eq_uniform_manifold_spectral_gap}
 \lambda_1(\mathcal S(\theta))\ge
 \lambda_{\mathcal S}^\star
 :=\frac{4}{3^{k+2}}
       \min_{1\le i\le N}\lambda_1(M_i)>0
 \qquad(\theta\in\Theta).
\end{equation}
The reverse comparison, obtained by reversing the preceding two estimates,
is
\[
 \lambda_1(\mathcal S(\theta))
 \le\frac{b^k}{a^{k+2}}\lambda_1(M_i).
\]
Thus
\[
 \Lambda_{\mathcal S}^\star
 :=\frac{b^k}{a^{k+2}}
       \max_{1\le i\le N}\lambda_1(M_i)<\infty
\]
is a common upper bound.  This argument uses tangential Jacobians rather
than Hausdorff convergence or unproved spectral continuity.

\paragraph{Uniform Neumann spectral stability of the tubes.}
For \(0<r<r_{\mathrm{tub}}\), let
\[
 \mathcal T_\theta(r):=
 \{x:\operatorname{dist}(x,\mathcal S(\theta))<r\}.
\]
Pull this tube back by the normal addition map
\(\Phi_\theta(u,n)=u+n\).  The normal connection splits the tangent space of
the radius-\(r\) normal disk bundle into horizontal and vertical parts.  If
\(v\) is horizontal and \(w\) vertical, the Weingarten identity gives
\[
 D\Phi_\theta(u,n)[v,w]=(I-A_n^\theta)v+w.
\]
Since the two terms on the right are tangent and normal, respectively, and
\(\|A_n^\theta\|_{\mathrm{op}}\le\mathsf C\|n\|\), after requiring
\(\mathsf Cr<1\) we have
\begin{align}
 (1-\mathsf Cr)^2(\|v\|^2+\|w\|^2)
 &\le\|D\Phi_\theta(u,n)[v,w]\|^2
 \le(1+\mathsf Cr)^2(\|v\|^2+\|w\|^2),
 \label{eq_uniform_poincare_tube_metric}\\
 (1-\mathsf Cr)^k
 &\le J_\theta(u,n)=\det(I-A_n^\theta)
 \le(1+\mathsf Cr)^k.
 \label{eq_uniform_poincare_tube_density}
\end{align}
Both estimates are intrinsic and uniform in \(\theta,u,n\).

We first prove a Poincar\'e inequality for the connection metric and product
density.  Normalize them to the probability measure
\[
 \mathrm d\nu^0_{\theta,r}(u,n)
 :=\frac{\mathrm d\mathcal M_\theta(u)\,\mathrm dn}
 {\mathcal M_\theta(\mathcal S(\theta))\,\omega_qr^q}.
\]
For a smooth function \(F\) on the normal disk bundle, define intrinsically
\[
 \overline F(u):=\frac1{\omega_qr^q}
     \int_{\substack{n\in N_u\mathcal S(\theta)\\\|n\|<r}}
          F(u,n)\,\mathrm dn.
\]
Conditional variance gives the exact decomposition
\begin{align*}
 \operatorname{Var}_{\nu^0_{\theta,r}}(F)
 &=\frac1{\mathcal M_\theta(\mathcal S(\theta))}
   \int_{\mathcal S(\theta)}
      \operatorname{Var}_{\mathrm{Unif}(B_{N_u}(0,r))}
         (F(u,\cdot))\,\mathrm d\mathcal M_\theta(u)\\
 &\quad+
   \operatorname{Var}_{\mathcal M_\theta/
      \mathcal M_\theta(\mathcal S(\theta))}(\overline F).
\end{align*}
The Euclidean-ball Poincar\'e inequality bounds the first term by
\[
 c_qr^2\int\|\nabla_{\mathrm{vert}}F\|^2
       \,\mathrm d\nu^0_{\theta,r},
\]
where \(c_q\) depends only on \(q\).  To control the second term, identify
nearby normal fibers by parallel transport in the metric normal connection.
Parallel transport is orthogonal, hence maps the fiber ball onto the fiber
ball and preserves its Lebesgue measure.  Differentiation followed by
Jensen's inequality therefore gives
\[
 \|\nabla_{\mathcal S}\overline F(u)\|^2
 \le\frac1{\omega_qr^q}
   \int_{\|n\|<r}\|\nabla_{\mathrm{hor}}F(u,n)\|^2\,\mathrm dn.
\]
Together with the base Poincar\'e inequality, this proves
\begin{align*}
 \operatorname{Var}_{\nu^0_{\theta,r}}(F)
 \le{}&c_qr^2\int\|\nabla_{\mathrm{vert}}F\|^2
             \,\mathrm d\nu^0_{\theta,r}\\
 &+\lambda_1(\mathcal S(\theta))^{-1}
       \int\|\nabla_{\mathrm{hor}}F\|^2
             \,\mathrm d\nu^0_{\theta,r}.
\end{align*}
Choose a common preliminary radius with
\(c_qr^2\le(\Lambda_{\mathcal S}^\star)^{-1}\).  Since
\(\lambda_1(\mathcal S(\theta))\le\Lambda_{\mathcal S}^\star\), the
right-hand side is bounded by
\(\lambda_1(\mathcal S(\theta))^{-1}\) times the connection-metric
Dirichlet energy.  The calculation extends from smooth functions to
\(H^1\) by local trivialization and the usual density of smooth functions on
the compact disk bundle with boundary.

We now transfer this inequality through \(\Phi_\theta\).  The metric
comparison \eqref{eq_uniform_poincare_tube_metric}, the density comparison
\eqref{eq_uniform_poincare_tube_density}, and the variational formula give
the explicit lower estimate
\[
 \lambda_1^{\mathrm N}(\mathcal T_\theta(r))
 \ge
 \frac{(1-\mathsf Cr)^k}{(1+\mathsf Cr)^{k+2}}
 \lambda_1(\mathcal S(\theta)).
\]
For clarity, the factors arise by bounding the variance with the largest
density and the Dirichlet energy with the smallest density and inverse-metric
factor.  Testing with functions pulled back from the base similarly gives
\[
 \lambda_1^{\mathrm N}(\mathcal T_\theta(r))
 \le
 \frac{(1+\mathsf Cr)^k}{(1-\mathsf Cr)^{k+2}}
 \lambda_1(\mathcal S(\theta)).
\]
For \(\mathsf Cr\le1/(k+2)\), Bernoulli's inequality and convexity of
\(s\mapsto(1+s)^{-(k+2)}\) give
\[
 \frac{(1-\mathsf Cr)^k}{(1+\mathsf Cr)^{k+2}}
 \ge(1-k\mathsf Cr)(1-(k+2)\mathsf Cr)
 \ge1-(2k+2)\mathsf Cr.
\]
Set \(B_{\mathrm{spec}}:=(2k+2)\mathsf C\).  Since the upper rational
factor also converges to one as \(r\downarrow0\), there is a common
\(r_{\mathrm{spec}}>0\), with \(r_{\mathrm{spec}}<r_{\mathrm{tub}}\),
such that, for \(0<r\le r_{\mathrm{spec}}\),
\begin{equation}\label{eq_uniform_tube_neumann_gap}
 \lambda_1^{\mathrm N}(\mathcal T_\theta(r))
 \ge(1-B_{\mathrm{spec}}r)\lambda_1(\mathcal S(\theta))
 \ge\frac12\lambda_{\mathcal S}^\star,
\end{equation}
and also
\(\lambda_1^{\mathrm N}(\mathcal T_\theta(r))
 \le2\Lambda_{\mathcal S}^\star\).
Here we reduced \(r_{\mathrm{spec}}\) so that
\(B_{\mathrm{spec}}r_{\mathrm{spec}}\le1/2\).  We reduce it once more, if
necessary, so that
\begin{equation}\label{eq_uniform_poincare_tube_auxiliary_smallness}
 \frac{M_\Delta C_\star}{2r_{\mathrm{spec}}^2}
 \ge2e^{4LC_\star}\Lambda_{\mathcal S}^\star.
\end{equation}
This last reduction is used only in the Lyapunov combination below.  The
entire argument used the normal connection and fiberwise Lebesgue measure;
it assumes neither orientability nor a global normal frame or a trivial
normal bundle.

\paragraph{The uniform Lyapunov estimate.}
Let
\[
 W_{\theta,\lambda}(x):=
       \exp\!\left(\frac{h_\theta(x)}{2\lambda}\right),
 \qquad
 \mathcal L_{\theta,\lambda}
 :=-\partial_xh_\theta\cdot\nabla+\lambda\Delta.
\]
Since \(h_\theta\ge0\), one has
\(W_{\theta,\lambda}\ge1\).
Direct calculation gives
\begin{equation}\label{eq_uniform_poincare_generator_identity}
 \frac{\mathcal L_{\theta,\lambda}W_{\theta,\lambda}}
      {\lambda W_{\theta,\lambda}}
 =\frac{\Delta_xh_\theta}{2\lambda}
  -\frac{\|\partial_xh_\theta\|^2}{4\lambda^2}.
\end{equation}
Put
\[
 U_{\theta,\lambda}
 :=\mathcal T_\theta(\sqrt{C_\star\lambda}).
\]
We record the four-region estimate underlying Lemma~3 of the cited source.
First, if \(\|x\|\ge2R_0\), then
\(\operatorname{dist}(x,\mathcal S(\theta))\ge\|x\|/2\).  The tail EB and
tail Laplacian bound, together with
\(\lambda\le\nu_{\mathsf{eb}}^2/(64M_\Delta)\), imply from
\eqref{eq_uniform_poincare_generator_identity}
\[
 \frac{\mathcal LW}{\lambda W}
 \le-\frac{\nu_{\mathsf{eb}}^2R_0^2}{128\lambda^2}.
\]
Second, on
\(\mathcal U_{\mathcal S}\setminus U_{\theta,\lambda}\), the local EB and
the definition \(C_\star=4M_\Delta/\nu_{\mathrm{loc}}^2\) give
\[
 \frac{\mathcal LW}{\lambda W}
 \le\frac{M_\Delta}{2\lambda}
  -\frac{\nu_{\mathrm{loc}}^2C_\star\lambda}{4\lambda^2}
 =-\frac{M_\Delta}{2\lambda}.
\]
Third, on \(\mathcal U_-\), negative curvature gives
\[
 \frac{\mathcal LW}{\lambda W}
 \le-\frac{d\mu_-}{2\lambda}.
\]
Finally, on \(\mathcal K\), if that set is nonempty, the gradient gap and
\(\lambda\le g_0^2/(4M_\Delta)\) give
\[
 \frac{\mathcal LW}{\lambda W}
 \le\frac{M_\Delta}{2\lambda}
    -\frac{g_0^2}{4\lambda^2}
 \le-\frac{M_\Delta}{2\lambda}
 \le-\frac{d\mu_-}{2\lambda}.
\]
These regions cover \(\mathbb R^d\setminus U_{\theta,\lambda}\): outside
the radius-\(2R_0\) ball use the first region; inside it, use successively
\(\mathcal U_{\mathcal S}\), \(\mathcal U_-\), and their compact
complement \(\mathcal K\).  The containment
\(U_{\theta,\lambda}\subset\mathcal U_{\mathcal S}\) follows from
\(C_\star\lambda\le r_{\mathcal S}^2\).

Consequently, with
\[
 \sigma:=\min\left\{
   \frac{\nu_{\mathsf{eb}}^2R_0^2}{128\lambda^2},
   \frac{d\mu_-}{2\lambda}\right\},
 \qquad
 b:=\sigma+\frac{M_\Delta}{2\lambda},
\]
the four estimates give the drift bound outside \(U_{\theta,\lambda}\).
On \(U_{\theta,\lambda}\), the global Laplacian bound and the nonpositive
gradient term in \eqref{eq_uniform_poincare_generator_identity} give
\[
 \frac{\mathcal L_{\theta,\lambda}W_{\theta,\lambda}}
      {\lambda W_{\theta,\lambda}}
 \le\frac{M_\Delta}{2\lambda}=-\sigma+b.
\]
Thus, on all of \(\mathbb R^d\),
\[
 \frac{\mathcal L_{\theta,\lambda}W_{\theta,\lambda}}
      {\lambda W_{\theta,\lambda}}
 \le-\sigma+b\mathbf 1_{U_{\theta,\lambda}}.
\]
The restriction
\(\lambda\le\nu_{\mathsf{eb}}^2R_0^2/(64d\mu_-)\) selects the branch
\[
 \sigma=\frac{d\mu_-}{2\lambda}.
\]

It remains to combine this drift with the truncated Poincar\'e inequality.
On \(U_{\theta,\lambda}\), global smoothness and nearest-point projection
give
\[
 0\le h_\theta(x)\le\frac L2C_\star\lambda,
\]
so the oscillation of \(h_\theta/\lambda\) is at most
\(LC_\star/2\), and in particular at most
\(\overline C:=4LC_\star\).  The Holley--Stroock comparison used in
Lemma~4 of the cited source therefore gives, for the spectral gap
\(\rho_{\theta,\lambda}^{U}\) of the Gibbs measure conditioned on
\(U_{\theta,\lambda}\),
\[
 e^{-4LC_\star}\lambda_1^{\mathrm N}(U_{\theta,\lambda})
 \le\rho_{\theta,\lambda}^{U}
 \le e^{4LC_\star}\lambda_1^{\mathrm N}(U_{\theta,\lambda}).
\]
If \(\sqrt{C_\star\lambda}\le r_{\mathrm{spec}}\), the upper tube estimate
and \eqref{eq_uniform_poincare_tube_auxiliary_smallness} imply
\(\rho_{\theta,\lambda}^{U}\le b\).

We now give the precise Lyapunov-to-Poincar\'e step.  This is the normalized
form of the argument underlying
\citep[Proposition~1]{gong2024poincare}.  It is included because the drift
inequality established immediately above is naturally expressed after
division by \(W_{\theta,\lambda}\).

Abbreviate
\[
    \mu:=\mu_\theta^\lambda,
    \qquad U:=U_{\theta,\lambda},
    \qquad W:=W_{\theta,\lambda},
    \qquad \rho_U:=\rho_{\theta,\lambda}^{U}.
\]
For a smooth compactly supported function \(\varphi\), put
\[
    c_U:=\int_U\varphi\,\mathrm d\mu_U,
    \qquad F:=\varphi-c_U,
\]
where \(\mu_U:=\mu(\,\cdot\mid U)\).  Symmetry of
\(\mathcal L_{\theta,\lambda}\) in \(L^2(\mu)\) and integration by parts
give
\begin{align*}
\int F^2\frac{-\mathcal L_{\theta,\lambda}W}{\lambda W}
      \,\mathrm d\mu
&=\int\left\langle
      \nabla\!\left(\frac{F^2}{W}\right),\nabla W
    \right\rangle\,\mathrm d\mu\\
&=\int\|\nabla F\|^2\,\mathrm d\mu
  -\int\left\|\nabla F-\frac{F}{W}\nabla W\right\|^2
      \,\mathrm d\mu\\
&\le\int\|\nabla F\|^2\,\mathrm d\mu.
\end{align*}
The integration-by-parts identity is justified first with a spatial cutoff
and then by letting the cutoff radius tend to infinity; global gradient
smoothness and the derived quadratic tail bound dominate the resulting
terms.  The normalized drift inequality above is equivalent to
\[
    \sigma
    \le
    \frac{-\mathcal L_{\theta,\lambda}W}{\lambda W}
    +b\mathbf 1_U.
\]
Consequently,
\begin{align*}
\operatorname{Var}_{\mu}(\varphi)
&\le\int F^2\,\mathrm d\mu\\
&\le\frac1\sigma\int\|\nabla\varphi\|^2\,\mathrm d\mu
     +\frac b\sigma\int_U F^2\,\mathrm d\mu.
\end{align*}
Since \(c_U\) is the \(\mu_U\)-mean of \(\varphi\), the conditioned
Poincar\'e inequality yields
\[
    \int_U F^2\,\mathrm d\mu
    =\mu(U)\operatorname{Var}_{\mu_U}(\varphi)
    \le\frac1{\rho_U}\int_U\|\nabla\varphi\|^2\,\mathrm d\mu.
\]
It follows that
\[
    \operatorname{Var}_{\mu}(\varphi)
    \le
    \frac{b+\rho_U}{\sigma\rho_U}
    \int\|\nabla\varphi\|^2\,\mathrm d\mu,
\]
and hence
\[
    \rho_{\theta,\lambda}
    \ge
    \frac{\sigma}{b+\rho_{\theta,\lambda}^{U}}
    \rho_{\theta,\lambda}^{U}.
\]
Approximation by smooth compactly supported functions extends this inequality
to \(H^1(\mu)\).  Using \(\rho_{\theta,\lambda}^{U}\le b\), the selected
branch \(\sigma=d\mu_-/(2\lambda)\), the definition
\(b=\sigma+M_\Delta/(2\lambda)\), and the lower Holley--Stroock comparison
now gives
\begin{align*}
\rho_{\theta,\lambda}
&\ge\frac12\frac{d\mu_-}{d\mu_-+M_\Delta}
      \rho_{\theta,\lambda}^{U}\\
&\ge\frac12\frac{d\mu_-}{d\mu_-+M_\Delta}
      e^{-4LC_\star}
      \lambda_1^{\mathrm N}(U_{\theta,\lambda}).
\end{align*}
Thus the common Lyapunov constants are obtained without a uniform
third-derivative bound.

\paragraph{Common temperature and final constant.}
Define
\begin{equation}\label{eq_uniform_poincare_temperature}
 \lambda_{\mathrm{PI}}^\star
 :=\min\left\{
   \frac{\nu_{\mathsf{eb}}^2}{64M_\Delta},
   \frac{r_{\mathcal S}^2}{C_\star},
   \frac{g_0^2}{4M_\Delta},
   \frac{\nu_{\mathsf{eb}}^2R_0^2}{64d\mu_-},
   \frac{r_{\mathrm{spec}}^2}{C_\star}
 \right\}>0.
\end{equation}
Every entry is positive and depends only on the common constants constructed
above.  If
\(\theta\in\Theta\) and
\(0<\lambda\le\lambda_{\mathrm{PI}}^\star\), then
\eqref{eq_uniform_tube_neumann_gap} gives
\[
 \lambda_1^{\mathrm N}(U_{\theta,\lambda})
 \ge\frac12\lambda_{\mathcal S}^\star.
\]
Hence
\[
 \rho_{\theta,\lambda}
 \ge\rho_{\mathrm{PI}}^\star
 :=\frac14\frac{d\mu_-}{d\mu_-+M_\Delta}
       e^{-4L_{g,xx}C_\star}\lambda_{\mathcal S}^\star>0.
\]
With
\(C_{\mathrm{PI}}^\star:=(\rho_{\mathrm{PI}}^\star)^{-1}\), this is
exactly \eqref{eq_uniform_poincare_gap}.  Thus the quantifiers are
\[
 \exists\lambda_{\mathrm{PI}}^\star>0,\ 
         C_{\mathrm{PI}}^\star<\infty
 \quad\forall\theta\in\Theta
 \quad\forall\,0<\lambda\le\lambda_{\mathrm{PI}}^\star,
\]
as asserted.
\end{proof}

\section{Proofs of Results in \Cref{sec:3}}\label[appendix]{append_proof_sec_3}
This appendix provides proofs for the geometric and measure-approximation results in \Cref{sec:3}. We first establish a uniform comparison between the volume of small geodesic balls on $\mathcal{S}(\theta)$ and Euclidean balls, and then use Wasserstein stability to control the superquantile functional under perturbations of the underlying measure.
\subsection{Proof of \Cref{lemm_comparison_volume_geod_ball_euc_ball}}\label[appendix]{append_proof_lemm_VaR_approx}
\begin{proof}\label{proof_lemm_comparison_volume_geod_ball_euc_ball}
Fix \(\theta\in\Theta\) and \(p\in\mathcal S(\theta)\), and write
\(M:=\mathcal S(\theta)\).
\paragraph{Reduction and proof plan.}
The quantity to be estimated is the volume of a small geodesic ball.  Below
the injectivity radius, the exponential map sends a ball in \(T_pM\) onto
that geodesic ball, so the change-of-variables formula writes its volume as
the integral of the Jacobian \(J_{\theta,p}\) of \(\exp_p\).  Consequently,
the lemma will follow once we find one radius range, independent of
\(\theta\) and \(p\), on which
\[
    0<\mathsf c_L\le J_{\theta,p}(v)\le\mathsf c_H<\infty.
\]
The proof establishes these pointwise Jacobian bounds as follows.
\begin{enumerate}[label=\emph{Step \arabic*.},leftmargin=*,itemsep=0.4ex]
 \item \emph{Choose a common radius range.}
 The uniform curvature and injectivity-radius bounds provide a range on which
 the same construction works for every \(\theta\) and \(p\).  The
 one-dimensional case is dealt with separately because it has no angular
 directions.

 \item \emph{Reduce the volume estimate to a Jacobian estimate.}
 On this radius range, define \(J_{\theta,p}\) by pulling the Riemannian
 volume density back through \(\exp_p\).  The change-of-variables formula then
 expresses the desired volume as the integral of \(J_{\theta,p}\) over the
 tangent-space ball.

 \item \emph{Express that Jacobian through angular Jacobi fields.}
 At \(v=t\omega\), the radial derivative of \(\exp_p\) has length one.  The
 remaining \(\dimk-1\) derivatives describe what happens when the direction
 \(\omega\) is varied; these are transverse Jacobi fields.  Placing them in a
 matrix \(\mathcal Y(t)\) gives
 \[
     J_{\theta,p}(t\omega)
     =t^{-(\dimk-1)}\lvert\det\mathcal Y(t)\rvert.
 \]
 Thus only the angular stretching represented by \(\mathcal Y\) remains to
 be bounded.

 \item \emph{Use curvature to bound the relative growth of the Jacobi
 fields.}
 The matrix \(\mathcal Y\) satisfies
 \(\mathcal Y''+\mathcal R\mathcal Y=0\), and the sectional-curvature bounds
 give \(-\mathsf C_0I\preceq\mathcal R\preceq\mathsf C_0I\).
 The matrix \(S:=\mathcal Y'\mathcal Y^{-1}\) is the analogue
 of the scalar logarithmic derivative \(y'/y\): bounds on \(S\) integrate to
 bounds on \(\|\mathcal Y(t)z\|\).  We first show that \(S\) is symmetric and
 satisfies the first-order equation \(S'+S^2+\mathcal R=0\).  Its small-time
 expansion starts a comparison with the logarithmic growth rates in the
 constant positive- and negative-curvature models.  An
 \(\varepsilon\)-perturbation makes the initial inequalities strict, and a
 first-crossing argument proves that they persist.

 \item \emph{Recover bounds on the Jacobi fields and the Jacobian.}
 Integrating the bounds on \(S\) gives
 \[
 \frac{\sin(\sqrt{\mathsf C_0}t)}{\sqrt{\mathsf C_0}}\|z\|
 \le \|\mathcal Y(t)z\|
 \le
 \frac{\sinh(\sqrt{\mathsf C_0}t)}{\sqrt{\mathsf C_0}}\|z\|.
 \]
 These are singular-value bounds for \(\mathcal Y(t)\).  Multiplying them
 bounds its determinant, and the formula in Step~3 then bounds
 \(J_{\theta,p}(t\omega)\).

 \item \emph{Integrate the resulting Jacobian bounds.}
 Steps~3--5 give constants \(\mathsf c_L,\mathsf c_H\) with
 \(\mathsf c_L\le J_{\theta,p}\le\mathsf c_H\).  Inserting these inequalities
 into the change-of-variables formula from Step~2 proves the desired
 comparison with the Euclidean volume \(\omega_{\dimk}r^{\dimk}\).
\end{enumerate}

\noindent\emph{Step 1: uniform geometric inputs and the
one-dimensional case.}
By \Cref{assum_compact_opt_set},
\[
    -\mathsf C_0\le \operatorname{Sec}_M\le \mathsf C_0,
    \qquad
    \operatorname{inj}_M(p)\ge r_0,
\]
where \(\dimk\), \(\mathsf C_0\), and \(r_0\) are independent of \(\theta\).

If \(\dimk=1\), every geodesic ball with \(0<r<r_0\) is an
open interval of length \(2r=\omega_1r\).  The conclusion follows with
\(c_1=\mathsf c_L=\mathsf c_H=1\).  The case \(\dimk=0\) is excluded by the
standing assumption \(\dimk\ge1\).  Hence it remains only to consider
\(\dimk\ge2\).

\noindent\emph{Step 2: reduction of the volume estimate
to a Jacobian estimate.}
If \(0<r<r_0\), the normal-ball characterization gives a
\(\mathcal C^1\) diffeomorphism
\[
    \exp_p:
    \mathbb B_{T_pM}(0;r)
    \longrightarrow
    \mathbb B_{\mathsf g_\theta}(p;r).
\]
Its standard proof uses only minimizing geodesics and the
Gauss lemma, both available for the present \(\mathcal C^2\) metric; see
\citep[Section~5.9.2]{petersen2006riemannian}.
Using Riemannian
densities, which do not require orientability, define \(J_{\theta,p}(v)>0\) by
\[
    \exp_p^*(\ud\mathcal M_\theta)(v)
    =J_{\theta,p}(v)\,\ud v,
\]
where \(\ud v\) is Euclidean Lebesgue measure on \(T_pM\).  The change of
variables formula yields
\begin{equation}\label{eq_geod_vol_ball_normal}
    \mathcal M_\theta\bigl(\mathbb B_{\mathsf g_\theta}(p;r)\bigr)
    =
    \int_{\mathbb B_{T_pM}(0;r)}J_{\theta,p}(v)\,\ud v.
\end{equation}
Thus it remains to bound \(J_{\theta,p}\) above and below
by positive constants that do not depend on \(\theta\), \(p\), or \(v\) in
the chosen radius range.

\noindent\emph{Step 3: expression of the Jacobian through
angular Jacobi fields.}
As established in the proof of \Cref{assum_compact_opt_set}, \(M\) is
\(\mathcal C^3\), so its induced metric is \(\mathcal C^2\), its curvature is
continuous, and its exponential map is \(\mathcal C^1\).  Fix an arbitrary unit vector
\(\omega\in T_pM\), choose an orthonormal basis
\(\omega,e_1,\ldots,e_{\dimk-1}\), and set
\[
 \gamma(t):=\exp_p(t\omega),
 \qquad
 Y_i(t):=D(\exp_p)_{t\omega}[t e_i].
\]
Here \(D\) denotes the differential defined in
\Cref{append_additional_notations}.  Since \(T_pM\) is a vector space, its tangent
space at \(t\omega\) is canonically identified with \(T_pM\).  Under this
identification,
\[
 D(\exp_p)_{t\omega}:T_pM\longrightarrow T_{\gamma(t)}M.
\]
Equivalently,
\[
 Y_i(t)
 =\left.\frac{\ud}{\ud s}\right|_{s=0}
   \exp_p\bigl(t\omega+s\,t e_i\bigr).
\]
Let \(E_1(t),\ldots,E_{\dimk-1}(t)\) be the parallel transports of
\(e_1,\ldots,e_{\dimk-1}\) along \(\gamma\).  The fields \(Y_i\) are
orthogonal to \(\dot\gamma\), and we define their component matrix by
\[
    \mathcal Y_{ji}(t):=
    \left\langle Y_i(t),E_j(t)\right\rangle .
\]
In the same frame, define the curvature matrix
\[
    \mathcal R_{ij}(t)
    :=
    \left\langle
      R^M\!\bigl(E_j(t),\dot\gamma(t)\bigr)\dot\gamma(t),
      E_i(t)
    \right\rangle ,
\]
where \(R^M\) is the Riemann curvature tensor of \(M\), with the convention
appearing in the Jacobi equation.  The Jacobi equation and the initial
conditions now read
\[
 \begin{gathered}
 \mathcal Y''+\mathcal R\mathcal Y=0,
 \qquad \mathcal Y(0)=0,
 \qquad \mathcal Y'(0)=I_{\dimk-1},\\
 -\mathsf C_0I_{\dimk-1}\preceq\mathcal R(t)
 \preceq \mathsf C_0I_{\dimk-1},
 \end{gathered}
\]
where \(\mathcal R\) is continuous and symmetric.  Indeed, for \(z\ne0\),
\[
 \frac{z^{\mathsf T}\mathcal R(t)z}{\|z\|^2}
\]
is the sectional curvature of the plane spanned by \(\dot\gamma(t)\) and
\(\sum_jz_jE_j(t)\); hence the sectional-curvature bounds give exactly the
displayed matrix inequalities.

The identities
\[
 D(\exp_p)_{t\omega}[\omega]=\dot\gamma(t),
 \qquad
 D(\exp_p)_{t\omega}[e_i]=t^{-1}Y_i(t)
\]
and the Gauss lemma show that the radial image has length one and is
orthogonal to all the angular images.  Therefore,
\begin{equation}\label{eq_geodesic_jacobian_from_jacobi_matrix}
    J_{\theta,p}(t\omega)
    =t^{-(\dimk-1)}\lvert\det\mathcal Y(t)\rvert.
\end{equation}
Thus the required pointwise bounds on \(J_{\theta,p}\) reduce to bounds on
the singular values of \(\mathcal Y(t)\).

{\noindent\emph{Step 4: comparison of the angular stretching using
curvature.}}
We next compare the growth of \(\mathcal Y\) with the constant-curvature
models.  Since \(\mathcal Y(t)/t\to I_{\dimk-1}\), the matrix is invertible
for all sufficiently small \(t>0\).  Define its first singular time---the
first time at which some angular direction collapses---by
\[
    \tau:=\inf\{t>0:\det\mathcal Y(t)=0\},
    \qquad \inf\varnothing:=+\infty.
\]
Then \(\tau>0\).  On \((0,\tau)\), define the
matrix logarithmic derivative
\(S:=\mathcal Y'\mathcal Y^{-1}\).  To see that it is symmetric, set
\[
 \mathcal W(t)
 :=\mathcal Y(t)^{\mathsf T}\mathcal Y'(t)
   -\mathcal Y'(t)^{\mathsf T}\mathcal Y(t).
\]
The Jacobi equation and the symmetry of \(\mathcal R\) give
\[
 \mathcal W'
 =\mathcal Y^{\mathsf T}\mathcal Y''
   -\mathcal Y''^{\mathsf T}\mathcal Y
 =-\mathcal Y^{\mathsf T}\mathcal R\mathcal Y
   +\mathcal Y^{\mathsf T}\mathcal R^{\mathsf T}\mathcal Y
 =0.
\]
Because \(\mathcal W(0)=0\), we have
\(\mathcal Y^{\mathsf T}\mathcal Y'
=\mathcal Y'^{\mathsf T}\mathcal Y\), and multiplication by
\(\mathcal Y^{-\mathsf T}\) and \(\mathcal Y^{-1}\) shows that
\(S^{\mathsf T}=S\).  This symmetry is what allows us to compare \(S\) in
the semidefinite order.  Moreover,
\[
 S'
 =\mathcal Y''\mathcal Y^{-1}
   -\mathcal Y'\mathcal Y^{-1}\mathcal Y'\mathcal Y^{-1}
 =-\mathcal R-S^2.
\]
Thus \(S\) satisfies the matrix Riccati equation
\(S'+S^2+\mathcal R=0\), to which the semidefinite-order comparison below
applies.

Because \(S\) is not defined at \(t=0\), we first determine its small-time
asymptotics.  Integrating the Jacobi equation once and twice gives
\[
 \mathcal Y'(t)
 =I_{\dimk-1}-\int_0^t\mathcal R(s)\mathcal Y(s)\,\ud s,
 \qquad
 \mathcal Y(t)
 =tI_{\dimk-1}
  -\int_0^t(t-s)\mathcal R(s)\mathcal Y(s)\,\ud s.
\]
The initial conditions imply
\(\mathcal Y(s)=sI_{\dimk-1}+o(s)\); continuity of \(\mathcal R\) therefore
gives \(\mathcal R(s)\mathcal Y(s)=s\mathcal R(0)+o(s)\).  Using
\(\int_0^t s\,\ud s=t^2/2\) and
\(\int_0^t(t-s)s\,\ud s=t^3/6\), we obtain
\[
 \mathcal Y'(t)
 =I_{\dimk-1}-\frac{t^2}{2}\mathcal R(0)+o(t^2),
 \qquad
 \mathcal Y(t)
 =tI_{\dimk-1}-\frac{t^3}{6}\mathcal R(0)+o(t^3).
\]
Factoring \(t\) from \(\mathcal Y(t)\) and expanding its inverse gives
\[
 \mathcal Y(t)^{-1}
 =t^{-1}\left[I_{\dimk-1}
       +\frac{t^2}{6}\mathcal R(0)+o(t^2)\right]
\]
and hence
\[
 S(t)=t^{-1}I_{\dimk-1}-\frac{t}{3}\mathcal R(0)+o(t)
 \qquad(t\to0^+).
\]
All remainders are in operator norm.  This expansion initializes the matrix
comparison at sufficiently small positive times.  Fix
\(0<T<\min\{\tau,\pi/\sqrt{\mathsf C_0}\}\).  For every sufficiently small
\(\varepsilon>0\), one has
\(T<\pi/\sqrt{\mathsf C_0+\varepsilon}\); put
\[
 a_\varepsilon(t):=\sqrt{\mathsf C_0+\varepsilon}
       \cot\!\bigl(\sqrt{\mathsf C_0+\varepsilon}\,t\bigr),
 \qquad
 b_\varepsilon(t):=\sqrt{\mathsf C_0+\varepsilon}
       \coth\!\bigl(\sqrt{\mathsf C_0+\varepsilon}\,t\bigr).
\]
These are the logarithmic growth rates of the constant positive- and
negative-curvature model solutions.  Direct differentiation gives
\[
 a_\varepsilon'+a_\varepsilon^2
    +(\mathsf C_0+\varepsilon)=0,
 \qquad
 b_\varepsilon'+b_\varepsilon^2
    -(\mathsf C_0+\varepsilon)=0.
\]
The perturbation by \(\varepsilon\) supplies the strict margin needed in the
first-crossing argument.
The scalar comparison functions satisfy
\[
 a_\varepsilon(t)
 =t^{-1}-\frac{\mathsf C_0+\varepsilon}{3}t+O(t^3),
 \qquad
 b_\varepsilon(t)
 =t^{-1}+\frac{\mathsf C_0+\varepsilon}{3}t+O(t^3).
\]
Subtracting these expansions from the expansion of \(S\) gives
\[
 \begin{aligned}
 S(t)-a_\varepsilon(t)I
 &=\frac{t}{3}
   \bigl[(\mathsf C_0+\varepsilon)I-\mathcal R(0)\bigr]+o(t),\\
 b_\varepsilon(t)I-S(t)
 &=\frac{t}{3}
   \bigl[(\mathsf C_0+\varepsilon)I+\mathcal R(0)\bigr]+o(t).
 \end{aligned}
\]
The curvature bound
\(-\mathsf C_0I\preceq\mathcal R(0)\preceq\mathsf C_0I\) implies that
both matrices in square brackets are bounded below by \(\varepsilon I\).
For fixed \(\varepsilon>0\), the operator norms of the \(o(t)\) remainders
are at most \(\varepsilon t/6\) when \(t>0\) is sufficiently small.  Therefore,
\[
 S(t)-a_\varepsilon(t)I
 \succeq\frac{\varepsilon t}{6}I\succ0,
 \qquad
 b_\varepsilon(t)I-S(t)
 \succeq\frac{\varepsilon t}{6}I\succ0.
\]
These
inequalities cannot fail at a first eigenvalue crossing: for a unit null
vector \(v\) there, the scalar quadratic form is nonnegative just before the
crossing and zero at the crossing, so its left derivative is nonpositive.
The Riccati equations instead give, respectively,
\[
 v^{\mathsf T}(S-a_\varepsilon I)'v
   =\mathsf C_0+\varepsilon-v^{\mathsf T}\mathcal Rv\ge\varepsilon,
 \qquad
 v^{\mathsf T}(b_\varepsilon I-S)'v
   =\mathsf C_0+\varepsilon+v^{\mathsf T}\mathcal Rv\ge\varepsilon,
\]
contradicting a first loss of positive semidefiniteness on \((0,T]\).  Letting
\(\varepsilon\to0^+\) gives
\[
 \sqrt{\mathsf C_0}\cot(\sqrt{\mathsf C_0}t)I
 \preceq S(t)\preceq
 \sqrt{\mathsf C_0}\coth(\sqrt{\mathsf C_0}t)I
 \qquad(0<t\le T).
\]
{\noindent\emph{Step 5: recovery of the Jacobi-field and Jacobian bounds.}}
For \(z\ne0\), integrate
\[
 \frac{\ud}{\ud t}\log\|\mathcal Y(t)z\|
 =\frac{\langle S(t)\mathcal Y(t)z,\mathcal Y(t)z\rangle}
        {\|\mathcal Y(t)z\|^2},
\]
and using \(\mathcal Y(t)/t\to I\), we obtain
\[
 \frac{\sin(\sqrt{\mathsf C_0}t)}{\sqrt{\mathsf C_0}}\|z\|
 \le\|\mathcal Y(t)z\|
 \le\frac{\sinh(\sqrt{\mathsf C_0}t)}{\sqrt{\mathsf C_0}}\|z\|.
\]
Because \(T\) was arbitrary, this holds for
\(0<t<\min\{\tau,\pi/\sqrt{\mathsf C_0}\}\).  If
\(\tau<\pi/\sqrt{\mathsf C_0}\), a nonzero
\(z\in\ker\mathcal Y(\tau)\), together with
the lower bound as \(t\to\tau^-\), gives the contradiction
\(0=\|\mathcal Y(\tau)z\|
\ge \sin(\sqrt{\mathsf C_0}\tau)\|z\|/\sqrt{\mathsf C_0}>0\).  Thus the
singular-value bound holds throughout
\(0<t<\pi/\sqrt{\mathsf C_0}\).  In particular,
\(\det\mathcal Y(t)>0\) there because
\(\mathcal Y(t)/t\to I\), and multiplication of the singular-value bounds
gives
\[
 \left(\frac{\sin(\sqrt{\mathsf C_0}t)}{\sqrt{\mathsf C_0}}\right)^{\dimk-1}
 \le\det\mathcal Y(t)
 \le
 \left(\frac{\sinh(\sqrt{\mathsf C_0}t)}{\sqrt{\mathsf C_0}}\right)^{\dimk-1}.
\]
Combining this determinant estimate with
\eqref{eq_geodesic_jacobian_from_jacobi_matrix} gives the required pointwise
comparison
\begin{equation}\label{eq_geodesic_ball_radial_jacobian_comparison}
 \left(
 \frac{\sin(\sqrt{\mathsf C_0}t)}{\sqrt{\mathsf C_0}t}
 \right)^{\dimk-1}
 \le J_{\theta,p}(t\omega)
 \le
 \left(
 \frac{\sinh(\sqrt{\mathsf C_0}t)}{\sqrt{\mathsf C_0}t}
 \right)^{\dimk-1}
\end{equation}
for \(0<t<\min\{r_0,\pi/\sqrt{\mathsf C_0}\}\).  Thus the upper curvature
bound gives the lower spherical-model density, whereas the lower curvature
bound gives the upper hyperbolic-model density.
No curvature derivative is used.  For the classical smooth analogue, see
\cite[Chapter~6, Theorem~27, p.~175]{petersen2006riemannian}; the argument
above covers the present \(\mathcal C^2\) metric.

\noindent\emph{Step 6: integration of the Jacobian bound.}
Set
\[
    c_1:=1,
    \qquad
    \mathsf c_L
    :=\left(\frac{\sin(1/\sqrt2)}{1/\sqrt2}\right)^{\dimk-1},
    \qquad
    \mathsf c_H
    :=\left(\frac{\sinh(1/\sqrt2)}{1/\sqrt2}\right)^{\dimk-1}.
\]
If
\[
    0<t<\min\left\{r_0,\frac1{\sqrt{2\mathsf C_0}}\right\},
\]
then \(0<\sqrt{\mathsf C_0}t<1/\sqrt2<\pi\).  Since
\(s\mapsto\sin(s)/s\) decreases on
\((0,\pi)\), whereas \(s\mapsto\sinh(s)/s\) increases on \((0,\infty)\),
the preceding Jacobian comparison implies
\[
    0<\mathsf c_L
    \le J_{\theta,p}(t\omega)
    \le\mathsf c_H<\infty,
    \qquad
    \mathsf c_L\le1\le\mathsf c_H.
\]
Since \(D(\exp_p)_0\) is the identity on \(T_pM\), one
also has \(J_{\theta,p}(0)=1\); hence the same bounds include the center of
the tangent-space ball.
Integrating this estimate in \eqref{eq_geod_vol_ball_normal} and using
\[
    \Leb_{T_pM}\bigl(\mathbb B_{T_pM}(0;r)\bigr)
    =\omega_{\dimk} r^{\dimk}
\]
gives, for every
\(0<r<\min\{r_0,c_1/\sqrt{2\mathsf C_0}\}\),
\[
    \mathsf c_L\,\omega_{\dimk}r^{\dimk}
    \le
    \mathcal M_\theta\bigl(\mathbb B_{\mathsf g_\theta}(p;r)\bigr)
    \le
    \mathsf c_H\,\omega_{\dimk}r^{\dimk}.
\]

Thus, also when \(\dimk\ge2\), the constants
\(c_1,\mathsf c_L,\mathsf c_H\) depend only on the common dimension \(\dimk\) and
are independent of \(\theta\), \(p\), and \(r\), as required.
\end{proof}

\subsection{A \(\mathcal C^1\) example showing polynomial-exponent sharpness}\label[appendix]{appendix_sharpness_sqr_rate}
\phantomsection
\label[appendix]{append_sharpness_sqr_rate}
We construct one pair \((f,g)\) satisfying the assumptions of
\Cref{th_CVaR_unif_approx} for which no polynomial rate with exponent larger
than \(1/\dimk\) is possible.  Retain the lower-level objective \(g\) from
the sphere family in \Cref{ex:sphere_family}.  At \(\theta=0\), its minimizer manifold is the unit
sphere \(\mathcal S(0)=\mathbb S^{d-1}\), so \(\dimk=d-1\).  Fix
\(p\in\mathbb S^{d-1}\) and define
\[
\psi(s):=
\begin{cases}
0, & s=0,\\[1mm]
\displaystyle\frac{s}{\log(1/s)}, & 0<s\le e^{-2},\\[3mm]
\displaystyle\frac{e^{-2}}2+\frac34(s-e^{-2}), & s>e^{-2}.
\end{cases}
\]
The last two formulas agree at \(e^{-2}\), as do their derivatives, whose common
value is \(3/4\).  On \(0<s\le e^{-2}\),
\[
\psi'(s)
=\frac{1}{\log(1/s)}+\frac{1}{\log^2(1/s)}
\in(0,3/4],
\]
and \(\psi'(s)\to0\) as \(s\to0^+\).  Thus \(\psi\) is \(\mathcal C^1\),
strictly increasing, and globally Lipschitz.  Define,
for every \(\theta\in\Theta\) and \(x\in\mathbb R^d\),
\[
    f(\theta,x):=-\psi(\|x-p\|)
\]
and note that, for \(x\ne p\),
\[
\nabla_xf(\theta,x)
=-\psi'(\|x-p\|)\frac{x-p}{\|x-p\|}.
\]
This gradient converges to zero as \(x\to p\), so setting
\(\nabla_xf(\theta,p)=0\) shows that \(f\) is \(\mathcal C^1\).  It is
independent of \(\theta\), globally Lipschitz, and has at most linear growth.
Consequently, together with the sphere family, it satisfies all assumptions
used in \Cref{th_CVaR_unif_approx}.  On \(\mathbb S^{d-1}\), it has the
unique maximum \(f(0,p)=F_{\max}(0)=0\).

The lower-level objective at \(\theta=0\) is rotationally invariant, and its
normal Hessian eigenvalue is the constant \(2/(1+a^2)\) on the unit sphere.
The density in \eqref{eqn_density_limiting_measure} is therefore constant, so
\(\mu_0^0\) is normalized surface measure.  The mass of a geodesic cap is
continuous and strictly increasing with its radius on \((0,\pi)\).  Hence,
for sufficiently small \(\delta>0\), there is a unique
\(r_\delta\in(0,\pi)\) satisfying
\[
\mu_0^0\bigl(\mathbb B_{\mathsf g_0}(p;r_\delta)\bigr)=\delta.
\]
If \(Y\sim\mu_0^0\) and
\(r=\operatorname{dist}_{\mathsf g_0}(Y,p)\), then
\(\|Y-p\|=2\sin(r/2)\).  Since \(\psi\) is strictly increasing, the cap
\(\mathbb B_{\mathsf g_0}(p;r_\delta)\) is the upper \(\delta\)-tail of
\(f(0,Y)\).  Its boundary has zero surface measure, so there is no atom at
the corresponding quantile.  More explicitly, with
\[
q_\delta:=-\psi\!\left(2\sin\frac{r_\delta}{2}\right),
\]
one has
\(\mathbb P(f(0,Y)>q_\delta)
=\mathbb P(f(0,Y)\ge q_\delta)=\delta\).  Thus \(q_\delta\) minimizes the
variational superquantile objective, and substitution in
\eqref{eq_def_sqr_loss} gives
\[
\widetilde F_{\SQR}(0)
=\frac1\delta
  \mathbb E\!\left[
      f(0,Y)\,
      \mathbf 1_{\mathbb B_{\mathsf g_0}(p;r_\delta)}(Y)
  \right].
\]

Let \(c_{\dimk}>0\) denote the normalizing constant in geodesic polar
coordinates on \(\mathbb S^{\dimk}\), namely the ratio of the surface areas
of \(\mathbb S^{\dimk-1}\) and \(\mathbb S^{\dimk}\).  Then
\[
\delta
=c_{\dimk}\int_0^{r_\delta}(\sin r)^{\dimk-1}\,dr.
\]
For small \(r_\delta\), chosen in particular so that
\(2\sin(r_\delta/2)\le e^{-2}\), comparison of
\(2\sin(r/2)\) and \(\sin r\) with
\(r\) gives constants \(0<a_0\le A_0<\infty\), independent of \(\delta\),
such that
\[
a_0r_\delta^{\dimk}
\le\delta\le
A_0r_\delta^{\dimk}.
\]
The same comparisons and the definition of \(\psi\) give
\[
F_{\max}(0)-\widetilde F_{\SQR}(0)
=\frac{c_{\dimk}}{\delta}
  \int_0^{r_\delta}
  \psi\!\left(2\sin\frac r2\right)
  (\sin r)^{\dimk-1}\,dr.
\]
For the upper bound on this integral, use
\(2\sin(r/2)\le r\) on \([0,r_\delta]\).  For the lower bound, restrict the
integral to \([r_\delta/2,r_\delta]\), where both sine factors are bounded
below by fixed multiples of \(r_\delta\).  It follows that there are
constants \(0<b_1\le B_1<\infty\), independent of sufficiently small
\(\delta\), for which
\[
b_1\frac{r_\delta^{\dimk+1}}{\log(1/r_\delta)}
\le
\int_0^{r_\delta}
  \psi\!\left(2\sin\frac r2\right)
  (\sin r)^{\dimk-1}\,dr
\le
B_1\frac{r_\delta^{\dimk+1}}{\log(1/r_\delta)}.
\]
After division by \(\delta\), the cap-volume bounds therefore imply, after
introducing constants \(0<a_1\le A_1<\infty\),
\[
a_1\frac{r_\delta}{\log(1/r_\delta)}
\le
F_{\max}(0)-\widetilde F_{\SQR}(0)
\le
A_1\frac{r_\delta}{\log(1/r_\delta)}.
\]
Combining this estimate with the two-sided comparison between \(\delta\) and
\(r_\delta^{\dimk}\), and using
\(\log(1/r_\delta)=\dimk^{-1}\log(1/\delta)+O(1)\), gives constants
\(0<a_2\le A_2<\infty\) such that
\[
a_2\frac{\delta^{1/\dimk}}{\log(1/\delta)}
\le
F_{\max}(0)-\widetilde F_{\SQR}(0)
\le
A_2\frac{\delta^{1/\dimk}}{\log(1/\delta)}
\]
for all sufficiently small \(\delta>0\).
Consequently, for every \(\alpha>1/\dimk\),
\[
    \frac{F_{\max}(0)-\widetilde F_{\SQR}(0)}
         {\delta^\alpha}
    \longrightarrow+\infty
    \qquad(\delta\to0^+).
\]
Thus no estimate of order \(O(\delta^\alpha)\) with
\(\alpha>1/\dimk\) can hold over the present class, which proves the
polynomial-exponent sharpness claimed in the main body.
The logarithmic factor is why this conclusion concerns the polynomial
exponent rather than a matching constant-order lower bound.  Since \(d\ge2\)
in the sphere family is arbitrary, the construction covers every intrinsic
dimension \(k\ge1\) by choosing \(d=k+1\).

\subsection{Proof of the Gibbs perturbation lemma}\label[appendix]{append_proof_Theo_func_approx_cvar}
We next provide the proof of \Cref{Theo_func_approx_cvar} which is essential for our approximation result in \Cref{rmk_func_approx_cvar_gibbs}.
\begin{proof}
Fix \(\theta\in\Theta\), \(\lambda>0\), and \(\delta\in(0,1)\).  By
\Cref{lemma_density_limiting_measure}, \(\mu_\theta^\lambda\) is a probability
measure and \(\mu_\theta^0\) is supported on the compact set
\(\mathcal S(\theta)\).  In particular, \(\mu_\theta^0\) has a finite first
moment.  To verify the same property for \(\mu_\theta^\lambda\), set
\[
    h_\theta(x):=g(\theta,x)-g^\star(\theta),
    \qquad
    B_\theta:=\max_{u\in\mathcal S(\theta)}\|u\|.
\]
For
\(\|x\|\ge\max\{\mathsf D_{\mathsf{qg}},2B_\theta\}\),
\eqref{eq_tail_EB_implies_QG} gives
\[
    h_\theta(x)
    \ge \frac{\mu_{\mathsf{qg}}}{2}
       \operatorname{dist}^2\bigl(x,\mathcal S(\theta)\bigr)
    \ge \frac{\mu_{\mathsf{qg}}}{8}\|x\|^2.
\]
Consequently,
\(\int_{\mathbb R^d}\|x\|e^{-h_\theta(x)/\lambda}\,\mathrm dx<\infty\),
so \(\mu_\theta^\lambda\) also has a finite first moment.  Since
\(f(\theta,\cdot)\) is globally Lipschitz, it is integrable under both
measures.

For either \(\nu\in\{\mu_\theta^\lambda,\mu_\theta^0\}\), the function
\[
    \Phi_\nu(\beta)
    :=\beta+\frac1\delta
       \int_{\mathbb R^d}(f(\theta,x)-\beta)_+\,\nu(\mathrm dx)
\]
is finite and continuous; indeed, it is
\((1+\delta^{-1})\)-Lipschitz in \(\beta\).  Moreover,
\(\Phi_\nu(\beta)\ge\beta\) as
\(\beta\to+\infty\), while
\[
    \Phi_\nu(\beta)
    \ge \frac1\delta\int f(\theta,x)\,\nu(\mathrm dx)
       +\left(1-\frac1\delta\right)\beta
    \longrightarrow+\infty
    \qquad(\beta\to-\infty).
\]
Thus both variational objectives attain their minima.

For fixed \(\beta\), let
\(G_\beta(x):=(f(\theta,x)-\beta)_+\).  The positive-part map is
\(1\)-Lipschitz, so \(G_\beta\) is \(L_{f,1}\)-Lipschitz, uniformly in
\(\beta\).  Hence, for every coupling
\(\pi\in\Pi(\mu_\theta^\lambda,\mu_\theta^0)\),
\[
\begin{aligned}
&\bigl|\phi_{\lambda,\delta}(\theta,\beta)
       -\phi_{0,\delta}(\theta,\beta)\bigr|\\
&\quad\le\frac1\delta
   \int_{\mathbb R^d\times\mathbb R^d}
   |G_\beta(x)-G_\beta(y)|\,\pi(\mathrm dx,\mathrm dy)\\
&\quad\le\frac{L_{f,1}}\delta
   \int_{\mathbb R^d\times\mathbb R^d}
   \|x-y\|\,\pi(\mathrm dx,\mathrm dy).
\end{aligned}
\]
Taking the infimum over all such couplings yields
\begin{equation}\label{eq_00121}
\sup_{\beta\in\mathbb R}
\bigl|\phi_{\lambda,\delta}(\theta,\beta)
      -\phi_{0,\delta}(\theta,\beta)\bigr|
\le
\frac{L_{f,1}}\delta
\mathbb W_1\!\left(\mu_\theta^\lambda,\mu_\theta^0\right).
\end{equation}
Writing the right-hand side of \eqref{eq_00121} as
\(C_{\theta,\lambda,\delta}\), the two pointwise inequalities implied by
\eqref{eq_00121}, followed by minimization over \(\beta\), give
\begin{equation}\label{eq_0026}
\begin{aligned}
\min_{\beta}\phi_{0,\delta}(\theta,\beta)
&\le
\min_{\beta}\phi_{\lambda,\delta}(\theta,\beta)
+C_{\theta,\lambda,\delta},\\
\min_{\beta}\phi_{\lambda,\delta}(\theta,\beta)
&\le
\min_{\beta}\phi_{0,\delta}(\theta,\beta)
+C_{\theta,\lambda,\delta}.
\end{aligned}
\end{equation}
Therefore,
\begin{equation}\label{eq_00035}
\left|
\min_{\beta}\phi_{\lambda,\delta}(\theta,\beta)
-\min_{\beta}\phi_{0,\delta}(\theta,\beta)
\right|
\le
\frac{L_{f,1}}\delta
\mathbb W_1\!\left(\mu_\theta^\lambda,\mu_\theta^0\right).
\end{equation}
The two minima in \eqref{eq_00035} are precisely
\(\tilde F_{\SQG}(\theta)\) and \(\tilde F_{\SQR}(\theta)\), respectively.
\end{proof}

\section{Proofs of Results in \Cref{sec:6}}
This appendix collects technical results for the oracle analysis in \Cref{sec:6}. We begin by deriving moment and tail bounds for the Gibbs measure; these bounds are then used to prove \Cref{lem_3}. Next, we establish accuracy guarantees for the PSGD-based inner estimator. Finally, we present the complete proof of \Cref{thm_cvg_alg_zeroth_goldstein}, separating the standard zeroth-order argument from our new control of the minima-selection error.
\subsection{Proof of \Cref{lem_3}}\label[appendix]{append_proof_lemm_lb_beta_min_rsq}

We start by presenting the following technical lemma, which is an essential prerequisite for the proof of \Cref{lem_3}. It establishes a technical result ensuring the finiteness of the moments of the Gibbs measure $\gibbs_{\theta}^{\lambda}$.
\begin{lem}\label{lem_gibbs_mass_outside}
Suppose that \Cref{mu_PL_assum,assum_sample_comp} hold, and let \(\mathsf D\)
be fixed as in \eqref{eq_working_radius_D}.  For every
\(n\in\mathbb N\cup\{0\}\), there exist constants \(\mathsf W_n<\infty\) and
\(\lambda_{\mathrm{tail},n}>0\), independent of \(\theta\), such that
\begin{equation}\label{eq_uniform_gibbs_polynomial_tail}
\mathbb E_{X\sim\gibbs_\theta^\lambda}
\left[\|X\|^n\mathbf 1_{\{X\notin\mathbb B_d(0;2\mathsf D)\}}\right]
\le\mathsf W_n\lambda^{1-d/2}
\exp\!\left(-\frac{c_{\mathrm{tail}}}{\lambda}\right)
\end{equation}
for every \(\theta\in\Theta\) and
\(0<\lambda\le\lambda_{\mathrm{tail},n}\), where
\[
    c_{\mathrm{tail}}:=\frac{\mu_{\mathsf{qg}}\mathsf D^2}{2}.
\]
\end{lem}

\begin{proof}
Fix \(\theta\), and set
\[
h_\theta:=g(\theta,\cdot)-g^\star(\theta),
\qquad
\widehat Z_\lambda(\theta):=\int e^{-h_\theta(x)/\lambda}\,\mathrm dx.
\]
The shift does not change the Gibbs measure.  Choose
\(u_\theta\in\mathcal S(\theta)\).  The descent lemma gives
\(h_\theta(x)\le L_{g,xx}\|x-u_\theta\|^2/2\), so, for
\(0<\lambda\le1\),
\begin{equation}\label{eq_tail_crude_shifted_normalizer}
\widehat Z_\lambda(\theta)
\ge e^{-L_{g,xx}/2}\operatorname{vol}(\mathbb B_d(0;1))\lambda^{d/2}
=:c_Z\lambda^{d/2}.
\end{equation}
If \(\|x\|\ge2\mathsf D\), then
\(\operatorname{dist}(x,\mathcal S(\theta))\ge\|x\|/2\).
The estimate \eqref{eq_tail_EB_implies_QG} therefore gives
\begin{equation}\label{eq_tail_radial_lower}
    h_\theta(x)\ge\frac{\mu_{\mathsf{qg}}}{8}\|x\|^2.
\end{equation}

Let \(s:=(n+d)/2\ge1\), and note that the surface area of the Euclidean
unit sphere is \(\lvert\mathbb S^{d-1}\rvert=2\pi^{d/2}/\Gamma(d/2)\).
Denote the upper incomplete
gamma function by
\(\Gamma(s,y):=\int_y^\infty t^{s-1}e^{-t}\,\mathrm dt\).
For \(\|x\|\ge2\mathsf D\), \eqref{eq_tail_radial_lower} first gives the
pointwise bound
\[
 \|x\|^n e^{-h_\theta(x)/\lambda}
 \le
 \|x\|^n
 \exp\!\left(-\frac{\mu_{\mathsf{qg}}\|x\|^2}{8\lambda}\right).
\]
To integrate this radial upper bound, write \(x=r\omega\), where
\(r=\|x\|\ge2\mathsf D\) and \(\omega\in\mathbb S^{d-1}\).  If
\(\mathrm d\sigma(\omega)\) denotes surface-area measure on
\(\mathbb S^{d-1}\), then the polar-coordinate volume element is
\(\mathrm dx=r^{d-1}\mathrm dr\,\mathrm d\sigma(\omega)\).  Consequently,
\begin{align*}
\int_{\|x\|\ge2\mathsf D}\|x\|^ne^{-h_\theta(x)/\lambda}\,\mathrm dx
&\le
\int_{2\mathsf D}^{\infty}\int_{\mathbb S^{d-1}}
r^n
\exp\!\left(-\frac{\mu_{\mathsf{qg}}r^2}{8\lambda}\right)
r^{d-1}\,\mathrm d\sigma(\omega)\,\mathrm dr\\
&=
\lvert\mathbb S^{d-1}\rvert
\int_{2\mathsf D}^{\infty}
r^{n+d-1}
\exp\!\left(-\frac{\mu_{\mathsf{qg}}r^2}{8\lambda}\right)
\,\mathrm dr.
\end{align*}
Set \(y:=\mu_{\mathsf{qg}}r^2/(8\lambda)\).  Since
\[
 r^{n+d-1}\,\mathrm dr
 =\frac12
 \left(\frac{8\lambda}{\mu_{\mathsf{qg}}}\right)^s
 y^{s-1}\,\mathrm dy,
 \qquad
 r=2\mathsf D
 \ \Longleftrightarrow\
 y=\frac{c_{\mathrm{tail}}}{\lambda},
\]
where the last equality uses
\(c_{\mathrm{tail}}=\mu_{\mathsf{qg}}\mathsf D^2/2\).
Substitution in the preceding radial integral gives
\begin{align}
\int_{\|x\|\ge2\mathsf D}\|x\|^ne^{-h_\theta(x)/\lambda}\,\mathrm dx
&\le\frac{\lvert\mathbb S^{d-1}\rvert}{2}
\left(\frac{8\lambda}{\mu_{\mathsf{qg}}}\right)^s
\Gamma\!\left(s,\frac{c_{\mathrm{tail}}}{\lambda}\right).
\label{eq_tail_incomplete_gamma}
\end{align}
For \(s\ge1\) and \(z>0\) satisfying \(z\ge2(s-1)\),
\begin{equation}\label{eq_incomplete_gamma_bound}
    \Gamma(s,z)\le2z^{s-1}e^{-z}.
\end{equation}
Indeed, writing \(t=z+u\) and using
\((1+u/z)^{s-1}\le e^{(s-1)u/z}\) bounds the remaining integral by
\(\int_0^\infty e^{-u/2}\,\mathrm du=2\).
Set
\[
\lambda_{\mathrm{tail},n}
:=\min\left\{1,\frac{c_{\mathrm{tail}}}{\max\{n+d-2,1\}}\right\}.
\]
Then \eqref{eq_incomplete_gamma_bound} applies for
\(0<\lambda\le\lambda_{\mathrm{tail},n}\), and
\[
\int_{\|x\|\ge2\mathsf D}\|x\|^ne^{-h_\theta(x)/\lambda}\,\mathrm dx
\le\lvert\mathbb S^{d-1}\rvert
\frac{8\cdot4^{s-1}}{\mu_{\mathsf{qg}}}
\mathsf D^{n+d-2}\lambda
e^{-c_{\mathrm{tail}}/\lambda}.
\]
Division by \eqref{eq_tail_crude_shifted_normalizer} proves
\eqref{eq_uniform_gibbs_polynomial_tail}, for example with
\[
\mathsf W_n
:=\frac{\lvert\mathbb S^{d-1}\rvert
8\cdot4^{(n+d)/2-1}\mathsf D^{n+d-2}}
{\mu_{\mathsf{qg}}c_Z}.
\]
Every displayed constant is independent of \(\theta\).
\end{proof}

We then prove \Cref{lem_3}, which states that every global minimizer of
\(\phi_{\lambda,\delta}(\theta,\cdot)\) lies in a common bounded interval.
\begin{proof}[Proof of \Cref{lem_3}]\label{proof_lem_3}
Fix \(\theta\in\Theta\), \(0<\delta\le1/2\), and
\(0<\lambda\le\lambda_\beta\).  Let
\[
 X\sim\mu_\theta^\lambda,\qquad Z:=f(\theta,X),
\]
and take
\(\beta^\star\in\argmin_{\beta\in\mathbb R}
\phi_{\lambda,\delta}(\theta,\beta)\).
Finiteness and attainment follow from the convex variational analysis in
\Cref{Theo_func_approx_cvar}.

Because the hinge term is nonnegative,
\begin{equation}\label{eq_029001}
 \beta^\star
 \le \phi_{\lambda,\delta}(\theta,\beta^\star)
 =\min_{\beta\in\mathbb R}\phi_{\lambda,\delta}(\theta,\beta).
\end{equation}
The perturbation estimate in \Cref{Theo_func_approx_cvar} and
\Cref{prop_W1_Gibbs_uniform} give
\begin{align}
\min_{\beta}\phi_{\lambda,\delta}(\theta,\beta)
&\le \min_{\beta}\phi_{0,\delta}(\theta,\beta)
 +\frac{L_{f,1}}{\delta}
   \mathbb W_1(\mu_\theta^\lambda,\mu_\theta^0)\nonumber\\
&\le \min_{\beta}\phi_{0,\delta}(\theta,\beta)
 +\frac{L_{f,1}K_0}{\delta}\sqrt\lambda .
\label{eq_00048}
\end{align}
Choosing \(\beta=F_{\max}(\theta)\) in the limiting objective makes its
hinge term zero, since \(f(\theta,Y)\le F_{\max}(\theta)\) for
\(Y\sim\mu_\theta^0\).  Thus
\[
 \min_\beta\phi_{0,\delta}(\theta,\beta)
 \le F_{\max}(\theta)\le L_{f,0}(\mathsf D).
\]
Combining this with \eqref{eq_029001}--\eqref{eq_00048} yields
\begin{equation}\label{eq_0101050}
 \beta^\star
 \le L_{f,0}(\mathsf D)
 +\frac{L_{f,1}K_0}{\delta}\sqrt\lambda .
\end{equation}

For the lower bound, the one-sided derivatives of the convex hinge
expectation give, without any atomlessness assumption,
\[
 \partial_\beta\phi_{\lambda,\delta}(\theta,\beta)
 =
 \left[
 1-\frac{\mathbb P(Z\ge\beta)}{\delta},
 1-\frac{\mathbb P(Z>\beta)}{\delta}
 \right].
\]
The condition \(0\in
\partial_\beta\phi_{\lambda,\delta}(\theta,\beta^\star)\) therefore implies
\[
 \mathbb P(Z>\beta^\star)\le\delta
 \le\mathbb P(Z\ge\beta^\star),
 \qquad
 \mathbb P(Z\le\beta^\star)\ge1-\delta\ge\frac12.
\]
By the definition of \(\lambda_{\mathrm{mass}}\) and \(\lambda_\beta\),
\Cref{lem_gibbs_mass_outside} gives
\[
 \mu_\theta^\lambda\!\left(\mathbb B_d(0;2\mathsf D)\right)\ge\frac34.
\]
Consequently, inclusion--exclusion yields
\[
\begin{aligned}
&\mu_\theta^\lambda\!\left(
 \{x:f(\theta,x)\le\beta^\star\}
 \cap\mathbb B_d(0;2\mathsf D)\right)\\
&\hspace{4cm}\ge\frac12+\frac34-1=\frac14.
\end{aligned}
\]
The intersection is nonempty; choose \(x_0\) in it.  Then
\[
 \beta^\star\ge f(\theta,x_0)\ge-L_{f,0}(2\mathsf D).
\]
Together with \eqref{eq_0101050}, this proves
\eqref{eq_bound_beta_minimizer}.  All thresholds and constants used here are
common in \(\theta\), whereas \(\mathrm B(\delta,\lambda)\) retains its
displayed dependence on \(\delta\) and \(\lambda\).
\end{proof}

\subsection{Proof of \Cref{lem_5_2}}\label[appendix]{append_proofs_lem_5_2}
This subsection proves the mean-square accuracy bound for the PSGD estimator
of $\widetilde F_{\SQG}(\theta)$, including the conditional stochastic-
subgradient and validation estimates used below.
\begin{proof}
Fix the displayed parameters in \Cref{lem_5_2}, let
$X\sim\mu_\theta^\lambda$, and write $Z:=f(\theta,X)$.
Dominated convergence applied to the hinge gives the
left and right derivatives
\[
D_-\phi_{\lambda,\delta}(\theta,\beta)
=1-\frac{\mathbb P(Z\ge\beta)}{\delta},
\qquad
D_+\phi_{\lambda,\delta}(\theta,\beta)
=1-\frac{\mathbb P(Z>\beta)}{\delta}.
\]
For a finite convex function of one real variable, the subdifferential is the
closed interval between its left and right derivatives.  Hence
\[
\partial_\beta\phi_{\lambda,\delta}(\theta,\beta)
=
\left[
1-\frac{\mathbb P(Z\ge\beta)}{\delta},
1-\frac{\mathbb P(Z>\beta)}{\delta}
\right].
\]
The endpoints are correctly ordered because
\(\mathbb P(Z\ge\beta)\ge\mathbb P(Z>\beta)\).  Its length is
\(\mathbb P(Z=\beta)/\delta\), so it collapses to a singleton exactly when
there is no atom at \(\beta\).
Choose the measurable sample-hinge subgradient
\[
s(\beta,Z):=-\mathbf 1_{\{Z>\beta\}}
\in\partial_\beta(Z-\beta)_+,
\qquad
G(\beta,X):=1+\delta^{-1}s(\beta,Z).
\]
Thus the choice at the atom $\{Z=\beta\}$ is zero, and
\[
\mathbb E[G(\beta,X)]
=1-\frac{\mathbb P(Z>\beta)}{\delta}
\in\partial_\beta\phi_{\lambda,\delta}(\theta,\beta).
\]
Moreover,
\[
|G(\beta,X)|\le1+\delta^{-1}=:\sigma,
\qquad
\mathbb E\exp\!\left(\frac{|G(\beta,X)|^2}{\sigma^2}\right)\le e.
\]
No atomlessness is used.  Since the hinge is $1$-Lipschitz, for all
$\beta_1,\beta_2\in\mathbb R$,
\begin{align}
\left|\phi_{\lambda,\delta}(\theta,\beta_1)
-\phi_{\lambda,\delta}(\theta,\beta_2)\right|
&\le \left(1+\delta^{-1}\right)|\beta_1-\beta_2|.
\label{eq_00110}
\end{align}

We next establish the high-probability PSGD estimate used in the statement.
Assume first that $\mathrm B>0$, set
$\eta:=2\mathrm B/(\sigma\sqrt L)$, and let
\[
\mathcal F_l
:=\sigma\!\left(\widetilde X^{(0)},\ldots,
\widetilde X^{(l-1)}\right),
\quad
G_l:=G(\beta^l,\widetilde X^{(l)}),
\quad
g_l:=\mathbb E[G_l\mid\mathcal F_l].
\]
Conditional independence gives
$g_l\in\partial_\beta\phi_{\lambda,\delta}(\theta,\beta^l)$.
By \Cref{lem_3}, an exact minimizer $\beta^\star$ lies in
$[-\mathrm B,\mathrm B]$.
Nonexpansiveness of the projection gives, for each \(l\),
\[
|\beta^{l+1}-\beta^\star|^2
\le |\beta^l-\beta^\star|^2
-2\eta G_l(\beta^l-\beta^\star)+\eta^2\sigma^2.
\]
Moreover, convexity and
\(g_l\in\partial_\beta\phi_{\lambda,\delta}(\theta,\beta^l)\) imply
\[
\phi_{\lambda,\delta}(\theta,\beta^l)
-\phi_{\lambda,\delta}(\theta,\beta^\star)
\le g_l(\beta^l-\beta^\star)
=G_l(\beta^l-\beta^\star)
 +(g_l-G_l)(\beta^l-\beta^\star).
\]
Rearranging the first display, substituting it into the second, and summing
from \(l=0\) to \(L-1\) yield
\begin{align*}
\frac1L\sum_{l=0}^{L-1}
\bigl[\phi_{\lambda,\delta}(\theta,\beta^l)
-\phi_{\lambda,\delta}(\theta,\beta^\star)\bigr]
&\le
\frac{(2\mathrm B)^2}{2\eta L}
+\frac{\eta\sigma^2}{2}\\
&\quad+
\frac1L\sum_{l=0}^{L-1}
(g_l-G_l)(\beta^l-\beta^\star).
\end{align*}
The last summands are martingale differences and have absolute value at most
$4\mathrm B\sigma$.  Azuma--Hoeffding and convexity at
$\widehat\beta=L^{-1}\sum_{l=0}^{L-1}\beta^l$ therefore imply, with
conditional probability at least $1-\Delta$,
\[
\phi_{\lambda,\delta}(\theta,\widehat\beta)
-\min_\beta\phi_{\lambda,\delta}(\theta,\beta)
\le
\frac{\mathrm B\sigma}{\sqrt L}
\left(2+4\sqrt{2\log(1/\Delta)}\right).
\]
For every \(\tau>0\) and \(0<\Delta\le1/2\), the sufficient
condition
\[
L\ge64\,\mathrm B^2\sigma^2\tau^{-2}\log^2(2/\Delta)
\]
thus gives
\[
\mathbb P\!\left(
\phi_{\lambda,\delta}(\theta,\widehat\beta)
-\min_\beta\phi_{\lambda,\delta}(\theta,\beta)>\tau
\,\middle|\,\theta\right)
\le\Delta.
\]
This is precisely \eqref{eq_whp_gaurante_PSGD} with
\(\widehat\epsilon=\tau\).
For comparison with
\citep[Proposition~2]{lan2012validation}, take
\(\omega(\beta)=\tfrac12(\beta-\beta^0)^2\) as the distance-generating
function on \([ -\mathrm B,\mathrm B]\).  Its prox map is Euclidean
projection and its minimizer is the prescribed initial point \(\beta^0\).
The feasible set is compact and convex, the objective is convex,
\(M_\ast=\sigma\) satisfies the source's exponential-moment condition, and
the source scaling constant
\(2\mathrm B/(\mathrm B+|\beta^0|)\in[1,2]\) makes its constant stepsize
equal to \(2\mathrm B/(\sigma\sqrt L)\).  Its weighted output is the
arithmetic average because the stepsize is constant.  Thus the cited
hypotheses are matched; the preceding projection--martingale derivation is
nevertheless self-contained.  If $\mathrm B=0$, the feasible interval is
$\{0\}$ and the
optimization gap is identically zero.

We now bound the estimator error.  By $|a+b|^2\le2|a|^2+2|b|^2$,
\begin{align}
\MoveEqLeft[2]
\mathbb E\!\left[
|\widetilde\psi(\theta)-\widetilde F_{\SQG}(\theta)|^2
\,\middle|\,\theta\right]\nonumber\\
&\le
2\underbrace{\mathbb E\!\left[
|\phi_{\lambda,\delta}(\theta,\widehat\beta;\widehat X_I)
-\phi_{\lambda,\delta}(\theta,\widehat\beta)|^2
\,\middle|\,\theta\right]}_{(i)}\nonumber\\
&\quad+
2\underbrace{\mathbb E\!\left[
\left(\phi_{\lambda,\delta}(\theta,\widehat\beta)
-\min_\beta\phi_{\lambda,\delta}(\theta,\beta)\right)^2
\,\middle|\,\theta\right]}_{(ii)}.
\label{eq_append_22}
\end{align}
Here $\widetilde F_{\SQG}(\theta)=
\min_\beta\phi_{\lambda,\delta}(\theta,\beta)$.

For (i), freshness of the validation batch implies that, conditionally on
$(\theta,\widehat\beta)$, its observations remain i.i.d. with law
$\mu_\theta^\lambda$.  Consequently,
\begin{align}
(i)
&=
\mathbb E\!\left[
\frac{1}{\delta^2|I|}
\operatorname{Var}\!\left(
(f(\theta,X)-\widehat\beta)_+
\,\middle|\,\theta,\widehat\beta\right)
\,\middle|\,\theta\right].
\label{eq_append_24}
\end{align}
Let $X'$ be an independent copy of $X$, also independent of
$\widehat\beta$.  For $Y=(f(\theta,X)-\widehat\beta)_+$ and the corresponding
$Y'$, the hinge contraction gives
\[
\begin{aligned}
\operatorname{Var}(Y\mid\theta,\widehat\beta)
&=\frac12\mathbb E[(Y-Y')^2\mid\theta,\widehat\beta]\\
&\le\frac12\mathbb E\!\left[
 (f(\theta,X)-f(\theta,X'))^2
 \,\middle|\,\theta,\widehat\beta\right]\\
&\le2\mathbb E[f(\theta,X)^2\mid\theta]
\le4\mathbb E[f(\theta,X)^2\mid\theta].
\end{aligned}
\]
Define
\[
\overline{\mathsf W}_{2\mathsf n_f}
:=\mathsf W_{2\mathsf n_f}
\sup_{0<t\le\lambda_{\mathrm{tail},2\mathsf n_f}}
t^{1-d/2}e^{-c_{\mathrm{tail}}/t}<\infty.
\]
The supremum is finite because the displayed function extends
continuously to \(t=0\) with value zero.
\begin{align}
4\mathbb E[f(\theta,X)^2]
&\le4\mathbb E[(C_f\|X\|^{\mathsf n_f}+D_f)^2]\nonumber\\
&\le8C_f^2\mathbb E[\|X\|^{2\mathsf n_f}]+8D_f^2\nonumber\\
&\le \widetilde C_f:=8C_f^2
    \left((2\mathsf{D})^{2\mathsf n_{f}}+\overline{\mathsf W}_{2\mathsf n_f}\right)
    +8D_{f}^{2},\label{eq_append_54}
\end{align}
where the first inequality follows from \Cref{assum_sample_comp}, and the last
uses
\Cref{lem_gibbs_mass_outside} under
\(0<\lambda\le\lambda_{\mathrm{tail},2\mathsf n_f}\).  Substitution into
\eqref{eq_append_24} gives
\begin{align}
(i)\le\frac{\widetilde C_f}{\delta^2|I|}.
\label{eq_append_25}
\end{align}

For (ii), \eqref{eq_00110}, the projection of every iterate, and
\Cref{lem_3} give the deterministic range bound
\[
0\le
Q:=\phi_{\lambda,\delta}(\theta,\widehat\beta)
-\min_\beta\phi_{\lambda,\delta}(\theta,\beta)
\le2\sigma\mathrm B.
\]
Set
\(\Delta:=\min\{1/2,
\epsilon_v^2/(32\sigma^2\mathrm B^2+\epsilon_v^2)\}\).  Then
\(0<\Delta\le1/2\) and
\(4\sigma^2\mathrm B^2\Delta\le\epsilon_v^2/8\).  Applying
\eqref{eq_whp_gaurante_PSGD} with target
\(\epsilon_v/(2\sqrt2)\) gives
\begin{align}
(ii)=\mathbb E[Q^2\mid\theta]
&\le\frac{\epsilon_v^2}{8}
    +4\sigma^2\mathrm B^2\Delta
\nonumber\\
&\le\frac{\epsilon_v^2}{4}.
\label{eq_append_23}
\end{align}
For $\mathrm B=0$, this conclusion holds directly because $Q=0$.
Finally, if
\[
|I|\ge4\widetilde C_f\delta^{-2}\epsilon_v^{-2},
\]
then \eqref{eq_append_25} gives $(i)\le\epsilon_v^2/4$.  Substitution of
\eqref{eq_append_23} and \eqref{eq_append_25} into \eqref{eq_append_22}
yields exactly
\[
\mathbb E\!\left[
|\widetilde\psi(\theta)-\widetilde F_{\SQG}(\theta)|^2
\,\middle|\,\theta\right]
\le2(i)+2(ii)\le\epsilon_v^2.
\]
Only now substitute $\delta=\epsilon_v^{\dimk}$.  The validation condition
becomes
\[
|I|\ge4\widetilde C_f\epsilon_v^{-2-2\dimk}.
\]
Moreover, $\sigma\le2\epsilon_v^{-\dimk}$, while the definition in
\Cref{lem_3} and $\lambda\le\lambda_{\mathrm{mom}}$ give
$\mathrm B\le C_{\mathrm B}\epsilon_v^{-\dimk}$ for a common finite
$C_{\mathrm B}$.  Hence
$\log(2/\Delta)=\mathcal O(\dimk[1+\log(\epsilon_v^{-1})])$.
Taking the maximum of \(1\) and the ceiling of the exact
lower bound in the statement gives a successful integer \(L\) satisfying
\[
L=\mathcal O\!\left(1+
\dimk^2\mathrm B^2\epsilon_v^{-2-2\dimk}
[1+\log^2(\epsilon_v^{-1})]\right).
\]
When the second term is at least one, this is the stated
asymptotic iteration order; when \(\mathrm B=0\), no PSGD update or training
sample is needed.
\end{proof}
\subsection{Proof of \Cref{thm_cvg_alg_zeroth_goldstein}}\label[appendix]{append_proof_thm_cvg_alg_zeroth_goldstein}
We now prove \Cref{thm_cvg_alg_zeroth_goldstein}. The argument follows the standard projected zeroth-order analysis for the smoothed objective, and the only additional step is to bound the error introduced by replacing $F_{\max}$ with $\tilde{\psi}$ using \Cref{eq_append_21} and \Cref{lem_5_2}.
\begin{proof}\label{proof_thm_cvg_alg_zeroth_goldstein}
Set
\[
    L_F:=L_{F_{\max},1}
    =L_{f,1}\left(1+\frac{2L_{g,\theta x}}{\mu}\right).
\]
By \Cref{lemm_hyp_blo_cont}, applied on the compact convex thickening,
$F_{\max}$ is $L_F$-Lipschitz there.  If $L_F=0$, then $F_{\max}$ is
constant on the connected thickening, $0\in\partial_\rho F_{\max}(\theta)$,
and the claimed residual is zero.  Hence assume $L_F>0$.

Because \(0<\rho<\bar\rho\), for every \(\theta\in\Theta\) the smoothing
ball \(\mathbb B_m(\theta;\rho)\) lies at distance at least
\(\bar\rho-\rho>0\) from the complement of
\(\Theta^{+\bar\rho}\).  Thus \(F_{\max}\) is locally Lipschitz on an
ambient open neighborhood of every such ball.  More precisely, set
\[
\Omega_\rho
:=\{z\in\mathbb R^m:
\operatorname{dist}(z,\Theta)<\bar\rho-\rho\}.
\]
For every \(z\in\Omega_\rho\) and \(v\in\mathbb B_m(0;1)\), one has
\(z+\rho v\in\operatorname{int}(\Theta^{+\bar\rho})\).  Hence the
uniform-ball convolution
\[
    \psi_\rho(z)
    :=\mathbb E_{V\sim\operatorname{Unif}(\mathbb B_m(0;1))}
      \bigl[F_{\max}(z+\rho V)\bigr],
      \qquad z\in\Omega_\rho,
\]
is well defined on an open neighborhood of \(\Theta\).
The standard boundary formula for ball convolution gives, for
$U\sim\operatorname{Unif}(\mathbb S^{m-1})$,
\begin{align*}
\nabla\psi_\rho(\theta)
&=\frac m\rho\,
  \mathbb E\bigl[F_{\max}(\theta+\rho U)U\bigr]\\
&=\frac{m}{2\rho}\,
  \mathbb E\!\left[
  \bigl(F_{\max}(\theta+\rho U)-F_{\max}(\theta-\rho U)\bigr)U
  \right].
\end{align*}
This is the ball-smoothing/sphere-estimator identity in
\citep[Supplement, Lemma~D.1]{lin2022gradient}; see also
\citep[Lemma~A.1]{liu2024zeroth}.

The convolution is $L_F$-Lipschitz, as in the localized form of
\citep[Proposition~2.3]{lin2022gradient}.  Moreover, for
$\theta,\vartheta\in\Theta$
and every unit vector $a$, the boundary formula and
$\mathbb E|\langle a,U\rangle|\le
(\mathbb E\langle a,U\rangle^2)^{1/2}=m^{-1/2}$ give
\[
\left|\left\langle a,
\nabla\psi_\rho(\theta)-\nabla\psi_\rho(\vartheta)
\right\rangle\right|
\le\frac{\sqrt m\,L_F}{\rho}\,\|\theta-\vartheta\|.
\]
Hence $\psi_\rho$ is continuously differentiable and
\[
    \operatorname{Lip}(\nabla\psi_\rho)
    \le L_{\psi,2}:=\frac{\sqrt m\,L_F}{\rho}.
\]

Finally, Rademacher's theorem applies on the preceding ambient neighborhood,
and the weak gradient of the convolution is the average of
$\nabla F_{\max}$ at differentiability points in
$\mathbb B_m(\theta;\rho)$.  Each such gradient belongs to the corresponding
ambient Clarke subdifferential.  The Clarke graph over this compact ball is
compact and bounded, so its convex hull is closed and contains the average.
Consequently,
\[
    \nabla\psi_\rho(\theta)\in\partial_\rho F_{\max}(\theta),
\]
which is the localized form of \citep[Theorem~3.1]{lin2022gradient}; compare
\citep[Lemma~3.8]{liu2024zeroth}.  This ambient inclusion remains valid when
$\theta\in\partial\Theta$; only the smoothing ball needs the strict thickening
margin.

We now analyze the projected updates with the prescribed common stepsize
$\eta=L_{\psi,2}^{-1}=\rho/(\sqrt m\,L_F)$.

Given the \( (L_{\psi,2}= \sqrt{m}L_{F_{\max},1}\rho^{-1}) \)-smoothness of \( \psi_{\rho}(\theta) \), we have
\begin{align}
\psi_{\rho}(\theta_{i+1}) &\leq \psi_{\rho}(\theta_{i}) + \langle \nabla \psi_{\rho}(\theta_{i}), \theta_{i+1} - \theta_{i} \rangle + \frac{L_{\psi,2}}{2} \| \theta_{i+1} - \theta_{i} \|^{2}
\end{align}
Substituting \( \theta_{i+1} - \theta_{i} = -\eta \mathcal{G}_{\Theta}(\theta_{i}, \tilde{\mathsf{g}}_{\rho}(\theta_{i}); \eta) \), we have
\begin{align}
\psi_{\rho}(\theta_{i+1})
&\leq \psi_{\rho}(\theta_{i})
- \eta \langle \nabla \psi_{\rho}(\theta_{i}),
\mathcal{G}_{\Theta}(\theta_{i}, \tilde{\mathsf{g}}_{\rho}(\theta_{i}); \eta) \rangle\nonumber\\
&\quad+ \frac{L_{\psi,2}}{2} \eta^{2}
\| \mathcal{G}_{\Theta}(\theta_{i}, \tilde{\mathsf{g}}_{\rho}(\theta_{i}); \eta) \|^{2}
\end{align}
Expanding the inner product and rearranging terms
\begin{align}
\psi_{\rho}(\theta_{i+1})&\le \psi_{\rho}(\theta_{i}) - \eta \langle \tilde{\mathsf{g}}_{\rho}(\theta_{i}), \mathcal{G}_{\Theta}(\theta_{i}, \tilde{\mathsf{g}}_{\rho}(\theta_{i}); \eta) \rangle + \frac{L_{\psi,2}}{2} \eta^{2} \| \mathcal{G}_{\Theta}(\theta_{i}, \tilde{\mathsf{g}}_{\rho}(\theta_{i}); \eta) \|^{2} \nonumber \\
&\quad - \eta \langle \nabla \psi_{\rho}(\theta_{i}) - \tilde{\mathsf{g}}_{\rho}(\theta_{i}), \mathcal{G}_{\Theta}(\theta_{i}, \tilde{\mathsf{g}}_{\rho}(\theta_{i}); \eta) \rangle
\end{align}
Further simplifying and using the projection-optimality inequality
$\| \mathcal{G}_{\Theta}(\theta, v; \eta) \|^{2}\le
\langle \mathcal{G}_{\Theta}(\theta, v; \eta),v\rangle$ for every
$\theta\in\Theta$ and $v\in\mathbb R^m$, we have
\begin{align}
\MoveEqLeft[1]\psi_{\rho}(\theta_{i+1})\le\psi_{\rho}(\theta_{i}) - \left[\eta - \frac{L_{\psi,2}}{2} \eta^{2}\right] \cdot \| \mathcal{G}_{\Theta}(\theta_{i}, \tilde{\mathsf{g}}_{\rho}(\theta_{i}); \eta) \|^{2}\nonumber\\
& -\eta \langle \nabla \psi_{\rho}(\theta_{i}) - \tilde{\mathsf{g}}_{\rho}(\theta_{i}),\mathcal{G}_{\Theta}(\theta_{i}, \nabla\psi_{\rho}(\theta_{i}); \eta) \rangle \nonumber \\
&\quad + \eta \langle \nabla \psi_{\rho}(\theta_{i}) - \tilde{\mathsf{g}}_{\rho}(\theta_{i}),\mathcal{G}_{\Theta}(\theta_{i}, \nabla\psi_{\rho}(\theta_{i}); \eta)-\mathcal{G}_{\Theta}(\theta_{i}, \tilde{\mathsf{g}}_{\rho}(\theta_{i}); \eta) \rangle\nonumber\\
&\le \psi_{\rho}(\theta_{i}) - \left[\eta - \frac{L_{\psi,2}}{2} \eta^{2}\right] \cdot \| \mathcal{G}_{\Theta}(\theta_{i}, \tilde{\mathsf{g}}_{\rho}(\theta_{i}); \eta) \|^{2}\nonumber \\
&-\eta \langle \nabla \psi_{\rho}(\theta_{i}) - \tilde{\mathsf{g}}_{\rho}(\theta_{i}),\mathcal{G}_{\Theta}(\theta_{i}, \nabla\psi_{\rho}(\theta_{i}); \eta) \rangle\nonumber\\
&\quad+ \eta \|\nabla \psi_{\rho}(\theta_{i}) - \tilde{\mathsf{g}}_{\rho}(\theta_{i})\|^{2},
\end{align}
where the last inequality uses Cauchy--Schwarz and
\[
\|\mathcal{G}_{\Theta}(\theta,v_{1};\eta)
-\mathcal{G}_{\Theta}(\theta,v_{2};\eta)\|
\le \|v_{1}-v_{2}\|.
\]
Define $\mathsf{g}_{\rho}(\theta):=b_u^{-1}
\sum_{t=1}^{b_u}\mathsf{g}_{\rho}(\theta,u_t)$, where
\[
\mathsf{g}_{\rho}(\theta,u_t)
:=\frac{m}{2\rho}
[F_{\max}(\theta+\rho u_t)-F_{\max}(\theta-\rho u_t)]u_t.
\]
Then
\begin{align}
   \MoveEqLeft[2]\psi_{\rho}(\theta_{i+1})\le \psi_{\rho}(\theta_{i}) - \left[\eta - \frac{L_{\psi,2}}{2} \eta^{2}\right] \cdot \| \mathcal{G}_{\Theta}(\theta_{i}, \tilde{\mathsf{g}}_{\rho}(\theta_{i}); \eta) \|^{2}\nonumber\\
   &-\eta \langle \nabla \psi_{\rho}(\theta_{i}) - \mathsf{g}_{\rho}(\theta_{i}),\mathcal{G}_{\Theta}(\theta_{i}, \nabla\psi_{\rho}(\theta_{i}); \eta) \rangle \nonumber \\
&-\eta \langle \mathsf{g}_{\rho}(\theta_{i}) - \tilde{\mathsf{g}}_{\rho}(\theta_{i}),\mathcal{G}_{\Theta}(\theta_{i}, \nabla\psi_{\rho}(\theta_{i}); \eta) \rangle + \eta \|\nabla \psi_{\rho}(\theta_{i}) - \tilde{\mathsf{g}}_{\rho}(\theta_{i})\|^{2}, 
\end{align}
Since $\eta=L_{\psi,2}^{-1}$, we have
\begin{align}
   \MoveEqLeft[4] \| \mathcal{G}_{\Theta}(\theta_{i}, \tilde{\mathsf{g}}_{\rho}(\theta_{i}); \eta) \|^{2}\le\frac{2[\psi_{\rho}(\theta_{i})-\psi_{\rho}(\theta_{i+1})]}{\eta} \nonumber \\
   &-2\langle \nabla \psi_{\rho}(\theta_{i}) - \mathsf{g}_{\rho}(\theta_{i}),\mathcal{G}_{\Theta}(\theta_{i}, \nabla\psi_{\rho}(\theta_{i}); \eta) \rangle\nonumber\\
&-2 \langle \mathsf{g}_{\rho}(\theta_{i}) - \tilde{\mathsf{g}}_{\rho}(\theta_{i}),\mathcal{G}_{\Theta}(\theta_{i}, \nabla\psi_{\rho}(\theta_{i}); \eta) \rangle + 2 \|\nabla \psi_{\rho}(\theta_{i}) - \tilde{\mathsf{g}}_{\rho}(\theta_{i})\|^{2}. 
\end{align}
Then by using Young's inequality ($\langle \av,\bv\rangle\le c\|\av\|^{2}/2+\|\bv\|^{2}/(2c)$ for the vectors $\av,\bv$ and positive constant $c$, we have
\begin{align}
   \MoveEqLeft[4] \| \mathcal{G}_{\Theta}(\theta_{i}, \tilde{\mathsf{g}}_{\rho}(\theta_{i}); \eta) \|^{2}\le\frac{2[\psi_{\rho}(\theta_{i})-\psi_{\rho}(\theta_{i+1})]}{\eta}\nonumber\\
   & -2\langle \nabla \psi_{\rho}(\theta_{i}) - \mathsf{g}_{\rho}(\theta_{i}),\mathcal{G}_{\Theta}(\theta_{i}, \nabla\psi_{\rho}(\theta_{i}); \eta) \rangle+4 \| \mathsf{g}_{\rho}(\theta_{i}) - \tilde{\mathsf{g}}_{\rho}(\theta_{i})\|^{2} \nonumber \\
&+4^{-1}\|\mathcal{G}_{\Theta}(\theta_{i}, \nabla\psi_{\rho}(\theta_{i}); \eta)\|^{2} + 2 \|\nabla \psi_{\rho}(\theta_{i}) - \tilde{\mathsf{g}}_{\rho}(\theta_{i})\|^{2}. 
\end{align}
By \citep[Supplement, Lemma~D.1]{lin2022gradient},
$\mathbb{E}_{u_{t}}[\mathsf{g}_{\rho}(\theta_{i})\mid\theta_{i}]
=\nabla\psi_{\rho}(\theta_{i})$. Taking expectation given \( \theta_{i} \)
\begin{align}
\MoveEqLeft[4] \mathbb{E}[\| \mathcal{G}_{\Theta}(\theta_{i}, \tilde{\mathsf{g}}_{\rho}(\theta_{i}); \eta) \|^{2}\mid \theta_{i}]\le\nonumber\\
&\frac{2[\psi_{\rho}(\theta_{i})-\mathbb{E}[\psi_{\rho}(\theta_{i+1})\mid\theta_{i}]]}{\eta}+4\,\mathbb{E}[\| \mathsf{g}_{\rho}(\theta_{i}) - \tilde{\mathsf{g}}_{\rho}(\theta_{i})\|^{2}\mid\theta_{i}]+\nonumber \\
&\quad+4^{-1}\mathbb{E}[\|\mathcal{G}_{\Theta}(\theta_{i}, \nabla\psi_{\rho}(\theta_{i}); \eta)\|^{2}\mid\theta_{i}] + 2 \mathbb{E}[\|\nabla \psi_{\rho}(\theta_{i}) - \tilde{\mathsf{g}}_{\rho}(\theta_{i})\|^{2}\mid\theta_{i}].\label{eq_append_65}
\end{align}
The term involving $\nabla \psi_{\rho}(\theta_{i})-\mathsf{g}_{\rho}(\theta_{i})$ vanishes because $\mathbb{E}_{u_{t}}[\mathsf{g}_{\rho}(\theta_{i})\mid\theta_{i}]=\nabla\psi_{\rho}(\theta_{i})$ and $\mathcal{G}_{\Theta}(\theta_{i}, \nabla\psi_{\rho}(\theta_{i}); \eta)$ is deterministic given $\theta_i$. Moreover, by $(a+b)^2\le 2(a^2+b^2)$,
\begin{align}
\MoveEqLeft[3]
\mathbb{E}[\|\nabla \psi_{\rho}(\theta_{i})
- \tilde{\mathsf{g}}_{\rho}(\theta_{i})\|^{2}\mid \theta_{i}]\nonumber\\
&\le 2\,\mathbb{E}[\|\nabla \psi_{\rho}(\theta_{i})
-\mathsf{g}_{\rho}(\theta_{i})\|^{2}\mid \theta_{i}]\nonumber\\
&\quad+2\,\mathbb{E}[\|\mathsf{g}_{\rho}(\theta_{i})
-\tilde{\mathsf{g}}_{\rho}(\theta_{i})\|^{2}\mid \theta_{i}].\nonumber
\end{align}
Substituting this into \eqref{eq_append_65}, we obtain
\begin{align}
\MoveEqLeft[4] \mathbb{E}[\| \mathcal{G}_{\Theta}(\theta_{i}, \tilde{\mathsf{g}}_{\rho}(\theta_{i}); \eta) \|^{2}\mid \theta_{i}] \le \nonumber\\
&\frac{2[\psi_{\rho}(\theta_{i})-\mathbb{E}[\psi_{\rho}(\theta_{i+1})\mid\theta_{i}]]}{\eta}
+4\,\mathbb{E}[\|\nabla \psi_{\rho}(\theta_{i})-\mathsf{g}_{\rho}(\theta_{i})\|^{2}\mid \theta_{i}] \nonumber\\
&\qquad +8\,\mathbb{E}[\| \mathsf{g}_{\rho}(\theta_{i}) - \tilde{\mathsf{g}}_{\rho}(\theta_{i})\|^{2}\mid\theta_{i}]
+4^{-1}\mathbb{E}[\|\mathcal{G}_{\Theta}(\theta_{i}, \nabla\psi_{\rho}(\theta_{i}); \eta)\|^{2}\mid\theta_{i}].\label{eq_append_65b}
\end{align}

\vspace{0.1cm}
In order to give a bound on $\mathbb{E}[\| \mathcal{G}_{\Theta}(\theta_{i}, \nabla \psi_{\rho}(\theta_{i}); \eta) \|^{2} |\theta_{i}] $, we have
\begin{align}
\MoveEqLeft[4]\mathbb{E}[\| \mathcal{G}_{\Theta}(\theta_{i}, \nabla \psi_{\rho}(\theta_{i}); \eta) \|^{2} |\theta_{i}] \leq 2 \mathbb{E}[\| \mathcal{G}_{\Theta}(\theta_{i}, \tilde{\mathsf{g}}_{\rho}(\theta_{i}); \eta) \|^{2} |\theta_{i}] \nonumber \\
& + 2 \mathbb{E}[\| \mathcal{G}_{\Theta}(\theta_{i}, \tilde{\mathsf{g}}_{\rho}(\theta_{i}); \eta) - \mathcal{G}_{\Theta}(\theta_{i}, \nabla \psi_{\rho}(\theta_{i}); \eta) \|^{2} |\theta_{i}]\nonumber\\
&\overset{(a)}{\le} 2 \mathbb{E}[\| \mathcal{G}_{\Theta}(\theta_{i}, \tilde{\mathsf{g}}_{\rho}(\theta_{i}); \eta) \|^{2} |\theta_{i}]+ 2 \mathbb{E}[\|\tilde{\mathsf{g}}_{\rho}(\theta_{i}) - \nabla \psi_{\rho}(\theta_{i}) \|^{2} |\theta_{i}]\nonumber\\
&\le 2 \mathbb{E}[\| \mathcal{G}_{\Theta}(\theta_{i}, \tilde{\mathsf{g}}_{\rho}(\theta_{i}); \eta) \|^{2} |\theta_{i}] \nonumber\\
&\qquad +4 \mathbb{E}[\|\nabla \psi_{\rho}(\theta_{i})-\mathsf{g}_{\rho}(\theta_{i}) \|^{2} |\theta_{i}]
+4 \mathbb{E}[\|\tilde{\mathsf{g}}_{\rho}(\theta_{i}) - \mathsf{g}_{\rho}(\theta_{i}) \|^{2} |\theta_{i}].\label{eq_append_66}
\end{align}
where (a) follows by $\|\mathcal{G}_{\Theta}(\theta, v_{1}; \eta)-\mathcal{G}_{\Theta}(\theta, v_{2}; \eta)\|\le \|v_{1}-v_{2}\|$. Incorporating \eqref{eq_append_65b} in \eqref{eq_append_66}, we get
\begin{align}
    \MoveEqLeft[2]\mathbb{E}[\| \mathcal{G}_{\Theta}(\theta_{i}, \nabla \psi_{\rho}(\theta_{i}); \eta) \|^{2} |\theta_{i}]\le \frac{4[\psi_{\rho}(\theta_{i})-\mathbb{E}[\psi_{\rho}(\theta_{i+1})\mid\theta_{i}]]}{\eta}\nonumber \\
&+12\,\mathbb{E}[\|\nabla \psi_{\rho}(\theta_{i})-\mathsf{g}_{\rho}(\theta_{i})\|^{2}\mid\theta_{i}]
+20\,\mathbb{E}[\| \mathsf{g}_{\rho}(\theta_{i}) - \tilde{\mathsf{g}}_{\rho}(\theta_{i})\|^{2}\mid\theta_{i}]\nonumber\\
&\qquad +2^{-1}\mathbb{E}[\|\mathcal{G}_{\Theta}(\theta_{i}, \nabla\psi_{\rho}(\theta_{i}); \eta)\|^{2}\mid\theta_{i}].
\end{align}
Finally, we get
\begin{align}\label{eq_append_0024}
\MoveEqLeft[4]\mathbb{E}[\| \mathcal{G}_{\Theta}(\theta_{i}, \nabla \psi_{\rho}(\theta_{i}); \eta) \|^{2} |\theta_{i}]\le\nonumber\\
&\frac{8[\psi_{\rho}(\theta_{i})-\mathbb{E}[\psi_{\rho}(\theta_{i+1})\mid\theta_{i}]]}{\eta}
+24\,\mathbb{E}[\|\nabla \psi_{\rho}(\theta_{i})-\mathsf{g}_{\rho}(\theta_{i})\|^{2}\mid\theta_{i}] \nonumber\\
&\qquad +40\,\mathbb{E}[\| \mathsf{g}_{\rho}(\theta_{i}) - \tilde{\mathsf{g}}_{\rho}(\theta_{i})\|^{2}\mid\theta_{i}].
\end{align}
Up to this point, the analysis is  standard for constrained zeroth-order methods and does not use the specific structure of $F_{\max}$.

\textit{Step 2: Controlling the directional and oracle errors.} We first bound the variance coming from the random directions.
Conditionally on the complete history before iteration
\(i\), the variables
$\{\mathsf{g}_{\rho}(\theta_i,u_{i,t})\}_{t=1}^{b_u}$ are i.i.d. and satisfy
$\mathbb{E}_{u_t}[\mathsf{g}_{\rho}(\theta_i,u_{i,t})\mid \theta_i]
=\nabla \psi_{\rho}(\theta_i)$. We therefore have
\begin{align}\label{eq_append_0025}
\MoveEqLeft[4]\mathbb{E}[\|\nabla \psi_{\rho}(\theta_i)-\mathsf{g}_{\rho}(\theta_i)\|^{2}\mid \theta_i]\nonumber\\
&= \frac{1}{b_u^{2}}\sum_{t=1}^{b_u}\mathbb{E}[\|\mathsf{g}_{\rho}(\theta_i,u_{i,t})-\nabla \psi_{\rho}(\theta_i)\|^{2}\mid \theta_i]\nonumber\\
&\le \frac{1}{b_u}\,\mathbb{E}[\|\mathsf{g}_{\rho}(\theta_i,u)\|^{2}\mid \theta_i]\nonumber\\
&\le \frac{m^{2}L_{F_{\max},1}^{2}}{b_u},
\end{align}
where we used
\[
\|\mathsf{g}_{\rho}(\theta_i,u)\|
= \frac{m}{2\rho}\,|F_{\max}(\theta_i+\rho u)-F_{\max}(\theta_i-\rho u)|
\le mL_{F_{\max},1}.
\]
For each direction and sign, put
\[
e_{i,t}^{\pm}
:=\widetilde\psi(\theta_i\pm\rho u_{i,t})
  -F_{\max}(\theta_i\pm\rho u_{i,t}).
\]
Condition first on $\theta_i$ and on all directions
$u_{i,1},\ldots,u_{i,b_u}$.  The training and validation samples at each
query are fresh and have the Gibbs law at that same perturbed parameter.
Thus \eqref{eq_append_21}, applied on $\Theta^{+\bar\rho}$ and followed by
the tower property, yields a finite constant $C_{\mathrm{or}}$, independent
of $i,t,\theta_i$ and the directions, such that
\[
\mathbb E[(e_{i,t}^{\pm})^2\mid
          \theta_i,u_{i,1},\ldots,u_{i,b_u}]
\le C_{\mathrm{or}}\epsilon_v^2.
\]
Moreover, Jensen's inequality gives the pointwise estimate
\[
\|\widetilde{\mathsf g}_\rho(\theta_i)-\mathsf g_\rho(\theta_i)\|^2
\le
\frac{m^2}{2\rho^2b_u}
\sum_{t=1}^{b_u}\bigl[(e_{i,t}^{+})^2+(e_{i,t}^{-})^2\bigr].
\]
Consequently,
\begin{align}
\mathbb E[\|\widetilde{\mathsf g}_\rho(\theta_i)
             -\mathsf g_\rho(\theta_i)\|^2\mid\theta_i]
\le C_{\mathrm{or}}\frac{m^2}{\rho^2}\epsilon_v^2.
\label{eq_append_0026}
\end{align}
Here \eqref{eq_append_21} is available whenever
\[
|I|\ge4\widetilde C_f\epsilon_v^{-2-2\dimk}
\]
and, for a common finite $C_{\mathrm{PSGD}}$,
\[
L\ge C_{\mathrm{PSGD}}\dimk^2\mathrm B^2
\epsilon_v^{-2-2\dimk}
\bigl[1+\log^2(\epsilon_v^{-1})\bigr].
\]
No unbiasedness of the approximate oracle is used.

Substituting \eqref{eq_append_0025}--\eqref{eq_append_0026} into
\eqref{eq_append_0024}, averaging over
$i=0,\ldots,N-1$, and telescoping gives
\begin{align*}
\frac1N\sum_{i=0}^{N-1}
\mathbb E\!\left[
\|\mathcal G_\Theta(\theta_i,\nabla\psi_\rho(\theta_i);\eta)\|^2
\right]
&\le
\frac{8[\psi_\rho(\theta_0)-\mathbb E\psi_\rho(\theta_N)]}{N\eta}
 +\frac{24m^2L_F^2}{b_u}\\
&\quad+40C_{\mathrm{or}}\frac{m^2}{\rho^2}\epsilon_v^2.
\end{align*}
Because $\psi_\rho$ is $L_F$-Lipschitz and every iterate lies in $\Theta$,
\[
\psi_\rho(\theta_0)-\mathbb E\psi_\rho(\theta_N)
\le L_F\operatorname{diam}(\Theta).
\]
With $\eta=\rho/(\sqrt mL_F)$, the preceding estimate becomes
\[
\frac1N\sum_{i=0}^{N-1}
\mathbb E\|\mathcal G_\Theta(
\theta_i,\nabla\psi_\rho(\theta_i);\eta)\|^2
\le
\frac{8\sqrt mL_F^2\operatorname{diam}(\Theta)}{N\rho}
+\frac{24m^2L_F^2}{b_u}
+40C_{\mathrm{or}}\frac{m^2}{\rho^2}\epsilon_v^2.
\]
The choice $\epsilon_v=\rho\epsilon/m$ controls the last term.  The
sufficient lower bounds on $N$ and $b_u$ in
\Cref{thm_cvg_alg_zeroth_goldstein} control the first two terms after
enlarging the fixed constants $C_N$ and $C_u$; for the first term, use
$m\ge1$ and $\rho<\bar\rho$.  Hence the right-hand side is at most a fixed
multiple of $\epsilon^2$.

Choose a deterministic index
\[
i_\star\in\argmin_{0\le i\le N-1}
\mathbb E\|\mathcal G_\Theta(
\theta_i,\nabla\psi_\rho(\theta_i);\eta)\|^2.
\]
The ambient Goldstein inclusion and Jensen's inequality imply
\begin{align*}
\mathbb E\!\left[
\min_{g\in\partial_\rho F_{\max}(\theta_{i_\star})}
\|\mathcal G_\Theta(\theta_{i_\star},g;\eta)\|
\right]
&\le
\mathbb E\|\mathcal G_\Theta(
\theta_{i_\star},\nabla\psi_\rho(\theta_{i_\star});\eta)\|\\
&\le C_{\mathrm{res}}\epsilon.
\end{align*}
This implication remains valid when $\theta_{i_\star}\in\partial\Theta$:
the strict margin is required for the ambient smoothing ball, whereas the
gradient mapping itself uses the prescribed projection onto $\Theta$.

Let \(J\sim\operatorname{Unif}\{0,\ldots,N-1\}\) be independent of the
algorithmic history and set \(\widehat\theta:=\theta_J\), as in
\Cref{alg: Stochastic-Zeroth-order}.  Conditioning on the history gives the
exact identity
\[
\mathbb E\!\left[
\|\mathcal G_\Theta(
\theta_J,\nabla\psi_\rho(\theta_J);\eta)\|^2\right]
=
\frac1N\sum_{i=0}^{N-1}
\mathbb E\!\left[
\|\mathcal G_\Theta(
\theta_i,\nabla\psi_\rho(\theta_i);\eta)\|^2\right].
\]
Thus Jensen's inequality and
\(\nabla\psi_\rho(\theta_J)\in\partial_\rho F_{\max}(\theta_J)\) yield
\[
\begin{aligned}
\mathbb E\!\left[
\min_{g\in\partial_\rho F_{\max}(\widehat\theta)}
\|\mathcal G_\Theta(\widehat\theta,g;\eta)\|\right]
&\le
\mathbb E\!\left[
\|\mathcal G_\Theta(
\theta_J,\nabla\psi_\rho(\theta_J);\eta)\|\right]\\
&\le C_{\mathrm{res}}\epsilon .
\end{aligned}
\]
This is the asserted guarantee for the iterate actually returned by the
algorithm, and drawing \(J\) uses no oracle call.

Finally, Algorithm~\ref{alg: Stochastic-Zeroth-order} uses at most
$2Nb_u(L+|I|)$ Gibbs samples; this count is exact when
\(\mathrm B>0\).  For the explicit ceiling choices in the
complexity paragraph following \Cref{thm_cvg_alg_zeroth_goldstein} (rather
than arbitrary larger admissible integers), and in the stated standard
small-accuracy regime in which the unrounded sufficient orders are at least
one, substitution gives
\[
2Nb_u(L+|I|)
=\mathcal O\!\left(
\frac{m^{2\dimk+5}
\left\{
\dimk^2\mathrm B^2
\left[1+\log^2\!\left(\frac{m}{\rho\epsilon}\right)\right]
+\widetilde C_f
\right\}}
{\epsilon^{2\dimk+6}\rho^{2\dimk+4}}
\right),
\]
with all polynomial exponents unchanged.
\end{proof}

\section{Computational Complexity of Gibbs Oracle}\label[appendix]{append_comp_complex_gibbs}
In this section, we modify the result of \Cref{sec:6} for the setting where the Gibbs samples are generated approximately by LMC rather than drawn exactly from the Gibbs distribution.

\subsection*{Proof of \Cref{prop_chewi2024logconcave_1}}
\phantomsection
\label[appendix]{append_proof_prop_chewi2024logconcave_1}
\begin{proof}
Let \((\pi_t)_{t\ge0}\) be the law of the Langevin diffusion with invariant
measure \(\pi\), initialized from \(\widehat\mu_0\), and put
\(\ell_{\mathrm R}:=\log(1/\varepsilon_{\mathrm R})\).  The proof has three
ingredients: order-three diffusion mixing, an order-four discretization
bound, and the weak triangle inequality relating them.

\paragraph{Diffusion mixing.}
Set \(r(t):=\mathbb R_3(\pi_t\Vert\pi)\).  The Poincar\'e part of the proof
of \cite[Theorem~2.2.25]{chewi2024logconcave}, applied at order three and
using \(C_{\mathrm{PI}}(\pi)\le\alpha^{-1}\), gives
\[
 r'(t)\le-\frac{4\alpha}{3}\bigl(1-e^{-r(t)}\bigr).
\]
This yields the two stages recorded in the cited theorem.  While \(r(t)\ge1\),
the elementary bound \(1-e^{-r(t)}\ge1/2\) gives
\(r'(t)\le-2\alpha/3\).  Hence the first time \(t_1\) at which
\(r(t_1)\le1\) satisfies
\[
 t_1\le\frac{3R_3}{2\alpha},
\]
with \(t_1=0\) if \(r(0)\le1\).  For \(t\ge t_1\), the bound
\(1-e^{-r}\ge r/2\) on \(0\le r\le1\) gives
\[
 r(t)\le\exp\!\left[-\frac{2\alpha}{3}(t-t_1)\right].
\]
Now choose \(C_T\ge6\) and set
\(T=C_T\alpha^{-1}(R_3+\ell_{\mathrm R})\).  Since
\(\ell_{\mathrm R}\ge\log2\), one has
\[
 T-t_1
 \ge\frac{3}{2\alpha}
       \log\!\left(\frac{4}{\varepsilon_{\mathrm R}^2}\right).
\]
Consequently,
\[
 \mathbb R_3(\pi_T\Vert\pi)=r(T)
 \le \frac{\varepsilon_{\mathrm R}^2}{4}.
\]

\paragraph{Order-three score transfer.}
By \cite[Lemma~6.2.7]{chewi2024logconcave}, there is a universal constant
\(C>0\) such that, for every \(s>0\),
\[
 \pi\!\left\{\|\nabla G\|\ge
 C\sqrt\beta(\sqrt d+s)\right\}\le e^{-s^2}.
\]
Apply \cite[Lemma~6.2.8]{chewi2024logconcave} at order three.  Data
processing for the diffusion semigroup gives
\(\mathbb R_3(\pi_t\Vert\pi)
 \le\mathbb R_3(\widehat\mu_0\Vert\pi)\le R_3\); hence, uniformly for
\(0\le t\le T\),
\[
 \pi_t\!\left\{\|\nabla G\|\ge
 C\sqrt\beta(\sqrt d+s)\right\}
 \le
 \exp\!\left[-\frac23(s^2-R_3)\right].
\]
We record explicitly the exponential-moment consequence used in the
discretization proof.  The last display implies that there are universal
constants \(c_1,C_1>0\) such that, with
\(r_0:=C_1\sqrt{\beta(d+R_3)}\),
\[
 \pi_t\{\|\nabla G\|\ge r\}
 \le e^{-c_1r^2/\beta},
 \qquad r\ge r_0,\quad 0\le t\le T.
\]
For \(Z:=\|\nabla G\|\) and \(a\ge0\), integration of the tail gives
\[
 \mathbb E_{\pi_t}e^{aZ^2}
 =1+\int_0^\infty 2ar e^{ar^2}\pi_t\{Z>r\}\,\ud r.
\]
Split this integral at \(r_0\).  If \(0\le a\le c_1/(2\beta)\), then
\[
 \mathbb E_{\pi_t}e^{aZ^2}
 \le e^{ar_0^2}
 +2a\int_{r_0}^\infty
 r e^{-(c_1/\beta-a)r^2}\,\ud r
 \le e^{C_1^2a\beta(d+R_3)}+C a\beta.
\]
Since \(d\ge1\), enlarging the universal constant gives
\[
 \log\mathbb E_{\pi_t}e^{a\|\nabla G\|^2}
 \le C a\beta(d+R_3),
 \qquad 0\le a\le \frac{c}{\beta},\quad 0\le t\le T.
\]
Thus the order-four discretization calculation needs only the order-three
warm start \(R_3\); no order-four initial divergence is required.

\paragraph{Discretization and the joint choice.}
Apply \cite[Theorem~6.2.1]{chewi2024logconcave} at order four and then use
data processing from path space to the endpoint.  At this fixed order, the
exponential-moment parameter in that theorem is a universal multiple of
\(\beta^2\bar h^2T\).  The restriction
\(\bar h\le c_h/(2\beta^{3/2}\sqrt T)\) places this parameter in the range
\(a\le c/\beta\) verified above, while
\(\bar h\le c_h/(4\beta^2T)\) is the theorem's other step restriction.
Substitution of the preceding exponential-moment bound therefore gives a
universal constant \(C_D\) such that
\[
 \mathbb R_4(\widehat\mu_n\Vert\pi_T)
 \le C_D\left[
 \beta^3\bar h^2T(d+R_3)+\beta^2d\bar hT
 \right],
\]
provided
\(\bar h\le c_h/(4\beta^2T)\) and
\(\bar h\le c_h/(2\beta^{3/2}\sqrt T)\), after decreasing the universal
constant \(c_h>0\) if necessary.

Now use the construction in the proposition:
\[
 h_0=c_h\min\left\{
 \frac{1}{4\beta^2T},
 \frac{1}{2\beta^{3/2}\sqrt T},
 \frac{\varepsilon_{\mathrm R}^2}{2\beta^2dT}
 \right\},
 \qquad
 n=\left\lceil\frac{T}{h_0}\right\rceil,
 \qquad
 \bar h=\frac{T}{n}.
\]
Then \(n\bar h=T\) and \(\bar h\le h_0\), so all three step restrictions
hold for the actual chain.  In particular,
\[
 \beta^2d\bar hT
 \le \frac{c_h}{2}\varepsilon_{\mathrm R}^2
\]
and
\[
 \beta^3\bar h^2T(d+R_3)
 \le
 \frac{c_h^2\varepsilon_{\mathrm R}^4}{4}
 \frac{d+R_3}{\beta d^2T}.
\]
Moreover,
\[
 \frac{d+R_3}{\beta d^2T}
 =\frac{\alpha}{C_T\beta}
   \frac{d+R_3}{d^2(R_3+\ell_{\mathrm R})}
 \le \frac{1}{C_T}
 \left(\frac{1}{\log 2}+1\right),
\]
where we used \(\alpha\le1\), \(\beta\ge1\), \(d\ge1\), and
\(\ell_{\mathrm R}\ge\log2\).  We may therefore choose \(c_h\) universally
small, and then choose one universal
\(\varepsilon_{\mathrm{src}}\in(0,1/2]\), so that
\[
 \mathbb R_4(\widehat\mu_n\Vert\pi_T)
 \le \frac{\varepsilon_{\mathrm R}^2}{2}
 \qquad
 (0<\varepsilon_{\mathrm R}\le\varepsilon_{\mathrm{src}}).
\]

Finally, \cite[Lemma~6.2.5]{chewi2024logconcave} with target order two and
interpolation parameter \(1/2\) gives the precise weak triangle inequality
\[
 \mathbb R_2(\widehat\mu_n\Vert\pi)
 \le \frac32\mathbb R_4(\widehat\mu_n\Vert\pi_T)
      +\mathbb R_3(\pi_T\Vert\pi)
 \le \varepsilon_{\mathrm R}^2.
\]

It remains to bound the iteration count for this same joint choice.  Since
\(n=\lceil T/h_0\rceil\),
\[
 n\le1+\frac{1}{c_h}\max\left\{
 4\beta^2T^2,
 2\beta^{3/2}T^{3/2},
 2\beta^2dT^2\varepsilon_{\mathrm R}^{-2}
 \right\}.
\]
Increasing \(C_T\), if necessary, ensures \(T\ge1\).  The assumptions
\(\beta\ge1\), \(d\ge1\), and
\(\varepsilon_{\mathrm R}\le1/2\) then show that the last member of the
maximum dominates the other two up to a universal factor and also absorbs
the leading \(1\).  Consequently,
\[
 n\le C\beta^2d\varepsilon_{\mathrm R}^{-2}T^2
 \le C_{\mathrm{Ch}}\beta^2\alpha^{-2}d
 \varepsilon_{\mathrm R}^{-2}
 \left[R_3^2+\log^2(1/\varepsilon_{\mathrm R})\right],
\]
where the last inequality uses \((a+b)^2\le2(a^2+b^2)\).  This proves the
claimed bound for the same pair \((\bar h,n)\) that gives the R\'enyi
accuracy.
\end{proof}

\begin{lem}[R\'enyi divergence for the \(\lambda\)-scaled Gaussian initialization]\label{lem_gaussian_lmc_initialization}
Let \(\mathcal P\subset\mathbb R^m\) be a compact parameter set on which
\Cref{mu_PL_assum,assum_sample_comp} and all constants used below hold
uniformly.  For \(\vartheta\in\mathcal P\), let
\(g^{\star}(\vartheta):=\min_{x\in\mathbb{R}^{d}}g(\vartheta,x)\) and
\[
\mu_{\vartheta}^{\lambda}(\ud x)
=
(Z_{\vartheta}^{\lambda})^{-1}
\exp(-g(\vartheta,x)/\lambda)\ud x,
\qquad
0<\lambda\le 1.
\]
Fix \(0<c<3/(2L_{g,2})\) and initialize the LMC chain from
\[
\nu_{0}^{\lambda}:=\mathcal{N}(0,c\lambda I_{d}).
\]
Then there exist constants \(C_{\mathrm{init}},C'_{\mathrm{init}}<\infty\),
independent of \(\vartheta\) and \(\lambda\), such that
\[
\sup_{\vartheta\in\mathcal P}
\mathbb{R}_{3}\!\left(
\nu_{0}^{\lambda}\,\middle\|\,\mu_{\vartheta}^{\lambda}\right)
\le
\frac{C_{\mathrm{init}}}{\lambda}
+C'_{\mathrm{init}}\bigl(1+\log(1/\lambda)\bigr)
.
\]
The constants may depend on \(d\) and the fixed problem
data on \(\mathcal P\); no dimension-free estimate is asserted.
\end{lem}
\begin{proof}
Fix \(\vartheta\in\mathcal P\).
Let \(p_{\lambda}\) denote the density of
\(\nu_{0}^{\lambda}\).
By \eqref{eq_tail_EB_implies_QG}, there is a constant
\(C_{Z}<\infty\), independent of \(\vartheta\) and
\(\lambda\in(0,1]\), such that
\[
Z_{\vartheta}^{\lambda}
=\int_{\mathbb{R}^{d}}\exp(-g(\vartheta,x)/\lambda)\ud x
\le C_{Z}\exp(-g^{\star}(\vartheta)/\lambda).
\]
The order-3 R\'enyi divergence satisfies
\[
\mathbb{R}_{3}\!\left(
\nu_{0}^{\lambda}\,\middle\|\,\mu_{\vartheta}^{\lambda}\right)
=
\frac12\log\!\left[
(Z_{\vartheta}^{\lambda})^{2}
\int_{\mathbb{R}^{d}}p_{\lambda}(x)^{3}
\exp(2g(\vartheta,x)/\lambda)\ud x
\right].
\]
Set \(a:=3/(2c)-L_{g,2}>0\). By \(L_{g,2}\)-smoothness in \(x\),
\[
g(\vartheta,x)
\le
g(\vartheta,0)+\langle \partial_x g(\vartheta,0),x\rangle
+\frac{L_{g,2}}{2}\|x\|^{2}.
\]
Since \(p_{\lambda}(x)=(2\pi c\lambda)^{-d/2}\exp(-\|x\|^{2}/(2c\lambda))\), completing the square gives
\[
\int_{\mathbb{R}^{d}}p_{\lambda}(x)^{3}
\exp(2g(\vartheta,x)/\lambda)\ud x
\le
C\lambda^{-d}
\exp\!\left(
\frac{2g(\vartheta,0)+
\|\partial_x g(\vartheta,0)\|^{2}/a}{\lambda}
\right),
\]
where \(C<\infty\) is independent of \(\vartheta\) and \(\lambda\).
Combining the last two displays,
\[
\mathbb{R}_{3}\!\left(
\nu_{0}^{\lambda}\,\middle\|\,\mu_{\vartheta}^{\lambda}\right)
\le
\frac{g(\vartheta,0)-g^{\star}(\vartheta)}{\lambda}
+\frac{\|\partial_x g(\vartheta,0)\|^{2}}{2a\lambda}
+\frac{d}{2}\log(1/\lambda)+C.
\]
The first two numerators are uniformly bounded over
\(\vartheta\in\mathcal P\) because \(\mathcal P\) is compact and
\(g\), \(\partial_x g\), and \(g^{\star}\) are continuous.  This proves the
claimed uniform bound.
\end{proof}
In Section~5 we apply
\Cref{lem_gaussian_lmc_initialization} with
\(\mathcal P=\Theta^{+\bar\rho}\), on which the lower-level constants are
assumed uniform.
Set
\[
 \lambda_{\mathrm{adm}}
 :=\min\{1,\lambda_{\mathrm{mom}},\lambda_{\mathrm{PI}}^\star\},
 \qquad
 \mathcal R_3(\lambda)
 :=\sup_{\vartheta\in\Theta^{+\bar\rho}}
 \mathbb R_3\!\left(
 \mathcal N(0,c\lambda I_d)\,\middle\|\,\mu_\vartheta^\lambda
 \right).
\]
The preceding estimate and the boundedness of
\(s(1+\log(1/s))\) on \((0,1]\) imply
\begin{equation}\label{eq_uniform_lmc_warm_factor}
 \mathfrak R_3^\star
 :=\sup_{0<\lambda\le\lambda_{\mathrm{adm}}}
 \lambda\,\mathcal R_3(\lambda)<\infty.
\end{equation}
This supremum is uniform in \(\vartheta\) and \(\lambda\), but it may depend
on the ambient dimension \(d\) and on the fixed problem constants.  In
particular, neither the displayed R\'enyi bound nor
\eqref{eq_uniform_lmc_warm_factor} is asserted to be dimension-free.

\subsection{Runtime Complexity of \Cref{alg:PSGD} by Using LMC}
We replace the exact samples
\(\{\widetilde X^{(l)}\}_{l=0}^{L-1}\) from
\(\mu_\theta^\lambda\) in \Cref{alg:PSGD} by endpoints
\(\{\widetilde X_l^n\}_{l=0}^{L-1}\) of independent, independently
initialized, \(n\)-step LMC chains targeting that Gibbs law.  Denote their
common law by \(\widehat\mu_{n,\theta}\).  For \(\mathrm B>0\), set
\[
    \sigma:=1+\delta^{-1},
    \qquad
    \eta:=\frac{2\mathrm B}{\sigma\sqrt L}.
\]
If \(\mathrm B=0\), the feasible interval is the singleton \(\{0\}\), the
optimization error vanishes, and the PSGD updates below are unnecessary.
\begin{algorithm}[t]
\caption{$\text{PSGD-LMC}(\theta,\beta_{0})$}
\textbf{Input:} $\theta$, $\beta_{0}\in[-\mathrm B,\mathrm B]$, $L$, $n$,
	the standard-potential LMC step $\bar h$
	(equivalently, the physical step $h=\bar h/\lambda$ in
	\eqref{eqn:langevin}), and
$\eta=2\mathrm B/(\sigma\sqrt L)$.
    \begin{algorithmic}[1]
        \FOR{$l$ = $0$ to $L-1$}
        \STATE $\beta^{l+1}=\text{Proj}_{[-\mathrm{B},\mathrm{B}]}\!\left(\beta^{l}-\eta\,\widetilde \partial_{\beta}\phivar(\theta,\beta^{l};\widetilde X^{n}_{l})\right)$\\
        \,[$\widetilde X_l^n\sim\widehat\mu_{n,\theta}$ are endpoints of independent LMC chains, independently initialized and independent of the preceding iterates.]
        \ENDFOR
        \STATE $\tilde\beta=\frac{1}{L}\sum_{l=0}^{L-1}\beta^{l}$
    \end{algorithmic}\label{alg:PSGD_LMC}
\end{algorithm}
The corresponding approximate objective is
\begin{align}
    \widetilde{\phi}_{\lambda,\delta}(\theta,\beta)
    := \beta + \frac{1}{\delta} \mathbb{E}_{X \sim \widehat{\mu}_{n,\theta}} \left[ \max\{ f(\theta, X) - \beta, 0 \} \right].
\end{align}
Thus \(\widehat\mu_{n,\theta}\) is a probability law, not an
empirical measure.
Since
\[
    \bigl| \tilde{\phi}_{\lambda,\delta}(\theta,\beta) - \phi_{\lambda,\delta}(\theta,\beta) \bigr|
\le \frac{L_{f,1}}{\delta}\,\mathbb{W}_1(\widehat{\mu}_{n,\theta}, \mu_{\theta}^{\lambda}),
\]
we can invoke the LMC Wasserstein bound from \Cref{section_gradient_oracle_complexity} to obtain
\begin{equation}\label{eq_083083}
\bigl| \tilde{\phi}_{\lambda,\delta}(\theta,\beta) - \phi_{\lambda,\delta}(\theta,\beta) \bigr|
\le \frac{L_{f,1}}{\delta}\,
C_W(C_{\mathrm{PI}}^\star)^{1/2}\varepsilon_{\mathrm R}.
\end{equation}
Setting
\(\varepsilon_{\mathrm R}:=\delta\epsilon_v\),
\(\delta=\epsilon_v^{\dimk}\), and
\(\lambda=\epsilon_v^{2(\dimk+1)}\), and choosing
\[
 (\bar h,n)
 =\bigl(
 \bar h_\star(\lambda,\varepsilon_{\mathrm R}),
 n_\star(\lambda,\varepsilon_{\mathrm R})
 \bigr)
\]
according to \eqref{eq_uniform_lmc_joint_policy}, gives
\(\mathbb R_2(\widehat\mu_{n,\theta}\|\mu_\theta^\lambda)
\le\varepsilon_{\mathrm R}^2\), provided
\(\varepsilon_{\mathrm R}\le\varepsilon_{\mathrm{Ch}}\).  The common
upper bound \eqref{eq_uniform_lmc_mixing_bound}, together with
\eqref{eq_uniform_lmc_warm_factor}, yields
\[
 n=n_\star(\lambda,\varepsilon_{\mathrm R})
 \le \mathfrak C_{\mathrm{mix}}d\,
       \epsilon_v^{-10(\dimk+1)}.
\]
Here \(\mathfrak C_{\mathrm{mix}}\) may depend on \(d\) and the fixed
problem constants.  With the Gaussian initialization in
\Cref{lem_gaussian_lmc_initialization}, the resulting Wasserstein estimate
therefore yields
\[
\bigl| \tilde{\phi}_{\lambda,\delta}(\theta,\beta) - \phi_{\lambda,\delta}(\theta,\beta) \bigr|
= \mathcal{O}(\epsilon_{v}).
\]
Set
\[
\begin{aligned}
 m_{\lambda,\delta}(\theta)
 &:=\min_{\beta\in\mathbb R}\phi_{\lambda,\delta}(\theta,\beta)
 =\min_{\beta\in[-\mathrm B,\mathrm B]}
       \phi_{\lambda,\delta}(\theta,\beta),\\
 \widetilde m_{\lambda,\delta}^{\,\mathrm B}(\theta)
 &:=\min_{\beta\in[-\mathrm B,\mathrm B]}
       \widetilde\phi_{\lambda,\delta}(\theta,\beta).
\end{aligned}
\]
The equality for \(m_{\lambda,\delta}\) follows from \Cref{lem_3}, whereas
\eqref{eq_083083} gives
\[
 \left|\widetilde m_{\lambda,\delta}^{\,\mathrm B}(\theta)
       -m_{\lambda,\delta}(\theta)\right|
 \le
 \sup_{\beta\in[-\mathrm B,\mathrm B]}
 \left|\widetilde\phi_{\lambda,\delta}(\theta,\beta)
       -\phi_{\lambda,\delta}(\theta,\beta)\right|
 =\mathcal O(\epsilon_v).
\]
No location claim is made for an unrestricted minimizer of
\(\widetilde\phi_{\lambda,\delta}\).
We then define the LMC–based estimator of $\tilde{F}_{\SQG}(\theta)$ by
\begin{equation}\label{eq:def_hat_psi_LMC}
    \hat{\psi}(\theta)
    :=
    \phi_{\lambda,\delta}(\theta,\tilde{\beta};\hat{X}^{n}_I)
    = \tilde{\beta}
      + \frac{1}{\delta|I|}
        \sum_{i\in I} \max\{f(\theta,\hat{X}^{n}_i) - \tilde{\beta},0\},
\end{equation}
where $\widehat X^{n}_I := \{\widehat X^{n}_i\}_{i\in I}$ are
i.i.d.\ endpoints with law $\widehat\mu_{n,\theta}$ obtained from fresh,
independently initialized LMC chains.  This validation family is independent
of every chain used by \Cref{alg:PSGD_LMC}, and hence independent of
$\widetilde\beta$.
We then restate \cref{lem_5_2} for \cref{alg:PSGD_LMC}.
\begin{lem}\label{lemm_E_1}
Suppose that the hypotheses of \Cref{lem_3,lem_5_2} hold uniformly when their
compact parameter set \(\Theta\) is replaced by \(\Theta^{+\bar\rho}\); use
their conclusions with the corresponding constants constructed on this
thickening.  Let \(\dimk\ge1\) and choose
\[
  \begin{gathered}
  0<\epsilon_v\le
  \min\left\{2^{-1/\dimk},
  \varepsilon_{\mathrm{Ch}}^{1/(\dimk+1)}\right\},
  \qquad \delta:=\epsilon_v^{\dimk},\\
  \lambda:=\epsilon_v^{2(\dimk+1)}
  \le\min\{\lambda_{\mathrm{mom}},\lambda_{\mathrm{PI}}^\star\}.
  \end{gathered}
\]
Thus the R\'enyi target
\(\varepsilon_{\mathrm R}:=\delta\epsilon_v
=\epsilon_v^{\dimk+1}\) lies in the range
\eqref{eq_uniform_lmc_accuracy_threshold}.  For every
\(\vartheta\in\Theta^{+\bar\rho}\), let
\(\mathrm B=\mathrm B(\delta,\lambda)\) be given by \Cref{lem_3}, take
\(\beta^0\in[-\mathrm B,\mathrm B]\), and set
\(\sigma:=1+\delta^{-1}\).  Let \(I\) be a nonempty finite index set and let
\(L\) be a positive integer.  Choose the common LMC policy
\[
 (\bar h,n)
 =\bigl(
 \bar h_\star(\lambda,\delta\epsilon_v),
 n_\star(\lambda,\delta\epsilon_v)
 \bigr)
\]
from \eqref{eq_uniform_lmc_joint_policy}.  Initialize every chain from
\(\mathcal N(0,c\lambda I_d)\), using the fixed
\(0<c<3/(2L_{g,2})\) in
\Cref{lem_gaussian_lmc_initialization}.  When \(\mathrm B>0\), run
\Cref{alg:PSGD_LMC} with stepsize
\(2\mathrm B/(\sigma\sqrt L)\), mutually independent training chains, and a
fresh validation family independent of all training chains.  If
\(\mathrm B=0\), set \(\widetilde\beta=0\); no training chain is needed.

There are constants \(\widetilde C_f\ge0\),
\(C_{\mathrm{PSGD}}>0\), and
\(C_{\mathrm{orc}}^{\mathrm{LMC}}>0\), depending only on the fixed problem
data on the thickening and independent of
\(\vartheta,\epsilon_v,\delta,\lambda,L,I,n\), and \(\bar h\), such that
\[
 \mathbb E\!\left[
 \left|\widehat\psi(\vartheta)
       -\widetilde F_{\mathrm{SQG}}(\vartheta)\right|^2
 \,\middle|\,\vartheta\right]
 \le C_{\mathrm{orc}}^{\mathrm{LMC}}\epsilon_v^2
\]
whenever
\[
 |I|\ge4\widetilde C_f\epsilon_v^{-2\dimk-2},
 \qquad
 L\ge\max\left\{1,C_{\mathrm{PSGD}}\dimk^2\mathrm B^2
       \epsilon_v^{-2-2\dimk}
       \left[1+\log^2(\epsilon_v^{-1})\right]\right\}.
\]
Moreover, \eqref{eq_uniform_lmc_mixing_bound} gives the genuine upper bound
\[
 n=n_\star(\lambda,\delta\epsilon_v)
 \le\mathfrak C_{\mathrm{mix}}d\,
       \epsilon_v^{-10-10\dimk}.
\]
If \(L\) is chosen as the maximum of \(1\) and the ceiling of the sufficient
expression displayed above, then
\[
 L n
 \le C\,\mathfrak C_{\mathrm{mix}}d
 \frac{1+\dimk^2\mathrm B^2
       [1+\log^2(\epsilon_v^{-1})]}
      {\epsilon_v^{12+12\dimk}},
\]
where \(C<\infty\) is independent of \(\vartheta\) and \(\epsilon_v\).
Here \(Ln\) bounds the LMC work in the PSGD training phase.  A complete
value-oracle evaluation, including validation, uses at most
\(n(L+|I|)\) LMC steps.
The fixed constant \(\mathfrak C_{\mathrm{mix}}\) may depend on \(d\).
\end{lem}
\begin{proof}
Fix \(\theta\in\Theta^{+\bar\rho}\).  The proof has three parts: we first
control PSGD for the LMC endpoint law, then transfer that bound to the exact
Gibbs objective, and finally control the fresh validation batch.

\noindent\textit{Optimization for the endpoint law.}\par
When \(\mathrm B>0\), let
\(\widetilde\beta\) be the output of \Cref{alg:PSGD_LMC} at \(\theta\), and let
\(\mathcal F_l\) be the sigma-field generated by
\(\beta^0,\widetilde X_0^n,\ldots,\widetilde X_{l-1}^n\).
Then \(\beta^l\) is \(\mathcal F_l\)-measurable and, because every training
chain is freshly and independently initialized,
\[
 \operatorname{Law}(\widetilde X_l^n\mid\mathcal F_l)
 =\widehat\mu_{n,\theta}.
\]
Consequently, the atom-robust hinge selection from the proof of
\Cref{lem_5_2} has conditional expectation in
\(\partial_\beta\widetilde\phi_{\lambda,\delta}
(\theta,\beta^l)\) and absolute value at most
\(\sigma=1+\delta^{-1}\).  Repeating the projection--martingale argument used
to prove \eqref{eq_whp_gaurante_PSGD} gives the corresponding estimate for
\(\widetilde\phi_{\lambda,\delta}\) and its restricted minimum.
When \(\mathrm B=0\), we use \(\widetilde\beta=0\), and the corresponding
optimization gap is identically zero.

Set
\(\Delta:=\min\{1/2,
\epsilon_v^2/(32\sigma^2\mathrm B^2+\epsilon_v^2)\}\).  Then
\(0<\Delta\le1/2\),
\(4\sigma^2\mathrm B^2\Delta\le\epsilon_v^2/8\).  For
\(\mathrm B>0\), use the corresponding estimate with target
\(\epsilon_v/(2\sqrt2)\).  Its exact sufficient iteration condition is
\begin{equation}\label{eq:LMC_L_rate}
 L\ge 512\,\mathrm B^2\sigma^2
       \epsilon_v^{-2}\log^2(2/\Delta),
\end{equation}
and under this condition the bound gives
\[
 \mathbb P\!\left(
 \widetilde\phi_{\lambda,\delta}(\theta,\widetilde\beta)
 -\widetilde m_{\lambda,\delta}^{\,\mathrm B}(\theta)
 >\frac{\epsilon_v}{2\sqrt2}
 \,\middle|\,\theta\right)\le\Delta.
\]
Both \(\widetilde\beta\) and a minimizer defining
\(\widetilde m_{\lambda,\delta}^{\,\mathrm B}\) lie in
\([-\mathrm B,\mathrm B]\); hence their objective gap is bounded by
\(2\sigma\mathrm B\).
Splitting according to the preceding event yields
\begin{equation}\label{eq:LMC_phi_tilde_quad}
\begin{aligned}
 \mathbb E\!\left[
 \left(
 \widetilde\phi_{\lambda,\delta}(\theta,\widetilde\beta)
 -\widetilde m_{\lambda,\delta}^{\,\mathrm B}(\theta)
 \right)^2
 \,\middle|\,\theta\right]
 &\le\frac{\epsilon_v^2}{8}
    +4\sigma^2\mathrm B^2\Delta\\
 &\le\frac{\epsilon_v^2}{4}.
\end{aligned}
\end{equation}
For \(\delta=\epsilon_v^{\dimk}\), the estimates used
at the end of the proof of \Cref{lem_5_2} give
\(\sigma\le2\epsilon_v^{-\dimk}\) and
\(\log(2/\Delta)\le
C\dimk[1+\log(\epsilon_v^{-1})]\).  Hence there is a sufficiently large
fixed \(C_{\mathrm{PSGD}}>0\) such that the simpler condition
\[
 L\ge C_{\mathrm{PSGD}}\dimk^2\mathrm B^2
 \epsilon_v^{-2-2\dimk}
 \left[1+\log^2(\epsilon_v^{-1})\right]
\]
implies \eqref{eq:LMC_L_rate}.  If \(\mathrm B=0\),
\eqref{eq:LMC_phi_tilde_quad} holds directly because the optimization gap is
zero.

\noindent\textit{Transfer to the exact Gibbs objective.}\par
Recall that
\(\widetilde F_{\mathrm{SQG}}(\theta)=m_{\lambda,\delta}(\theta)\).  The exact
restricted-minimum decomposition is
\[
\begin{aligned}
 \phi_{\lambda,\delta}(\theta,\widetilde\beta)
 -m_{\lambda,\delta}(\theta)
 &\le
 \left|\phi_{\lambda,\delta}(\theta,\widetilde\beta)
 -\widetilde\phi_{\lambda,\delta}(\theta,\widetilde\beta)\right|\\
 &\quad+
 \left[
 \widetilde\phi_{\lambda,\delta}(\theta,\widetilde\beta)
 -\widetilde m_{\lambda,\delta}^{\,\mathrm B}(\theta)
 \right]\\
 &\quad+
 \left|\widetilde m_{\lambda,\delta}^{\,\mathrm B}(\theta)
 -m_{\lambda,\delta}(\theta)\right|.
\end{aligned}
\]
Define
\[
 C_{\mathrm{pert}}
 :=L_{f,1}C_W(C_{\mathrm{PI}}^\star)^{1/2}.
\]
Because \(\varepsilon_{\mathrm R}=\delta\epsilon_v\),
\eqref{eq_083083} and the restricted-minimum estimate following it give,
respectively,
\[
 \left|\phi_{\lambda,\delta}(\theta,\widetilde\beta)
       -\widetilde\phi_{\lambda,\delta}
        (\theta,\widetilde\beta)\right|
 \le C_{\mathrm{pert}}\epsilon_v,
 \qquad
 \left|\widetilde m_{\lambda,\delta}^{\,\mathrm B}(\theta)
       -m_{\lambda,\delta}(\theta)\right|
 \le C_{\mathrm{pert}}\epsilon_v.
\]
Combining these deterministic bounds with
\eqref{eq:LMC_phi_tilde_quad} and
\((a+b+c)^2\le3(a^2+b^2+c^2)\) gives
\begin{equation}\label{eq:LMC_phi_true_quad}
 \mathbb E\!\left[
 \left(
 \phi_{\lambda,\delta}(\theta,\widetilde\beta)
 -\widetilde F_{\mathrm{SQG}}(\theta)
 \right)^2
 \,\middle|\,\theta\right]
 \le C_2\epsilon_v^2,
 \qquad C_2:=6C_{\mathrm{pert}}^2+\frac34,
\end{equation}
Here
\((\bar h,n)=(\bar h_\star,n_\star)\) is the pair from
\eqref{eq_uniform_lmc_joint_policy}, and
\[
 n=n_\star(\lambda,\delta\epsilon_v)
 \le\mathfrak C_{\mathrm{mix}}d\,
      \epsilon_v^{-10-10\dimk}
\]
by \eqref{eq_uniform_lmc_mixing_bound} and
\eqref{eq_uniform_lmc_warm_factor}.

\noindent\textit{Validation error.}\par
The finite-second-moment verification in
\Cref{section_gradient_oracle_complexity} ensures that the following
optimal couplings exist.
For each \(i\in I\), choose a \(W_2\)-optimal coupling of the fresh validation
endpoint \(\widehat X_i^n\sim\widehat\mu_{n,\theta}\) with an auxiliary
\(X_i\sim\mu_\theta^\lambda\), taking the coupled pairs independently across
\(i\) and independently of all PSGD--LMC training chains.  Thus they are
independent of \(\widetilde\beta\).  Define the auxiliary ideal-validation estimator
\[
 \psi_I^{\mathrm{id}}(\theta)
 :=\widetilde\beta+\frac1{\delta|I|}
   \sum_{i\in I}(f(\theta,X_i)-\widetilde\beta)_+.
\]
Using \((a+b+c)^2\le3(a^2+b^2+c^2)\), we have
\begin{align}
\MoveEqLeft[4]\mathbb E\!\left[
\left|\widehat\psi(\theta)-\widetilde F_{\mathrm{SQG}}(\theta)\right|^2
\,\middle|\,\theta\right]\nonumber\\
&\le
3\,\underbrace{\mathbb E\!\left[
\left|\widehat\psi(\theta)-\psi_I^{\mathrm{id}}(\theta)\right|^2
\,\middle|\,\theta\right]}_{(1)}
+3\,\underbrace{\mathbb E\!\left[
\left|\psi_I^{\mathrm{id}}(\theta)
-\phi_{\lambda,\delta}(\theta,\widetilde\beta)\right|^2
\,\middle|\,\theta\right]}_{(2)}
\nonumber\\
&\quad+
3\,\underbrace{\mathbb E\!\left[
\left|\phi_{\lambda,\delta}(\theta,\widetilde\beta)
-\widetilde F_{\mathrm{SQG}}(\theta)\right|^2
\,\middle|\,\theta\right]}_{(3)}.
\label{eq:LMC_three_terms}
\end{align}

For the first term, the hinge contraction and
\(\|f(\theta,\cdot)\|_{\mathrm{Lip}}\le L_{f,1}\) give, with
\[
D_i:=(f(\theta,\widehat X_i^n)-\widetilde\beta)_+
     -(f(\theta,X_i)-\widetilde\beta)_+,
\]
\begin{align}\label{eq_008787}
\mathbb E\!\left[
\left|\widehat\psi(\theta)-\psi_I^{\mathrm{id}}(\theta)\right|^2
\,\middle|\,\theta\right]
&\le\frac1{\delta^2|I|}
\sum_{i\in I}
\mathbb E[|D_i|^2\mid\theta]\nonumber\\
&\le
\frac{L_{f,1}^2}{\delta^2}
\mathbb W_2^2(\widehat\mu_{n,\theta},\mu_\theta^\lambda).
\end{align}
The common LMC policy at
\(\varepsilon_{\mathrm R}=\delta\epsilon_v\) gives
\(\mathbb R_2(\widehat\mu_{n,\theta}\|\mu_\theta^\lambda)
\le\delta^2\epsilon_v^2\).  The Poincar\'e--transport bound in
\Cref{section_gradient_oracle_complexity} therefore yields
\[
 \mathbb W_2^2(\widehat\mu_{n,\theta},\mu_\theta^\lambda)
 \le C_W^2C_{\mathrm{PI}}^\star\delta^2\epsilon_v^2,
\]
uniformly in \(\theta\).  Here
\(n=n_\star(\lambda,\delta\epsilon_v)\) obeys the upper bound displayed in
the statement.  Substitution into \eqref{eq_008787} gives
\begin{equation}\label{eq:LMC_term1}
\mathbb E\!\left[
\left|\widehat\psi(\theta)-\psi_I^{\mathrm{id}}(\theta)\right|^2
\,\middle|\,\theta\right]
\le C_3\epsilon_v^2,
\qquad C_3:=L_{f,1}^2C_W^2C_{\mathrm{PI}}^\star.
\end{equation}

For the second term, conditional on \((\theta,\widetilde\beta)\), the
\(X_i\)'s remain i.i.d. with law \(\mu_\theta^\lambda\).  The calculation in
\eqref{eq_append_54}--\eqref{eq_append_25} therefore gives
\begin{equation}\label{eq:LMC_term2}
\mathbb E\!\left[
\left|\psi_I^{\mathrm{id}}(\theta)
-\phi_{\lambda,\delta}(\theta,\widetilde\beta)\right|^2
\,\middle|\,\theta\right]
\le\frac{\widetilde C_f}{\delta^2|I|}.
\end{equation}
With \(\delta=\epsilon_v^{\dimk}\) and
\(|I|\ge4\widetilde C_f\epsilon_v^{-2\dimk-2}\), this gives
\[
 \mathbb E\!\left[
 \left|\psi_I^{\mathrm{id}}(\theta)
 -\phi_{\lambda,\delta}(\theta,\widetilde\beta)\right|^2
 \,\middle|\,\theta\right]
 \le\frac{\epsilon_v^2}{4}.
\]

For the third term, \eqref{eq:LMC_phi_true_quad} gives
\begin{equation}\label{eq:LMC_term3}
\mathbb{E}\bigl[|\phi_{\lambda,\delta}(\theta,\tilde{\beta})-\tilde{F}_{\SQG}(\theta)|^2 \mid \theta\bigr]
\le C_{2}\epsilon_{v}^2.
\end{equation}

Substituting \eqref{eq:LMC_term1}--\eqref{eq:LMC_term3} into
\eqref{eq:LMC_three_terms} yields
\[
 \mathbb E\!\left[
 |\widehat\psi(\theta)-\widetilde F_{\mathrm{SQG}}(\theta)|^2
 \,\middle|\,\theta\right]
 \le C_{\mathrm{orc}}^{\mathrm{LMC}}\epsilon_v^2,
 \qquad
 C_{\mathrm{orc}}^{\mathrm{LMC}}
 :=3\left(C_3+\frac14+C_2\right).
\]
This constant is independent of \(\theta,\delta,\lambda\), and
\(\epsilon_v\) on the coupled range.

For the work bound, choose \(L\) exactly as in the statement, namely as the
maximum of one and the ceiling of its displayed sufficient expression.
The elementary inequality \(\lceil x\rceil\le1+x\), together with
\(0<\epsilon_v<1\), then gives
\[
 L\le C\left\{1+
 \dimk^2\mathrm B^2\epsilon_v^{-2-2\dimk}
 \left[1+\log^2(\epsilon_v^{-1})\right]\right\}.
\]
Multiplying this upper bound by
\(n=n_\star(\lambda,\delta\epsilon_v)
\le\mathfrak C_{\mathrm{mix}}d\epsilon_v^{-10-10\dimk}\) gives
\[
L\times n
\le C\mathfrak C_{\mathrm{mix}}
\frac{d\left\{1+\dimk^{2}\,\mathrm{B}^{2}
\left[1+\log^{2}(\epsilon_{v}^{-1})\right]\right\}}
{\epsilon_{v}^{12+12\dimk}}.
\]
Thus \(Ln\) is the training-chain work stated in the lemma; validation adds
\(n|I|\) LMC steps.
This completes the proof.
\end{proof}

\begingroup
\subsection{Proof of \Cref{theorem_gradient_complexity}}\label[appendix]{append_e_3}
Algorithm~\ref{alg:zroth-order-LMC} replaces every exact Gibbs sample in
\Cref{alg: Stochastic-Zeroth-order} by the endpoint of a query-specific LMC
chain and replaces PSGD by PSGD--LMC.  The exact directional-estimator
analysis is therefore unchanged.  The proof below establishes the one new
ingredient---a uniform mean-square bound for the LMC-based value oracle---and
then applies the projected-descent argument of
\Cref{thm_cvg_alg_zeroth_goldstein}.

\endgroup
\begingroup
\makeatletter
\def\fs@ruled{\def\@fs@cfont{\bfseries}\let\@fs@capt\floatc@ruled
  \def\@fs@pre{\begingroup\normalcolor\hrule height.8pt depth0pt\endgroup\kern2pt}\def\@fs@post{\kern2pt\begingroup\normalcolor\hrule\endgroup\relax}\def\@fs@mid{\kern2pt\begingroup\normalcolor\hrule\endgroup\kern2pt}\let\@fs@iftopcapt\iftrue}
\makeatother
\begin{algorithm}[t]
\caption{Projected Stochastic Zeroth-order with LMC (PSZO-LMC)}\label{alg:zroth-order-LMC}
\textbf{Input:}
$\mathrm B=\mathrm B(\delta,\lambda)$ from \Cref{lem_3},
$\beta_{0}\in[-\mathrm B,\mathrm B]$, $\theta_{0}\in\Theta$,
$\lambda>0$, $0<\delta\le1/2$, an oracle accuracy $\epsilon_v>0$ satisfying
$\delta\epsilon_v\le\varepsilon_{\mathrm{Ch}}$, $0<\rho<\bar\rho$, a fixed
$c\in(0,3/(2L_{g,2}))$ with every LMC chain initialized from
$\mathcal N(0,c\lambda I_d)$, the outer stepsize $\eta>0$, the common pair
\[
(\bar h,n)=\bigl(
\bar h_\star(\lambda,\delta\epsilon_v),
n_\star(\lambda,\delta\epsilon_v)\bigr),
\]
and positive integers $N,b_{u},L,|I|$.
\begin{algorithmic}[1]
        \FOR{$k$ = $0$ to $N-1$}
\STATE Draw fresh directions
         \(u_{k,1},\ldots,u_{k,b_u}\) independently and uniformly from
         \(\mathbb S^{m-1}\), independently of all randomness generated
         before iteration \(k\).
         \STATE For each \(t\in[1:b_u]\), perform the following query-specific operations.
         \STATE $\vartheta_{k,t}^{\pm}:=\theta_k\pm\rho u_{k,t}$,
         $\quad\beta_{k,t}^{\pm}:=\text{\rm PSGD-LMC}(\vartheta_{k,t}^{\pm},\beta_0)$, using independent training-chain families.
\STATE At each current query
         \(\vartheta_{k,t}^{\pm}\), independently initialize fresh validation
         chains and use independent LMC noises to generate
         $\widehat X_{I,k,t}^{n,\pm}\sim
         (\widehat\mu_{n,\vartheta_{k,t}^{\pm}})^{\otimes |I|}$.
         Use distinct validation chains for every
         \(t\) and sign, sharing no initialization or LMC noise with earlier
         chains or with the corresponding PSGD--LMC training chains.
         \STATE $\widehat\psi(\vartheta_{k,t}^{\pm})
         :=\phi_{\lambda,\delta}
         (\vartheta_{k,t}^{\pm},\beta_{k,t}^{\pm};
          \widehat X_{I,k,t}^{n,\pm})$.
        \STATE $\hat{\mathsf{g}}_{\rho}(\theta_{k})=\frac{m}{2\rho}\cdot\frac{1}{b_{u}}\sum_{t=1}^{b_{u}}\left[\hat\psi(\theta_{k}+\rho u_{k,t})-\hat\psi(\theta_{k}-\rho u_{k,t})\right]u_{k,t}$
        \STATE $\theta_{k+1}=\operatorname{Proj}_{\Theta}
        (\theta_{k}-\eta\hat{\mathsf{g}}_{\rho}(\theta_{k}))$
        \ENDFOR
\STATE Draw
        $J\sim\operatorname{Unif}\{0,\ldots,N-1\}$ independently of the
        entire algorithmic history.
        \STATE \textbf{Return}
        $\widehat\theta:=\theta_J$.
    \end{algorithmic}
\end{algorithm}
\endgroup

\begingroup
\begin{proof}
The case \(L_F=0\) was settled in the theorem statement, so assume
\(L_F>0\).  We first recall the projected zeroth-order estimate, then bound
the LMC-oracle error at every actual query, and finally telescope and count
the oracle calls.

\textit{Step 1: Projected zeroth-order estimate.}
\begingroup
As in the proof of \Cref{thm_cvg_alg_zeroth_goldstein}, use
\[
\psi_\rho(\theta)
:=\mathbb E_{V\sim\operatorname{Unif}(\mathbb B_m(0;1))}
  [F_{\max}(\theta+\rho V)].
\]
The sphere two-point estimator is unbiased for the gradient of this
ball-smoothed function.  Because $0<\rho<\bar\rho$, every smoothing ball is
contained in $\operatorname{int}(\Theta^{+\bar\rho})$.  When $L_F>0$, use
the same constant outer step
\[
\eta:=\frac{\rho}{\sqrt m\,L_F}=L_{\psi,2}^{-1}.
\]
The smoothness and projected gradient-mapping analysis in
\Cref{append_proof_thm_cvg_alg_zeroth_goldstein} yield
\endgroup
\begingroup
(see \cref{eq_append_0024} in \Cref{append_proof_thm_cvg_alg_zeroth_goldstein})
\begin{align}
\label{eq:E2_grad_map_base}
\mathbb{E}\big[\|\mathcal{G}_{\Theta}(\theta_i,\nabla\psi_{\rho}(\theta_i);\eta)\|^2
\mid \theta_i\big]
&\le
\frac{8\big(\psi_{\rho}(\theta_i)-\mathbb{E}[\psi_{\rho}(\theta_{i+1})\mid\theta_i]\big)}
     {\eta}\nonumber\\
&\qquad\qquad+
24\,\mathbb{E}\big[\|\nabla\psi_{\rho}(\theta_i)-\mathsf{g}_{\rho}(\theta_i)\|^2
\mid \theta_i\big]\nonumber\\
&\qquad\qquad+
40\,\mathbb{E}\big[\|\mathsf{g}_{\rho}(\theta_i)-\hat{\mathsf{g}}_{\rho}(\theta_i)\|^2
\mid \theta_i\big],
\end{align}
\endgroup
where
\[
\mathsf{g}_{\rho}(\theta_i)
:= \frac{1}{b_{u}}\sum_{t=1}^{b_{u}} \mathsf{g}_{\rho}(\theta_i,u_{i,t}),
\qquad
\mathsf{g}_{\rho}(\theta,u)
:= \frac{m}{2\rho}\big[F_{\max}(\theta+\rho u)-F_{\max}(\theta-\rho u)\big]u,
\]
and $\hat{\mathsf{g}}_{\rho}(\theta_i)$ is the stochastic gradient estimator
used in Algorithm~\ref{alg:zroth-order-LMC}.
Up to \eqref{eq:E2_grad_map_base}, the analysis is identical to
Theorem~\ref{thm_cvg_alg_zeroth_goldstein}.  Indeed, conditionally on
\(\theta_i\), the fresh directions are i.i.d., each directional estimator has
mean \(\nabla\psi_\rho(\theta_i)\), and its norm is at most \(mL_F\).
Therefore
\[
\mathbb{E}\big[\|\nabla\psi_{\rho}(\theta_i)-\mathsf{g}_{\rho}(\theta_i)\|^2\mid\theta_i\big]
\le \frac{m^{2}L_F^{2}}{b_u}.
\]
The only new task is to control
$\mathbb{E}[\|\hat{\mathsf{g}}_{\rho}(\theta)-\mathsf{g}_{\rho}(\theta)\|^2\mid\theta]$
when $\hat{\mathsf{g}}_{\rho}$ is built from the LMC-based oracle $\hat{\psi}$.

\textit{Step 2: Querywise LMC-oracle error.}
Fix \(i\in\{0,\ldots,N-1\}\) and abbreviate
\(\theta:=\theta_i\) and \(u_t:=u_{i,t}\) for \(t=1,\ldots,b_u\).
{Because \(0<\rho<\bar\rho\), the ball
\(\mathbb B_m(\theta;\rho)\) lies in
\(\operatorname{int}(\Theta^{+\bar\rho})\), and every
\(\theta\pm\rho u_t\) belongs to that thickening.
The assumptions and all constants used in
\Cref{lemm_E_1} are uniform on that thickening, and each oracle evaluation
uses training and validation chains targeting the Gibbs law at that same
parameter.}
The proxy gradient used in the projected update is
\begin{equation}
\label{eq:E2_grad_def}
\hat{\mathsf{g}}_{\rho}(\theta)
=
\frac{m}{2\rho}\cdot \frac{1}{b_{u}}
\sum_{t=1}^{b_{u}}
\big[\hat{\psi}(\theta+\rho u_{t})-\hat{\psi}(\theta-\rho u_{t})\big]u_{t}.
\end{equation}
Here PSGD--LMC returns the threshold \(\widetilde\beta\), and
\(\widehat\psi\) is formed from this threshold and the independent validation
family, as specified in Algorithm~\ref{alg:zroth-order-LMC}.
We need to bound
\(\mathbb E[\|\hat{\mathsf g}_\rho(\theta)
-\mathsf g_\rho(\theta)\|^2\mid\theta]\) in
\eqref{eq:E2_grad_map_base}.
\begingroup
For \(t=1,\ldots,b_u\), let
\[
 e_t^\pm
 :=\hat\psi(\theta\pm\rho u_t)
   -F_{\max}(\theta\pm\rho u_t).
\]
Condition on \(\theta\) and on the realized directions
\(u_1,\ldots,u_{b_u}\).  The \(2b_u\) query parameters are then fixed, and
the corresponding training and validation families are mutually independent
and target their query-specific Gibbs laws.  Pointwise,
\[
\begin{aligned}
 \|\hat{\mathsf g}_\rho(\theta)-\mathsf g_\rho(\theta)\|^2
 &=\left\|
 \frac{m}{2\rho b_u}\sum_{t=1}^{b_u}
 (e_t^+-e_t^-)u_t
 \right\|^2\\
 &\le
 \frac{m^2}{2\rho^2b_u}
 \sum_{t=1}^{b_u}\bigl[(e_t^+)^2+(e_t^-)^2\bigr].
\end{aligned}
\]
Therefore
\begin{align}
\MoveEqLeft[2]
\mathbb E\!\left[
\|\hat{\mathsf g}_\rho(\theta)-\mathsf g_\rho(\theta)\|^2
\middle|\theta,u_1,\ldots,u_{b_u}\right]\nonumber\\
&\le
\frac{m^2}{2\rho^2b_u}
\sum_{t=1}^{b_u}
\left\{
\mathbb E[(e_t^+)^2\mid\theta,u_1,\ldots,u_{b_u}]
+\mathbb E[(e_t^-)^2\mid\theta,u_1,\ldots,u_{b_u}]
\right\}.
\label{eq:E2_g_diff}
\end{align}
\endgroup

{Because \(\epsilon_v=\rho\epsilon/m\), the choices of \(|I|\)
and \(L\) in \Cref{theorem_gradient_complexity} are exactly the sufficient
orders required by \Cref{lemm_E_1}; the theorem also uses that lemma's common
joint pair \((\bar h,n)\).  Thus \Cref{lemm_E_1} gives the following bound
with \(C_E:=C_{\mathrm{orc}}^{\mathrm{LMC}}\), which is independent of the
query parameter and the accuracy variables: for every
\(\vartheta\in\Theta^{+\bar\rho}\),}
\begin{equation}
\label{eq:E2_psi_hat_vs_Fsqg}
{{\mathbb{E}\big[
|\hat{\psi}(\vartheta)-\tilde{F}_{\SQG}(\vartheta)|^{2}
\mid\vartheta\big]}
\le C_E\epsilon_{v}^{2},}
\end{equation}
{
{Moreover, with \(\epsilon_v=\rho\epsilon/m\) and the admissible
parameter range in \Cref{thm_cvg_alg_zeroth_goldstein}, the uniform
approximation result \Cref{rmk_func_approx_cvar_gibbs} gives}}
\begin{equation}
\label{eq:E2_Fsqg_vs_Fmax}
{{
\big|\tilde{F}_{\SQG}(\vartheta)-F_{\max}(\vartheta)\big|
\le C_{\mathrm{app}}\epsilon_{v}}
\qquad {\text{for all } \vartheta\in\Theta^{+\bar\rho}.}}
\end{equation}
\begingroup
Combining \eqref{eq:E2_psi_hat_vs_Fsqg} and \eqref{eq:E2_Fsqg_vs_Fmax} and
using $(a+b)^{2}\le 2(a^{2}+b^{2})$, we obtain
\begin{align}
\MoveEqLeft[2]
\mathbb{E}\!\left[
|\hat{\psi}(\theta)-F_{\max}(\theta)|^{2}\,\middle|\,\theta\right]
\nonumber\\
&\le
2\,\mathbb{E}\!\left[
|\hat{\psi}(\theta)-\tilde{F}_{\SQG}(\theta)|^{2}
\,\middle|\,\theta\right]\nonumber\\
&\quad+2\,|\tilde{F}_{\SQG}(\theta)-F_{\max}(\theta)|^{2}
\nonumber\\
&\le 2\bigl(C_E+C_{\mathrm{app}}^2\bigr)\epsilon_v^2
 =:C_{\mathrm{orc}}\epsilon_v^2.
\label{eq:E2_psi_hat_vs_Fmax}
\end{align}
\endgroup
\begingroup
Apply \eqref{eq:E2_psi_hat_vs_Fmax} conditionally at each of the fixed
query parameters \(\theta\pm\rho u_t\) in \eqref{eq:E2_g_diff}.
The resulting conditional bound is independent of the directions.  Taking
expectation over the directions therefore yields
\begin{equation}
\label{eq:E2_g_diff_final}
\mathbb{E}\!\left[
\|\hat{\mathsf{g}}_{\rho}(\theta)-\mathsf{g}_{\rho}(\theta)\|^{2}
\middle|\theta\right]
\le
C_{\mathrm{orc}}\frac{m^{2}}{\rho^{2}}\epsilon_{v}^{2}.
\end{equation}
\endgroup

\textit{Step 3: Convergence rate and choice of parameters.}
\begingroup
Plugging \eqref{eq:E2_g_diff_final} into
\eqref{eq:E2_grad_map_base}, using
$\eta=\rho/(\sqrt m\,L_F)$, and summing
$i=0,\dots,N-1$ as in the proof of
\Cref{thm_cvg_alg_zeroth_goldstein} gives the explicit estimate
\begin{align}
\MoveEqLeft[2]
\frac{1}{N}\sum_{i=0}^{N-1}
\mathbb{E}\!\left[
\|\mathcal{G}_{\Theta}(\theta_i,\nabla\psi_{\rho}(\theta_i);\eta)\|^{2}
\right]
&\le
\frac{8\sqrt m\,L_F^2\operatorname{diam}(\Theta)}{N\rho}
+\frac{24m^2L_F^2}{b_u}\nonumber\\
&\quad
+40C_{\mathrm{orc}}\frac{m^2}{\rho^2}\epsilon_v^2.
\label{eq:LMC_outer_average}
\end{align}
Indeed, the first term follows by telescoping and using
\[
 \psi_\rho(\theta_0)-\mathbb E\psi_\rho(\theta_N)
 \le L_F\operatorname{diam}(\Theta).
\]

With the choices in \Cref{theorem_gradient_complexity} and
\(\epsilon_v=\rho\epsilon/m\), the three terms on the right-hand side of
\eqref{eq:LMC_outer_average} are bounded, respectively, by
\[
 \frac{8L_F^2\operatorname{diam}(\Theta)\bar\rho}{C_N}\epsilon^2,
 \qquad
 \frac{24L_F^2}{C_u}\epsilon^2,
 \qquad
 40C_{\mathrm{orc}}\epsilon^2.
\]
Thus the average in \eqref{eq:LMC_outer_average} is at most
\(C_\star\epsilon^2\), where
\[
 C_\star
 :=\frac{8L_F^2\operatorname{diam}(\Theta)\bar\rho}{C_N}
   +\frac{24L_F^2}{C_u}+40C_{\mathrm{orc}}.
\]
\begingroup
Let \(J\sim\operatorname{Unif}\{0,\ldots,N-1\}\) be the independent index
drawn by \Cref{alg:zroth-order-LMC}, and let
\(\widehat\theta:=\theta_J\).  Because \(J\) is drawn after all outer
iterations and independently of their history, conditioning on that history
gives the exact identity
\[
\mathbb E\!\left[
 \|\mathcal G_\Theta(
 \theta_J,\nabla\psi_\rho(\theta_J);\eta)\|^2\right]
=\frac1N\sum_{i=0}^{N-1}
 \mathbb E\!\left[
 \|\mathcal G_\Theta(
 \theta_i,\nabla\psi_\rho(\theta_i);\eta)\|^2\right].
\]
Jensen's inequality and the pointwise inclusion
\(\nabla\psi_\rho(\theta_J)\in
\partial_\rho F_{\max}(\theta_J)\) give
\[
\begin{aligned}
\mathbb E\!\left[
 \min_{g\in\partial_\rho F_{\max}(\widehat\theta)}
 \|\mathcal G_\Theta(\widehat\theta,g;\eta)\|\right]
&\le
\mathbb E\!\left[
 \|\mathcal G_\Theta(
 \theta_J,\nabla\psi_\rho(\theta_J);\eta)\|\right]\\
&\le \sqrt{C_\star}\,\epsilon.
\end{aligned}
\]
Thus one may take
\(C_{\mathrm{res}}^{\mathrm{LMC}}:=\sqrt{C_\star}\), so the residual is
controlled at the iterate actually returned by the
algorithm; drawing \(J\) requires no oracle call.
\endgroup
\endgroup

\textit{Step 4: Oracle counts.}
{Each outer iteration of Algorithm~\ref{alg:zroth-order-LMC} uses
\(2b_u\) calls to \(\widehat\psi\).  In the present \(L_F>0\) branch,
\(L_{f,1}>0\), and the definition in \Cref{lem_3} gives
\(\mathrm B\ge L_{f,1}K_0\sqrt\lambda/\delta>0\).  Thus each call uses
exactly \(L\) training chains and \(|I|\) fresh validation chains, all of
length \(n=n_\star\).  The exact counts are therefore
\[
 Q_f=2Nb_u(L+|I|),
 \qquad
 Q_g=n_\star Q_f.
\]
When \(L_F=0\), both counts are zero, as stated in the theorem.}
\begingroup
\begingroup
Choose \(N,b_u,|I|\), and \(L\) at the ceiling values stated in
\Cref{theorem_gradient_complexity}, and choose
\[
 (\bar h,n)=
 \bigl(\bar h_\star(\lambda,\delta\epsilon_v),
       n_\star(\lambda,\delta\epsilon_v)\bigr).
\]
Since \(\lceil a\rceil\le1+a\) for \(a\ge0\), the choices in the theorem
satisfy
\[
\begin{aligned}
N&\le1+C_Nm\epsilon^{-2}\rho^{-2},\\
b_u&\le1+C_um^2\epsilon^{-2},\\
|I|&\le1+4\widetilde C_fm^{2\dimk+2}
       (\rho\epsilon)^{-2\dimk-2},\\
L&\le1+C_{\mathrm{PSGD}}\dimk^2m^{2\dimk+2}\mathrm B^2
       (\rho\epsilon)^{-2\dimk-2}
       \left[1+\log^2\!\left(\frac{m}{\rho\epsilon}\right)\right].
\end{aligned}
\]
Consequently, throughout the complete admissible range,
\[
\begin{aligned}
Q_f
&\le2\left(1+C_Nm\epsilon^{-2}\rho^{-2}\right)
       \left(1+C_um^2\epsilon^{-2}\right)\\
&\quad\times\Bigg[2
 +4\widetilde C_fm^{2\dimk+2}(\rho\epsilon)^{-2\dimk-2}\\
&\hspace{21mm}
 +C_{\mathrm{PSGD}}\dimk^2m^{2\dimk+2}\mathrm B^2
  (\rho\epsilon)^{-2\dimk-2}
  \left[1+\log^2\!\left(\frac{m}{\rho\epsilon}\right)\right]
 \Bigg],\\
Q_g&\le n_\star Q_f.
\end{aligned}
\]
These bounds hold on the complete admissible range and are upper bounds for
the stated ceiling choices, not for arbitrary larger integers.  Under the
three explicit size conditions in \Cref{theorem_gradient_complexity}, the
ceiling factors are bounded by universal multiples of their unrounded
expressions.  This gives the displayed bound for \(Q_f\).

Furthermore, \eqref{eq_uniform_lmc_mixing_bound},
\eqref{eq_uniform_lmc_warm_factor}, and
\(\lambda=(\delta\epsilon_v)^2\) give
\[
 n_\star
 \le\mathfrak C_{\mathrm{mix}}d\,
       \epsilon_v^{-10(\dimk+1)}
 =\mathfrak C_{\mathrm{mix}}d\,
       m^{10(\dimk+1)}(\rho\epsilon)^{-10(\dimk+1)}.
\]
Multiplying the \(Q_f\) bound by this estimate gives exactly the displayed
bound for \(Q_g\).  Hence the asserted powers of
\(m,\rho,\epsilon\), and \(\dimk\) follow without changing any exponent.
The factor
\(\mathfrak C_{\mathrm{mix}}\) is independent of the parameter and accuracy
variables, but it may itself depend on \(d\); this calculation therefore
does not assert a complete polynomial ambient-dimension bound.
\endgroup
\endgroup
\end{proof}

\endgroup

\section{Compute-matched Comparisons under a Fixed Wall-clock Budget}
\label{app:fixed_runtime}
This appendix complements the iteration-based results in \Cref{sec:experiment} by comparing methods under a fixed wall-clock time budget.
As emphasized in \Cref{sec:experiment}, PSZO-MinSel invokes a Gibbs--superquantile oracle implemented via PSGD and parallel-chain Langevin Monte Carlo, so the cost per outer iteration can differ substantially from gradient-based baselines.

\textit{Timing protocol and selection.}
Toy-example methods are implemented in JAX, while data hyper-cleaning methods are implemented in PyTorch; all experiments are executed on the same machine described in \Cref{sec:experiment}.
To obtain reliable timings, we exclude one-time compilation by running a short warm-up pass before timing and synchronize the device when recording timestamps.
Runtimes are measured as wall-clock seconds elapsed from the start of the timed run; for JAX experiments we time a JIT-compiled step function and record timestamps with device synchronization to avoid asynchronous timing artifacts.
As in \Cref{sec:experiment}, for each method we select the best-so-far iterate within a budget $T$,
\(
\hat{\theta}(T) \in \arg\min\{\mathcal{E}_n:\ t_n\le T\},
\)
where $\mathcal{E}_n$ is $F(\theta_n)$ for the toy example and the validation loss in \eqref{eq:hypercleaning} for data hyper-cleaning.
All tables below report performance evaluated at $\hat{\theta}(T)$.

\subsection{Toy example: fixed-runtime results}
\label{app:fixed_runtime_toy}

\textit{Per-iteration timing.}
\Cref{tab:time_per_iter_toy} reports the mean wall-clock time per outer iteration (post warm-up) for the toy example.
This highlights the substantial per-iteration cost gap between PSZO-MinSel and gradient-based baselines.

\begin{table}[t]
\centering
\begin{tabular}{c|ccccc}
\hline
$d$ & V-PBGD & G-PBGD & RHG & IAPTT-GM & \textbf{PSZO-MinSel}\\
\hline
2   & 0.022 & 0.018 & 0.026 & 0.027 & 6.74 \\
5   & 0.021 & 0.018 & 0.027 & 0.032 & 7.90 \\
10  & 0.022 & 0.018 & 0.036 & 0.045 & 8.10 \\
20  & 0.021 & 0.017 & 0.045 & 0.046 & 9.70 \\
\hline
\end{tabular}
\caption{Mean wall-clock time per outer iteration (ms) on the toy example with $\dimk=d-1$.}
\label{tab:time_per_iter_toy}
\end{table}

\subsubsection{Optimistic formulation: dimension sweep under fixed time}
\label{app:fixed_runtime_dim_sweep}

\Cref{tab:dim_sweep_time} replicates \Cref{tab:dim_sweep} under a fixed wall-clock budget $T_{\mathrm{toy}}=10$ seconds and reports the best-so-far absolute error $|\hat{\theta}(T_{\mathrm{toy}})-\theta^\star|$.

\begin{table}[t]
\centering
\begin{tabular}{c|ccccc}
\hline
$d$ & V-PBGD & G-PBGD & RHG & IAPTT-GM & \textbf{PSZO-MinSel} \\ \hline
2   & 0.510 & 0.264 & 0.313 & 0.094 & 0.016$\pm$0.014 \\
5   & 0.203 & 0.213 & 0.235 & 0.113 & 0.034$\pm$0.017 \\
10  & 0.245 & 0.278 & 0.260 & 0.180 & 0.053$\pm$0.013 \\
20  & 0.333 & 0.201 & 0.322 & 0.190 & 0.083$\pm$0.006 \\
\hline
\end{tabular}
\caption{Compute-matched analogue of \Cref{tab:dim_sweep}: best-so-far absolute error $|\hat{\theta}(T_{\mathrm{toy}})-\theta^\star|$ under a fixed wall-clock budget $T_{\mathrm{toy}}=10$ s.}
\label{tab:dim_sweep_time}
\end{table}

\noindent
Under equal time, PSZO-MinSel retains a clear advantage across dimensions, consistent with the iteration-based conclusions in \Cref{sec:experiment}.

\subsubsection{Optimistic formulation: constant manifold dimension under fixed time}

\Cref{tab:dim_fixed_time} reports the original fixed-\(\dimk\) implementation under the same budget
$T_{\mathrm{toy}}=10$ seconds. As noted in \Cref{sec:experiment}, that implementation omitted the stabilizing term in
the inactive coordinates and must be rerun for the corrected objective.

\begin{table}[t]
\centering
\begin{tabular}{c|ccccc}
\hline
$d$ & V-PBGD & G-PBGD & RHG & IAPTT-GM & \textbf{PSZO-MinSel} \\ \hline
5   & 0.105 & 0.235 & 0.202 & 0.137 & 0.029$\pm$0.038 \\
10  & 0.105 & 0.235 & 0.237 & 0.170 & 0.0208$\pm$0.022 \\
20  & 0.104 & 0.238 & 0.237 & 0.170 & 0.021$\pm$0.014 \\
30  & 0.104 & 0.254 & 0.237 & 0.155 & 0.018$\pm$0.011 \\
\hline
\end{tabular}
\caption{Compute-matched analogue of \Cref{tab:dim_fixed}: best-so-far absolute error $|\hat{\theta}(T_{\mathrm{toy}})-\theta^\star|$ under $T_{\mathrm{toy}}=10$ s, with fixed $\dimk=1$.}
\label{tab:dim_fixed_time}
\end{table}

\noindent
The reported values are stable across \(d\), but they do not establish intrinsic-dimension dependence because they
were produced by the unstabilized objective described above.

\subsection{Data hyper-cleaning: fixed-runtime results}
\label{app:fixed_runtime_hyperclean}

\paragraph{Per-iteration timing.}
\Cref{tab:time_per_iter_hyper} reports the mean wall-clock time per outer iteration (post warm-up) on MNIST for representative models.

\begin{table}[t]
\centering
\begin{tabular}{c|ccccc}
\hline
Model & V-PBGD & G-PBGD & RHG & IAPTT-GM & \textbf{PSZO-MinSel} \\
\hline
Linear      & 0.015 & 0.016 & 15.6 & 0.068 & 13.5 \\
2-layer MLP & 0.022 & 0.026 & 18.4 & 0.087 & 17.1 \\
\hline
\end{tabular}
\caption{Mean wall-clock time per outer iteration (s) on MNIST data hyper-cleaning (post warm-up).}
\label{tab:time_per_iter_hyper}
\end{table}

\noindent
RHG and PSZO-MinSel have much higher per-outer costs in this task, making fixed-runtime evaluation particularly important.

\paragraph{Fixed-runtime performance.}
\Cref{tab:hypercleaning-results-time} replicates \Cref{tab:hypercleaning-results} under a fixed budget $T_{\mathrm{hc}}=1800$ seconds.
For each method, we select the iterate achieving the lowest validation loss within $T_{\mathrm{hc}}$ and report the corresponding test accuracy and F1-score.

\begin{table}[t]
\centering
\caption{Compute-matched analogue of \Cref{tab:hypercleaning-results} under a fixed wall-clock budget $T_{\mathrm{hc}}=1800$ seconds (pollute rate $p\in\{40,60,80\}\%$): test accuracy (\%) and F1 score (\%), evaluated at the best-so-far iterate within the time budget (mean $\pm$ 95\% CI over 10 runs).}
\label{tab:hypercleaning-results-time}
\tiny
\resizebox{\textwidth}{!}{
\begin{tabular}{cc|cc|cc}
\toprule
\multirow{2}{*}{$p$} & \multirow{2}{*}{Method}
& \multicolumn{2}{c|}{Linear model}
& \multicolumn{2}{c}{2-layer MLP} \\
\cmidrule(lr){3-4}\cmidrule(lr){5-6}
& & Test acc. & F1 & Test acc. & F1 \\
\midrule
\multirow{5}{*}{40}
& RHG        & $60.50 \pm 0.32$ & $45.5\pm 0.4$ & $60.54 \pm 0.29$ & $42.3\pm 0.4$ \\
& IAPTT-GM   & $\mathbf{90.4\pm 0.45}$   & $80.5\pm0.5$  & $\mathbf{93.76\pm 0.33}$  & $76.9\pm 0.3$ \\
& G-PBGD     & $83.55\pm 0.57$  & $71.3\pm 0.4$ & $82.64\pm 0.62$  & $72.3\pm 0.4$ \\
& V-PBGD     & $89.43\pm 0.23$  & $\mathbf{83.2\pm 0.3}$ & $86.56\pm 0.10$  & $\mathbf{80.2\pm 0.2}$ \\
& \textbf{PSZO-MinSel}
             & $88.5\pm 0.51$   & $78.2\pm 0.4$ & $90.3\pm 0.44$   & $79.2\pm 0.2$ \\
\cmidrule(lr){1-6}
\multirow{5}{*}{60}
& RHG        & $62.74 \pm 0.69$ & $33.8\pm 0.4$ & $63.34 \pm 0.42$ & $35.5\pm 0.3$ \\
& IAPTT-GM   & $79.42\pm 0.38$  & $45.4\pm0.3$  & $78.82\pm 0.45$  & $50.5\pm0.4$ \\
& G-PBGD     & $69.48\pm 0.42$  & $40.2\pm 0.2$ & $66.65\pm 0.42$  & $46.5\pm0.6$ \\
& V-PBGD     & $73.14\pm 0.25$  & $48.7\pm 0.3$ & $72.63\pm 0.28$  & $52.4\pm0.3$ \\
& \textbf{PSZO-MinSel}
             & $\mathbf{80.50\pm 0.30}$  & $\mathbf{57.3\pm 0.4}$ & $\mathbf{82.50\pm 0.40}$  & $\mathbf{56.5\pm 0.4}$ \\
\cmidrule(lr){1-6}
\multirow{5}{*}{80}
& RHG        & $50.2 \pm 0.4$   & $24.5\pm 0.4$ & $54.5 \pm 0.5$   & $20.9\pm 0.3$ \\
& IAPTT-GM   & $65.4\pm 0.54$   & $30.1\pm0.4$  & $64.54\pm 0.42$  & $32.9\pm 0.2$ \\
& G-PBGD     & $67.5\pm 0.33$   & $30.9\pm 0.8$ & $62.60\pm 0.34$  & $30.4\pm0.4$ \\
& V-PBGD     & $64.23\pm 0.22$  & $32.5\pm 0.2$ & $61.22\pm 0.23$  & $30.9\pm 0.4$ \\
& \textbf{PSZO-MinSel}
             & $\mathbf{70.34\pm 0.45}$  & $\mathbf{36.7\pm 0.5}$ & $\mathbf{68.43\pm 0.38}$  & $\mathbf{39.1\pm 0.4}$ \\
\bottomrule
\end{tabular}
}
\end{table}

\noindent
At low noise ($p=40\%$), PSZO-MinSel is slightly below the strongest baselines but remains comparable under the same runtime budget; as the pollute rate increases ($p=60\%,80\%$), PSZO-MinSel achieves the best performance, consistent with the benefits of trajectory-based selection under heavier corruption.

\color{black}
 
\bibliographystyle{plainurl}
\bibliography{reference}
\end{document}